\documentclass{amsart}

\usepackage{amsfonts,amssymb,verbatim,amsmath,amsthm,latexsym,textcomp,amscd}
\usepackage{latexsym,amsfonts,amssymb,epsfig,verbatim}
\usepackage{amsmath,amsthm,amssymb,latexsym,graphics,textcomp}
\usepackage{graphicx}
\usepackage{color}
\usepackage{url}
\usepackage{stmaryrd}

\usepackage{tikz}
\usepackage{mathtools}
\usetikzlibrary{positioning}
\usepackage[normalem]{ulem}
\newcommand\greensout{\bgroup\markoverwith
	{\textcolor{blue}{\rule[.5ex]{2pt}{.5pt}}}\ULon}
\usepackage{xcolor}
\begin{document}

\newtheorem{theorem}{Theorem}[section]
\newtheorem{prop}[theorem]{Proposition}
\newtheorem{lemma}[theorem]{Lemma}
\newtheorem{cor}[theorem]{Corollary}
\newtheorem{defn}[theorem]{Definition}
\newtheorem{notation}[theorem]{Notation}
\newtheorem{convention}[theorem]{Convention}
\newtheorem{conj}[theorem]{Conjecture}
\newtheorem{claim}[theorem]{Claim}
\newtheorem{rem}{Remark}
\newtheorem{example}{Example}
\newcommand{\map}{\rightarrow}
\newcommand{\C}{\mathcal C}
\newcommand\AAA{{\mathcal A}}
\def\AA{\mathcal A}

\def\L{{\mathcal L}}
\def\al{\alpha}
\def\A{{\mathcal A}}

\newcommand\BB{{\mathcal B}}
\newcommand\DD{{\mathcal D}}
\newcommand\EE{{\mathcal E}}
\newcommand\FF{{\mathcal F}}
\newcommand\GG{{\mathcal G}}
\newcommand\GB{{\mathbb G}}
\newcommand\HH{{\mathcal H}}
\newcommand\II{{\mathcal I}}
\newcommand\JJ{{\mathcal J}}
\newcommand\KK{{\mathcal K}}
\newcommand\LL{{\mathbb L}}
\newcommand\LS{{\mathcal L}} 
\newcommand\MM{{\mathcal M}}
\newcommand\NN{{\mathbb N}}
\newcommand\OO{{\mathcal O}}
\newcommand\PP{{\mathcal P}}
\newcommand\QQ{{\mathbb Q}}
\newcommand\QQQ{{\mathcal Q}}
\newcommand\RR{{\mathbb R}}
\newcommand\SSS{{\Sigma}}
\newcommand\TT{{\mathcal T}}
\newcommand\UU{{\mathcal U}}
\newcommand\VV{{\mathcal V}}
\newcommand\WW{{\mathcal W}}
\newcommand\XX{{\mathcal X}}
\newcommand\YY{{\mathcal Y}}
\newcommand\ZZ{{\mathbb Z}}
\newcommand\hhat{\widehat}
\newcommand{\Lam}[1]{\ensuremath{\partial^{(2)}_{#1}}}

\newcommand{\secref}[1]{Section~\ref{#1}}
\newcommand{\thmref}[1]{Theorem~\ref{#1}}
\newcommand{\lemref}[1]{Lemma~\ref{#1}}
\newcommand{\rmkref}[1]{Remark~\ref{#1}}
\newcommand{\propref}[1]{Proposition~\ref{#1}}
\newcommand{\corref}[1]{Corollary~\ref{#1}}
\newcommand{\probref}[1]{Problem~\ref{#1}}
\newcommand{\eqnref}[1]{~{\textrm(\ref{#1})}}
\newcommand{\boldG}{\mbox{{\bf G}}}

\def\Ga{\Gamma}
\def\Z{\mathbb Z}

\def\diam{\operatorname{diam}}
\def\dist{\operatorname{dist}}
\def\hull{\operatorname{Hull}}

\def\length{\operatorname{length}}
\newcommand\RED{\textcolor{red}}
\newcommand\GREEN{\textcolor{green}}
\newcommand\BLUE{\textcolor{blue}}
\def\mini{\scriptsize}

\def\acts{\curvearrowright}
\def\embed{\hookrightarrow}

\def\ga{\gamma}
\newcommand\la{\lambda}
\newcommand\eps{\epsilon}
\def\geo{\partial_{\infty}}
\def\bhb{\bigskip\hrule\bigskip}

\title[Pullbacks of metric bundles]{Pullbacks of metric bundles and Cannon-Thurston maps}

\author{Swathi Krishna}
\author{Pranab Sardar}
\address{Indian Institute of Science Education and Research, Mohali, Punjab 140306, India.}
\date{\today}

\begin{abstract}
Metric (graph) bundles were defined by Mj and Sardar in \cite{pranab-mahan}. In this paper, we introduce the notion of morphisms and pullbacks of metric (graph) bundles.
Given a metric (graph) bundle $X$ over $B$ where $X$ and all the fibers are uniformly (Gromov) hyperbolic and nonelementary,
and a Lipschitz quasiisometric embedding $i: A\map B$ we show that the pullback $i^{*}X$ is
hyperbolic and the map $i^{*}: i^{*}X\map X$ admits a continuous boundary extension, i.e. the Cannon-Thurston (CT) map 
$\partial i^{*}:\partial( i^{*}X) \map \partial X$.  As an application of our theorem, we show that given a short exact 
sequence of nonelementary hyperbolic groups $1\map N\map G\stackrel{\pi}{\map} Q\map 1$ and a finitely generated quasiisometrically
embedded subgroup $Q_1<Q$, $G_1:=\pi^{-1}(Q_1)$ is hyperbolic and the inclusion $G_1\map G$ admits the CT map $\partial G_1\map \partial G$.
We then derive several interesting properties of the CT map.
\end{abstract}

\maketitle

\tableofcontents


\section{Introduction}

Given a hyperbolic group $G$ and a hyperbolic subgroup $H$ a natural question to ask is 
if the inclusion $H\map G$ always extends continuously to $\partial H\map \partial G$ (see \cite[Q 1.19]{bestvinahp}).
This question was posed by Mahan Mitra (Mj) motivated by the seminal article of Cannon and Thurston (see \cite{CT}).
In \cite{CT} the authors found the first instance of this phenomenon where $H$ is not quasiisometrically embedded in $G$.
It follows from their work that if $G=\pi_1(M)$ where $M$ is a closed hyperbolic $3$-manifold fibering over a circle and $H=\pi_1(S)$
with $S$ (an orientable closed surface of genus at least $2$) being the fiber, then the boundary extension $\partial H\map \partial G$ exists. 
More generally, one may ask for a pair of (Gromov) hyperbolic metric spaces $Y\subset X$ if there is a 
continuous extension of the inclusion $Y\map X$ to $\partial Y\map \partial X$. Such an extension is by definition
unique (see Definition \ref{ct-defn}) when it exists and is popularly known 
as the {\em Cannon-Thurston map} or {\em `CT map'} for short in Geometric Group Theory. The above question of Mahan Mitra (Mj)
has motivated numerous works. The reader is referred to \cite{mahan-icm} for a detailed history of the problem. 
Although the general question for groups has been answered in the negative recently by Baker and Riley (\cite{baker-riley})
there are many interesting questions to be answered in this context. In this paper, we pick up the following.

{\bf Question.} {\em Suppose $1\map N\map G\stackrel{\pi}{\map} Q\map 1$ is a short exact sequence of hyperbolic groups. 
Suppose $Q_1<Q$ is quasiisometrically embedded and $G_1=\pi^{-1}(Q_1)$.
Then does the inclusion $G_1<G$ admit the CT map? }

\smallskip
It follows by the results of \cite{pranab-mahan} that $G_1$ is hyperbolic (see Remark 4.4, \cite{pranab-mahan}) 
so that the question makes sense. In this paper, we answer the above question affirmatively.
However, we reformulate this question in terms of metric (graph) bundles as defined in \cite{pranab-mahan} (see section $3$ of this paper) and obtain the following more general result. One is referred to Lemma \ref{defn: bary map} and the discussion following it for the definition of barycenter map. Coarsely surjective maps are introduced in Definition \ref{defn: 1.2}(3).

{\bf Theorem \ref{CT for mgbdl}.} {\em Suppose $\pi:X\map B$ is a metric (graph) bundle such that \\
(1) $X$ is hyperbolic and \\
(2) all the fibers are uniformly hyperbolic and nonelementary, i.e. there are $\delta\geq 0$ and $R\geq 0$
such that any fiber $F$ is $\delta$-hyperbolic and the barycenter map $\partial^3_s F\map F$ is $R$-coarsely surjective. }

Suppose $i:A\map B$ is a Lipschitz, quasiisometric embedding
and $\pi_Y:Y\map A$ is the pullback bundle under $i$ (see Definition \ref{pullback-defn}). Then $i^{*}:Y\map X$ admits the CT map.

\smallskip
There are two main sources of examples of metric graph bundles mentioned in this paper where the above theorem can be applied. The first one is that
of short exact sequences of groups.

{\bf Theorem \ref{main application}.} 
{\em Suppose $1\map N\map G\stackrel{\pi}{\map} Q\map 1$ is a short exact sequence of hyperbolic groups. Suppose $Q_1<Q$ is 
quasiisometrically embedded and $G_1=\pi^{-1}(Q_1)$.
Then $G_1$ is a hyperbolic group and the inclusion $G_1<G$ admits the CT map. }

\smallskip
We note that special cases of Theorem \ref{CT for mgbdl} and Theorem \ref{main application}, namely when $A$ is a point and $Q_1=(1)$ respectively, were already known.
See Theorem 5.3 in \cite{pranab-mahan} and Theorem 4.3 in \cite{mitra-ct}. Another context where Theorem \ref{CT for mgbdl} applies is that of complexes of
hyperbolic groups. We refer to Section \ref{example: complex of groups}  for relevant definitions.

Suppose $\mathcal Y$ is a finite simplicial complex and $ \GB(\mathcal Y)$
is a developable complex of nonelementary hyperbolic groups over $\mathcal Y$. 
Suppose that for all face $\sigma$ of $\YY$, $G_{\sigma}$ is a nonelementary hyperbolic group
and for any two faces $\sigma \subset \tau$ the corresponding homomorphism $G_{\tau}\map G_{\sigma}$ is an isomorphism
onto a finite index subgroup of $G_{\sigma}$. Suppose that the fundamental group of the complex of groups, $G$ say, 
is hyperbolic. Suppose we have a {\em good} subcomplex $\YY_1\subset \YY$  i.e. one for which the following two conditions
are satisfied. \\ (1) The natural homomorphism $\pi_1(\GG, \YY_1)\map \pi_1(\GG,\YY)$ is injective.

Let $G_1=\pi_1(\GG, \YY_1)$. 
Suppose $G_1$ and $G$ are both endowed with word metrics with respect to some finite generating sets.
Let $\hat{G}$ and $\hat{G}_1$ be the coned off spaces a la Farb (\cite{farb-relhyp}), obtained by coning off all
the face groups in $G$ and $G_1$ respectively.\\
(2) Then the induced map $\hat{G}_1\map \hat{G}$ of the coned off spaces is a quasiisometric embedding. 
With these hypotheses we have:

{\bf Theorem \ref{CT for complex of groups}.} 
{\em The group $G_1$ is hyperbolic and the inclusion $G_1\map G$ admits the CT map.}

Particularly interesting cases to which the above theorem applies are obtained in \cite{min} and \cite{mahan-pritam}. There
graphs of groups are considered where all the vertex and edge groups are either surface groups (\cite{min}) or free groups of rank $\geq 3$ (\cite{mahan-pritam}) respectively.

Next, we explore properties of the Cannon-Thurston map $\partial Y\map \partial X$ proved in Theorem \ref{CT for mgbdl}. 
Suppose $F$ is a fiber of the bundle $Y$ over $A$. Then there is a CT map for the
inclusions $i_{F,X}: F\map X$ and $i_{F,Y}: F\map Y$, and the map $i^{*}:Y\map X$. Since $\partial i_{F,X}=\partial i^{*}\circ \partial i_{F,Y}$, if $\alpha, \beta\in \partial F$ are identified under
$\partial i_{F,X}$ then under $\partial i^{*}$ the points $\partial i_{F,Y}(\alpha)$ and   $\partial i_{F,Y}(\beta)$ are identified too. 
It turns out that a sort of `converse' of this is also true.

\smallskip
 {\bf Theorem \ref{CT leaf in Y}.} {\em Suppose we have the hypotheses of Theorem \ref{CT for mgbdl} and also that the fibers of the
bundle are proper metric spaces.
Suppose $\gamma$ is a (quasi)geodesic line in $Y$ such that $\gamma(\infty)$
and $\gamma(-\infty)$ are identified by the CT map $\partial i^{*}: \partial Y\map \partial X$. Then $\pi_Y(\gamma)$ is bounded.
In particular, given any fiber $F$ of the metric bundle, $\gamma$ is at a finite Hausdorff distance from a quasigeodesic line of $F$.}

\smallskip
On the other hand as an immediate application of Theorem \ref{CT leaf in Y} (in fact, see Corollary \ref{thm:CT surjective}
and Proposition \ref{bundle boundary}) we get the following:

{\bf Theorem.} {\em Suppose we have the hypotheses of Theorem \ref{CT for mgbdl} and also that the fibers of the bundle are proper metric spaces.
Let $F$ be the fiber over a point $b\in A$. Then the CT map $\partial i_{F,X}:\partial F\map \partial X$
is surjective if and only if the CT maps $\partial i_{F, Y_{\xi}}:\partial F\map \partial Y_{\xi}$ are surjective 
for all $\xi\in \partial B$ where $Y_{\xi}$ is the pullback of a (quasi)geodesic ray in $B$ asymptotic to $\xi$.

In particular $\partial i_{F,Y}:\partial F\map \partial Y$ is surjective if $\partial i_{F,X}: \partial F\map \partial X$ is surjective.
}

\smallskip
Following Mitra (\cite{mitra-endlam}) we define the {\em Cannon-Thurston lamination} $\Lam{X}(F)$ to be 
$\{(z_1,z_2)\in \partial F\times \partial F:z_1\neq z_2,\, \partial i_{F,X}(z_1)=\partial i_{F,X}(z_2)\}$ and following Bowditch
(\cite[Section 2.3]{bowditch-stacks}) we define for any point $\xi\in \partial B$ a subset of this lamination denoted by $\Lam{\xi,X}(F)$  or simply $\Lam{\xi}(F)$
when $X$ is understood, where $(z_1, z_2)\in \Lam{\xi,X}(F)$ if and only if
$\partial i_{F,X}(z_1)= \partial i_{F,X}(z_2)=\tilde{\gamma}(\infty)$ where $\tilde{\gamma}$ is a quasiisometric 
lift in $X$ of a (quasi)geodesic ray $\gamma$ in $B$ converging to $\xi$. 
If $(z_1,z_2)\in \Lam{\xi,X}(F)$ and $\alpha$ is a (quasi)geodesic line in $F$ connecting $z_1,z_2$ then $\alpha$ is referred to be a {\em leaf} of the lamination
$\Lam{\xi,X}(F)$. Leaves are assumed to be uniform quasigeodesics in the following theorem using Proposition \ref{visibility}.

\smallskip
{\bf Theorem.} (See Lemma \ref{CT lamination} through Lemma \ref{CT property final}.) {\bf (Properties of $\Lam{X}(F)$)}

{\em $(1)$ $\Lam{X}(F)= \coprod_{\xi\in \partial B} \Lam{\xi,X}(F)$.

$(2)$ $\Lam{X}(F)$ and $\Lam{\xi,X}(F)$ are all closed subsets of $\partial^{(2)} F$
where $\partial^{(2)} F=\{(z_1, z_2)\in \partial F\times \partial F:z_1\neq z_2\}.$

$(3)$ The leaves of $\Lam{\xi_1, X}(F), \Lam{\xi_2, X}(F)$ are coarsely transverse to each other for all $\xi_1\neq \xi_2\in \partial B$:

Given $\xi_1\neq \xi_2\in \partial B$ and $D>0$  there exists $R>0$  such that if $\gamma_i$ is leaf of $\Lam{\xi_i,X}(F)$, $i=1,2$ then $\gamma_1\cap N_D(\gamma_2)$ has 
diameter less than $R$.

$(4)$ If $\xi_n\map \xi$ in $\partial B$ and $\alpha_n$ is a leaf of $\Lam{\xi_n, X}(F)$ for all $n\in \NN$  which converge to a
geodesic line $\alpha$ then $\alpha$ is a leaf of $\Lam{\xi,X}(F)$.

$(5)$ $\Lam{\xi,X}(F)=\Lam{\xi,Y}(F)$ for all $\xi \in \partial A$ if we have the hypothesis of Theorem \ref{CT for mgbdl}. }

\smallskip
Finally, we also prove the following interesting property of the CT lamination.

{\bf Theorem \ref{main application2}.} {\em
Suppose $X$ is a metric (graph) bundle over $B$ satisfying the hypotheses of Theorem \ref{CT for mgbdl} such
that $X$ is a proper metric space. Let $F=F_b$ where $b\in B$.
Suppose $\partial F$ is not homeomorphic to a dendrite and also the CT map $\partial F\map \partial X$
is surjective. 

Then for all $\xi\in \partial B$ we have $\Lam{\xi, X}(F)\neq \emptyset$.
}

This applies in particular to the examples of short exact sequence of hyperbolic groups and the complexes of
hyperbolic groups mentioned in Theorem \ref{main application} and Theorem \ref{CT for complex of groups}
above.

\medskip

\noindent 
{\bf Outline of the paper:} In section 2 we recall basic hyperbolic geometry, Cannon-Thurston maps, etc. In section 3 we recall the
basics of metric (graph) bundles and we introduce morphisms of bundles, pullbacks. Here we prove the existence of pullbacks under suitable assumptions.
In section 4 we mainly recall the machinery of \cite{pranab-mahan} and we prove a few elementary results. Section 5 is devoted to the proof of the 
main theorem. In section 6 we derive applications of the main result and we mention some related results.

\smallskip

\section{Hyperbolic metric spaces}
In this section, we remark on the notation and convention to be followed in the rest of the paper
and we put together basic definitions and results about hyperbolic metric spaces. 
We begin with some basic notions from large scale geometry.
Most of these are quite standard, e.g. see \cite{gromov-hypgps}, \cite{GhH}. We have used 
\cite{pranab-mahan} where all the basic notions can be quickly found in one place.

{\bf Notation, convention and some metric space notions.} 
One is referred to \cite[Chapter I.1, I.3]{bridson-haefliger}
for the definitions and basic facts about geodesic metric spaces, metric graphs and length spaces.\\
(0) For any set $A$, $Id_A$ will denote the {\em identity} map $A\map A$. If $A\subset B$ then we denote by
$i_{A,B}:A\map B$ the inclusion map of $A$ into $B$.\\
(1) If $x\in X$ and $A\subset X$ then $d(x,A)$ will denote
$\inf\{d(x,y):y\in A\}$ and will be referred to as the {\em distance of $x$ from $A$}. 
For $D\geq 0$ and $A\subset X$,  
$N_D(A):=\{x\in X: \, d(x,a)\leq D \, \mbox{ for some} \, a\in A\}$ will be called
the $D$-{\em neighborhood} of $A$ in $X$. 
For $A,B\subset X$ we shall denote by $d(A,B)$ the quantity $\inf\{d(x,B):x\in A\}$ and by $Hd(A,B)$
the quantity $\inf\{D>0: A\subset N_D(B), B\subset N_D(A)\}$ and will refer to it as the {\em Hausdorff distance} of $A,B$.\\
(2) If $X$ is a length space we consider only subspaces $Y\subset X$ such that the induced length
metric on $Y$ takes values in $[0,\infty)$, or equivalently for any pair of points in $Y$ there is a rectifiable
path in $X$ joining them which is contained in $Y$. We shall refer to such subsets as {\em rectifiably path connected}.
If $\gamma$ is a rectifiable path in $X$ then $l(\gamma)$ will denote the length of $\gamma$.\\
(3) All graphs are connected for us.
If $X$ is a metric graph then $\VV(X)$ will denote the set of vertices of $X$. Generally, we shall write $x\in X$
to mean $x\in \VV(X)$. In metric graphs (see \cite[Chapter I.1]{bridson-haefliger}) all the edges are assumed to have 
length $1$.  In a graph $X$ the paths are assumed to be a sequence of vertices. In other words, these are maps
$I\cap \mathbb Z\map X$ where $I$ is a closed interval in $\RR$ with end points in $ \mathbb Z\cup \{\pm \infty\}$. We shall
informally write this as $\alpha: I\map X$ and sometimes refer to it as a {\em dotted path} for emphasis. 
Length of such a path $\alpha: I\map X$ is defined to be $l(\alpha)=\sum d(\alpha(i), \alpha(i+1))$ where the sum is
taken over all $i\in \ZZ$ such that $i, i+1\in I$.
If $\alpha: [0,n]\map X$ and $\beta: [0, m]\map X$ are two paths
with $\alpha(n)=\beta(0)$ then their concatenation $\alpha*\beta$ will be the path $[0,m+n]\map X$ defined by 
$\alpha*\beta(i)=\alpha(i)$ if $i\in [0,n]$ and $\alpha*\beta(j)=\beta(j-n)$ if $j\in [n,m+n]$.\\
(4) If $X$ is a geodesic metric space and $x,y\in X$ then we shall use $[x,y]_X$ or simply $[x,y]$ to denote a
geodesic segment joining $x$ to $y$. This applies in particular to metric graphs.  For $x,y,z\in X$ we 
shall denote by $\Delta xyz$ some geodesic triangle with vertices $x,y,z$.  \\
(5) If $X$ is any metric space then for all $A\subset X$, $diam(A)$ will denote the diameter of $A$.

\subsection{Basic notions from large scale geometry}
Suppose $X$, $Y$ are any two metric spaces and $k\geq 1, \, \epsilon \geq 0, \epsilon'\geq 0$. 

\begin{defn}{\em (\cite[Definition 1.1.1]{pranab-mahan})}\label{defn: 1.2}
\begin{enumerate}
\item A map $\phi:X\rightarrow Y$ is said to be {\em metrically proper} if there is an increasing function
$f:[0,\infty)\map [0,\infty)$ with $\lim_{t\map \infty} f(t)=\infty$ 
such that for any $x,y\in X$ and $R\in [0,\infty)$, $d_Y(\phi(x),\phi(y))\leq R$
implies $d_X(x,y)\leq f(R)$.  In this case we say that $\phi$ is proper {\em as measured by $f$}.
\item A subset $A$ of a metric space $X$ is said to be $r$-{\em dense} in $X$ for some $r\geq 0$ if $N_r(A)=X$.
\item Suppose $A$ is a set. A map $\phi:A \rightarrow Y$ is said to be $\epsilon$-{\em coarsely surjective}
if $\phi(A)$ is $\epsilon$-dense in $Y$. We will say that it is coarsely surjective if it is $\epsilon$-coarsely
surjective for some $\epsilon\geq 0$.
\item A map $\phi: X\rightarrow Y$ is said to be  {\em coarsely $(\epsilon, \epsilon')$-Lipschitz} if for every $x_1,x_2\in X$, we have
$d(\phi(x_1),\phi(x_2))\leq \epsilon d(x_1,x_2) + \epsilon' $. A coarsely $(\epsilon, \epsilon)$-Lipschitz map
will be simply called a {\em coarsely $\epsilon$-Lipschitz} map. A map $\phi$ is {\em coarsely Lipschitz} if it is coarsely
$\epsilon$-Lipschitz for some $\epsilon \geq 0$.
\item $(i)$ A map $\phi:X\rightarrow Y$ is said to be
a $(k,\epsilon)$-{\em quasiisometric embedding} if for every $x_1,x_2\in X$, one has 
$$ -\epsilon + d(x_1,x_2)/k \leq d(\phi(x_1),\phi(x_2))\leq \epsilon+kd(x_1,x_2).$$ 
A map $\phi:X\rightarrow Y$ will simply be referred to as
a quasiisometric embedding if it is a $(k,\epsilon)$-quasiisometric embedding for some $k\geq 1$, $\epsilon\geq 0$.
A $(k,k)$-quasiisometric embedding will  be referred to as a $k$-{\em quasiisometric embedding}. 

$(ii)$ A map $\phi:X\rightarrow Y$ is a $(k,\epsilon)$-{\em quasiisometry}
(resp. $k$-{\em quasiisometry}) if it is a $(k,\epsilon)$-quasiisometric embedding (resp. 
$k$-quasiisometric embedding) and moreover, it is $D$-coarsely surjective for some $D\geq 0$.

$(iii)$ A $(k,\epsilon)$-{\em quasigeodesic} (resp. a $k$-{\em quasigeodesic}) in a metric space $X$ is 
a $(k,\epsilon)$-quasiisometric embedding (resp. a $k$-quasiisometric embedding) $\gamma:I\rightarrow X$, where
$I\subseteq \mathbb R$ is an interval. 

We recall that a $(1,0)$-quasigeodesic is called a {\em geodesic}.
 
If $I=[0,\infty)$, then $\gamma$ will be called a {\em quasigeodesic ray}. If $I=\RR$, then we
call it a {\em quasigeodesic line}. One similarly defines a {\em geodesic ray} and a {\em geodesic line}.
 We refer to the constant(s) $k$ (and $\epsilon$) as {\em quasigeodesic constant(s).}

Quasigeodesics in a metric graph $X$ will be maps $I\cap \ZZ\map X$, informally written as
$I\map X$ where $I$ is a closed interval in $\RR$.
\begin{comment}\item Suppose $I\subset \RR$ is a closed interval with end points in $\ZZ\cup \{\pm \infty\}$ and $J=\ZZ\cap I$ with the restricted metric
from $\RR$. Then a $(k,\epsilon)$-qi embedding $\alpha: J\map X$ will be called a {\em dotted $(k,\epsilon)$-quasigeodesic}.
If $I$ is a finite interval, say $I=[m,n]$ then $\sum_{i=m}^{n-1} d_X(\alpha(i), \alpha(i+1))$ will be called the
{\em length} of the dotted quasigeodesic $\alpha$. 

Similarly we define {\em dotted $k$-quasigeodesic, dotted $k$-quasigeodesic ray, dotted $k$-quasigeodesic line} etc.
\item A map $\psi: Y\rightarrow X$ is said to be an $\epsilon$-{\em coarse left (or right) inverse} of a map $\phi: X\rightarrow Y$ if 
for all $x\in X$ one has $d_X(\psi\circ \phi(x),x)\leq \epsilon$ (respectively for all $y\in Y$, $d_X(\phi\circ \psi(y),y)\leq \epsilon$).
\end{comment}
\item Suppose $\phi,\phi':X\map Y$ are two maps and $\epsilon\geq 0$.

(i) We define $d(\phi,\phi')$ to be the quantity $\sup\{d_Y(\phi(x),\phi'(x)): x\in X\}$ provided the supremum
exists in $\RR$; otherwise we write $d(\phi, \phi')=\infty$.

(ii) A map $\psi:Y\map X$ is called an $\epsilon$-coarse left (right) inverse of $\phi$ if $d(\psi\circ \phi, Id_X)\leq \epsilon$
(resp. $d(\phi\circ \psi, Id_Y)\leq \epsilon$).

If $\psi$ is both an $\epsilon$-coarse left and right inverse then it is simply called an $\epsilon$-{\em coarse inverse} of $\phi$.
\item Suppose $S$ is any set. A map $f:S\map X$ satisfying some properties $\mathcal P_1, \cdots, \mathcal P_k$
will be called {\em coarsely unique} if for any other map $g:S\map X$ with properties $\mathcal P_1, \cdots, \mathcal P_k$
there is a constant $D$ such that $d(f,g)\leq D$.

\end{enumerate}
\end{defn}
The definition (7) above is taken from \cite{pranab-mahan}. See the definition following Lemma 2.9 there. 
In places where this definition will be used the properties may not be explicitly stated but they will be clear
from the context. If $S$ is finite then we talk about a finite subset of $X$ to be coarsely unique, e.g. see the
remark following Lemma \ref{qc last lemma}.

{\bf Remark on terminology:} 
(1) All the above definitions are about certain properties of maps and in each
case some parameters are involved. 

(i) When the parameters are not important or they are clear from the context then we say that the map
has the particular property without explicit mention of the parameters, e.g. `$\phi:X\map Y$ is metrically
proper' if $\phi$ is metrically proper as measured by some function.

(ii) When we have a set of pairs of metric spaces and a map between each pair possessing the
same property with the same parameters then we say that the set of maps `uniformly' have the property, e.g.
{\em uniformly metrically proper, uniformly coarsely Lipschitz, uniform qi embeddings, uniform approximate nearest point
projection} etc. 

(2) We often refer to quasiisometric embeddings as `qi embedding' and quasiisometry as `qi'.

The following gives a characterization of quasiisometry to be used in the discussion on metric bundles.
\begin{lemma}\label{elem-lemma1} 
{\em (\cite[Lemma 1.1]{pranab-mahan})}
\begin{enumerate}
\item 
For every $K_1,K_2\geq 1$ and $D\geq 0$ there are $K_{\ref{elem-lemma1}}=K_{\ref{elem-lemma1}}(K_1,K_2,D)$,
such that the following hold.

A $K_1$-coarsely Lipschitz map with a $K_2$-coarsely Lipschitz, $D$-coarse inverse is a 
$K_{\ref{elem-lemma1}}$-quasiisometry.

\item  Given $K\geq 1$, $\epsilon \geq 0$  and $R\geq 0$ there are constants $C_{\ref{elem-lemma1}}=C_{\ref{elem-lemma1}}(K,\epsilon, R)$ and 
$D_{\ref{elem-lemma1}}=D_{\ref{elem-lemma1}}(K,\epsilon, R)$ such that the following holds:

Suppose $X,Y$ are any two metric spaces and $f:X\map Y$ is a $(K,\epsilon)$-quasiisometry which is $R$-coarsely surjective.
Then there is a $(K_{\ref{elem-lemma1}},C_{\ref{elem-lemma1}})$-quasiisometric $D_{\ref{elem-lemma1}}$-coarse inverse of $f$.
\end{enumerate}
\end{lemma}

The following lemmas follow from simple calculations and hence we omit their proofs.
\begin{lemma}\label{qi composition}
(1) Suppose we have a sequence of maps $X\stackrel{f}{\map} Y\stackrel{g}{\map}Z$ where $f,g$ are coarsely $L_1$-Lipschitz
and $L_2$-Lipschitz respectively. Then $g\circ f$ is  coarsely $(L_1L_2, L_1L_2+L_2)$-Lipschitz.

(2) Suppose $f:X\map Y$ is a $(K_1,\epsilon_1)$-qi embedding and $g:Y\map Z$ is a $(K_2,\epsilon_2)$-qi embedding. 
Then $g\circ f:X\map Z$ is a $(K_1K_2,K_2\epsilon_1+\epsilon_2)$-qi embedding. 

Moreover, if $f$ is $D_1$-coarsely surjective and $g$ is $D_2$-coarsely surjective then $g\circ f$ is
$(K_2D_1+\epsilon_2+D_2)$-coarsely surjective.

In particular, the composition of finitely many quasiisometries is a quasiisometry.
\end{lemma}

\begin{lemma} \label{lemma: qi graph}
Suppose $X'$ is any connected graph and $r>0$. Suppose $X$ is another graph obtained from $X'$ by introducing
some new edges to $X'$ where $e=[v,w]$ is an edge in $X$ but not in $X'$ implies $d_{X'}(v,w)\leq r$.
Then the inclusion map $X'\map X$ is a quasiisometry.
\end{lemma}


The following lemma appears in \cite[Section 1.5]{ps-kap} in a somewhat different form. We include a proof for the sake of completeness.
\begin{lemma}\label{quasigeod criteria}
Let $X$ be any metric space, $x,y\in X$, $\gamma$ be a (dotted) $k$-quasigeodesic joining $x,y$ and 
$\alpha:I\map X$ is a (dotted) coarsely
$L$-Lipschitz path joining $x,y$. Suppose moreover, $\alpha$ is a proper embedding as measured by a function 
$f:[0,\infty)\map [0,\infty)$ and that $Hd(\alpha, \gamma)\leq D$ for some $D\geq 0$.
Then $\alpha$ is (dotted) $K_{\ref{quasigeod criteria}}=K_{\ref{quasigeod criteria}}(k,f,D,L)$-quasigeodesic in $X$.
\end{lemma}

\proof Suppose $\gamma$ is defined on an interval $J$.
Let $a,b\in I$. Then we have $$ \,\,d(\alpha(a), \alpha(b))\leq L|a-b|+L \longrightarrow (1)$$ since $\alpha$ is coarsely
$L$-Lipschitz. Now let $a',b'\in J$ be such that $d(\alpha(a), \gamma(a')\leq D$ and $d(\alpha(b), \gamma(b'))\leq D$.
Let $R=d(\alpha(a), \alpha(b))$. Then by triangle inequality $d(\gamma(a'), \gamma(b'))\leq 2D+R$.
Since $\gamma$ is a $k$-quasigeodesic we have $-k+ |a'-b'|/k\leq d(\gamma(a'), \gamma(b'))\leq 2D+R$.
Hence, $|a'-b'|\leq k(2D+R)+k^2$. Without loss of generality suppose $a'\leq b'$. Consider the sequence
of points $a'_0=a', a'_1,\cdots, a'_n=b'$ in $J$ such that $a'_{i+1}=1+a'_i$ for $0\leq i\leq n-2$ and
$a'_n-a'_{n-1}\leq 1$. We note that $n\leq 1+k(2D+R)+k^2$.
Let $a_i\in I$ be such that $d(\gamma(a'_i), \alpha(a_i))\leq D$, $0\leq i\leq n$ where $a_0=a, a_n=b$.
Once again by triangle inequality we have 
$$d(\alpha(a_i), \alpha(a_{i+1}))\leq 2D+d(\gamma(a'_i), \gamma(a'_{i+1}))\leq 2D+2k$$
for $0\leq i\leq n-1$ since $\gamma$ is a $k$-quasigeodesic. This implies $|a_i-a_{i+1}|\leq f(2D+2k)$
since $\alpha$ is a proper embedding as measured by $f$. Hence,
$$|a-b|\leq \sum_{i=0}^{n-1} |a_i-a_{i+1}|\leq nf(2D+k)\leq (1+k(2D+R)+k^2)f(2D+2k).$$
Thus we have 
$$\,\,-\frac{1+2kD+k^2}{k}+\frac{1}{kf(2D+2k)} |a-b|\leq R=d(\alpha(a),\alpha(b)). \longrightarrow (2)$$
Hence, by (1) and (2) we can take
$$K_{\ref{quasigeod criteria}}=1+2D+k+kf(2D+2k)+L. \qed$$

The following lemma is implicit in the proof  \cite[Proposition 2.10]{pranab-mahan}. 
The proof of this lemma being immediate we omit it. 
\begin{lemma}\label{proving Lipschitz}
Suppose $X$ is a length space and $Y$ is any metric space. Let $f:X\map Y$
be any map. Then $f$ is coarsely $C$-Lipschitz for some $C\geq 0$ if
for all $x_1, x_2\in X$, $d_X(x_1,x_2)\leq 1$ implies
$d_Y(f(x_1),f(x_2))\leq C$.
\end{lemma}

\begin{rem}
We spend quite some time restating some results proved in \cite{pranab-mahan} in the generality of length spaces
since the main result in our paper is about length spaces. For instance (1) the existence of pullback of metric
bundles to be defined below is unclear within the category of geodesic metric spaces; and (2) we observe that for 
the definition of Cannon-Thurston maps the assumption of (Gromov) hyperbolic {\em geodesic} metric spaces is rather 
restrictive and unnecessary. 
\end{rem}
\noindent
In a length metric space geodesics may not exist joining a pair of points. However, we still have the following.
\begin{lemma}\label{length space lemma1}
Suppose $X$ is a length space.
(1) Given any $\epsilon>0$, any pair of points of $X$ can be joined by a continuous,
rectifiable, arc length parameterized path which is a $(1,\epsilon)$-quasigeodesic.

(2) Any pair of points of $X$ can be joined by a dotted $1$-quasigeodesic.
\end{lemma}

{\bf Metric graph approximation to a length space}

\smallskip
Given any length space $X$, we define a metric graph $Y$ as follows. We take the vertex set $V(Y)=X$. We join $x,y\in X$ by an
edge (of length $1$) if and only if $d_X(x,y)\leq 1$. We let $\psi_X:X\map \mathcal V(Y)\subset Y$ be the identity map. Let 
$\phi_X: Y\map X$ be  defined to be the inverse of $\psi_X$ on $\VV(Y)$ and for any point $y$ in the interior of an edge $e$ of $Y$
we define $\phi_X(y)$ to be one of the end points of the edge $e$. The following hold.

\begin{lemma}\label{length space qi to graph}\cite[Lemma 1.32]{ps-kap}
(1) $Y$ is a (connected) metric graph. (2) The maps $\psi_X$ and $\phi_X|_{\VV(Y)}$ are coarsely $1$-surjective, 
$(1,1)$-quasiisometries. (3)
The map $\phi_X$ is a $(1,3)$-quasiisometry and it is a $1$-coarse inverse of $\psi_X$.
\end{lemma}
\begin{comment}
\proof  
Given $x,y$, as in the proof of Lemma \ref{length space lemma1}(1), we join them by an arc length parameterized path 
$\gamma:[0,l]\map X$ such that $l\leq d(x,y)+\epsilon$ where $\epsilon>0$ is chosen in such a way that $d(x,y)<m+1$ where $m$
is the non-negative integer determined by $m\leq d(x,y)<m+1$. Since $l(\gamma)<m+1$ it follows that
$d_Y(x,y)\leq m+1\leq d_X(x,y)+1$. Suppose $x,y\in X$ such that $d_Y(x,y)= n$. Let $x=x_0,x_1,\ldots, x_n=y$ the consecutive vertices 
on a geodesic in $Y$ joining $x,y$. Then we know that $d_X(x_i, x_{i+1})\leq 1$. Thus $d_X(x,y)\leq \sum^n_{i=1}d_X(x_{i-1},x_i)\leq n$.
Thus we get $d_X(x,y)\leq d_Y(x,y)\leq d_X(x,y)+1$. This proves the first statement of the lemma. 

Finally, it is clear that $N_1(\psi_X(X))=Y$ and hence $\psi_X$ is coarse $1$-surjective. The remaining part of the proof follows
from a simple calculation and so we omit the proof.\qed	
\end{comment}
\begin{rem}
We shall refer to the space $Y$ constructed in the proof of the above lemma as the 
{\em (canonical) metric graph approximation} to $X$.
We also preserve the notations $\psi_X$ and $\phi_X$ to be used in this context only.

\end{rem}

\begin{defn} {\em Gromov inner product:} Let $X$ be any metric space and let $p,x,y\in X$.
Then the {\em Gromov inner product of $x,y$ with respect to $p$} is defined to
be the number $\frac{1}{2}(d(p,x)+d(p,y)-d(x,y))$. It is denoted by $(x.y)_p$.
\end{defn}

\begin{lemma}\label{gromov product step1}
Suppose $X$ is a length space and $x_1,x_2,x_3\in X$. Let $\gamma_{ij}, i< j, 1\leq i,j\leq 3$ 
denote $(1,1)$-quasigeodesics joining the respective pairs of points $x_i, x_j$. Suppose there are
points $w_1\in \gamma_{23}, w_2\in \gamma_{13}$ and $w_3\in \gamma_{12}$ such that
$d(w_1, w_i)\leq R$ for some $R\geq 0$, $i=2,3$. Then $|(x_2.x_3)_{x_1}-d(x_1,w_1)|\leq 3+2R$.
\end{lemma}
\proof By triangle inequality we have $|d(x_2,w_1)-d(x_2,w_2)|\leq R$, $|d(x_3,w_1)-d(x_3,w_2)|\leq R$,
$|d(x_1,w_1)-d(x_1,w_i)|\leq R$, $i=2,3$. Since the $\gamma_{ij}$'s
are $(1,1)$-quasigeodesics it is easy to see that
$d(x_1,w_3)+d(w_3,x_2)\leq d(x_1,x_2)+3$, $d(x_1,w_2)+d(w_2,x_3)\leq d(x_1,x_3)+3$
and $d(x_2,w_1)+d(w_1,x_3)\leq d(x_2,x_3)+3$. It then follows by a simple calculation that
$$2d(x_1,w_1)-6-4R\leq d(x_1,x_2)+d(x_1,x_3)-d(x_2,x_3)\leq 2d(x_1,w_1)+3+4R.$$
Hence, we have $|(x_2.x_3)_{x_1}-d(x_1,w_1)|\leq 3+2R$.
\qed

\begin{defn}
\begin{enumerate}
\item Suppose $X$ is a length space and $Y_1,Y_2, Z$ are nonempty subsets of $X$. We say that $Z$ {\em coarsely disconnects}
$Y_1,Y_2$ in $X$ if  (i) $Y_i\setminus Z\neq \emptyset$, $i=1,2$ and (ii) for all $K\geq 1$ there is $R\geq 0$ 
such that the following holds: For any $y_i\in Y_i$, $i=1,2$ and any $K$-quasigeodesic $\gamma$ in $X$ joining 
$y_1, y_2$ we have $\gamma \cap N_R(Z)\neq \emptyset$.

\item Suppose $Y, Z\subset X$, $Y_1,Y_2\subset Y$. We say that $Z$ {\em coarsely bisects} $Y$ into $Y_1,Y_2$ in $X$ if 
$Y=Y_1\cup Y_2$ and $Z$ coarsely disconnects $Y_1, Y_2$ in $X$.

\item Suppose $\{X_i\}$ is a collection of length spaces and there are nonempty sets $Y_i, Z_i\subset X_i$, 
$Y^{+}_i, Y^{-}_i\subset Y_i$ such that $Y_i=Y^{+}_i\cup Y^{-}_i$, $Y^{+}_i\setminus Z_i\neq \emptyset$, and 
$Y^{-}_i\setminus Z_i\neq \emptyset$ for all $i$. We say that $Z_i$'s {\em uniformly coarsely bisect} $Y_i$'s 
into $Y^{+}_i$'s, and $Y^{-}_i$'s if for all $K\geq 1$ there is $R=R(K)\geq 0$ with the following property:
For any $i$, and for any $x^{+}_i\in Y^{+}_i, x^{-}_i\in Y^{-}_i$ and any $K$-quasigeodesic 
$\gamma_i\subset X_i$ joining $x^{\pm}_i$ we have $N_R(Z_i)\cap \gamma_i\neq \emptyset$.
\end{enumerate}
\end{defn}
We note that the first part of the above definition implies $Y_1\cap Y_2\subset N_{R(1)}(Z)$. Moreover one
would like to impose the condition that $Y_i\setminus Z$ are of infinite diameter. Keeping the
application we have in mind we do not assume that. 

\begin{defn}{\em (Approximate nearest point projection)}
(1) Suppose $X$ is any metric space, $A\subset X$, and $x\in X$. Given $\epsilon \geq 0$ and $y\in A$
we say that $y$ is an $\epsilon$-{\em approximate nearest point projection} of $x$ on $A$ if
for all $z\in A$ we have $d(x,y)\leq d(x,z)+\epsilon$.

(2) Suppose $X$ is any metric space, $A\subset X$ and $\epsilon \geq 0$. An $\epsilon$-{\em approximate
nearest point projection map} $f: X\map A$ is a map such that $f(a)=a$ for all $a\in A$ and 
$f(x)$ is an $\epsilon$-approximate nearest point projection of $x$ on $A$ for all $x\in X\setminus A$.

For $\epsilon=0$ an $\epsilon$-approximate nearest point projection is simply referred to as a {\em nearest
point projection}. A nearest point projection map from $X$ onto a subset $A$ will be denoted by $P_{A,X}:X\map A$
or simply $P_A:X\map A$ when there is no possibility of confusion.
\end{defn}
We note that given a metric space $X$ and $A\subset X$ a nearest point projection map $X\map A$
may not be defined in general but an $\epsilon$-approximate nearest point projection map $X\map A$ exists
by axiom of choice for all $\epsilon>0$.

\begin{lemma}\label{lem: approx proj}
Suppose $X$ is a metric space and $A\subset X$. Suppose $y\in A$ is an $\epsilon$-approximate
nearest point projection of $x\in X$. Suppose $\alpha:I\map X$ is a $(1,1)$-quasigeodesic joining
$x,y$. Then $y$ is an $(\epsilon+3)$-approximate nearest point of $x'$ on $A$ for all $x'\in \alpha$.
\end{lemma}
\proof Suppose $z\in A$ is any point. Then we know that $d(x,y)\leq d(x,z)+\epsilon$.
Since $\alpha$ is a $(1,1)$-quasigeodesic it is easy to see that $d(x,x')+d(x',y)\leq d(x,y)+3$.
Hence, $d(x,x')+d(x',y)\leq d(x,z)+3+\epsilon$ which in turn implies that
$d(x',y)\leq d(x,z)-d(x,x')+3+\epsilon\leq d(x',z)+\epsilon+3$. Hence, $y$ is an $(\epsilon+3)$-approximate
nearest point projection of $x'$ on $A$. 
\qed

\begin{cor}\label{projection on geodesic}
Suppose $X$ is any metric space and $x,y,z\in X$. Suppose $\alpha$, $\beta$ are $(1,1)$-quasigeodesics
joining $x,y$ and $y,z$ respectively. If $y$ is an $\epsilon$-approximate nearest point projection of
$x$ on $\beta$ then $\alpha*\beta$ is $(3,3+\epsilon)$-quasigeodesic. 
\end{cor}
\proof Let $x'\in \alpha$ and $y'\in \beta$. Let $\beta'$ denote the segment of $\beta$ from $y$ to $y'$.
Then $y$ is an $\epsilon$-approximate nearest point projection of $x$ on $\beta'$ too. Hence, by the
previous lemma $y$ is an $(\epsilon+3)$-approximate nearest point projection of $x'$ on $\beta'$.
Without loss of generality, suppose $\alpha(a)=x'$, $\alpha(a+m)=y$, $\beta(0)=y$, and $\beta(n)=y'$.
Now, $d(x',y)\leq d(x',y')+\epsilon+3$. Hence $d(y,y')\leq d(x',y')+d(x',y)\leq 2d(x',y')+\epsilon+3$.
Since $\alpha, \beta$ are both $(1,1)$-quasigeodesics it follows that 
$m-1\leq d(x',y)\leq d(x',y')+\epsilon+3$ and $n-1\leq d(y,y')\leq 2d(x',y')+\epsilon+3$. Adding these
we get $m+n-2\leq 3d(x',y')+2\epsilon+6$. On the other hand, $d(x',y')\leq d(x',y)+d(y, y')\leq m+n+2$.
Putting everything together we get
$$\frac{1}{3}(m+n)-\frac{2\epsilon+8}{3}\leq d(x',y')\leq (m+n)+2$$
from which the corollary follows immediately.
\qed


\subsection{Rips hyperbolicity vs Gromov hyperbolicity}
This subsection gives a quick introduction to some basic notions and results about hyperbolic metric spaces. One is
referred to  \cite{gromov-hypgps}, \cite{GhH}, \cite{Shortetal} for more details. The following definition of 
hyperbolic metric spaces is due to E. Rips and hence we refer to this as the Rips hyperbolicity.

\begin{defn} 
$(1)$ Suppose $\Delta x_1x_2x_3$ is a geodesic triangle in a metric space $X$
and $\delta\geq 0$, $K\geq 0$. We say that the triangle $\Delta x_1x_2x_3$ is $\delta$-{\em slim} if any side of the triangle is contained in the
$\delta$-neighborhood of the union of the remaining two sides.

$(2)$ Let $\delta\geq 0$ and $X$ be a geodesic metric space. We say that $X$ is $\delta$-{\em hyperbolic} 
{\em (}in the sense of Rips{\em )} if all  geodesic triangles in $X$ are $\delta$-slim.

A geodesic metric space is said to be  {\em (Rips) hyperbolic} if it is $\delta$-hyperbolic in the sense of Rips
for some $\delta\geq 0$.
\end{defn}

 However,  in this paper
we need to deal with length spaces a lot which a priori need not be geodesic. The following definition is more relevant in that case.

\begin{defn}{\em (Gromov hyperbolicity)}
Suppose $X$ is any metric space, not necessarily geodesic and $\delta \geq 0$.

(1)  Let $p\in X$. We say that the Gromov inner product on $X$ with respect to $p$, i.e. 
the map $X\times X\map \RR$ defined by $(x,y)\mapsto (x.y)_p$,
is $\delta$-{\em hyperbolic} if $$(x.y)_p\geq \,\mbox{min}\{(x.z)_p, (y.z)_p\}-\delta$$
for all $x,y,z\in X$.

$(2)$ The metric space $X$ is called $\delta$-{\em hyperbolic} in the sense of Gromov if
the Gromov inner product on $X$ is $\delta$-hyperbolic with respect to any point of $X$.

A metric space is called {\em (Gromov) hyperbolic} if it is $\delta$-hyperbolic 
in the sense of Gromov for some $\delta\geq 0$.
\end{defn}



 However, it is a standard fact that for geodesic metric spaces the two concepts are equivalent. See 
\cite[Section 6.3C]{gromov-hypgps}, or \cite[Proposition 1.22, Chapter III.H]{bridson-haefliger} for instance.
In this subsection we observe an analog of Rips hyperbolicity in Gromov hyperbolic length spaces using the
next two lemmas.

The following lemma is a crucial property of Rips hyperbolic metric spaces.

\begin{lemma}\label{stab-qg}  
{\em (Stability of quasigeodesics in a Rips hyperbolic space, \cite{GhH})} 
For all $\delta\geq 0$ and $k\geq 1$, $\epsilon\geq 0$ there
is a constant $D_{\ref{stab-qg}}=D_{\ref{stab-qg}}(\delta, k, \epsilon)$ such that
the following holds:

Suppose $Y$ is a geodesic metric space $\delta$-hyperbolic in the sense of Rips. Then
the Hausdorff distance between a geodesic and a $(k, \epsilon)$-quasigeodesic 
joining the same pair of end points is less than or equal to $D_{\ref{stab-qg}}$.
\end{lemma}

One is referred to \cite[Theorem 3.18, Theorem 3.20]{vaisala} for a proof of the following lemma. 
\begin{lemma}\label{qi vs gromov hyp}
Suppose $X$ is a metric space which is $\delta$-hyperbolic in the sense of Gromov. If $f:X\map Y$ is a $R$-coarsely surjective,
$(1,C)$-quasiisometry then $Y$ is $D=D_{\ref{qi vs gromov hyp}}(\delta, R,C)$-hyperbolic in the sense of Gromov.
\end{lemma}

Using metric graph approximations to length spaces (Lemma \ref{length space qi to graph}) and the fact that for geodesic
metric spaces Gromov hyperbolicity implies Rips hyperbolicity we obtain the following three corollaries.

\begin{cor}{\em (Stability of quasigeodesics in a Gromov hyperbolic space)} \label{cor: stab-qg}
Given $\delta\geq 0, k\geq 1, \epsilon\geq 0$ there is $D=D_{\ref{cor: stab-qg}}(\delta, k, \epsilon)$
such that the following holds.

\noindent
Suppose $X$ is metric space which is $\delta$-hyperbolic in the sense of Gromov.
Then given $(k, \epsilon)$-quasigeodesics $\gamma_i$, $i=1,2$ with the same end points
we have $Hd(\gamma_1, \gamma_2)\leq D$.
\end{cor}

\begin{cor}\label{slim iff gromov}{\em (Analog of Rips hyperbolicity for length spaces)}
Suppose $X$ is a length space. If $X$ is $\delta$-hyperbolic in the sense of Gromov then for all $K\geq 1$, $\epsilon\geq 0$
all $(K,\epsilon)$-quasigeodesic triangles in $X$ are $D_{\ref{slim iff gromov}}=D_{\ref{slim iff gromov}}(\delta,K,\epsilon)$-slim.

Conversely if all $(K,\epsilon)$-quasigeodesic triangles in $X$ are $R$-slim for some $R\geq 0$ and for some sufficiently
large $K, \epsilon$ then $X$ is 
$\lambda_{\ref{slim iff gromov}}=\lambda_{\ref{slim iff gromov}}(R,K,\epsilon)$-hyperbolic in the sense of Gromov.
\end{cor}

Slimness of triangles immediately implies slimness of polygons:
\begin{cor}\label{slim finite polygons} {\em (Slimness of polygons)}
Suppose that $X$ is a length space. If $X$ is $\delta$-hyperbolic in the sense of Gromov then for all $K\geq 1$, $\epsilon\geq 0$
all $(K,\epsilon)$-quasigeodesic $n$-gons in $X$ are 
$(n-2)D_{\ref{slim iff gromov}}=(n-2)D_{\ref{slim iff gromov}}(\delta,K,\epsilon)$-slim.
\end{cor}

\begin{convention}
For the rest of the paper a $\delta$-hyperbolic (or simply hyperbolic) space will refer either to
 (1) a $\delta$-hyperbolic (resp. hyperbolic) space in the sense of Rips if it is a geodesic metric
space or (2) a $\delta$-hyperbolic (resp. hyperbolic) space in the sense of Gromov if it is not a 
geodesic metric space. However, in this case the space will be assumed to be a length space.
The constant $\delta$ will be referred to as the {\em hyperbolicity constant} for the space involved.

\end{convention}


\subsection{Quasiconvex subspaces of hyperbolic spaces}

\begin{defn}
Let $X$ be a hyperbolic geodesic metric space and let $A\subseteq X$. For $K\geq 0$,
we say that $A$ is $K$-{\em quasiconvex} in $X$ if any geodesic with end points
in $A$ is contained in $N_K(A)$. \\
If $X$ is a Gromov hyperbolic length space and $A\subset X$ then we will say that $A$ is
$K$-{\em quasiconvex} if any $(1,1)$-quasigeodesic joining
a pair of points of $A$ is contained in $N_K(A)$.\\
A subset $A\subset X$
is said to be {\em quasiconvex} if it is $K$-quasiconvex for some $K\geq 0$.
\end{defn}

The following lemma relates quasiconvexity with qi embedding. It is straightforward  and is proved in the context of geodesic metric spaces in
\cite[Chapter 1, section 1.11]{ps-kap}. Hence we skip the proof.
\begin{lemma}\label{qc vs qi emb}
(1) Given $\delta\geq 0$ and $k\geq 0$ there are constants $D=D(\delta, k)$ and $K=K(\delta, k)$  such that the following holds:

Suppose $X$ is a $\delta$-hyperbolic metric space and $A\subset X$ is $k$-quasiconvex.
Then $N_D(A)$ is path connected and with respect to the induced path metric on $N_D(A)$ from $X$ 
the inclusion map $N_D(A)\map X$ is a $K$-qi embedding.

(2) Suppose $X$ is a hyperbolic metric space and $Y$ is a quasiconvex subset. Suppose $Y$ is path connected and with respect to the
induced path metric on $Y$ from $X$ the inclusion map $Y\map X$ is metrically proper. Then the inclusion map is a qi embedding.
\end{lemma}

In this subsection, in a Gromov hyperbolic setting, we prove a number of results about quasiconvex sets analogous to those 
in \cite[Section 1.2]{pranab-mahan} which were proved in a Rips hyperbolic setting. The importance of the following 
lemma for this paper can be hardly exaggerated.

\begin{lemma}\label{subqc-elem} {\em (Projection on a quasiconvex set)}
Let $X$ be a $\delta$-hyperbolic metric space, $U\subset X$ is a $K$-quasiconvex set and $\epsilon \geq 0$.
Suppose $y\in U$ is an $\epsilon$-approximate nearest point projection of a point $x\in X$ on $U$. 
Let $z\in U$. Suppose $\alpha$ is a (dotted) $k$-quasigeodesic joining $x$ to $y$ and $\beta$ is a (dotted) $k$-quasigeodesic
joining $y$ to $z$. Then $\alpha*\beta$ is a (dotted) 
$K_{\ref{subqc-elem}}=K_{\ref{subqc-elem}}(\delta, K,k,\epsilon)$-quasigeodesic in $X$. 

In particular, if $\gamma$ is $k$-quasigeodesic joining $x,z$ then $y$ is contained
in the $D_{\ref{subqc-elem}}(\delta,K,k,\epsilon)$-neighborhood of $\gamma$.
\end{lemma}

\proof Without loss of generality we shall assume that $X$ is a $\delta$-hyperbolic length space.
Suppose $\beta_1$ is a $(1,1)$-quasigeodesic in $X$ joining $y,z$. Since $U$ is $K$-quasiconvex it is
clear that $y$ is an $(\epsilon+K)$-approximate nearest point projection of $x$ on $\beta_1$. Hence, if
$\alpha_1$ is a $(1,1)$-quasigeodesic joining $x,y$ then $\alpha_1*\beta_1$ is a $(3,3+\epsilon+K)$-quasigeodesic
in $X$ by Corollary \ref{projection on geodesic}. By stability of quasigeodesics
$Hd(\alpha, \alpha_1)\leq D_{\ref{cor: stab-qg}}(\delta, k, \epsilon)$, and
$Hd(\beta, \beta_1)\leq D_{\ref{cor: stab-qg}}(\delta, k, \epsilon)$.
Hence, $Hd(\alpha*\beta, \alpha_1*\beta_1)\leq D_{\ref{cor: stab-qg}}(\delta, k, \epsilon)$.
By Lemma \ref{quasigeod criteria} it is enough to show now that $\gamma=\alpha*\beta$ is uniformly properly
embedded. 
Let 
$\gamma_1=\alpha_1*\beta_1$ and $R=D_{\ref{cor: stab-qg}}(\delta, k, \epsilon)$.
Suppose $\alpha:[0,l]\map X$ with $\alpha(0)=x, \alpha(l)=y$ and $\beta:[0,m]\map X$
with $\beta(0)=y, \beta(m)=z$. Let $s\leq t\in [0,l+m] $ and $d(\gamma(s), \gamma(t))\leq D$ for some $D\geq 0$.
We need find a constant $D_1$ such that $t-s\leq D_1$  where $D_1$ depends on $\delta, k, K$ and $D$ only. 
However, if $s,t\in [0,l]$ or $s,t\in [l,l+m]$ then we have $-k+(t-s)/k\leq D$ since both $\alpha, \beta$ are 
$k$-quasigeodesics. Hence, in that case $t-s\leq k^2+kD$. Suppose $s\in [0,l)$ and $t\in (l,m]$. In this
case $\gamma(s)=\alpha(s), \gamma(t)=\beta(t-l)$. Let $x'\in \alpha_1, y'\in\beta_1$ be such that 
$d(x', \gamma(s))\leq R$ and $d(y',\gamma(t))\leq R$. Then $d(x',y')\leq 2R+D$. 
Suppose $\gamma_1(s')=x', \gamma_1(t')=y', \gamma_1(u)=y$
where $s'\leq u\leq t'$. Since $\gamma_1$ is a $(3,3+\epsilon+K)$-quasigeodesic
we have $|s'-t'|\leq 3(3+\epsilon+K)+3d(x',y')\leq 3(3+\epsilon+K)+3(2R+D)$. 
It follows that $|s'-u|$ and $|u-t'|$ are 
both at most $3(3+\epsilon+K)+3(2R+D)=9+3\epsilon+3K+6R+3D$. Hence, $d(x',y), d(y,y')$ are both at most
$3(9+3\epsilon+3K+6R+3D)+3+\epsilon+K=30+10\epsilon+10K+18R+9D=D'$, say. 
Hence, $d(\gamma(s),y)$, $d(y,\gamma(t))$ are both at most $R+D'$.
Since $\alpha, \beta$
are $k$-quasigeodesics it follows that $l-s$ and $t-l$ are both at most
$k^2+k(R+D')$. Hence, $t-s\leq 2(k^2+k(R+D'))$. Hence, we can take $D_1=2k^2+2kR+2kD'$. This completes the proof
of the existence of $K_{\ref{subqc-elem}}$.

Clearly one can set
$D_{\ref{subqc-elem}}(\delta,K,k,\epsilon)=D_{\ref{cor: stab-qg}}(\delta, K_{\ref{subqc-elem}}(\delta,K,k,k),K_{\ref{subqc-elem}}(\delta,K,k,k))$.
\qed 

\begin{cor}\label{gluing quasigeodesics}
Suppose $X$ is a $\delta$-hyperbolic metric space and $\alpha$ is a $k$-quasigeodesic in $X$ with an end point $y$.
Suppose $x\in X$ and $y$ is an $\epsilon$-approximate nearest point projection of $x$ on $\alpha$. Suppose
$\beta$ is a $k$-quasigeodesic joining $x$ to $y$. Then $\beta*\alpha$ is a 
$K_{\ref{gluing quasigeodesics}}(\delta, k,\epsilon)$-quasigeodesic.
\end{cor}
\proof We briefly indicate the proof. One first notes by stability of quasigeodesics that images of uniform quasigeodesics
are uniformly quasiconvex. Then one applies the preceding lemma.
\qed

The following corollary easily follows from Lemma \ref{subqc-elem} and
Lemma \ref{lem: approx proj}. For instance, the proof is similar to that of \cite[Lemma 1.32]{pranab-mahan}.

\begin{cor}\label{nested qc sets}{\em (Projection on nested quasiconvex sets)}
Suppose $X$ is a $\delta$-hyperbolic metric space and $V\subset U$ are two $K$-quasiconvex subsets of $X$.
Suppose $x\in X$ and $x_1\in U$, $x_2\in V$ are $\epsilon$-approximate nearest point projection of
$x$ on $U$ and $V$ respectively. Suppose $x_3$ is an $\epsilon$-approximate nearest point projection of
$x_1$ on $V$. Then $d(x_2,x_3)\leq D_{\ref{nested qc sets}}(\delta, K, \epsilon)$.

In particular, for any two $\epsilon$-approximate nearest point projections $x_1,x_2$ of $x$ on $U$
we have $d(x_1,x_2)\leq D_{\ref{nested qc sets}}(\delta, K, \epsilon)$.
\end{cor}

\begin{cor}\label{cor: lip proj}
Given $\delta\geq0, K\geq 0, \epsilon\geq 0$ there are constants $L=L_{\ref{cor: lip proj}}(\delta, K,\epsilon)$,
$D=D_{\ref{cor: lip proj}}(\delta, K,\epsilon)$ and $R=R_{\ref{cor: lip proj}}(\delta, K,\epsilon)$ such that the following hold:

(1) Suppose $X$ is a $\delta$-hyperbolic metric space and $U$ is a $K$-quasiconvex subset of $X$.
Then for all $\epsilon\geq 0$ any $\epsilon$-approximate nearest point projection map $P:X\map U$
is coarsely $L$-Lipschitz.

(2) Suppose $V$ is another $K$-quasiconvex subset of $X$ and $v_1,v_2\in V$ and $u_i=P(v_i)$, $i=1,2$.
If $d(u_1,u_2)\geq D$ then $u_1, u_2\in N_R(V)$. 

In particular, if the diameter of $P(V)$ is at least $D$ then $d(U,V)\leq R$.
\end{cor}

\proof (1) Suppose $x,y\in X$ with $d(x,y)\leq 1$. Then $P(x)$ is an
$(\epsilon+1)$-approximate nearest point projection of $y$ on $U$. Hence, by Corollary
\ref{nested qc sets} we have $d(P(x), P(y))\leq D_{\ref{nested qc sets}}(\delta, K, \epsilon+1)$.
Hence, we may take $L_{\ref{cor: lip proj}}(\delta, K,\epsilon)=D_{\ref{nested qc sets}}(\delta, K, \epsilon+1)$
by Lemma \ref{proving Lipschitz}.

(2) Consider the quadrilateral formed by $(1,1)$-quasigeodesics
joining the pairs $(u_1, u_2), (u_2,v_2), (v_2,v_1)$ and $(v_1,u_1)$. This is 
$2D_{\ref{slim iff gromov}}(\delta,1,1)$-slim by Corollary \ref{slim finite polygons}. Let
$\delta'=2D_{\ref{slim iff gromov}}(\delta,1,1)$. Suppose no point of the side $v_1v_2$ is contained
in a $\delta'$-neighborhood of the side $u_1u_2$. Then there are two points say $x_1, x_2\in v_1v_2$
such that $x_i\in N_{\delta'}(u_iv_i)$, $i=1,2$ and $d(x_1, x_2)\leq 2$. Hence there are points
$y_i\in u_iv_i$, $i=1,2$ such that $d(y_1, y_2)\leq 2+2\delta'$. However, $u_i$ is an
$(\epsilon+3)$-approximate nearest point projection of $y_i$ on $U$ by Lemma \ref{lem: approx proj}.
Hence, by the first part of the Corollary \ref{cor: lip proj} we have 
$d(u_1, u_2)\leq L_{\ref{cor: lip proj}}(\delta, K,\epsilon+3)+(2+2\delta')L_{\ref{cor: lip proj}}(\delta, K,\epsilon+3)$.
Hence, if the diameter of $P(V)$ is bigger than
$D= L_{\ref{cor: lip proj}}(\delta, K,\epsilon+3)+(2+2\delta')L_{\ref{cor: lip proj}}(\delta, K,\epsilon+3)$
then there is a point $x\in v_1v_2$ and $y\in u_1u_2$ such that $d(x,y)\leq \delta'$. Since
$U$ is $K$-quasiconvex we have thus $x\in N_{K+\delta'}(U)$. Thus we may choose
$R=K+\delta'$. \qed

The second part of the above corollary is implied in Lemma 1.35 of \cite{pranab-mahan} too.
 The next lemma roughly says that the
nearest point projection of a quasigeodesic on a quasiconvex set is close to a quasigeodesic.

\begin{lemma}\label{trivial lemma}
Given $K\geq 0, R\geq 0, \delta\geq 0$ there is a constant $D=D_{\ref{trivial lemma}}(R,K,\delta)$ such that the
following holds:\\
Suppose $X$ is a $\delta$-hyperbolic metric space and $A$ is a $K$-quasiconvex subset of $X$.
Suppose $x,y\in X$ and $\bar{x}, \bar{y}\in A$ respectively are their $1$-approximate nearest point projections on 
$A$. Let $[x,y], [\bar{x},\bar{y}]$ denote $1$-quasigeodesics in $X$ joining $x,y$ and $\bar{x},\bar{y}$ respectively.
Suppose $z\in [x,y]$ and $\bar{z}$ is a $1$-approximate nearest point projection of $z$ on $A$ and $d(z,\bar{z})\leq R$. 
Then $d(z, [\bar{x},\bar{y}])\leq D$.
\end{lemma}

\proof 
By Corollary \ref{slim finite polygons} quadrilaterals in $X$ formed by $1$-quasigeodesics are
$2D_{\ref{slim iff gromov}}(\delta,1,1)$-slim. Hence, there is 
$z'\in [x,\bar{x}]\cup [\bar{x}, \bar{y}]\cup [y,\bar{y}]$ such that $d(z,z')\leq 2D_{\ref{slim iff gromov}}(\delta,1,1)$.
If $z'\in [\bar{x},\bar{y}]$ then we are done. Suppose not. Without loss of generality
let us assume that $z'\in [x,\bar{x}]$. Then $d(z',A)\leq d(z,z')+d(z,A)\leq 2D_{\ref{slim iff gromov}}(\delta,1,1)+R$.
Since $\bar{x}$ is a $1$-approximate nearest point projection of $x$ on $A$, $\bar{x}$ is a $4$-approximate nearest point
projection of $z'$ on $A$ by Lemma \ref{lem: approx proj}. Hence, by Corollary \ref{cor: lip proj}, 
$d(\bar{x},\bar{z})\leq L_{\ref{cor: lip proj}}(\delta, K,4)d(z', z)\leq L_{\ref{cor: lip proj}}(\delta, K,4)(2D_{\ref{slim iff gromov}}(\delta,1,1)+R)$. But $d(z,\bar{z})\leq R$. Hence, $d(z, \bar{x})\leq R+ L_{\ref{cor: lip proj}}(\delta, K,4)(2D_{\ref{slim iff gromov}}(\delta,1,1)+R)$. Thus we can take 
$D_{\ref{trivial lemma}}(R,K,\delta)=\max\{2D_{\ref{slim iff gromov}}(\delta,1,1),R+ L_{\ref{cor: lip proj}}(\delta, K,4)(2D_{\ref{slim iff gromov}}(\delta,1,1)+R)\}$. \qed

The following lemma asserts that quasiconvexity and nearest point projections are preserved under qi embeddings.
\begin{lemma}\label{qc in subspace}
Suppose $X$ is a $\delta$-hyperbolic metric graph and $Y\subset X$ is a connected sub-graph such that
the inclusion $(Y,d_Y)\map (X,d_X)$ is a $k$-qi embedding. Suppose $A\subset Y$ is $K$-quasiconvex in $Y$.
Then the following holds.

(1) $A$ is $K_{\ref{qc in subspace}}(\delta,k,K)$-quasiconvex in $X$.

(2) For any $x\in Y$ if $x_1,x_2\in A$ are the nearest point projections of $x$ on $A$
in $Y$ and $X$ respectively then $d_Y(x_1,x_2)\leq D_{\ref{qc in subspace}}(\delta,k,K)$.
\end{lemma}
\proof (1) Suppose $x,y\in A$ and let $\alpha , \beta$ be geodesics joining $x,y$
in $Y$ and $X$ respectively. Since, $Y$ is $k$-qi embedded $\alpha$ is a $(k, k)$-quasigeodesic
in $X$ by Lemma \ref{qi composition}. Hence, by stability of quasigeodesics 
$Hd(\alpha,\beta)\leq D_{\ref{stab-qg}}(\delta, k, k)$. However, $A$ being $K$-quasiconvex in $Y$,
$\alpha\subset N_K(A)$ in $Y$ and hence in $X$ as well. Thus $\beta\subset N_{K+D_{\ref{stab-qg}}(\delta, k,k)}(A)$
in $X$. Hence, we can take $K_{\ref{qc in subspace}}(\delta,k,K)=K+D_{\ref{stab-qg}}(\delta, k, k)$.

(2) Suppose $K_1=K_{\ref{qc in subspace}}(\delta,k,K)$. Then $x_2\in N_D([x,x_1]_X)$ in $X$ where
$D=D_{\ref{subqc-elem}}(\delta, K_1,1,1)$. We have $Hd([x,x_1]_Y, [x,x_1]_X)\leq D_{\ref{stab-qg}}(\delta,k,k)$
by stability of quasigeodesics. Thus there is a point $x'_2\in [x,x_1]_Y$ such that
$d_X(x_2,x'_2)\leq D+D_{\ref{stab-qg}}(\delta,k,k)=D_1$, say. Then $d_Y(x_2,x'_2)\leq k(D_1+k)$ since
$Y$ is $k$-qi embedded in $X$. Since $x_1$ is a nearest point projection
of $x$ on $A$ in $Y$, it is also a nearest point projection of $x'_2$ on $A$ in $Y$. Hence,
$d_Y(x'_2,x_1)\leq d_Y(x'_2,x_2)\leq k(D_1+k)$. Hence, $d_Y(x_1,x_2)\leq 2k(D_1+k)$ by triangle inequality.
Thus we can take $D_{\ref{qc in subspace}}(\delta,k,K)=2k(D_1+k)$.
\qed

\begin{defn} Suppose $X$ is a $\delta$-hyperbolic metric space
and $A, B$ are two quasiconvex subsets. Let  $R>0$. 
We say that $A, B$ are {\em mutually $R$-cobounded, or simply  $R$-cobounded,} if the set of all $1$-approximate nearest point
projections of the points of $A$ on $B$ has a diameter at most $R$ and vice versa.

When the constant $R$ is understood or is not important we just say that $A, B$ are cobounded.
\end{defn}
The following corollary is an immediate consequence of Corollary \ref{cor: lip proj}(2).

\begin{cor}\label{cobounded cor}{\em (\cite[Lemma 1.35]{pranab-mahan})}
Given $\delta\geq 0, k\geq 0$ there are constants $D=D_{\ref{cobounded cor}}(\delta, k)$ and $R=R_{\ref{cobounded cor}}(\delta, k)$
such that the following holds.

Suppose $X$ is a $\delta$-hyperbolic metric space and $A, B\subset X$ are two $k$-quasiconvex subsets.
If $d(A,B)\geq D$ then $A,B$ are mutually $R$-cobounded.
\end{cor}

The following proposition and its proof are motivated by an analogous result due to Hamenstadt (\cite[Lemma 3.5]{hamenst-word}).
See also \cite[Corollary 1.52]{pranab-mahan}. Before we state the proposition let us explain the set-up.

$(P\,0)$ Suppose $X$ is a $\delta$-hyperbolic metric graph and $Y\subset X$ is a $K$-quasiconvex subgraph, for some $\delta \geq 0, K\geq 0$.
Suppose $I$ is an interval in $\RR$ with end points in $\mathbb Z\cup \{\infty, -\infty\}$ and $\Pi:Y\map I$
is a map such that $I\cap \ZZ\subset \Pi(Y)$. Let $Y_i:=\Pi^{-1}(i)$ for all $i\in I\cap \mathbb Z$ and 
$Y_{ij}=\Pi^{-1}([i,j])$ for all $i,j\in I\cap \mathbb Z$ with $i<j$ such that the following hold.

$(P\,1)$ All the sets $Y_i$ and $Y_{ij}$, $i,j\in I$, $i<j$ are $K$-quasiconvex in $X$.

$(P\,2)$ $Y_{i}$ uniformly coarsely bisects $Y$ into $Y^{-}_i:=\Pi^{-1}((-\infty, i]\cap I)$ and $Y^{+}_i:=\Pi^{-1}([i,\infty)\cap I)$
for all $i\in I$. Let $R\geq 0$ be such that any geodesic in $Y$ joining $Y^{+}_i$ and
$Y^{-}_i$ passes through $N_R(Y_i)$ for all $i\in I\cap \ZZ$.

$(P\,3)$ $d(Y_{ii+1}, Y_{jj+1})>2K+1$ for all $i,j\in I$ if $j+1\in I$ and $i+1<j$.

$(P\, 4)$ There is $D\geq 0$ such that the sets $Y_i$ and $Y_j$ are $D$-cobounded in $X$ for all $i,j\in I\cap \ZZ$ with
$i<j$ unless $j=i+1$ and $i,j$ are the end points of $I$.

The proposition below is about a description of uniform quasigeodesics in $X$ joining points of $Y$.
\begin{prop}\label{hamenstadt}
Given $\delta\geq 0, K\geq 0, D\geq 0$, $\lambda\geq 1$, $\epsilon\geq 1$ and $R\geq 0$ there are 
$\lambda'=\lambda_{\ref{hamenstadt}}(\delta, K, D, \lambda, \epsilon, R)\geq 1$ and
$\mu_{\ref{hamenstadt}}= \mu_{\ref{hamenstadt}}(\delta, K, D, \epsilon, R)\geq 0$
such that the following holds.

Suppose we have the aforementioned hypotheses $(P\, 0)$, $(P\, 1)$, $(P\,2)$, $(P\,3)$ and $(P\,4)$.
Suppose $m,n \in I\cap \ZZ$ and  $y\in Y_m, y'\in Y_n$. Suppose $y_i\in Y$, $m\leq i\leq n$ are defined as follows:
$y_m=y$, $y_{i+1}$ is an $\epsilon$-approximate nearest point projection of $y_i$ on $Y_{i+1}$ for $m\leq i\leq n-1$. 
Suppose $\alpha_i\subset Y_{ii+1}$ is a $\lambda$-quasigeodesic in $X$ joining
$y_i$ and $y_{i+1}$, $m\leq i\leq n-1$ and $\beta$ is a $\lambda$-quasigeodesic joining $y_n$ and $y'$.

Then the concatenation of the
all the $\alpha_i$'s and $\beta$ is a $\lambda'$-quasigeodesic in $X$ joining $y, y'$.
Moreover, each $y_i$ is an $\mu_{\ref{hamenstadt}}$-approximate nearest point projection of $y$ on $Y_i$ for $m+2\leq i\leq n$.
\end{prop}
\proof The proof is broken into the following three claims. In course of the proof we shall denote the concatenation of
the $\alpha_i$'s and $\beta$ by $\alpha$.

{\bf Claim 1:} Suppose $x\in Y^{-}_i$ for some $i$. Let $\bar{x}$ be an $\epsilon$-approximate
nearest point projection of $x$ on $Y_i$. Then $\bar{x}$ is an $\epsilon'$-approximate nearest point projection of
$x$ on $Y^{+}_i$ where $\epsilon'$ depends only on $\epsilon$ and the parameters $\delta, D, K$ and $R$.

{\em Proof of Claim 1:}
Suppose $x'$ is a $1$-approximate nearest point projection of $x$ on $Y^{+}_i$. 
Since $Y^{+}_i$ is $K$-quasiconvex $[x,x']*[x',\bar{x}]$ is a $K_{\ref{subqc-elem}}(\delta, K, 1, 1)$-quasigeodesic
by Lemma \ref{subqc-elem}. Let $k_1=K_{\ref{subqc-elem}}(\delta, K, 1, 1)$.
Then by stability of quasigeodesics there is a point $z\in [x,\bar{x}]$ such that
$d(x',z)\leq D_{\ref{stab-qg}}(\delta, k_1)=D_1$, say. We claim that $z$ is uniformly close to $Y_i$.
Since $Y^{-}_i$ is $K$-quasiconvex there is a point $w\in Y^{-}_i$ such that $d(z,w)\leq K$. It follows that
$d(w, x')\leq D_1+K$. By $(P\, 2)$, there is a point $z_1\in [w,x']$ 
such that $d(z_1, Y_i)\leq R$. Since, $d(z_1,w)\leq d(w,x')\leq D_1+K$
and $d(w,z)\leq K$ it follows by triangle inequality that $d(z, Y_i)\leq 2K+D_1+R$.
Now, by Lemma \ref{lem: approx proj} $\bar{x}$ is an $(\epsilon+3)$-approximate nearest point projection of 
$z$ on $Y_i$. Hence, $d(x',\bar{x})\leq d(x',z)+d(z,\bar{x})\leq D_1+\epsilon+3+d(z,Y_i)$. 
It follows that $\epsilon'=3+\epsilon+2K+2D_1+R$ works.

{\em Note: We shall use $D_1$ again in the proof of Claim $3$ to denote the same constant as in the proof of Claim $1$ above.}

{\bf Claim 2.} Next we claim that for all $m+2\leq i\leq n-1$ there is uniformly bounded set $A_i\subset Y_i$ such 
that $\epsilon$-nearest point projection of any point of $Y^{-}_j$, $j<i$ on $Y_i$ is contained in $A_i$.

{\em Proof of Claim 2:} Consider any $Y_i$, $m+2\leq i\leq n-1$.
Let $B_i\subset Y_i$ be the set of all $1$-approximate nearest point projections of points of $Y_{i-1}$ on $Y_i$ in $X$.
Then the diameter of $B_i$ is at most $D$ by $(P\, 4)$. Suppose $x\in Y^{-}_j$, $j<i$. Let $x_1,x_2$ be respectively
$\epsilon$-approximate nearest point projections of $x$ on $Y_{i-1}$ and $Y_i$ respectively. Let
$x_3$ be an $\epsilon$-nearest point projection of $x_1$ on $Y_i$. Now, by Step 1
$x_1$ is an $\epsilon'$-approximate nearest point projection of $x$ on $Y^{+}_{i-1}$ and
$x_2$, $x_3$ are $\epsilon'$-approximate nearest point projection of $x$ and $x_1$ respectively
on $Y^{+}_i$. Therefore, by the first part of Corollary \ref{nested qc sets} we have 
$d(x_2,x_3)\leq D_{\ref{nested qc sets}}(\delta, K, \epsilon')$. 
However, if $x'_1\in B_i$ is a $1$-approximate nearest point projection of $x_1$ on $Y_i$ then by  
the second part of the Corollary \ref{nested qc sets} we have 
$d(x_3, B_i)\leq d(x_3,x'_1)\leq D_{\ref{nested qc sets}}(\delta, K, \epsilon)$ since $\epsilon \geq 1$.
Hence, $d(x_2,B_i)\leq 2D_{\ref{nested qc sets}}(\delta, K, \epsilon)$. Therefore,
we can take $A_i=N_{2D_{\ref{nested qc sets}}(\delta, K, \epsilon)}(B_i)\cap Y_i$.

\smallskip
Let $r= \sup_{m+2\leq i\leq n-1} \{diam(A_i)\}.$
We note that $r\leq D+ 2D_{\ref{nested qc sets}}(\delta, K, \epsilon)$. 

{\bf Claim 3.} Finally we claim that (1) $\alpha$ is contained in a uniformly small neighborhood of a geodesic joining $y,y'$
and (2) $\alpha$ is uniformly properly embedded in $X$. 

We note that the proposition follows from Claim 3 using Lemma \ref{quasigeod criteria}.

{\em Proof of Claim 3:}
Suppose $x,x'\in \alpha$, $\Pi(x)<\Pi(x')$. Choose smallest $k,l$ such that 
$x\in \alpha\cap Y_{kk+1}, x'\in \alpha\cap Y_{ll+1}$, where $m\leq k\leq l\leq n$. 
Let $\gamma$ be a geodesic in $X$ joining $x,x'$.

(1) 
It is enough to show that the segment of $\alpha$ joining $x$ to $x'$ is contained in a uniformly small
neighborhood of $\gamma$. Hence, without loss of generality $k<l$.
Due to Corollary \ref{slim finite polygons} it is enough to prove that the points 
$y_i$, $k+1\leq i\leq l-1$ are contained in a uniformly small neighborhood of $\gamma$
in order to show that the segment of $\alpha$ joining $x$ to $x'$ is contained in a uniformly small
neighborhood of $\gamma$. (We note that the path $\alpha_{n-1}*\beta$ is a 
$D_{\ref{subqc-elem}}(\delta, K, \lambda, \epsilon)$-quasigeodesic joining $y_{n-1}$ and $y'$.)
For this first we note that $x$ is on $\alpha_k$. 
Let $\gamma_k$ be a geodesic joining $y_k,y_{k+1}$. Then by stability of quasigeodesics 
there is a point $x_1\in \gamma_k$ such that
$d(x_1, x)\leq  D_{\ref{cor: stab-qg}}(\delta, \lambda,\lambda)$. Since $y_{k+1}$ is an $\epsilon$-approximate
nearest point projection of $y_k$ on $Y_{k+1}$, by Lemma \ref{lem: approx proj} $y_{k+1}$
is an $(\epsilon+3)$-approximate nearest point projection of $x_1$ on $Y_{k+1}$.
Hence, $y_{k+1}$ is an $(\epsilon+3+D_{\ref{cor: stab-qg}}(\delta, \lambda,\lambda))$-approximate nearest point
projection of $x_1$ on $Y_{k+1}$. Let $\epsilon_1=\epsilon+3+D_{\ref{cor: stab-qg}}(\delta, \lambda,\lambda)$. 
By Step 1 $y_{k+1}$ is an $\epsilon'_1$-nearest point projection of $x$ on $Y^{+}_{k+1}$
where $\epsilon'_1=3+\epsilon_1+2D_1++2K+R$.
Now the concatenation of a geodesic joining $y_{k+1}$ to $x'$ with the segment of $\alpha$ from $x$ to $y_{k+1}$
is a uniform quasigeodesic  by Lemma \ref{subqc-elem}. Thus by Corollary \ref{cor: stab-qg} $y_{k+1}$ is uniformly
close to $\gamma$. On the other hand by Step 2 $y_i$ is an $(\epsilon+r)$-approximate nearest point projection
of $x$ on $Y_i$ and hence an $(\epsilon+r)'$-approximate nearest point projection on $Y^{+}_i$ for all 
$k+2\leq i \leq l-1$. Hence, again by Lemma \ref{subqc-elem} and Corollary \ref{cor: stab-qg}
$y_i$ is within a uniformly small neighborhood of $\gamma$. This proves (1).


(2)
Suppose $L=\sup \{d(y_i,\gamma): k+1\leq i\leq l-1\}$. Suppose $x,x'\in \alpha$ as above with
$d(x, x')\leq N$.  Once again, without loss of generality $k<l$. We claim that $l\leq k+N$. 
To see this consider two adjacent vertices $v_i, v_{i+1}$ on $\gamma$. If $v_i\in N_K(Y_{ss+1})$ and 
$v_{i+1}\in N_K(Y_{tt+1})$ with $s<t$ then by the $(P\,3)$ we have $t=s+1$. The claim follows from this.
Suppose $\alpha(s_k)=x$, $\alpha(s_i)=y_i$ for $k+1\leq i\leq l-1$ and $\alpha(s_l)=x'$.
We note that $d(\alpha(s_i),\alpha(s_{i+1}))\leq N+2L$ for $k\leq i\leq l-1$.
Since $l-k\leq N$ and since the segments of $\alpha$ joining $\alpha(s_i),\alpha(s_{i+1})$, $k\leq i\leq l-1$ are 
uniform quasigeodesics we are done.

For the second part of the proposition we have already noticed that $y_i$ is an $(\epsilon+r)$-approximate
nearest point projection of any point $Y^{-}_j$, in particular of $y$, on $Y_i$ for all $j<i$, $m+2\leq i\leq n-1$. 
On the other hand, $y_{n-1}$ is an $(\epsilon+r)'=(\epsilon+r+3+2D_1+2K+R)$-approximate nearest point projection
of $y$ on $Y^{+}_{n-1}$. Hence, by Corollary \ref{nested qc sets} if $y'_n$ is a $1$-approximate point projection of
$y$ on $Y_n\subset Y^{+}_{n-1}$ then $d(y'_n, y_n)\leq D_{\ref{nested qc sets}}(\delta, K, (\epsilon+r)')$. Thus $y_n$
is an $(1+D_{\ref{nested qc sets}}(\delta, K, (\epsilon+r)'))$-approximate nearest point projection
of $y$ on $Y_n$.
\qed

\begin{lemma}\label{gromov product meaning}
Given $\delta\geq 0,  k\geq 1, \epsilon\geq 0$ there is a constant $D=D_{\ref{gromov product meaning}}(\delta, k, \epsilon)$ such that
the following is true.

Suppose $X$ is a $\delta$-hyperbolic metric space.
Suppose $x_1,x_2,p\in X$ and $\alpha$ is a $(k,\epsilon)$-quasigeodesic in $X$ joining $x_1,x_2$.
Then $|(x_1.x_2)_p-d(p,\alpha)|\leq D$.
\end{lemma}


\proof Without loss generality we shall assume that $X$ is a length space $\delta$-hyperbolic in the
sense of Gromov. Let $w\in \alpha$ be a $1$-approximate nearest point projection of $p$ on $\alpha$.
Let $\beta_1, \beta_2$ be $(1,1)$-quasigeodesics joining the pairs of points $(x_1,p), (x_2,p)$ respectively.
Let $\gamma$ be a $(1,1)$-quasigeodesic joining $p,w$ and let $\alpha'$ be a $(1,1)$-quasigeodesic joining $x_1,x_2$.
Let $C=D_{\ref{cor: stab-qg}}(\delta, k, \epsilon+1)$. Now, by Corollary \ref{cor: stab-qg} $Hd(\alpha, \alpha')\leq C$
and $\alpha$ is $C$-quasiconvex.
Let $\alpha_1$ be the portion of $\alpha$ from $x_1$ to $w$ and let $\alpha_2$ be the portion of $\alpha$ from
$w$ to $x_2$. Then $\alpha_1*\gamma$, $\alpha_2*\gamma$ are 
$K=K_{\ref{subqc-elem}}(\delta, C,k+\epsilon, k+\epsilon)$-quasigeodesics. Hence by Corollary \ref{cor: stab-qg}
$Hd(\beta_i, \alpha_i*\gamma)\leq D_{\ref{cor: stab-qg}}(\delta, K, K)$. Let $w_i\in \beta_i$
be such that $d(w,w_i)\leq D_{\ref{cor: stab-qg}}(\delta, K, K)$. 
Since $Hd(\alpha, \alpha')\leq C$, there is a point $w'\in \alpha'$ such that
$d(w,w')\leq C$. Hence, $d(w',w_i)\leq C+D_{\ref{cor: stab-qg}}(\delta, K, K)=R$, say.
Now by Lemma \ref{gromov product step1} $|(x_1.x_2)_p-d(p,w')|\leq 3+2R$.
It follows that $|(x_1.x_2)_p-d(p,w)|\leq 3+2R+C$. Since $w$ is a $1$-approximate nearest point projection of $p$
on $\alpha$ we have for all $z\in \alpha$, $d(p,w)\leq d(p,z)+1$. Thus $|d(p, \alpha)-d(p,w)|\leq 1$.
Hence, $|(x_1.x_2)_p-d(p,\alpha)|\leq 4+2R+C$.
\qed


\subsection{Boundaries of hyperbolic spaces and CT maps}

Given a hyperbolic metric space, there are the following three standard ways to define 
a boundary. Some of the results in this subsection are mentioned without proof. 
One may refer to \cite{bridson-haefliger} and \cite{Shortetal} for details.
\begin{defn}
\begin{enumerate}
\item {\bf Geodesic boundary.}
Suppose $X$ is a (geodesic) hyperbolic metric space. Let $\mathcal G$ denote the set of all
geodesic rays in $X$. The {\em geodesic boundary} $\partial X$ of $X$ is defined to be $\mathcal G/\sim$
where $\sim$ is the equivalence relation on $\mathcal G$ defined by setting $\alpha\sim \beta$ 
iff $Hd(\alpha,\beta)< \infty$.

\item {\bf Quasigeodesic boundary.}
Suppose $X$ is a hyperbolic metric space in the sense of Gromov. Let $\mathcal Q$ be the set of all
quasigeodesic rays in $X$. Then the {\em quasigeodesic boundary} $\partial_q X$ is defined to be
$\mathcal Q/\sim$ where $\sim$ is defined as above.

\item {\bf Gromov boundary or sequential boundary.}
Suppose $X$ is a hyperbolic metric space in the sense of Gromov and $p\in X$.
Let $\mathcal S$ be the set of all sequences $\{x_n\}$ in $X$ such that $\lim_{i,j\map \infty}(x_i.x_j)_p=\infty$. 
All such sequences are said to converge to infinity. On $\mathcal S$ we define an equivalence relation where 
$\{x_n\}\sim \{y_n\}$ if and only if $\lim_{i,j\map \infty}(x_i.y_j)_p=\infty$ for some (any) base point $p\in X$.
The {\em Gromov boundary or the sequential boundary} $\partial_s X$ of $X$, as a set, is defined to be $\mathcal S/\sim$.
\end{enumerate}
\end{defn}

{\bf Notation and convention.} (1)
The equivalence class of a geodesic ray or a quasigeodesic ray $\alpha$ in $\partial X$ or
$\partial_q X$ is denoted by $\alpha(\infty)$. It is customary to fix a base point and require that all the
rays start from there to define $\partial X$ and $\partial_q X$ but it is not essential.

(2) If $\alpha$ is a (quasi)geodesic ray with $\alpha(0)=x$, $\alpha(\infty)=\xi$ then we say that $\alpha$ joins $x$
to $\xi$. We use $[x,\xi)$ to denote any (quasi)geodesic ray joining $x$ to $\xi$
when the parametrization of the (quasi)geodesic ray is not important or is understood. 

(3) If $\alpha$ is a quasigeodesic line with $\alpha(\infty)=\xi_1, \alpha(-\infty)=\xi_2\in \partial_q X$
then we say that $\alpha$ joins $\xi_1, \xi_2$. We denote by $(\xi_1, \xi_2)$ any quasigeodesic line
joining $\xi_1, \xi_2$ when the parameters of the quasigeodesic are understood.

(4) If $\xi=[\{x_n\}]\in \partial_s X$ then  we write $x_n\map \xi$ or $\xi=\lim_{n\map \infty} x_n$ and say that the 
sequence $\{x_n\}$ converges to $\xi$.

(5) We shall denote by $\widehat{X}$ the set $X\cup\partial_s X$.

The following lemma and proposition summarizes all the basic properties of the boundary of hyperbolic spaces that
we will need in this paper. 

\begin{lemma}\label{bdry lemma 1}{\em (\cite[Theorem 11.108]{kap-drutu-book})}
Let $X$, $Y$ be hyperbolic metric spaces.\\	
(1) Given a qi embedding $\phi: X\map Y$ we have an injective map $\partial \phi: \partial_s X \map \partial_s Y$.

(2) (i) If $X\stackrel{\phi}{\map} Y\stackrel{\psi}{\map} Z$ are qi embeddings then 
$\partial (\psi\circ \phi)= \partial \psi \circ \partial \phi$

(ii) $\partial (Id_X)$ is the identity map on $\partial_s X$. 

(iii) A qi induces a bijective boundary map.
\end{lemma}

The following proposition relates the three definitions of boundaries.

\begin{prop}\label{visibility}
(1) For any metric space $X$ the inclusion $\mathcal G\map \mathcal Q$ induces an injective map $\partial X\map \partial_q X$.

(2) Given a quasigeodesic ray $\alpha$, $\lim_{n\map \infty} \alpha(n)$ is well defined and $\alpha\sim \beta$ implies
$\lim_{n\map \infty}\alpha(n)=\lim_{n\map \infty}\beta(n)$. This induces an injective map $\partial_q X\map \partial _s X$.

(3) If $X$ is a proper geodesic hyperbolic metric space then the map $\partial X\map \partial_q X$ is a bijection.

(4) The map $\partial_q X\map \partial_s X$ is a bijection for all Gromov hyperbolic length spaces. 

In fact, given $\delta \geq 0$ there is a constant $k_{\ref{visibility}}=k_{\ref{visibility}}(\delta)$ 
such that given any $\delta$-hyperbolic length space $X$, any 
pair of points $x,y\in \widehat{X}$ can be joined by a $k_{\ref{visibility}}$-quasigeodesic.
\end{prop}

\proof (1), (2), (3) are standard. See \cite[Chapter III.H]{bridson-haefliger} for instance.

(4) is proved for geodesic metric spaces in Section 2 of \cite{pranab-mahan}. See Lemma 2.4 there.
The same result for a general length space then is a simple consequence of the existence of a metric graph approximation
of a length space and the preceding lemma.
\qed

\begin{lemma}\label{ideal triangles are slim}{\bf (Ideal triangles are slim)}
Suppose $X$ is a $\delta$-hyperbolic metric space in the sense of Rips or Gromov.
Suppose $x,y,z\in \widehat{X}$ and we have three $k$-quasigeodesics joining each pair of points from $\{x,y,z\}$.
Then the triangle is $R=R_{\ref{ideal triangles are slim}}(\delta, k)$-slim.

In particular, if $\gamma_1, \gamma_2$ are two $k$-quasigeodesic rays with $\gamma_1(0)=\gamma_2(0)$
and $\gamma_1(\infty)=\gamma_2(\infty)$ then $Hd(\gamma_1, \gamma_2)\leq R$. 
\end{lemma}
The proof of the above lemma is pretty standard and hence we omit it. However, slimness of ideal triangles
immediately implies slimness of ideal polygons:

\begin{cor}\label{ideal polygons are slim}{\bf (Ideal polygons are slim)}
Suppose $X$ is a $\delta$-hyperbolic metric space in the sense of Rips or Gromov.
Suppose $x_1, x_2,\ldots, x_n\in \widehat{X}$ are $n$ points and we have $n$ $k$-quasigeodesics joining 
pairs of points $(x_1,x_2),(x_2,x_3),\ldots,(x_{n-1},x_n)$ and $(x_n,x_1)$. 
 Then this $n$-gon
is $R=R_{\ref{ideal polygons are slim}}(\delta, k,n)$-slim, i.e. every side is contained
in $R$-neighborhood of the union of the remaining $n-1$ sides.
\end{cor}

The following lemma gives a geometric interpretation for sequential boundary in terms of quasigeodesics.
\begin{lemma}\label{lem: bdry defn}
Let $x\in X$ be any point.
Suppose $\{x_n\}$ is any sequence of points in $X$ and $\beta_{m,n}$ is a $k$-quasigeodesic
joining $x_m$ to $x_n$ for all $m,n\in \NN$. Suppose $\alpha_n$ is a $k$-quasigeodesic joining $x$ to $x_n$.
Then

(1) $\{x_n\}\in \mathcal S$ if and only if $\lim_{m,n\map\infty} d(x, \beta_{m,n})=\infty$
if and only if there is a constant $D$ such that for all $M>0$ there is $N>0$ with 
$Hd(\alpha_m\cap B(x;M), \alpha_n\cap B(x;M))\leq D$ for all $m,n\geq N$.

(2) Suppose moreover $\xi\in \partial_s X$ and $\gamma_n$ is a $k$-quasigeodesic in $X$ joining
$x_n$ to $\xi$ for all $n\in \NN$ and $\alpha$ is a $k$-quasigeodesic joining $x$ to $\xi$. 

Then $x_n\map \xi$ if and only if $d(x, \gamma_n)\map \infty$ iff there is constant $D>0$ such that for all
$M>0$ there is $N>0$ with $Hd(\alpha\cap B(x;M), \alpha_n\cap B(x;M))\leq D$
for all $n\geq N$.
\end{lemma}
We skip the proof of this lemma. In fact,
the first statement of the lemma is an easy consequence of Lemma \ref{gromov product meaning}
and stability of quasigeodesics. The second statement is a simple consequence of Lemma
\ref{gromov product meaning}, stability of quasigeodesics and the Lemma \ref{ideal triangles are slim}.


The following lemma is proved in section 2 of \cite{pranab-mahan} (see Lemma 2.7 and Lemma 2.9 there)
for hyperbolic geodesic metric spaces. The same statements are true for length spaces too.
To prove it for length spaces one just takes a metric graph approximation. Since the proof is straightforward
we omit it.

\begin{lemma}{\bf (Barycenters of ideal triangles)}\label{defn: bary map}
Given $\delta\geq 0$ there is $r_0\geq 0$ such that for any $\delta$-hyperbolic length space $X$, 
any three distinct points $x,y,z\in \widehat{X}$ and any three $k_{\ref{visibility}}(\delta)$-quasigeodesics 
joining $x,y,z$ in pairs there is a point $x_0\in X$ such that $N_{r_0}(x_0)$ intersects all the three quasigeodesics.

We refer to a point with this property to be a {\em barycenter} of the ideal triangle $\Delta xyz$. There is a constant $L_0$ 
such that if $x_0, x_1$ are two barycenters of $\Delta xyz$ then $d(x_0, x_1)\leq L_0$.
\end{lemma}

Thus we have a coarsely well-defined map $\partial^3_s X\map X$. We shall refer to this map as the
{\bf barycenter map}. It is a standard fact that for a non-elementary hyperbolic group $G$ if $X$ is a Cayley
graph of $G$ then the barycenter map $\partial^3_s X\map X$ is coarsely surjective and vice versa. 
If $X$ is a hyperbolic metric space such that the barycenter map for $X$ is coarsely surjective then $X$ will
be called a {\bf nonelementary} hyperbolic space.
In section $4$ and $5$ we deal with spaces with this property. 

The following lemma is clear. For instance, we can apply the proof of 
\cite[Lemma 2.9]{pranab-mahan}.

\begin{lemma}\label{lem: barycenter map}
Barycenter maps being coarsely surjective is a qi invariant property among hyperbolic length spaces.
\end{lemma}


\subsubsection{Topology on $\partial_s X$ and Cannon-Thurston maps}
\begin{defn}
(1) If $\{\xi_n\}$ is a sequence of points in $\partial_s X$, we say that $\{\xi_n\}$ converges to $\xi\in \partial_s X$
if the following holds: Suppose $\xi_n=[\{x^n_k\}_k]$ and $\xi=[\{x_k\}]$. Then 
$\lim_{n\map \infty}(\liminf_{i,j\map \infty}(x_i.x^n_j)_p)=\infty$.


(2) A subset $A\subset \partial_s X$ is said to be closed if for any sequence $\{\xi_n\}$  in $A$, $\xi_n\map \xi$ implies
$\xi\in A$.
\end{defn} 
The definition of convergence that we have stated here is equivalent to the one stated in \cite{Shortetal}.
Moreover, that the convergence mentioned above is well-defined follows from \cite{Shortetal} and hence
we skip it. The next two lemmas give a geometric meaning of the convergence.
\smallskip

\begin{lemma}\label{lem: lim defn}
Given $k\geq 1$ and $\delta\geq 0$ there are constants 
$D=D_{\ref{lem: lim defn}}(k,\delta)$, $L=L_{\ref{lem: lim defn}}(k,\delta)$ and 
$r=r_{\ref{lem: lim defn}}(k,\delta)$ with the following properties:

Suppose $\alpha, \beta$ are two $k$-quasigeodesic rays starting from a point $x\in X$ such that 
$\alpha(\infty)\neq \beta(\infty)$ and $\gamma$ is a $k$-quasigeodesic line joining $\alpha(\infty)$ and $\beta(\infty)$.
Then the following hold:

(1) There exists $N\in \NN$ such that $|(\alpha(m).\beta(n))_x-d(x,\gamma)|\leq D$ for
all $m,n\geq N$.

In particular, $|\liminf_{m,n\map \infty} (\alpha(m).\beta(n))_x-d(x,\gamma)|\leq D$.

(2) Suppose $R=d(x,\gamma)$ then $Hd(\alpha\cap B(x;R-r), \beta\cap B(x;R-r))\leq L$. 
\end{lemma}
\proof (1) Since $\alpha(\infty)\neq \beta(\infty)$ by Lemma \ref{ideal triangles are slim} there is $N\in \NN$ such that
for all $m,n\geq N$, $\alpha(m)\in N_{R_{\ref{ideal triangles are slim}}}(\gamma)$ and
$\beta(n)\in N_{R_{\ref{ideal triangles are slim}}}(\gamma)$. Let $x_m, y_n\in \gamma$ be such that
$d(x_m, \alpha(m))\leq R_{\ref{ideal triangles are slim}}$ and $d(y_n, \beta(n))\leq R_{\ref{ideal triangles are slim}}$.
Then by joining $x_m, \alpha(m)$ and $y_n, \beta(n)$ and applying Corollary \ref{slim finite polygons}
we see that Hausdorff distance between any $(1,1)$-quasigeodesic joining $\alpha(m), \beta(n)$, say $c_{m,n}$ and the portion
of $\gamma$ between $x_m, y_n$ is at most 
$R_{\ref{ideal triangles are slim}}+2D_{\ref{slim iff gromov}}(\delta,k,k)$.
It is clear that for large enough $N$, $d(x, \gamma)$ is the same as the distance of $x$ and the segment of $\gamma$
between $x_m ,y_n$ if $m,n\geq N$. Thus for such $m,n$ we have
$|d(x, c_{m,n})-d(x,\gamma)|\leq R_{\ref{ideal triangles are slim}}+2D_{\ref{slim iff gromov}}(\delta,k,k)$.
But by Lemma \ref{gromov product meaning},
$|(\alpha(m).\beta(n))_x-d(x, c_{m,n}|\leq D_{\ref{gromov product meaning}}(\delta, k, k)$.
Hence, $|(\alpha(m).\beta(n))_x-d(x,\gamma)|\leq R_{\ref{ideal triangles are slim}}+2D_{\ref{slim iff gromov}}(\delta,k,k)+
D_{\ref{gromov product meaning}}(\delta, k, k)$ for all large $m,n$. 

(2) To see this we take a $1$-approximate nearest point projection, say $z$, of $x$ on $\gamma$. 
 Let $xz$ denote a $1$-quasigeodesic joining $x, z$. Then by Corollary 
\ref{gluing quasigeodesics} concatenation of $xz$ and the portions of $\gamma$ joining $z$ to $\gamma(\pm \infty)$ respectively are both
$K_{\ref{gluing quasigeodesics}}(\delta, k,k)$-quasigeodesics. Call them $\alpha'$ and $\beta'$ respectively.
Note that $\alpha(\infty)=\alpha'(\infty)$ and $\beta(\infty)=\beta'(\infty)$. Let 
$K=\max\{k, K_{\ref{gluing quasigeodesics}}(\delta, k,\epsilon)\}$.
Then by the last part of Lemma \ref{ideal triangles are slim} 
it follows that $z\in N_{r}(\alpha)\cap N_{r}(\beta)$ where $r=R_{\ref{ideal triangles are slim} }(\delta, K)$.
Suppose $x'\in \alpha, y'\in \beta$ are such that $d(z,x')\leq r$ and $d(y', z)\leq r$.
By Lemma \ref{slim iff gromov} the Hausdorff distance between $xz$ and the portions of $\alpha$ from $x$ to $x'$
and the portion of $\beta$ from $x$ to $y'$ are each at most $D_{\ref{slim iff gromov}}(\delta, k,k)+r$.
Thus these segments of $\alpha$ and $\beta$ are at a Hausdorff distance at most $L=2D_{\ref{slim iff gromov}}(\delta, k,k)+2r$
from each other. This completes the proof. \qed

\begin{lemma}\label{convergence explained}
Let $x\in X$ be any point. Suppose $\{\xi_n\}$ is any sequence of points in $\partial_s X$. 
Suppose $\beta_{m,n}$ is a $k$-quasigeodesic line joining $\xi_m$ to $\xi_n$ for all $m, n\in \NN$ 
and $\alpha_n$ is a $k$-quasigeodesic ray joining $x$ to $\xi_n$
for all $n\in \NN$. Then

(1) $\lim_{m,n\map \infty} d(x,\beta_{m,n})= \infty$ iff there is a constant $D=D(k,\delta)$ such that for all $M>0$ 
there is $N>0$ with $Hd(\alpha_m\cap B(x;M), \alpha_n\cap B(x;M))\leq D$ for all $m,n\geq N$ and in this case
$\{\xi_n\}$ converges to some point of $\partial_s X$.

(2) Suppose moreover $\xi\in \partial_s X$, $\gamma_n$ is a $k$-quasigeodesic ray in $X$ joining $\xi_n$ to $\xi$
for all $n$, and $\alpha$ is a $k$-quasigeodesic ray joining $x$ to $\xi$. Then 
$\xi_n\map \xi$ iff $d(x, \gamma_n)\map \infty$ iff there is constant $D'=D'(k,\delta)$ such that for all
$M>0$ there is $N>0$ with $Hd(\alpha\cap B(x;M), \alpha_n\cap B(x;M))\leq D$
for all $n\geq N$. In this case $\lim_{m,n\map \infty} d(x,\beta_{m,n})= \infty$.
\end{lemma}
\proof (1) 
The `iff' part is an immediate consequence of Lemma \ref{lem: lim defn}. We prove the last part.
Let $n_i$ be an increasing sequence in $\NN$ such that
for all $m,n\geq n_i$ we have $Hd(\alpha_m\cap B(x;i), \alpha_n\cap B(x;i))\leq D$.
Let $y_i$ be a point of $\alpha_{n_i}\cap B(x;i)$ such that $d(x, y_i)+1\geq sup\{d(x,y):x\in \alpha_{n_i}\cap B(x;i)\}$.
We claim that $y_i$ converges to a point of $\partial_s X$. Clearly $d(x,y_i)\map \infty$.
Given $i\leq j\in \NN$ we have $d(y_i,\alpha_n)\leq D$ and $d(y_j,\alpha_n)\leq D$ for all $n\geq n_j$.
By slimness of polygons we see that any $(1,1)$-quasigeodesic joining $y_i,y_j$ is uniformly close to $\alpha_n$.
It follows that $\lim_{i,j\map \infty} (y_i.y_j)_x=\infty$. Let $\xi=[\{y_n\}]$. It is clear that
$\xi_n\map \xi$.

(2) Both iff statements are immediate from Lemma \ref{lem: lim defn}. The last part follows from slimness 
of ideal triangle since $d(x, \gamma_n)\map \infty$.
\qed

\begin{cor}\label{cor: convergence explained}
Suppose $\{x_n\}$ is a sequence of points in $\widehat{X}$ such that $\{x_n\}\subset X$ or
$\{x_n\}\subset \partial_s X$. Suppose $x_n\map \xi\in \partial_s X$ and $\gamma_n$ is a $k$-quasigeodesic 
joining $x_n$ to $\xi$ for each $n$. Let $y_n\in \gamma_n$ such that $d(x,y_n)\map \infty$.
Then $\lim_{n\map \infty} y_n=\xi$.
\end{cor}

\begin{defn}{\em({\bf Cannon-Thurston map}, \cite{mitra-trees})} \label{ct-defn}
If $f:Y\map X$ is any map of hyperbolic metric spaces then we say that {\em the Cannon-Thurston} or {\bf the CT map} exists
for $f$ or that $f$ admits the CT map if $f$ gives rise to a continuous map $\partial{f}:\partial{Y}_s\map \partial{X}_s$ in the following sense:

Given any $\xi\in \partial_s Y$ and any sequence of points $\{y_n\}$ in $Y$ converging to $\xi$, the sequence $\{f(y_n)\}$ 
converges to a definite point of $\partial_s X$ independent of the $\{y_n\}$ and the resulting map $\partial{f}:\partial_s{Y}\map \partial_s{X}$ is continuous. 
\end{defn}
Generally, one assumes that the map $f$ is a proper embedding but for the sake of the definition it is unnecessary.
We note that the CT map is unique when it exists. The following lemma gives a sufficient condition for the existence
of CT maps.


\begin{lemma}{\em({\bf Mitra's criterion}, \cite[Lemma 2.1]{mitra-trees})}\label{CT existence}
Suppose $X$, $Y$ are geodesic hyperbolic metric spaces and $f:Y\map X$ is a metrically proper map. Then $f$ admits 
the CT map if the following holds:

$(\ast)$ Let $y_0\in Y$. There exists a function $\tau:\RR_{\geq 0}\map \RR_{\geq 0}$, with the property that $\tau(n)\map \infty$ as 
$n\map \infty$ such that for all geodesic segments $[y_1,y_2]_Y$ in $Y$ lying outside the $n$-ball around $y_0\in Y$, 
any geodesic segment $[f(y_1),f(y_2)]_X$ in $X$ joining the pair of points $f(y_1),f(y_2)$ lies outside the $\tau(n)$-ball 
around $f(y_0)\in X$.
\end{lemma}

\begin{rem}\label{rem: CT criteria}
{\em 
(1) The main set of examples where Lemma \ref{CT existence} applies comes from taking $Y$ to be a rectifiably path connected subspace of
a hyperbolic space $X$ with induced length metric and the map $f$ is assumed to be the inclusion map. One also considers
the orbit map $G\map X$ where $G$ is a hyperbolic group acting properly by isometries on a hyperbolic metric space $X$.
In these examples, the map $f$ is coarsely Lipschitz as well as metrically proper. The proof of the lemma by Mitra also
assumes that $X$, $Y$ are proper geodesic metric spaces and Mitra considered the geodesic boundaries. 
However, these conditions are not necessary as the following lemma and examples show. 

(2) The proof of Lemma \ref{CT existence} by Mitra only checks that the map is a well-defined extension of $f$ rather than it is continuous.
However, with very little effort the condition $(\ast)$ can be shown to be sufficient for the well-definedness as well as 
the continuity of the CT map.

(3) One can easily check that the condition $(\ast)$ is also necessary provided $X,Y$ are proper hyperbolic spaces and $f$ is  
coarsely Lipschitz and metrically proper.} 
\end{rem}

The following lemma is the main tool for the proof of our theorem of Cannon-Thurston map. We shall refer to this
as Mitra's lemma.

\begin{lemma}\label{mitra-lemma}
Suppose $X,Y$ are length spaces hyperbolic in the sense of Gromov, and $f:Y\map X$ is any map. Let $p\in Y$.

$(\ast\ast)$ Suppose for all $N>0$ there is $M=M(N)>0$ such that $N\map \infty$ implies $M\map \infty$ with the following property:
For any $y_1, y_2\in Y$, any $(1,1)$-quasigeodesic $\alpha$ in $Y$ joining $y_1,y_2$ and any
$(1,1)$-quasigeodesic $\beta$ in $X$ joining $f(y_1),f(y_2)$, $B(p,N)\cap \alpha=\emptyset$ 
implies $B(f(p), M)\cap\beta=\emptyset$. 

Then the CT map exists for $f:Y\map X$.
\end{lemma}
\proof
Suppose $\{y_n\}$ is any sequence in $Y$. Suppose $\alpha_{i,j}$ is a $(1,1)$-quasigeodesic in $Y$ joining $y_i,y_j$
and suppose $\gamma_{i,j}$ is a $(1,1)$-quasigeodesic in $X$ joining $f(y_i), f(y_j)$. Then by Lemma \ref{gromov product meaning}
$\lim_{i,j\map \infty}(y_i.y_j)_p=\infty$ if and only if $\lim_{i,j\map \infty} d(p,\alpha_{i,j})=\infty$ and
$\lim_{i,j\map \infty}(f(y_i).f(y_j))_{f(p)}=\infty$ if and only if $\lim_{i,j\map\infty}d_X(f(p), \gamma_{i,j})=\infty$.
On the other hand by $(\ast\ast)$ $\lim_{i,j\map \infty} d(p,\alpha_{i,j})=\infty$ implies
$\lim_{i,j\map\infty}d_X(f(p), \gamma_{i,j})=\infty$. Thus $\{y_n\}$ converges to a point of $\partial_s Y$ implies
$\{f(y_n)\}$ converges to a point of $\partial_s X$.
The same argument shows that if $\{y_n\}$ and $\{z_n\}$ are two sequences in $Y$ representing the same point of 
$\partial_s Y$ then $\{f(y_n)\}$ and $\{f(z_n)\}$ also represent the same point of $\partial_s X$. 
Thus we have a well-defined map $\partial f: \partial_s Y\map \partial_s X$.

Now we prove the continuity of the map.
We need to show that if $\xi_n\map \xi$ in $\partial_s Y$ then $\partial f(\xi_n)\map \partial f(\xi)$. 
Suppose $\xi_n$ is represented by the class of $\{y^n_k\}_k$ and $\xi$ is the equivalence
class of $\{y_k\}$. Then $$\lim_{n\map \infty}(\liminf_{i,j\map \infty}(y^n_i.y_j)_p)=\infty.$$ 
By Lemma \ref{gromov product meaning}
then we have  $$\lim_{n\map \infty}(\liminf_{i,j\map \infty}d(p,\alpha^n_{i,j})=\infty$$ for any $(1,\epsilon)$-quasigeodesic
$\alpha^n_{i,j}$ in $Y$ joining $y^n_i$ and $y_j$. By $(\ast)$ then we have 
$$\lim_{n\map \infty}(\liminf_{i,j\map \infty}d(f(p),\gamma^n_{i,j})=\infty$$ 
where $\gamma^n_{i,j}$ is any $(1,\epsilon)$-quasigeodesic in $X$ joining $f(y^n_i), f(y_j)$. This in turn implies
by Lemma \ref{gromov product meaning} that
$$\lim_{n\map \infty}(\liminf_{i,j\map \infty}(f(y^n_i).f(y_j))_{f(p)})=\infty.$$ 
Therefore, $\partial f(\xi_n)\map \partial f(\xi)$ as was required. \qed

\noindent
{\bf Examples and remarks:}
\begin{enumerate}
\item Suppose $f:\RR_{\geq 0}\map \RR_{\geq 0}$ is the function $f(x)=e^x-1$. Then $f$ is not coarsely Lipschitz but $f$ admits the CT map.

\item One can easily cook up an example along the line of the above example where metric properness is also violated but the CT map exists
as we see in the example below. We will see another interesting example in Corollary \ref{cor: ct to base}.

\item The condition $(\ast)$ in the above lemma is also not necessary in general for the existence of the CT map. Here is an example
in which both metric properness and $(\ast)$ fail to hold but nevertheless the CT map exists. Suppose $X$ is a tree 
built in two steps. First we have a star, i.e. a tree with one central vertex on which end points of finite intervals are glued where
the lengths of the intervals are unbounded. Then two distinct rays are glued to each vertex of the star other than the central vertex. 
Suppose $Y$ is obtained by collapsing the central star in $X$ to a point and $f$ is the quotient map. Then clearly the CT map
exists but $(\ast)$ is violated.
\end{enumerate}

The following lemma is very standard and hence we skip mentioning its proof.

\begin{lemma}{\em({\bf Functoriality of CT maps})}\label{CT properties}
(1) Suppose $X,Y,Z$ are hyperbolic metric spaces and $f:X\map Y$ and $g:Y\map Z$ admit the CT maps. Then
so does $g\circ f$ and $\partial (g\circ f)=\partial g \circ \partial f$.

(2) If $i:X\map X$ is the identity map then it admits the CT map $\partial i$  which is the identity map on $\partial_s X$

(3) If two maps $f,h:X\map Y$ are at a finite distance admitting the CT maps then they induce the same CT map.

(4) Suppose $f:X\map Y$ is a qi embedding of hyperbolic length spaces. 
Then $f$ admits the CT map
$\partial f:\partial_s X\map \partial_s Y$ which is a homeomorphism onto the image.

If $f$ is a quasiisometry then $\partial f$ is a homeomorphism. In particular, the action by left multiplication
of a hyperbolic group $G$ on itself induces an action of $G$ on $\partial G$ by homeomorphisms.
\end{lemma}

\subsubsection{Limit sets}
\begin{defn}
Suppose $X$ is a hyperbolic metric space and $A\subset X$. Then the {\em limit set} of $A$ in $X$ is
the set $\Lambda_X(A)=\{\lim_{n\map \infty} a_n\in \partial_s X: \{a_n\}\, \mbox{is a sequence in}\, A\}$.
\end{defn}
When $X$ is understood then the limit set of $A\subset X$ will be denoted simply by $\Lambda(A)$.
In this subsection, we collect some basic results on limit sets that we need in Section 6 of the
paper. In each case, we briefly indicate the proofs for the sake of completeness. The following is
straightforward.

\begin{lemma}\label{hausdorff limset}
Suppose $X$ is a hyperbolic metric space and $A, B\subset X$ with $Hd(A,B)<\infty$.
Then $\Lambda(A)= \Lambda(B)$.
\end{lemma}

\begin{lemma}\label{coarse separation vs limit set}
Suppose $X$ is a hyperbolic metric space and $Y\subset X$. Suppose $Z\subset Y$ coarsely bisects $Y$
in $X$ into $Y_1, Y_2$ where $Z\subset Y_1\cap Y_2$. Then $\Lambda(Y_1)\cap \Lambda(Y_2)=\Lambda(Z)$.
\end{lemma}
\proof This is a straightforward consequence of Lemma \ref{gromov product meaning}. \qed

\begin{lemma}\label{limset lem}
Suppose $X$ is a $\delta$-hyperbolic metric space and $A\subset X$ is $\lambda$-quasiconvex.
Suppose $\xi\in \Lambda(A)$ and $\gamma$ is a $K$-quasigeodesic ray converging to $\xi$.
Then there are $N\in \NN$ and $D=D_{\ref{limset lem}}(\delta, \lambda, K)>0$
such that $\gamma(n)\in N_D(A)$ for all $n\geq N$. 
\end{lemma}
\proof Rather than explicitly computing the constants we indicate how to obtain them.
Suppose $\{x_n\}$ is a sequence in $A$ such that $x_n\map \xi$.
Let $y_1\in \gamma$ be a $1$-approximate nearest point projection of $x_1$ on $\gamma$.
Let $\alpha_1$ denote a $(1,1)$-quasigeodesic joining $x_1,y_1$.
Then the concatenation, say $\gamma_1$, of $\alpha_1$ and the segment of $\gamma$ from $y_1$ to $\xi$ is a uniform 
quasigeodesic by Corollary \ref{gluing quasigeodesics}. For all $m>1$, let $y_m$ denote a $1$-approximate nearest
point projection of $x_m$ on $\gamma_1$. Then $y_m$ is contained in $\gamma_1$ for all large $m$.
However, once again by Corollary \ref{gluing quasigeodesics} the concatenation of the portion of $\gamma_1$
between $x_1, y_m$ and a $1$-quasigeodesic joining $x_m,y_m$ is a uniform quasigeodesic. Now it follows
by stability of quasigeodesics that the segment of $\gamma_1$ between $y_1,y_m$ is contained
in a uniformly small neighborhood of $A$ since $A$ is quasiconvex.
\qed

\begin{lemma}\label{CT-limset}
Suppose $X,Y$ are hyperbolic metric spaces, and $f:Y\map X$ is any metrically proper map.
Suppose that the CT map exists for $f$. Then we have $\Lambda(f(Y))=\partial f(\partial Y)$
in each of the following cases:

(1) $Y$ is a proper metric space.

(2) $f$ is a qi embedding.
\end{lemma}
\proof (1) It is clear that $\partial f(\partial Y)\subset \Lambda(f(Y))$.
Suppose $y_n$ is any sequence such that $f(y_n)\map \xi$ for some $\xi\in \partial_s X$.
Since $f$ is proper $\{y_n\}$ is an unbounded sequence. Since $Y$ is a proper length space
it is a geodesic metric space by Hopf-Rinow theorem (see \cite{bridson-haefliger}, Proposition 3.7, Chapter I.3).
Now it is a standard fact that any unbounded sequence in a proper geodesic metric space has a subsequence converging
to a point of the Gromov boundary of the space. Since $Y$ is proper, we have a subsequence $\{y_{n_k}\}$
of $\{y_n\}$ such that $y_{n_k}\map \eta$ for some $\eta\in \partial_s Y$. It is clear that 
$\partial f(\eta)=\xi$. Hence $\Lambda(f(Y))\subset \partial f(\partial Y)$.

(2) Let $y\in Y$ and $x=f(y)$. Suppose $\{y_n\}$ is a sequence of points in $Y$ such that 
$\lim_{m,n\map \infty} (f(y_m).f(y_n))_x=\infty$ and $\eta=[\{f(y_m)\}]$. Then by Lemma \ref{gromov product meaning}
for any $1$-quasigeodesic $\beta_{m,n}$ in $X$ joining $f(y_m),f(y_n)$ for all $m,n\in \NN$, we have $\lim_{m,n\map \infty}d_X(x, \beta_{m,n})=\infty$.
Since $f$ is a qi embedding if $\alpha_{m,n}$ is a $1$-quasigeodesic in $Y$ joining $y_m, y_n$ for all $m,n\in \NN$ then
$f(\alpha_{m,n})$ are uniform quasigeodesics in $X$. Hence, by stability of quasigeodesics in $X$ we have
$Hd(f(\alpha_{m,n}), \beta_{m,n})<D$ for some constant $D\geq 0$. Thus $\lim_{m,n\map \infty}d_X(x, f(\alpha_{m,n}))=\infty$.
Since $f$ is a qi embedding and $x=f(y)$ it follows that $\lim_{m,n\map \infty}d_Y(y, \alpha_{m,n})=\infty$. Therefore,
$\lim_{m,n\map \infty} (y_m.y_n)_y=\infty$ again by Lemma \ref{gromov product meaning}. Hence, if $\xi=[\{y_n\}]$ then $\partial f(\xi)=\eta$.
\qed

\begin{lemma}{\em (Projection of boundary points on quasiconvex sets)}\label{qc last lemma}
Given $\delta\geq 0$ and $k\geq 0$ there is a constant $R=R_{\ref{qc last lemma}}(\delta, k)$ such that
the following holds:

Suppose $X$ is a $\delta$-hyperbolic metric space, $A\subset X$ is $k$-quasiconvex and $\xi\in \partial X\setminus \Lambda(A)$.
Then there is a point $x\in A$ with the following property: Suppose $\{x_n\}$ is any sequence where $x_n\map \xi$. Then there is an $N>0$
 such that for all $n\geq N$ we have $P_A(x_n)\in A\cap B(x,R)$.
\end{lemma}
\proof Suppose $\{x_n\}, \{y_n\}$ are two sequences in $X$ such that $x_n\map \xi$ and $y_n\map \xi$. 
Let $\alpha_{m,n}$ be a $1$-quasigeodesic in $X$ joining $x_m,y_n$ for all $m,n\in \NN$. Let $P_A:X\map A$ be a $1$-approximate 
nearest point projection on $A$.

{\bf Claim:} There is a constant $R_0>0$ depending only on $\delta$ and $k$ and there is $N>0$ such that
$diam(P_A(\alpha_{m,n}))\leq R_0$ for all $m,n\geq N$.

We first note that $\lim_{m,n\map \infty} d(A,\alpha_{m,n})=\infty$. In fact,
if this is not the case then there is $r>0$ such that for all $N>0$ there are $m_N,n_N\geq N$ with $d(A,\alpha_{m_N,n_N})\leq r$. 
In that case let $a_N\in A$ be such that $d(a_N, \alpha_{m_N,n_N})\leq r$. It is then clear that $a_N\map \xi$ by Lemma \ref{lem: bdry defn}(1),
contradicting the hypothesis that $\xi\not\in \Lambda(A)$. 
By stability of quasigeodesics, any $1$-quasigeodesic is uniformly quasiconvex in $X$ and $A$ is given to be $k$-quasiconvex.
Hence, by Corollary \ref{cobounded cor} there are constants $D_0,R_0$ such that $d(A, \alpha_{m,n})>D_0$ implies that 
$diam(P_A(\alpha_{m,n}))\leq R_0$. Since, $\lim_{m,n\map \infty} d(A,\alpha_{m,n})=\infty$ there is $N>0$ such that
$d(A,\alpha_{m,n})>D_0$ for all $m,n\geq N$. This proves the existence of $N$ and $R_0$.

Now, by specializing the claim to the case $\{x_n\}=\{y_n\}$ we have $N_0>0$ such that if $\beta_{m,n}$ is a $1$-quasigeodesic
joining $x_m,x_n$ then $diam(P_A(\beta_{m,n}))\leq R_0$ for all $m,n\geq N_0$. Let $x=P_A(x_{N_0})$. 
Now, given any sequence $\{x'_n\}$ in $X$ with $x'_n\map \xi$ by the claim there is $M>0$ such that for all
$m,n\geq M$, $d(P_A(x_m), P_A(x'_n))\leq R_0$. Hence, if $N=\max\{N_0, M\}$ then 
$d(x,P_A(x'_n))\leq d(x, P_A(x_N))+d(P_A(x_N), P_A(x'_n))\leq 2R_0$. Thus we can take $R=2R_0$.\qed

Since the point $x\in A$ in the above lemma is coarsely unique we shall call any such point to be {\em the nearest
point projection} of $\xi$ on $A$ and we shall denote it by $P_A(\xi)$.



\section{Metric bundles}

In this section, we recall necessary definitions and some elementary properties of the primary objects of
study in this paper namely, metric bundles and metric graph bundles from \cite{pranab-mahan}.
We make a minor modification (see Definition \ref{defn-mbdle2}) to the definition of a metric bundle but use the same
definition of metric graph bundles as in \cite{pranab-mahan}.

\subsection{Basic definitions and properties.}
\begin{defn}\label{defn-mbdle}{\bf (Metric bundles \cite[Definition 1.2]{pranab-mahan})}
Suppose $(X,d)$ and $(B, d_B)$ are geodesic metric spaces; let $c\geq 1$ and let 
$\eta:[0,\infty) \rightarrow [0,\infty)$ be a function.
We say that $X$ is an $(\eta,c)-$ {\bf metric bundle} over $B$ if there is a surjective $1$-Lipschitz
map $\pi:X\rightarrow B$ such that the following conditions hold:\\
{\bf (1)} For each point $z\in B$, $F_z:=\pi^{-1}(z)$ is a geodesic metric space
with respect to the path metric $d_z$ induced from $X$. The inclusion maps
$i: (F_z,d_z) \rightarrow X$ are uniformly metrically proper as measured by $\eta$. \\
{\bf (2)}  Suppose $z_1,z_2\in B$, $d_B(z_1,z_2)\leq 1$ and let $\gamma$ be
a geodesic in $B$ joining them. \\
Then for any point $z\in \gamma$ and $x\in F_z$ there is a path $\tilde{\gamma}:[0,1]\map \pi^{-1}(\gamma)\subset X$
of length at most $c$ such that $\tilde{\gamma}(0)\in F_{z_1}$, $\tilde{\gamma}(1)\in F_{z_2}$
and $x\in \tilde{\gamma}$.
\end{defn}

If $X$ is a metric bundle over $B$ in the above sense then we shall refer to it as a {\bf geodesic metric bundle}
in this paper. However, the above definition seems a little restrictive. Therefore, we propose the following. 

\begin{defn}\label{defn-mbdle2}{\bf (Length metric bundles)}
Suppose $(X,d)$ and $(B, d_B)$ are {\em length spaces}, $c\geq 1$ and we have a function
$\eta:[0,\infty) \rightarrow [0,\infty)$.
We say that $X$ is an $(\eta,c)-$ {\bf length metric bundle} over $B$ if there is a surjective $1$-Lipschitz
map $\pi:X\rightarrow B$ such that the following conditions hold:\\
{\bf (1)} For each point $z\in B$, $F_z:=\pi^{-1}(z)$ is a length space
with respect to the path metric $d_z$ induced from $X$. The inclusion maps
$i: (F_z,d_z) \rightarrow X$ are uniformly metrically proper as measured by $\eta$. \\
{\bf (2)}  Suppose $z_1,z_2\in B$, and let $\gamma$ be a {\em path} of length at most $1$ in $B$ joining them. \\
Then for any point $z\in \gamma$ and $x\in F_z$ there is a path $\tilde{\gamma}:[0,1]\map \pi^{-1}(\gamma)\subset X$
of length at most $c$ such that $\tilde{\gamma}(0)\in F_{z_1}$, $\tilde{\gamma}(1)\in F_{z_2}$
and $x\in \tilde{\gamma}$.
\end{defn}
Given length spaces $X$ and $B$ we will say that $X$ is a {\bf length metric bundle} over $B$ if $X$ is an
$(\eta,c)$-length metric bundle over $B$ in the above sense for some function $\eta:{\mathbb R}^+ \rightarrow {\mathbb R}^+$ 
and some constant $c\geq 1$.  

\begin{convention}
From now on whenever we speak of a metric bundle we mean a length metric bundle.
\end{convention}

\begin{defn}\label{defn-mgbdl}{\bf (Metric graph bundles \cite[Definition 1.5]{pranab-mahan})}
Suppose $X$ and $B$ are metric graphs. Let $\eta:[0,\infty) \rightarrow [0,\infty)$ be
a function. We say that $X$ is an $\eta$-{\bf metric graph bundle}
over $B$ if there exists a surjective simplicial map $\pi:X\rightarrow B$  such that:\\
$1.$ For each $b\in \mathcal{ V}(B)$, $F_b:=\pi^{-1}(b)$ is a connected subgraph of $X$ and the inclusion maps
$i: F_b\rightarrow X$ are 
uniformly metrically proper as measured by $\eta$ for the  path metrics $d_b$ induced on $F_b$.\\
$2.$ Suppose $b_1,b_2\in \mathcal{ V}(B)$ are adjacent vertices.
Then each vertex $x_1$ of $F_{b_1}$ is connected by an edge with a vertex in $F_{b_2}$.
\end{defn}

\begin{rem}
Since the map $\pi$ is simplicial it follows that it is $1$-Lipschitz.
\end{rem}


For a metric (graph) bundle the spaces $(F_z,d_z)$, $z\in B$ will
be referred to as {\bf fibers} and the $d_z$-distance between two points in $F_z$ will be referred to as their
{\bf fiber distance}. A geodesic in $F_z$ will be called a {\bf fiber geodesic}. The spaces $X$ and $B$
will be referred to as the {\em total space}  and the {\em base space} of the bundle respectively. 
By a statement of the form `$X$ is a metric bundle (resp. metric graph bundle)' we will mean that it 
is the total space of a metric bundle (resp. metric graph bundle).

Most of the results proved for geodesic metric bundles in \cite{pranab-mahan} have their analogs for length
metric bundles. We explicitly prove this phenomenon or provide sufficient arguments for all the results
needed for our purpose.

\begin{convention}
Very often in a lemma, proposition, corollary, or a theorem we shall omit explicit mention of some of the parameters
on which a constant may depend if the parameters are understood.
\end{convention}




\begin{defn}
Suppose $\pi: X\map B$ is a metric (graph) bundle. 

(1) Suppose $A\subset B$ and $k\geq 1$. A {\em $k$-qi section} over $A$ is a $k$-qi embedding
$s:A\map X$ (resp. $s:\VV(A)\map X$)such that $\pi \circ s=Id_A$ (resp. $\pi \circ s=Id_{\VV(A)}$) 
where $A$ has the restricted metric from $B$ and $Id_A$ (resp. $Id_{\VV(A)}$)
denotes the identity map on $A\map A$ (resp. $\VV(A)\map \VV(A)$).

(2) Given any metric space (resp. graph) $Z$ and any qi embedding $f: Z\map B$ (resp. $f:\VV(Z)\map \VV(B)$)
a {\em $k$-qi lift} of $f$ is a $k$-qi embedding $\tilde{f}:Z\map X$ (resp. $\tilde{f}:\VV(Z)\map \VV(X)$)
such that $\pi\circ \tilde{f}=f$.
\end{defn}

\begin{convention}
(1) Most of the time we shall refer to the image of a qi section (or a qi lift) to be the qi section (resp. the qi lift). \\
(2) Suppose $\gamma:I\map B$ is a (quasi)geodesic and $\tilde{\gamma}$ is a qi lift of $\gamma$. Let $b=\gamma(t)$ for some $t\in I$.
Then we will denote $\tilde{\gamma}(t)$ by $\tilde{\gamma}(b)$ also. \\ 
(3) In the context of a metric graph bundle $(X,B,\pi)$, when we talk about a point in $X$, $B$ or a fiber, we mean that the point is a vertex in the corresponding space.
\end{convention}

The following lemma is immediate from the definition of a metric (graph) bundle. Hence we briefly indicate its
proof.

\begin{lemma}{\em ( Path lifting lemma)}\label{lifting geodesics}
Suppose $\pi:X\map B$ is an $(\eta,c)$-metric bundle or an $\eta$-metric graph bundle.
\begin{enumerate}
\item Suppose $b_1, b_2\in B$. Suppose $\gamma:[0,L]\map B$ is a continuous, 
rectifiable, arc length parameterized path (resp. an edge path) in $B$ joining $b_1$ to $b_2$. 
Given any $x\in F_{b_1}$ there is a path $\tilde{\gamma}$ in $\pi^{-1}(\gamma)$
such that $l(\tilde{\gamma})\leq (L+1)c$ (resp $l(\tilde{\gamma})=L)$ joining $x$ to some point of $F_{b_2}$. 

In particular, in case $X$ is a metric graph bundle over $B$ any geodesic $\gamma$ of $ B$ can be lifted
to a geodesic starting from any given point of $\pi^{-1}(\gamma)$.

\item For any $k\geq 1$ and $\epsilon\geq 0$, any {\em dotted} 
$(k,\epsilon)$-quasigeodesic $\beta:[m,n]\map B$ has a lift
$\tilde{\beta}$ starting from any point of $F_{\beta(m)}$ such that the following hold, where we assume $c=1$
for metric graph bundles.

For all $i,j\in [m,n]$ we have $$-\epsilon+\frac{1}{k}|i-j|\leq d_X(\tilde{\beta}(i), \tilde{\beta}(j))\leq c\cdot(k+\epsilon+1)|i-j|.$$
In particular it is a $c\cdot(k+\epsilon+1)$-qi lift of $\beta$. Also we have $$l(\tilde{\beta})\leq ck(k+\epsilon+1)(\epsilon+d_B(b_1,b_2)).$$
\end{enumerate}
\end{lemma}
\proof (1)
We fix a sequence of points $0=t_0, t_1,\cdots, t_n=L$ in $[0,L]$ such that 
$l(\gamma |_{[t_i, t_{i+1}]})=1$ for $0\leq i< n-1$ and $l(\gamma |_{[t_{n-1}, t_{n}]})\leq 1$ for the metric bundle case. 
For the metric graph bundle $\gamma(t_i)$ are the consecutive vertices on $\gamma$, $0\leq i\leq L=n$. 
Now given any $x=:x_0\in F_{t_0}$ we can inductively construct a sequence of points $x_i\in F_{t_i}$, $0\leq i\leq n$ and a sequence of
paths $\alpha_i$ of length at most $c$ (resp. an edge) joining $x_i$ to $x_{i+1}$ for $0\leq i\leq n-1$. 
Concatenation of these paths gives a candidate for $\tilde{\gamma}$. 

The second statement for metric graph bundles follow because $\pi:X\map B$ is a $1$-Lipschitz map.

(2)
We construct a lift $\tilde{\beta}$ of $\beta$ starting from any point $x\in F_{\beta(m)}$ inductively as follows. 
We know that $d_B(\beta(i), \beta(i+1))\leq k+\epsilon$. Let $\beta_i$ be a path in $B$ joining $\beta(i)$ to $\beta(i+1)$
which is of length at most $k+\epsilon+1$ for $m\leq i\leq n-1$. We can then find a sequence of paths of length at most
$(k+\epsilon+1)\cdot c$ in $\pi^{-1}(\beta_i)$ (where $c=1$ for metric graph bundle) 
$m\leq i\leq n-1$ using the first part of the lemma such that $\beta_m$ starts at $x$ and $\beta_{i+1}$ starts at the
end point of $\beta_i$ for $m+1\leq i\leq n-1$. Let $x_i$ be the starting point of $\beta_i$ for $m\leq i\leq n-1$ and let
$x_n$ be the end point of $\beta_{n-1}$. Then we define $\tilde{\beta}$ by setting $\tilde{\beta}(i)=x_i$, $m\leq i\leq n$.

Clearly $d_X(\tilde{\beta}(i), \tilde{\beta}(j))\leq c\cdot(k+\epsilon+1)|i-j|$. Also, 
$d_B(\pi\circ\tilde{\beta}(i), \pi\circ\tilde{\beta}(j))=d_B(\beta(i), \beta(j))\leq d_X(\tilde{\beta}(i), \tilde{\beta}(j))$ since
$\pi$ is $1$-Lipschitz. Since $\beta$ is a dotted $(k, \epsilon)$ quasigeodesic, we have 
$-\epsilon+\frac{1}{k}|i-j|\leq d_B(\beta(i), \beta(j))$. This proves that
$$-\epsilon+\frac{1}{k}|i-j|\leq d_X(\tilde{\beta}(i), \tilde{\beta}(j))\leq c\cdot(k+\epsilon+1)|i-j|.$$
For the last part of (2) we see that 
$$l(\tilde{\beta})=\sum_{i=m}^{n-1} d_X(\tilde{\beta}(i), \tilde{\beta}(i+1))\leq \sum_{i=m}^{n-1} c\cdot(k+\epsilon+1)=(n-m)c\cdot(k+\epsilon+1).$$
On the other hand since $\beta$ is a $(k,\epsilon)$-quasigeodesic we have $-\epsilon+ \frac{1}{k}(n-m)\leq d_B(b_1, b_2)$.
The conclusion immediately follows from these two inequalities. 
\qed

The following corollary follows from the proof of Proposition 2.10 of \cite{pranab-mahan}. We include it for the sake of completeness.
\begin{cor}\label{path lifting remark}
Given  any metric (graph) bundle $\pi:X\map B$ and $b_1, b_2\in B$ 
we can define a map $\phi:F_{b_1}\map F_{b_2}$ such that 
$d_X(x, \phi(x))\leq 3c +3cd_B(b_1, b_2)$ (resp. $d(x, \phi(x))=d_B(b_1,b_2)$) for all $x\in F_{b_1}$.
\end{cor}

\proof  The statement about the metric graph bundle is trivially true by Lemma \ref{lifting geodesics} (1).
For the metric bundle case, fix a dotted $1$-quasigeodesic  $\gamma$ joining
$b_1$ to $b_2$. Then for all $x\in F_{b_1}$ fix for once and all a dotted lift  $\tilde{\gamma}$ as constructed in the proof of the
Lemma \ref{lifting geodesics} which starts from $x$ and set $\phi(x)=\tilde{\gamma}(b_2)$. The statement then follows from  Lemma \ref{lifting geodesics}(2).
\qed
\begin{rem}
For all $b_1, b_2\in B$ any  map $f:F_{b_1}\map F_{b_2}$ such that $d_X(x, f(x))\leq D$ for some constant
$D$ independent of $x$ will be referred to as a {\bf fiber identification map}.
\end{rem}

The proof of the first part of the following lemma is immediate from Corollary \ref{path lifting remark}
whereas the next two parts essentially follow from the proof of Proposition 2.10 of \cite{pranab-mahan}.
Hence we skip the proofs.

\begin{lemma}\label{fibers qi}
Suppose $\pi:X\map B$ is an $(\eta,c)$-metric bundle or an $\eta$-metric graph bundle and 
$R\geq 0$. Suppose $b_1, b_2\in B$. The we have the following.
\begin{enumerate}
\item $Hd(F_{b_1},F_{b_2})\leq 3c+ 3cd_B(b_1,b_2)$ (resp. $Hd(F_{b_1},F_{b_2})= d_B(b_1,b_2)$).

\item Suppose $\phi_{b_1b_2}:F_{b_1}\map F_{b_2}$ is a map such that for all $x\in F_{b_1}$, $d(x,\phi_{b_1b_2}(x))\leq R$
for all $x\in F_{b_1}$. 

Then $\phi_{b_1b_2}$ is a $K_{\ref{fibers qi}}=K_{\ref{fibers qi}}(R)$-quasiisometry which is $D_{\ref{fibers qi}}$-surjective. 

\item If $\psi_{b_1b_2}:F_{b_1}\map F_{b_2}$ is any other map such that $d(x,\psi_{b_1b_2}(x))\leq R'$
for all $x\in F_{b_1}$ then $d(\phi_{b_1b_2}, \psi_{b_1b_2})\leq \eta(R+R')$. 

In particular, the maps $\phi_{b_1b_2}$ are coarsely unique (see Definition \ref{defn: 1.2}(7)).
\end{enumerate}
\end{lemma}

In this lemma, we deliberately suppress the dependence of $K_{\ref{fibers qi}}$ on the parameter(s) of the bundle.

\begin{cor}\label{fibers unif qi}
Suppose $\pi:X\map B$ is a metric (graph) bundle and $b_1, b_2\in B$ (resp. $b_1, b_2\in \VV(B)$) such that $d_B(b_1, b_2)\leq R$.
Suppose $\phi_{b_1b_2}:F_{b_1}\map F_{b_2}$ is a fiber identification map as constructed in the proof of
Corollary \ref{path lifting remark}. Then $\phi_{b_1b_2}$ is a $K_{\ref{fibers unif qi}}=K_{\ref{fibers unif qi}}(R)$-quasiisometry.
\end{cor}
\proof By Corollary \ref{path lifting remark} $d_X(x, \phi_{b_1b_2}(x))\leq 3c +3cd_B(b_1,b_2)\leq 3c+3cR$ 
for all $x\in F_{b_1}$ (resp. $d_X(x, \phi_{b_1b_2}(x))=d_B(b_1,b_2)\leq R$ for all $x\in \VV(B)$).
Hence by Lemma \ref{fibers qi}(2) $\phi_{b_1b_2}$ is $K_{\ref{fibers unif qi}}=K_{\ref{fibers qi}}(3c+3cR)$-qi for the
metric bundle and $K_{\ref{fibers unif qi}}=K_{\ref{fibers qi}}(R)$-qi for the metric graph bundle case.
\qed

The following corollary is proved as a simple consequence of Lemma \ref{fibers qi} and Corollary
\ref{path lifting remark}. (See Corollary 1.14, and Corollary 1.16 of \cite{pranab-mahan}.)
Therefore, we skip the proof of it.

\begin{cor}\label{bdd-flaring}{\em (Bounded flaring condition)}
For all $k\in \mathbb R$, $k\geq 1$ there is a function $\mu_k:\NN \rightarrow \NN$
such that the following holds:

Suppose $\pi: X\map B$ is an $(\eta,c)$-metric bundle or an $\eta$-metric graph bundle.
Let $\gamma\subset B$ be a dotted $(1,1)$-quasigeodesic (resp. a geodesic) joining $b_1,b_2\in B$,
and let $\tilde{\gamma_1}$, $\tilde{\gamma_2}$ be two $k$-qi lifts of $\gamma$ in $X$.
Suppose $\tilde{\gamma}_i(b_1)=x_i\in F_{b_1}$ and $\tilde{\gamma}_i(b_2)=y_i\in F_{b_2}$, $i=1,2$.

Then  
\[
d_{b_2}(y_1,y_2)\leq \mu_k(N) \mbox{max}\{d_{b_1}(x_1,x_2),1\}.
\]
if $d_B(b_1,b_2)\leq N$.
\end{cor}
\begin{comment}
\proof By Remark \ref{path lifting remark} we may define a fiber identification map $\phi_{b_1b_2}:F_{b_1}\map F_{b_2}$
where $d_X(x, \phi_{b_1b_2}(x))\leq N$. By Lemma \ref{fibers qi}(2) $\phi_{b_1b_2}$ is a
$K_{\ref{fibers qi}}(N)$-quasiisometry. Hence, $d_{b_2}(\phi_{b_1b_2}(x_1), \phi_{b_1b_2}(y_1))\leq K_{\ref{fibers qi}}(N)(1+d_{b_1}(x_1, y_1))$.
On the other hand $d_X(x_1,x_2)\leq kd_B(b_1,b_2)+k\leq kN+k$ and similarly $d_X(x_1,x_2)\leq kN+k$. Hence,
$d_{b_2}(\phi_{b_1b_2}(x_1), x_2)\leq \eta(N+kN+k)$ and similarly $d_{b_2}(\phi_{b_1b_2}(y_1), y_2)\leq \eta(N+kN+k)$.
Thus $d_{b_2}(x_2,y_2)\leq d_{b_2}(\phi_{b_1b_2}(x_1), x_2)+ d_{b_2}(\phi_{b_1b_2}(y_1), y_2)+ d_{b_2}(\phi_{b_1b_2}(x_1), \phi_{b_1b_2}(y_1))$.
It follows that $d_{b_2}(x_2,y_2)\leq 2\eta(N+kN+k)+ K_{\ref{fibers qi}}(N)(1+d_{b_1}(x_1, y_1))$. Since $d_{b_1}(x_1, y_1))$ is a nonnegative
integer we can clearly take $\mu_k(N)= 2\eta(N+kN+k)+ K_{\ref{fibers qi}}(N)$. \qed
\end{comment}
In the rest of the paper, we will summarize the conclusion 
of Corollary \ref{bdd-flaring} by saying that a  metric (graph) bundle satisfies the {\bf bounded flaring condition}.

\begin{rem}{\em (Metric bundles in the literature)}
Metric (graph) bundles appear in several places in other people's work. In \cite[Section 2.1]{bowditch-stacks}
Bowditch defines {\em stacks of (hyperbolic) spaces} which can easily be shown to be quasiisometric to
metric graph bundles over an interval in $\RR$. Conversely, a metric (graph) bundle whose base
is an interval in $\RR$ is clearly a stack of spaces as per \cite[Section 2.1]{bowditch-stacks}.
In \cite{whyte-bundles} Whyte defines {\em coarse bundles} which are also
quasiisometric to metric graph bundles but with additional restrictions.
\end{rem}


\smallskip


\subsection{Some natural constructions of metric bundles} 
In this section, we discuss a few general constructions that produce metric (graph) bundles. 

\smallskip

\begin{defn}
$(1)$ {\bf (Metric bundle morphisms)}
Suppose $(X_i, B_i,\pi_i)$, $i=1,2$ are metric bundles. A morphism from $(X_1,B_1,\pi_1)$ to  $(X_2,B_2,\pi_2)$ (or simply
from $X_1$ to $X_2$ when there is no possibility of confusion) consists of 
a pair of coarsely $L$-Lipschitz maps $f:X_1\map X_2$ and $g:B_1\map B_2$ for some $L\geq 0$ such that $\pi_2\circ f=g\circ \pi_1$,
i.e. the following diagram (Figure \ref{pic}) is commutative.


\begin{figure}[ht]
	\centering
	\begin{tikzpicture}[node distance=2cm,auto]
	\node (A) {$X_1$}; 
	\node (B) [right of= A] {$X_2$};
	\node (C) [below=1cm of A] {$B_1$};
	\node (D) [below=1cm of B] {$B_2$};
	\draw [->] (A) to node {$f$} (B);
	\draw [->] (A) to node [swap] {$\pi_1$} (C);
	\draw [->] (B) to node {$\pi_2$} (D);
	\draw [->] (C) to node [swap]{$g$}(D);
	\end{tikzpicture}
	\caption{}\label{pic}
\end{figure}


$(2)${\bf (Metric graph bundle morphisms)}
Suppose $(X_i, B_i,\pi_i)$, $i=1,2$ are metric graph bundles. A morphism from $(X_1,B_1,\pi_1)$ to  $(X_2,B_2,\pi_2)$ (or simply
from $X_1$ to $X_2$ when there is no possibility of confusion) consists of 
a pair of coarsely $L$-Lipschitz maps $f:\VV(X_1)\map \VV(X_2)$ and $g:\VV(B_1)\map \VV(B_2)$ for some $L\geq 0$ such that 
$\pi_2\circ f=g\circ \pi_1$.

$(3)$ {\bf (Isomorphisms)} A morphism $(f,g)$ from a metric (graph) bundle $(X_1,B_1,\pi_1)$ to a metric (graph) bundle $(X_2, B_2, \pi_2)$
is called an isomorphism if there is a morphism $(f', g')$ from $(X_2, B_2, \pi_2)$ to $(X_1,B_1,\pi_1)$ such that $f'$
is a coarse inverse of $f$ and $g'$ is a coarse inverse of $g$.
\end{defn}
We note that for any  morphism $(f,g)$ from a metric (graph) bundle $(X_1,B_1,\pi_1)$ to a metric (graph) bundle $(X_2, B_2, \pi_2)$
we have $f(\pi^{-1}_1(b))\subset \pi^{-1}_2(g(b))$ for all $b\in B_1$.
We will denote by $f_{b}:\pi^{-1}_1(b)\map \pi^{-1}_2(g(b))$ the restriction of $f$ to $\pi^{-1}_1(b)$ for all $b\in B_1$.
We shall refer to these maps as the {\bf fiber maps} of the morphisms. We also note that in the case of
metric graph bundles coarse Lipschitzness is equivalent to Lipschitzness.

\begin{lemma}\label{metric bundle map}
Given $k\geq 1, K\geq 1$ and $L\geq 0$ there are constants $L_{\ref{metric bundle map}}, K_{\ref{metric bundle map}}$ such that the following hold.

Suppose $(f,g)$ is a morphism of metric (graph) bundles as in the definition above. Then the following hold:

(1) For all $b\in B_1$  the map $f_b:\pi^{-1}_1(b)\map \pi^{-1}_2(g(b))$ 
is coarsely $L_{\ref{metric bundle map}}$-Lipschitz with respect to
the induced length metric on the fibers.

(2)  Suppose $\gamma:I\map B_1$ is a dotted $(1,1)$-quasigeodesic (or simply a geodesic in the case of a metric graph
bundle) and suppose $\tilde{\gamma}$ is a $k$-qi lift of $\gamma$. If $g$ is a $K$-qi embedding
then $f\circ \tilde{\gamma}$ is a $K_{\ref{metric bundle map}}=K_{\ref{metric bundle map}}(k,K,L)$-qi lift of $g\circ \gamma$. 
\end{lemma}
\proof  We shall check the lemma only for the metric bundle case because for metric graph bundles the proofs are similar
and in fact easier. 

Suppose $\pi_i:X_i\map B_i$, $i=1,2$ are $(\eta_i, c_i)$-metric bundles.

(1) Let $b\in B_1$ and $x,y\in \pi^{-1}_1(b)$ be such that $d_b(x,y)\leq 1$. Since $f$ is coarsely $L$-Lipschitz,
$d_{X_2}(f(x),f(y))\leq L+Ld_{X_1}(x,y)\leq L+Ld_b(x,y)\leq 2L$. Now, the fibers of $\pi_2$ are uniformly properly embedded as measured by $\eta_2$.
Hence, $d_{g(b)}(f(x),f(y))\leq \eta_2(2L)$. Therefore, by Lemma \ref{proving Lipschitz} the fiber map
$f_b:\pi^{-1}_1(b)\map \pi^{-1}_2(g(b))$ is $\eta_2(2L)$-coarsely Lipschitz. Hence, $L_{\ref{metric bundle map}}=\eta_2(2L)$ will do.

(2) Let $\gamma_2=g\circ\gamma$ and $\tilde{\gamma}_2=f\circ \tilde{\gamma}$.
Then clearly, $\pi_2\circ \tilde{\gamma}_2=\gamma_2$ whence $\tilde{\gamma}_2$ is a lift of $\gamma_2$.
By Lemma \ref{qi composition}(1) $\tilde{\gamma}_2=f\circ \tilde{\gamma}$ is coarsely $(kL, kL+L)$-Lipschitz. Hence, for all $s,t\in I$ we have 
$$d_{X_2}(\tilde{\gamma}_2(s), \tilde{\gamma}_2(t))\leq kL|s-t|+(kL+L).$$

On the other hand, for $s,t \in I$ we have 
$$d_{X_2}(\tilde{\gamma}_2(s),\tilde{\gamma}_2(t))\geq d_{B_2}(\pi_2 \circ \tilde{\gamma}_2(s),\pi_2 \circ \tilde{\gamma}_2(t)) = d_{B_2}(\gamma_2(s),\gamma_2(t)).$$
However, by Lemma \ref{qi composition}(2) $\gamma_2=g\circ \gamma$ is a $(K, 2K)$-qi embedding. Hence, we have
$$d_{X_2}(\tilde{\gamma}_2(s),\tilde{\gamma}_2(t)) \geq d_{B_2}(\gamma_2(s),\gamma_2(t))\geq -2K+\frac{1}{K}|s-t|.$$
Therefore, it follows that $\tilde{\gamma}_2$ is a $K_{\ref{metric bundle map}}=\max\{2K, kL+L\}$-qi lift of $\gamma_2$. \qed

\smallskip
The following theorem characterizes isomorphisms of metric (graph) bundles.

\begin{theorem}\label{bundle isomorphism}
If $(f,g)$ is an isomorphism of metric (graph) bundles as in the above definition then the maps $f,g$ are quasiisometries 
and all the fiber maps are uniform quasiisometries.

Conversely, if the map $g$ is a qi and the fiber maps are uniform qi then $(f,g)$ is an isomorphism.
\end{theorem}
\proof We shall prove the theorem in the case of a metric bundle only. The proof in the case of a metric graph bundle is very similar
and hence we skip it.

If $(f,g)$ is an isomorphism then $f,g$ are qi by Lemma \ref{elem-lemma1}(1). We need to show that the fiber maps are quasiisometries.

Suppose $(f', g')$ is a coarse inverse of $(f,g)$ such that $d_{X_2}(f\circ f'(x_2), x_2)\leq R$ and
$d_{X_1}(f'\circ f(x_1), x_1)\leq R$ for all $x_1\in X_1$ and $x_2\in X_2$. It follows that for all $b_1\in B_1, b_2\in B_2$ we have 
$d_{B_1}(b_1, g'\circ g(b_1))\leq R$ and $d_{B_2}(b_2, g\circ g'(b_2))\leq R$ since the maps $\pi_1, \pi_2$ are $1$-Lipschitz. 
Suppose $f', g'$ are coarsely $L'$-Lipschitz. Let $L_1=\eta_2(2L)$ and $L_2=\eta_1(2L')$.
Then for all $u\in B_1$, $f_u:\pi^{-1}_1(u)\map \pi^{-1}_2(g(u))$ is coarsely $L_1$-Lipschitz and for all
$v\in B_2$, $f'_v:\pi^{-1}_2(v)\map \pi^{-1}_1(g'(v))$ is coarsely $L_2$-Lipschitz by Lemma \ref{metric bundle map}(1).

Let $b\in B_1$. To show that $f_b:\pi^{-1}_1(b)\map \pi^{-1}_2(g(b)$ is a uniform quasiisometry, it is enough by 
Lemma \ref{elem-lemma1}(1) to find a uniformly coarsely Lipschitz map $\pi^{-1}_2(g(b))\map \pi^{-1}_1(b)$ which is
uniform coarse inverse of $f_b$. We already know that $f'_{g(b)}$ is $L_2$-coarsely Lipschitz. 
Let $b_1=g'\circ g(b)$. We also noted that $d_{B_1}(b,b_1)\leq R$. Hence, it follows by Corollary 
\ref{path lifting remark} and Corollary \ref{fibers unif qi} 
that we have a $K_{\ref{fibers qi}}(R)$-qi $\phi_{b_1b}: \pi^{-1}_1(b_1)\map \pi^{-1}_1(b)$ such that 
$d_{X_1}(x, \phi_{b_1b}(x))\leq 3{c_1}+3{c_1}R$ for all $x\in \pi^{-1}_1(b_1)$. Let $h= \phi_{b_1b}\circ f'_{g(b)}$. 
We claim that $h$ is a uniformly coarsely Lipschitz, uniform coarse inverse of $f_b$. Since $f'_{g(b)}$ is 
$L_2$-coarsely Lipschitz and clearly $\phi_{b_1b}$ is $K_{\ref{fibers qi}}(R)$-coarsely Lipschitz, it follows by
Lemma \ref{qi composition}(1) that $h$ is $(L_2K_{\ref{fibers qi}}(R)+K_{\ref{fibers qi}}(R))$-coarsely Lipschitz.

Moreover, for all $x\in \pi^{-1}_1(b)$ we have
$d_{X_1}(x, h\circ f_b(x))\leq d_{X_1}(x, f'_{g(b)}\circ f_b(x))+d_{X_1}(f'_{g(b)}\circ f_b(x), h\circ f_b(x))\leq R+3{c_1}+3{c_1}R$.
Hence, 
$d_b(x, h\circ f_b(x))\leq \eta_1(R+3{c_1}+3{c_1}R)$. Let $y\in \pi^{-1}_2(g(b))$. Then {\small
$$d_{X_2}(y, f_b\circ h(y))=
d_{X_2}(y, f\circ \phi_{b_1b}\circ f'(y))\leq d_{X_2}(y, f\circ f'(y))+d_{X_2}(f\circ f'(y), f\circ \phi_{b_1b}\circ f'(y))$$
$$\leq R+L(3{c_1}+3{c_1}R)+L$$} since $d_{X_1}(f'(y), \phi_{b_1b}\circ f'(y))\leq 3{c_1}+3{c_1}R$. Hence, 
$d_{g(b)}(y, f_b\circ h(y))\leq \eta_2(R+L(3{c_1}+3{c_1}R)+L)$. Hence by Lemma \ref{elem-lemma1}(1) $f_b$ is a uniform qi.

Conversely, suppose all the fiber maps of the morphism $(f,g)$ are $(\lambda,\epsilon)$-qi which are $R$-coarsely surjective
and $g$ is a $(\lambda_1,\epsilon_1)$-qi which is $R_1$-surjective. Let $g'$ be a coarsely 
$(K,C)$-quasiisometric, $D$-coarse inverse of $g$ where
$K=K_{\ref{elem-lemma1}}(\lambda_1,\epsilon_1, R_1)$, $C=C_{\ref{elem-lemma1}}(\lambda_1,\epsilon_1, R_1)$ and 
$D=D_{\ref{elem-lemma1}}(\lambda_1,\epsilon_1, R_1)$.
For all $u\in B_1$ let $\bar{f}_u$ be a $D_1$-coarse inverse of $f_u:F_u\map F_{g(u)}$.
We will define a map $f':X_2\map X_1$ such that $(f',g')$ is morphism from $X_2$ to $X_1$
and $f'$ is a coarse inverse of $f$ as follows. 

For all $u\in B_2$
we define $f'_{u}: F_{u}\map F_{g'(u)}$ as the composition $\bar{f}_{g'(u)}\circ \phi_{u g(g'(u))}$
where $\phi_{u g(g'(u))}$ is a fiber identification map as constructed in the proof of Corollary \ref{path lifting remark}.
Collectively this defines $f'$. Now we shall check that $f'$ satisfies the desired properties.

(i) We first check that $(f',g')$ is a morphism.
It is clear from the definition that $\pi_1\circ f'=g'\circ \pi_2$. Hence we will be done by showing that $f'$ is coarsely Lipschitz.
By Lemma \ref{proving Lipschitz} it is enough to show that for all $u_2, v_2\in B_2$ and $x\in F_{u_2}, y\in F_{v_2}$ with 
$d_{X_2}(x,y)\leq 1$, $d_{X_1}(f'(x), f'(y))$ is uniformly small. We note that $d_{B_2}(u_2, v_2)\leq 1$. 
Let $u_1=g'(u_2)$ and $v_1=g'(v_2)$. Then $d_{B_1}(u_1, v_1)\leq K+C$, $d_{B_2}(u_2, g(u_1))\leq D$ and $d_{B_2}(v_2, g(v_1))\leq D$. 
This means $d_{X_2}(x,\phi_{u_2g(u_1)}(x))\leq 3Dc_2+3c_2$ and $d_{X_2}(y, \phi_{v_2g(v_1)})\leq 3Dc_2+3c_2$ by Lemma \ref{lifting geodesics} and
Corollary \ref{path lifting remark}. Hence, $d_{X_2}(\phi_{u_2g(u_1)}(x), \phi_{v_2g(v_1)}(y))\leq 1+6c_2+6Dc_2$.
Let $x_2=\phi_{u_2g(u_1)}(x)$, $y_2=\phi_{v_2g(v_1)}(y)$, $x_1=f'(x_2)=\bar{f}_{g(u_1)}(x_2)$ and
$y_1=f'(y_2)=\bar{f}_{g(v_1)}(y_2)$. Therefore, $d_{X_2}(x_2, y_2)\leq 1+6c_2+6Dc_2=R_2$, say and we want to show that $d_{X_1}(x_1, y_1)$ is
uniformly small.
Let $x'_2=f(x_1)=f_{u_1}(x_1), y'_2=f(y_1)=f_{v_1}(y_1)$. Then $d_{X_2}(x_2, x'_2)\leq D_1$ and $d_{X_2}(y_2, y'_2)\leq D_1$.
Hence, $d_{X_2}(x'_2, y'_2)\leq R_2+2D_1$.
Since $d_{B_1}(u_1, v_1)\leq K+C$ there is a point $y'_1\in F_{u_1}$ such that $d_{X_1}(x_1, y'_1)\leq (K+C)c_1+c_1$. Hence,
$d_{X_2}(x'_2,f(y'_1))\leq ((K+C)c_1+c_1).L+L$. Hence, 
$d_{X_2}(f(y'_1),y'_2)\leq d_{X_2}(f(y'_1),x'_2)+d_{X_2}(x'_2, y'_2)\leq ((K+C)c_1+c_1).L+L+2D_1+R_2$. This implies that
$d_{v_2}(f(y'_1), f(y_1))\leq \eta_2(((K+C)c_1+c_1).L+L+2D_1+R_2)=D_2$, say. Since $f_{v_1}$ is a $(\lambda, \epsilon)$-qi
we have $-\epsilon+\frac{1}{\lambda} d_{v_1}(y_1, y'_1)\leq D_2$. Hence, $d_{v_1}(y_1, y'_1)\leq (\epsilon+D_2)\lambda$.
Thus, $d_{X_1}(x_1, y_1)\leq d_{X_1}(x_1, y'_1)+d_{X_1}(y'_1, y_1)\leq (K+C)c_1+c_1+(\epsilon+D_2)\lambda$. 

(ii) We already know that $g'$ is a coarse inverse of $g$. Hence we will be done by checking that $f'$ is a coarse
inverse of $f$. We will check only that $d(f'\circ f, Id_{X_1})<\infty$ leaving the proof of 
$d(f\circ f', Id_{X_2})<\infty$ for the reader. Suppose $b\in B_1$ and $x\in \pi^{-1}_1(b)$. Then 
$f'(f(x))=\bar{f}_{g'\circ g(b)}\circ \phi_{g(b) g\circ g'(g(b))}\circ f_b(x)$. 
We want to show that $d_{X_1}(x, f'(f(x)))$ is uniformly small. Let $h=f_{g'\circ g(b)}\circ\bar{f}_{g'\circ g(b)}$.
Then $d_{X_2}(f(x),f(f'(f(x))))= d_{X_2}(f_b(x), h\circ \phi_{g(b) g\circ g'(g(b))}\circ f_b(x))\leq 
d_{X_2}(f_b(x), \phi_{g(b) g\circ g'(g(b))}(f_b(x)))+ d_{X_2}(\phi_{g(b) g\circ g'(g(b))}(f_b(x)), h\circ \phi_{g(b) g\circ g'(g(b))}\circ f_b(x))$. Now since, $d(g\circ g', Id_{B_2})\leq D$, $d_{X_2}(f_b(x), \phi_{g(b) g\circ g'(g(b))}(f_b(x)))\leq 3Dc_2+3c_2$.
Since $d(h, Id_{F_{g(g'(g(b)))}})\leq D_1$ we have
$d_{X_2}(\phi_{g(b) g\circ g'(g(b))}(f_b(x)), h\circ \phi_{g(b) g\circ g'(g(b))}\circ f_b(x))\leq D_1$.
Thus $d_{X_2}(f(x),f(f'(f(x))))\leq 3Dc_2+3c_2+D_1$. Hence, it is enough to show that $f$ is a proper embedding.
Here is how this is proved. Suppose $b,b'\in B$, $x\in \pi^{-1}_1(b)$ and $x'\in \pi^{-1}_1(b')$.
Suppose $d_{X_2}(f(x),f(x'))\leq N$ for some $N\geq 0$. This implies 
$d_{B_2}(g(b),g(b'))=d_{B_2}(\pi_2\circ f(x),\pi_2\circ f(x'))\leq N$. Since $g$ is a $(\lambda_1, \epsilon_1)$-qi
we have $-\epsilon_1+d_{B_1}(b,b')/\lambda_1\leq N$, i.e.  $d_{B_1}(b,b')\leq (N+\epsilon_1)\lambda_1=N_1$, say.
Hence by Corollary \ref{path lifting remark} there is a point $x''\in \pi^{-1}_1(b')$ such that
$d_{X_1}(x,x'')\leq 3N_1c_1+3c_1$. Since $f$ is coarsely $L$-Lipschitz we have $d_{X_2}(f(x),f(x''))\leq L(3N_1c_1+3c_1)+L$.
It follows that $d(f(x'),f(x''))\leq d(f(x'),f(x))+d(f(x), f(x''))\leq N+L(3N_1c_1+3c_1)+L=N_2$, say. 
Hence, $d_{g(b')}(f(x'),f(x''))\leq \eta_2(N_2)$. Since $f_{b'}$ is a $(\lambda,\epsilon)$-qi
we have $d_{X_1}(x',x'')\leq d_{b'}(x',x'')\leq \lambda(\epsilon+\eta_2(N_2))$.
Hence, $d_{X_1}(x,x')\leq d_{X_1}(x,x'')+d_{X_1}(x',x'')\leq 3N_1c_1+3c_1+ \lambda(\epsilon+\eta_2(N_2))$. This
completes the proof. \qed

\begin{defn}{\em (Subbundle)}
Suppose $(X_i, B,\pi_i)$, $i=1,2$ are metric (graph) bundles with the same base space $B$. 
We say that $(X_1, B,\pi_1)$ is subbundle of $(X_2, B,\pi_2)$ or simply
$X_1$ is a subbundle of $X_2$ if there is a metric (graph) bundle morphism $(f,g)$ from $(X_1, B,\pi_1)$ to $(X_2, B,\pi_2)$ such
that all the fiber maps $f_b$, $b\in B$ are uniform qi embeddings and $g$ is the identity map on $B$ (resp. on $\VV(B)$).
\end{defn}

The most important example of a subbundle that concerns us is that of ladders which we discuss in a later section.
The following gives another way to construct a metric (graph) bundle. We omit the proof since it is immediate.

\begin{lemma}{\em (Restriction bundle)}\label{restriction bundle}
Suppose $\pi:X\map B$ is a metric (graph) bundle and $A\subset B$ is a connected subset such that
any pair of points in $A$ can be joined by a path of finite length in $A$ (resp. $A$ is a connected subgraph).
Then the restriction of $\pi$ to $Y=\pi^{-1}(A)$ gives a metric (graph) bundle with the same
parameters as that of $\pi: X\map B$ where $A$ and $Y$ are
given the induced length metrics from $B$ and $X$ respectively.

Moreover, if $f:Y\map X$ and $g:A\map B$ are the inclusion maps then $(f,g):(Y,A)\map (X,B)$ is a morphism of metric
(graph) bundles.
\end{lemma}

\begin{defn}\label{pullback-defn} 
$(1)$ {\bf(Pullback of a metric bundle)} Given a metric bundle $(X,B,\pi)$ and a coarsely Lipschitz map $g:B_1\map B$
a pullback of $(X,B,\pi)$ under $g$ is a metric bundle $(X_1,B_1,\pi_1)$ together 
with a morphism $(f:X_1\map X, g:B_1\map B)$ such that the following universal property holds: Suppose $\pi_2:Y\map B_1$ is 
another metric bundle and $(f^Y, g)$ is a morphism from $Y$ to $X$. Then there is a coarsely unique morphism $(f', Id_{B_1})$ 
from $Y$ to $X_1$ making the following diagram commutative.

\begin{figure}[ht]
\centering
\begin{tikzpicture}[node distance=1.5cm,auto]
  \node (P) {$X_1$};
  \node (B) [right of=P] {$X$};
  \node (A) [below of=P] {$B_1$};
  \node (C) [below of=B] {$B$};
  \node (P1) [node distance=1cm, left of=P, above of=P] {$Y$};
  \draw[->] (P) to node {$f$} (B);
  \draw[->] (P) to node [swap] {$\pi_1$} (A);
  \draw[->] (A) to node [swap] {$g$} (C);
  \draw[->] (B) to node {$\pi$} (C);
  \draw[->, bend right] (P1) to node [swap] {$\pi_2$} (A);
  \draw[->, bend left] (P1) to node {$f^Y$} (B);
  \draw[->, dashed] (P1) to node {$f'$} (P);
\end{tikzpicture}
\caption{}\label{pullback defn figure}
\end{figure}

$(2)$ {\bf(Pullback of a metric graph bundle)}
In the case of a metric graph bundle, the diagram is replaced by one where we have the vertex sets instead of the whole spaces.
\end{defn}

The following lemma follows by a standard argument.

\begin{lemma}\label{pullback unique}
Suppose we have a metric bundle $(X,B,\pi)$ and a coarsely Lipschitz map $g:B_1\map B$ for which there are two pullbacks
i.e. metric bundles $(X_i,B_1,\pi_i)$ together with a morphisms $(f_i:X_i\map X, g:B_1\map B)$, $i=1,2$
satisfying the universal property of the Definition \ref{pullback-defn}. Then there is a coarsely unique metric (graph) bundle
isomorphism from $X_1$ to $X_2$.
\end{lemma}

With the above lemma in mind, in the context of Definition \ref{pullback-defn}, we say that 
$f:X_1\map X$ is the pullback of $X$ under $g:B_1\map B$ or simply $X_1$ is {\em the} pullback
of $X$ under $g$ when all the other maps are understood.

\begin{lemma}\label{pullback lemma}
Given $L\geq 0$ and functions $\phi_1, \phi_2:[0,\infty)\map [0,\infty)$ there is a function $\phi:[0,\infty)\map [0,\infty)$ such that
the following hold:

Suppose we have the following commutative diagram of maps between metric spaces satisfying the properties (1)-(3) below.

\begin{figure}[ht]
\centering
\begin{tikzpicture}[node distance=1.5cm,auto]
  \node (P) {$X_1$};
  \node (B) [right of=P] {$X$};
  \node (A) [below of=P] {$B_1$};
  \node (P1) [node distance=1cm, left of=P, above of=P] {$Y$};
  \draw[->] (P) to node {$f$} (B);
  \draw[->] (P) to node [swap] {$\pi_1$} (A);
  \draw[->, bend right] (P1) to node [swap] {$\pi_2$} (A);
  \draw[->, bend left] (P1) to node {$f^Y$} (B);
  \draw[->, dashed] (P1) to node {$f'$} (P);
\end{tikzpicture}\label{pic3.1}

\end{figure}

(1) All the maps (except possibly $f'$) are coarsely $L$-Lipschitz.

(2) If $d_{B_1}(b,b')\leq N$ then $Hd(\pi^{-1}_1(b), \pi^{-1}_1(b'))\leq \phi_1(N)$
for all $b,b'\in B_1$ and $N\in [0,\infty)$.

(3) The restrictions of $f$ on the fibers of $\pi_1$ are uniformly properly embedded
as measured by $\phi_2$.

Then $d_Y(y,y')\leq R$ implies $d_{X_1}(f'(y), f'(y'))\leq \phi(R)$ for all $y', y\in Y$ and $R\in [0,\infty)$. 
In particular, if $Y$ is a length space
or the vertex set of a connected metric graph with restricted metric then $f'$
is coarsely $\phi(1)$-Lipschitz.

Moreover, $f'$ is coarsely unique, i.e. there is a constant $D>0$ such that if 
$f'':Y\map X_1$ is another map making the above diagram commutative then $d(f', f'')\leq D$.
\end{lemma}
\proof Suppose $y,y'\in Y$ with $d_Y(y,y')\leq R$. Let $x=f'(y), x'=f'(y')$.
Then $d_{B_1}(\pi_1(x),\pi_1(x'))=d_{B_1}(\pi_2(y), \pi_2(y'))\leq LR+L$.
Let $b=\pi_2(y), b'=\pi_2(y')$. Then $Hd(\pi^{-1}_1(b), \pi^{-1}_1(b'))\leq \phi_1(LR+L)=R_1$, say.
Let $x'_1\in \pi^{-1}_1(b')$ be such that $d_{X_1}(x, x'_1)\leq R_1$.
Then $d_X(f(x), f(x'_1))\leq LR_1+L$. On the other hand $d_X(f(x),f(x'))=d_X(f^Y(y), f^Y(y'))\leq LR+L$.
By triangle inequality, we have $d_X(f(x'), f(x'_1))\leq LR+L+LR_1+L=2L+RL+R_1L$.
Hence, by the hypothesis (3) of the lemma $d_{X_1}(x', x'_1)\leq \phi_2(2L+RL+R_1L)$.
Thus $d_{X_1}(x, x')\leq d_{X_1}(x, x'_1)+d_{X_1}(x', x'_1)\leq R_1+ \phi_2(2L+RL+R_1L)$.
Hence, we may choose $\phi(t)=\phi_1(Lt+L)+\phi_2(2L+tL+L\phi_1(Lt+L))$. 

In case $Y$ is a length space or the vertex set of a connected metric graph it follows by Lemma \ref{proving Lipschitz} that $f'$ is coarsely $\phi(1)$-Lipschitz. 

Lastly, suppose $f'':Y\map X_1$ is another map making the diagram commutative. In particular we have
$f^Y=f\circ f'=f\circ f''$. Hence for all $y\in Y$ we have $f(f'(y))=f(f''(y))$. Since $\pi_1(f'(y))=\pi_1(f''(y))=\pi_2(y)$
by the hypothesis (3) of the lemma it follows that $d_{X_1}(f'(y), f''(y))\leq  \phi_2(0)$.
Hence $d(f', f'')\leq \phi_2(0)$. \qed

\begin{rem}\label{pullback lemma remark}
We note that the condition (2) of the lemma above holds in case $\pi_1:X_1\map B_1$ is a metric (graph) bundle.
\end{rem}

\begin{prop}{\bf (Pullbacks of metric bundles)}\label{pullback bundle prop}
Suppose $(X,B,\pi)$ is a metric bundle and $g:B_1\map B$ is a Lipschitz map. Then there is a pullback.

More precisely the following hold: Suppose $X_1$ is the set theoretic pullback with the induced length metric
from $X\times B_1$ and let $\pi_1:X_1\map B_1$ be the projection on the second coordinate and let $f:X_1\map X$
be the projection on the first coordinate.  Then (1) $\pi_1:X_1\map B_1$ is metric bundle and $f$ is a coarsely Lipschitz
map so that $(f,g)$ is a morphism from $X_1$ to $X$. (2) $f:X_1\map X$ is the metric bundle pullback of $X$ under
$g$. (3) All the fiber maps $f_b:\pi^{-1}_1(b)\map \pi^{-1}(g(b))$, $b\in B_1$ are isometries with respect to induced length metrics
from $X_1$ and $X$ respectively.
\end{prop}

\proof By definition  $X_1=\{(x,t)\in X\times B_1:g(t)=\pi(x)\}$. 
We put on it the induced length metric from $X\times B_1$. 
Let $\pi_1:X_1\map B_1$ be the restriction of the projection map $X\times B_1\map B_1$ to $X_1$. 
We first show that $X_1$ is a length space.
Suppose $g$ is $L$-Lipschitz. Let $(x,s),(y,t)\in X_1$. Let $\alpha$ be a rectifiable
path joining $s, t$ in $B_1$. Then $g\circ \alpha$ is a rectifiable path in $B$ of length at most $l(\alpha)L$. By Lemma \ref{lifting geodesics}
and Corollary \ref{path lifting remark} this path can be lifted to a rectifiable path in $X$ starting
from $x$ and ending at some point say $z$ in $F_t$ such that the length of the path is at most $3c+3cLl(\alpha)$. 
By construction this lift is contained in $X_1$. Finally we can join $(y,t), (z,t)$ by
a rectifiable path in $F_t$. This show that $(x,s)$ and $(y,t)$ can be joined in $X_1$ by a rectifiable path. 
This proves that $X_1$ is a length space. Now, since $\pi^{-1}_1(t)= \pi^{-1}(g(t))$ is uniformly properly embedded 
in $X$ for all $t\in B_1$ and $X$ is properly embedded in $X\times B_1$,
$\pi^{-1}_1(t)$ is uniformly properly embedded in $X_1$ for all $t\in B_1$. The same argument also shows that any path
in $B_1$ of length at most $1$ can be lifted
to a path of length at most $3c+3cL$ verifying the condition 2 of metric bundles.

Hence $(X_1,B_1,\pi_1)$ is a metric bundle. Let $f:X_1\map X$ be the restriction of the projection map $X\times B_1\map X$
to $X_1$. Clearly $f:X_1 \map X$ is a morphism of metric bundles. Finally, we check the universal property.
If there is a metric bundle $\pi_2:Y\map B_1$ and a morphism $(f^Y,g)$ from $Y$ to $X$ then there is a map $f':Y\map X_1$
making the diagram \ref{pullback defn figure} commutative since we are working with the set theoretic pullback. 
That $f'$ is a coarsely unique, coarsely Lipschitz map now follows from Lemma \ref{pullback lemma}. In fact, 
condition (2) of that lemma follows from Lemma \ref{fibers qi}(1) since $\pi_1:X_1\map B_1$ is a metric bundle and 
(3) follows because fibers of metric bundles are uniformly properly embedded and in this case 
the restriction of $f$, $\pi^{-1}_1(b)\map \pi^{-1}(g(b))\subset X$ is an isometry with respect to 
the induced path metric on  $\pi^{-1}_1(b)$ and $\pi^{-1}(g(b))$ for all $b\in B_1$.
\qed

\begin{cor}\label{pullback cor}
Suppose $(X,B,\pi)$ is a metric bundle and $g:B_1\map B$ is a Lipschitz map. Suppose $\pi_2:X_2\map B_1$
is an arbitrary metric bundle and $(f_2:X_2\map X, g)$ is a morphism of metric bundles. If $X_2$ is the 
pullback of $X$ under $g$ and $f_2:X_2\map X$ is the pullback map then for all $b\in B_1$ the fiber map
$(f_2)_b:\pi^{-1}_2(b)\map \pi^{-1}(g(b))$ is a uniform quasiisometry with respect to the induced
length metrics on the fibers of $\pi_2$ and $\pi$ respectively.
\end{cor}
\proof Suppose $X_1$ is the pullback of $X$ under $g$ as constructed in the proof of the proposition above.
Then the fiber maps $f_b:\pi^{-1}_1(b)\map \pi^{-1}(g(b))$ are isometries with respect to the induced
metrics on the fibers of $\pi_1$ and $\pi$ respectively. On the other hand by Lemma 
\ref{pullback unique} there is a coarsely unique metric bundle isomorphism $(h, Id)$ from $X_2$ to $X_1$
making the diagram \ref{pullback cor figure} below commutative.

\begin{figure}[ht]
\centering
\begin{tikzpicture}[node distance=1.5cm,auto]
  \node (P) {$X_1$};
  \node (B) [right of=P] {$X$};
  \node (A) [below of=P] {$B_1$};
  \node (C) [below of=B] {$B$};
  \node (P1) [node distance=1cm, left of=P, above of=P] {$X_2$};
  \draw[->] (P) to node {$f$} (B);
  \draw[->] (P) to node [swap] {$\pi_1$} (A);
  \draw[->] (A) to node [swap] {$g$} (C);
  \draw[->] (B) to node {$\pi$} (C);
  \draw[->, bend right] (P1) to node [swap] {$\pi_2$} (A);
  \draw[->, bend left] (P1) to node {$f_2$} (B);
  \draw[->, dashed] (P1) to node {$h$} (P);
\end{tikzpicture}
\caption{}\label{pullback cor figure}
\end{figure}
Now, by Theorem \ref{bundle isomorphism} the fiber maps $h_b:\pi^{-1}_2(b)\map \pi^{-1}_1(b)$ are 
uniform quasiisometries with respect to the induced length metrics on the fibers of $\pi_2$ and
$\pi_1$ respectively. Since $(f_2)_b=f_b\circ h_b$ for all $b\in B_1$ are done by Lemma
\ref{qi composition}(2). \qed

\begin{example} Suppose $(X,B,\pi)$ is a metric bundle and $B_1\subset B$ which is path connected and such that with 
respect to the path metric induced from $B$, $B_1$ is a length space. Let $X_1=\pi^{-1}(B_1)$ be endowed with the induced path metric
from $X$. Let $\pi_1:X_1\map B_1$ be the restriction of $\pi$ to $X_1$. Let $g:B_1\map B$ and $f:X_1\map X$ be the inclusion 
maps. It is clear that $(X_1,B_1,\pi_1)$ is a metric bundle and also that $X_1$ is the pullback of $g$. 
\end{example}

\begin{rem}
The notion of morphisms of metric bundles was implicit in the work of Whyte(\cite{whyte-bundles}).
Along the line of \cite{whyte-bundles}, one can define a more general notion of metric bundles by relaxing the hypothesis of 
length spaces. In that category of spaces, pullbacks should exist under any coarsely Lipschitz maps. 
However, we do not delve into it here.
\end{rem}

\begin{prop}{\bf (Pullbacks for metric graph bundles)}\label{pullback graph prop}
Suppose $(X,B,\pi)$ is an $\eta$-metric graph bundle, $B_1$ is a metric graph and $g:\VV(B_1)\map \VV(B)$ 
is a coarsely $L$-Lipschitz map for some constant $L\geq 1$. Then there is a pullback $\pi_1:X_1\map B_1$ of $g$
such that all the fiber maps $f_b:\pi^{-1}_1(b)\map \pi^{-1}(g(b))$, $b\in \VV(B_1)$ are isometries with respect to 
induced length metrics from $X_1$ and $X$ respectively.
\end{prop}
\proof The proof is a little long. Hence we break this into steps for the sake of clarity.

{\bf Step 1. Construction of $X_1$ and $\pi_1:X_1\map B_1$ and $f:\VV(X_1)\map \VV(X)$.}
We first construct a metric graph $X_1$, a candidate for the total space of the bundle. The vertex set of $X_1$ is the 
disjoint union of the vertex sets of $\pi^{-1}(g(b))$, $b\in \VV(B_1)$. There are two types of edges. First of all
for all $b\in \VV(B_1)$, we take all the edges appearing in $\pi^{-1}(g(b))$. In other words, the full subgraph $\pi^{-1}(g(b))$ 
is contained in $X_1$. Let us denote that by $F_{b}$. For all adjacent vertices $s,t\in B_1$ we introduce some other edges with 
one end point in $F_s$ and the other in $F_t$.
We note that $F_s, F_t\subset X_1$ are identical copies of $F_{g(s)}$ and $F_{g(t)}$ respectively. Let $f_s:F_s\map F_{g(s)}$ denote this identification.
Let $e$ be an edge joining $s,t$ and let $\alpha$ be a geodesic in $B$ joining $g(s), g(t)$. 
Now for each $x\in F_s$ we lift the path $\alpha$ starting from $f_s(x)$ isometrically by Lemma \ref{lifting geodesics}(1)
to say $\tilde{\alpha}$. For each such lift we join $x$ by 
an edge to $y\in V(F_t)$ if and only if $f_t(y)=\tilde{\alpha}(g(t))$. This completes the construction of $X_1$.
We note that $d_B(g(s), g(t))\leq 2L$ and hence $l(\tilde{\alpha})\leq 2L$ too. Now we define $f:\VV(X_1)\map \VV(X)$ by setting
$f(x)=f_{\pi_1(x)}(x)$ for all $x\in \VV(X_1)$. It is clear that this map is $2L$-Lipschitz.

{\bf Step 2. $\pi_1:X_1\map B_1$ is a metric graph bundle and $(f,g)$ is a morphism.}
We need to verify that the fibers are uniformly properly embedded in $X_1$ so that $X_1$ is a metric graph bundle.
Suppose $x,y\in F_s$ and $d_{X_1}(x,y)\leq D$. Let $\alpha$ be a (dotted) geodesic in $X_1$ joining $x,y$. Then
$f\circ \alpha$ is a (dotted) path of length at most $2LD$. Thus $d_X(f(x), f(y))\leq 2LD$. Since $X$ is an $\eta$-metric graph bundle
$d_{g(s)}(f(x),f(y))\leq \eta(2LD)$. Since $f$ is an isometry when restricted to $F_s$ we have $d_s(x,y)\leq \eta(2LD)$.
This proves that $X_1$ is a metric graph bundle over $B_1$. 

On the other hand, $f$ is $2L$-Lipschitz by step 1 and $g$ is coarsely $L$-Lipschitz by hypothesis. 
It is also clear that $\pi\circ f=g\circ \pi_1$ by the definition of $f$. 
Thus $(f,g)$ is a morphism of metric graph bundles from $X_1$ to $X$.

\begin{figure}[ht]
\centering
\begin{tikzpicture}[node distance=1.5cm,auto]
  \node (P) {$\VV(X_1)$};
  \node (B) [right of=P] {$\VV(X)$};
  \node (A) [below of=P] {$\VV(B_1)$};
  \node (C) [below of=B] {$\VV(B)$};
  \node (P1) [node distance=1cm, left of=P, above of=P] {$\VV(Y)$};
  \draw[->] (P) to node {$f$} (B);
  \draw[->] (P) to node [swap] {$\pi_1$} (A);
  \draw[->] (A) to node [swap] {$g$} (C);
  \draw[->] (B) to node {$\pi$} (C);
  \draw[->, bend right] (P1) to node [swap] {$\pi_2$} (A);
  \draw[->, bend left] (P1) to node {$f^Y$} (B);
  \draw[->, dashed] (P1) to node {$f'$} (P);
\end{tikzpicture}
\caption{}\label{pic3}
\end{figure}

{\bf Step 3. $X_1$ is a pullback.}
Now we check that $X_1$ is a pullback of $X$ under $g$. Suppose $\pi_2:Y\map B_1$ is a metric graph bundle and 
$(f^Y,g)$ is a morphism of metric graph bundles from $Y$ to $X$ where $f^Y$ is coarsely $L_1$-Lipschitz
 We need to find a coarsely unique, coarsely Lipschitz 
map $f':\VV(Y)\map \VV(X_1)$ such that $(f', Id)$ is a morphism from $Y$ to $X_1$ and the whole diagram 
\ref{pic3} is commutative where $Id:\VV(B_1)\map \VV(B_1)$ is the identity map.

{\bf The map $f'$:}
For all $s\in \VV(B_1)$ we define $f'$ on $\VV(\pi^{-1}_2(s))$ as the composition $f^{-1}_s\circ f^Y_s$.
Collectively these maps define $f'$. It is clear that $f'$ makes the whole diagram above commutative.

The rest of the argument follows from Lemma \ref{pullback lemma}.
In fact, condition (2) of that lemma follows from Lemma \ref{fibers qi}(1) since $\pi_1:X_1\map B_1$ is a metric graph bundle and 
(3) follows because fibers of metric graph bundles are uniformly properly embedded and in this case 
the restriction of $f$, $\pi^{-1}_1(b)\map \pi^{-1}(g(b))\subset X$ is an isometry with respect to 
the induced path metric on  $\pi^{-1}_1(b)$ and $\pi^{-1}(g(b))$ for all $b\in \VV(B_1)$.
\qed

\smallskip
The corollary below follows immediately from the proof of the above proposition.
\begin{cor} \label{cor: restriction bundle}
Suppose $\pi:X\map B$ is a metric graph bundle. 
Suppose $A$ is a connected subgraph of $B$. Let $g:A\map B$ denote the inclusion map. Let $X_A=\pi^{-1}(A)$,
$\pi_A$ be the restriction of $\pi$ and let $f:X_A\map X$ denote the inclusion map.
Then $X_A$ is the pullback of $X$ under $g$.
\end{cor}

The proof of the following corollary is similar to that of Corollary \ref{pullback cor} and hence
we omit the proof.
\begin{cor}\label{pullback graph cor}
Suppose $(X,B,\pi)$ is a metric graph bundle and $g:\VV(B_1)\map \VV(B)$ is a coarsely Lipschitz map. Suppose $\pi_2:X_2\map B_1$
is an arbitrary metric graph bundle and $(f_2:\VV(X_2)\map \VV(X), g)$ is a morphism of metric graph bundles. If $X_2$ is the 
pullback of $X$ under $g$ and $f_2:\VV(X_2)\map \VV(X)$ is the pullback map then for all $b\in \VV(B_1)$ the fiber map
$(f_2)_b:\VV(\pi^{-1}_2(b))\map \VV(\pi^{-1}(g(b)))$ is a uniform quasiisometry with respect to the induced
length metrics on the fibers of $\pi_2$ and $\pi$ respectively.
\end{cor}



\subsection{Some examples} 
In this section we discuss in detail two main sources of examples for metric graph bundles to which the
main theorem of this paper will be applied.

\subsubsection{Short exact sequence of groups}

\begin{example}\label{exact sequence example}
{\em Given a short exact sequence of finitely generated groups
$$1\map N\map G\stackrel{\pi}{\map} Q\map 1$$ we have a naturally associated metric graph bundle. This is the main
motivating example of metric graph bundles. We recall the definition from \cite[Example 1.8]{pranab-mahan} with a minor modification.

Suppose $H<Q$ is a finitely generated subgroup. Let $G_1=\pi^{-1}(H)$. We fix a generating set $S_N$ of $N$, a generating
set $S\supseteq S_N$ of $G$ such that $S$ contains a generating set $S_1$ of $G_1$, $S_N\subset S_1$ and $N\cap S=S_N$. 
Let $S_Q=\pi(S)\setminus \{1\}$
and $S_H=\pi(S_1)\setminus \{1\}$. Then we have a metric graph bundle $\pi:\Gamma(G, S)\map \Gamma(Q,S_Q)$.
Clearly $\Gamma(H, S_H)$ is a subgraph of $\Gamma(Q, S_Q)$ and $\Gamma(G_1, S_1)=\pi^{-1}(\Gamma(H,S_H))$.
Hence, by Corollary \ref{cor: restriction bundle} it follows that $\Gamma(G_1, S_1)$ is the pullback of $\Gamma(G,S)$ under the
inclusion $\Gamma(H, S_H)\hookrightarrow \Gamma(Q, S_Q)$.}
\end{example}

\subsubsection{Complexes of groups}\label{example: complex of groups}
For this example, we refer to \cite{haefliger-cplx}.
Suppose $\mathcal Y$ is a finite simplicial complex and $ \mbox{{\bf G}}(\mathcal Y)$
is a developable complex of groups defined over $\mathcal Y$. 
(See \cite[Definition 2.2]{haefliger-cplx}.) For any face $\sigma$ of $\mathcal Y$, let $K_{\sigma}$ be a 
$K(G_{\sigma},1)$-space. Then by \cite[Theorem 3.4.1]{haefliger-cplx} there
is a complex of spaces $p: \mathcal X\map \mathcal Y$ (compare with {\em good} complexes of spaces
due to Corson \cite{corson1}) which is a cellular aspherical realization (see \cite[Definition 3.3.4]{haefliger-cplx})
of the complex of groups $ \boldG(\mathcal Y)$ such that inverse image under $p$ of the
barycenter of each face $\sigma$ is $K_{\sigma}$. 
It follows from the construction of $\mathcal X$ that there is a continuous section $s$ of $p: \mathcal X\map \mathcal Y$
over the $1$-skeleton $\mathcal Y^{(1)}$ of $\mathcal Y$. We fix a maximal tree of $s(\mathcal Y^{(1)})$ and a base
vertex $v_0\in \mathcal Y^{(0)}$ in it. Let $G=\pi_1(\mathcal X, s(v_0))$. Thus for any $v\in \mathcal Y^{(0)}$
we have a natural injective homomorphism $\pi_1(X_v,s(v))\map G$. We identify the image of the same with $G_v$.
Next following Corson \cite{corson1} we take the universal cover
$\pi_{\XX}: \tilde{\mathcal X}\map \mathcal X$. We put a CW complex structure on $\tilde{\mathcal X}$
in the standard way so that $\pi_{\XX}$ is a cellular map. Then for all $y\in \mathcal Y$, we collapse each connected 
component of $(p\circ \pi_{\XX})^{-1}(y)$ to a point. Suppose $\BB$ is the quotient complex thus obtained
and let $q:\tilde{\XX}\map \BB$ be the quotient map. 
Then we note that there is a cellular map $\bar{\pi}_{\XX}:\BB\map \YY$ making the following diagram
commutative.

\begin{figure}[ht]
	\centering
	\begin{tikzpicture}[node distance=2cm,auto]
	\node (A) {$\tilde{\XX}$}; 
	\node (B) [right of= A] {$\BB$};
	\node (C) [below=1cm of A] {$\XX$};
	\node (D) [below=1cm of B] {$\YY$};
	\draw [->] (A) to node {$q$} (B);
	\draw [->] (A) to node [swap] {$\pi_{\XX}$} (C);
	\draw [->] (B) to node {$\bar{\pi}_{\XX}$} (D);
	\draw [->] (C) to node [swap]{$p$}(D);
	\end{tikzpicture}
	\caption{}\label{complex of groups}
\end{figure}

Now for our purpose, we shall also assume that
all the face groups $G_{\sigma}$ are finitely generated, the $0$-skeleton of each $K_{\sigma}$ 
is a point $x_{\sigma}$, the $1$-skeleton is a wedge of finitely many circles and the
developable complex of groups satisfies the qi condition as defined below.

\begin{defn} Suppose we have a developable complex of groups $(\GG, \YY)$.

(1) We say that it satisfies the {\em qi condition}
if for any faces $\sigma \subset \tau$ of $\YY$ the corresponding homomorphism $G_{\tau}\map G_{\sigma}$
is an isomorphism onto a finite index subgroup of $G_{\sigma}$.

(2) If all the face groups of $G_{\sigma}$ satisfies a group theoretic property $\mathcal P$ then we
shall say that $(\GG,\YY)$ is a developable complex of groups with property $\mathcal P$.
\end{defn}
For instance, we shall work in section 6 with the developable complexes of nonelementary hyperbolic groups.

However, we now aim to associate to the complex of groups a metric graph bundle as follows.
Let $X'=(p\circ \pi_{\XX})^{-1}(\YY^{(1)})^{(1)}$ and $B= \bar{\pi}_{\XX}^{-1}(\YY^{(1)})^{(1)}$
where we denote by $\mathcal Z^{(1)}$ the $1$-skeleton of any CW complex $\mathcal Z$. Now we construct a metric graph bundle
$\pi:X\map B$ as follows. For all $v\in \BB^{(0)}$ let $F_v:=q^{-1}(v)^{(1)}$. Suppose $v,w\in \BB^{(0)}$ are connected by an edge $e$.
We look at the subcomplex $\tilde{\XX}_{[v,w]}=q^{-1}([v,w])$. Let $v_0=\bar{\pi}_{\XX}(v), w_0=\bar{\pi}_{\XX}(w)$
and $e_0=\bar{\pi}_{\XX}(e)$. Then $\tilde{\XX}_{[v,w]}\subset \pi^{-1}_{\XX}(p^{-1}([v_0,w_0])$. However, we recall
from Haefliger \cite{haefliger-cplx} how $p^{-1}([v_0,w_0])\subset \XX$ is built from the spaces the $K_{v_0}, K_{w_0}$
and $K_{e_0}$. There are injective homomorphisms $G_{e_0}\map G_{v_0}, G_{e_0}\map G_{w_0}$. We choose cellular maps
$f_0: K_{e_0}\map K_{v_0}, f_1: K_{e_0}\map K_{w_0}$ such that the induced maps in the fundamental groups are those group
homomorphisms.  Then one glues $K_{e_0}\times [0,1]$ to $K_{v_0}\bigsqcup K_{w_0}$ by gluing $K_{e_0}\times \{0\}$
to $K_{v_0}$ and $K_{e_0}\times \{1\}$ to $K_{w_0}$ using the maps $f_0, f_1$ respectively. Let $m_0$
be the midpoint of $x_{e_0}\times [0,1]\subset p^{-1}([v_0,w_0])$ and let $m\in e$ be the midpoint of $e$. Then
through any $a\in q^{-1}(m)^{(0)}$ we lift  $x_{e_0}\times [0,1]$. The lift is a $1$-cell joining $a_v\in q^{-1}(v)^{(0)}$
to $a_w\in q^{-1}(w)^{(0)}$. Let us denote the map $a\mapsto a_v$ by $f_{e,v}$ and the map $a\mapsto a_w$ by $f_{e,w}$
 
\begin{lemma}\label{lemma: complex gps}
(1) The map $f_{e,v}: q^{-1}(m)^{(0)}\map q^{-1}(v)^{(0)}$ is uniformly coarsely surjective with respect to the
graph metric on $q^{-1}(m)^{(0)}, q^{-1}(v)^{(0)}$ coming from $q^{-1}(m)^{(1)}, q^{-1}(v)^{(1)}$ respectively.

(2) Similar statement holds for $f_{e,w}$.
\end{lemma}
\proof We will only prove (1) as the proof of (2) is similar. The group $G_{v}<G$ is isomorphic to $G_{v_0}$ and
$q^{-1}(v)$ is a universal cover of $K_{v_0}$ since the complex of groups is developable. The groups $G_v$ acts
properly discontinuously with quotient $K_{v_0}$. Since the action is cellular the action of $G_v$ on $q^{-1}(v)^{(1)}$
is simply transitive. Similarly the action of $G_m$ is simply transitive on $q^{-1}(m)^{(1)}$. We note that $G_m<G_v$
and the map $f_{e,v}$ is equivariant. It is also clear that $[G_v:G_m]=[G_{v_0}:G_{e_0}]$. Finally we note that
$q^{-1}(v)^{(1)}$ is naturally isometric to a Cayley graph of $G_v$ when $q^{-1}(v)^{(1)}$ is given graph metric
where each edge has length $1$. The lemma is immediate from this. \qed 

Let $R>0$ be such that $f_{e,v}$ is coarsely $R$-surjective for all $0$-cell $v$ and $1$-cell $e$ of $\BB$ where $e$ is incident on $v$.
Then we construct a graph $X$ from $X'$ by introducing new edges as follows. Given $v,w\in\BB^{(0)}$ connected by
an edge $e$ we join all $x\in q^{-1}(v)$ to $y\in q^{-1}(w)$ by an edge if there is $a\in q^{-1}(m_e)^{(0)}$ such that
$d(x,f_{e,v}(a))\leq R$, $d(y,f_{e,w}(a))\leq R$ where the distances are taken in the respective $1$-skeletons
of $q^{-1}(v)$ and $q^{-1}(w)$.

\begin{prop}\label{bundle from complx}
Suppose we identify $G$ as the group of deck transformation on the covering map $\pi_{\XX}:\tilde{\XX}\map \XX$.
Then we have the following:

(1) $G$ acts on $X$ and on $B$ through simplicial maps. The map $q$ is $G$-equivariant.

(2) The $G$-action is proper and cofinite on $X$ but it is only cofinite on $B$. Also $B/G$ is isomorphic to $\YY^{(1)}$.

(3) For all $v\in \YY^{(0)}$ and $\tilde{v}\in \bar{\pi}_{\XX}^{-1}(v)$, $G_{\tilde{v}}$ is 
a conjugate of $G_v$ in $G$. 

(4) The action of $G_{\tilde{v}}$ on $X_{\tilde{v}}=q^{-1}(\tilde{v})$ is proper and cocompact.
In fact the action on $V(X_{\tilde{v}})$ is transitive and on $E(X_{\tilde{v}})$ is cofinite. In particular
if the $G_v$ is hyperbolic for all $v\in \mathcal Y^{(0)}$ then for all $v\in \mathcal Y^{(0)}$ and
$\tilde{v}\in \bar{\pi}_{\XX}^{-1}(v)$,  $X_{\tilde{v}}$ is uniformly hyperbolic.

(5) $\pi:X\map B$ is a metric graph bundle.
\end{prop}
\proof The group $G$ acts through deck transformations of the covering map $\pi_{\XX}: \tilde{\XX}\map \XX$.
Hence it follows that $G$ permutes the connected components of $ (p\circ \pi_{\XX})^{-1}(y)$ for all
$y\in \YY$. The action is also simplicial. Hence, (1) follows from this. For (2) we note that the action of
$G$ on $X'$ is proper and cofinite. On the other hand, the inclusion map $X'\map X$ is a $G$-equivariant
quasiisometry by Lemma \ref{lemma: qi graph}.  Hence the $G$-action on $X$ is proper and cofinite.
Clearly, $B/G$ is isomorphic to $\YY^{(1)}$ whence the $G$-action on $B$ is cofinite.
(3) is a consequence of a basic
covering space argument using the $G$-equivariance of the map $q$. In (4) the properness follows from
the properness of the action of $G$ on $X'$. Cocompactness is due to the fact that $X'/G$ is finite.
The second part also follows from the nature of $K(G_v,1)$ used to construct $\XX$,
where $v=\bar{\pi}_{\XX}(\tilde{v})$. The last part follows from the second by Milnor-Scwarz lemma.
What remains is to prove (5). For all $\tilde{v}\in V(B)$, let $X_{\tilde{v}} =\pi^{-1}(\tilde{v})$.
Since $B/G$ is finite and the map $q$ is $G$-equivariant the $X_{\tilde{v}}$'s are uniformly properly embedded
in $X'$ iff for all $w\in V(B/G)$ there is one $\tilde{w}\in \bar{\pi}_{\XX}^{-1}(w)$ such that 
$X_{\tilde{w}}$ is uniformly properly embedded in $X'$. However, each inclusion $X_{\tilde{v}}\map X'$ is 
$G_{\tilde{v}}$-equivariant, the $G_{\tilde{v}}$ action on $X_{\tilde{v}}$ is proper and cocompact and
$G_{\tilde{v}}$ is a finitely generated subgroup of $G$. Since each finitely generated subgroup of a finitely
generated group is uniformly properly embedded it follows that $X_{\tilde{v}}$ is properly embedded in $X'$.
Since $X'$ is quasiisometric to $X$, it follows that $X_{\tilde{v}}$'s are properly embedded in $X$.
This verifies property (1) of metric graph bundles. Property (2) follows from Lemma \ref{lemma: complex gps}
and the construction of the new edges. \qed

\medskip
\noindent
{\bf Subcomplexes of groups}

In the above set-up we now assume further that we have a connected subcomplex $\YY_1\subset \YY$.
Let $\XX_1=p^{-1}(\YY_1)$. We shall assume that the base point $x_0\in \XX$ is contained in $\XX_1$ and a maximal tree
of $s(\YY^{(1)}_1)$ is chosen so that it is contained in the chosen maximal tree of $s(\YY^{(1)})$.
Suppose the inclusion $\XX_1\map \XX$ is $\pi_1$-injective. Then the restriction
$\boldG(\YY_1)$ of $\boldG(\YY)$ to $\YY_1$ is a developable complex of groups by \cite[Corollary 2.15]{bridson-haefliger}. 
Let $G_1=\pi_1(\XX_1, x_0)$. However, $\XX_1\map \YY_1$ is a complex of spaces which is a cellular aspherical realization
of the complex of groups $ \boldG(\mathcal Y_1)$. Hence, we can build a metric graph bundle $\pi_1:X_1\map B_1$
as described in Proposition \ref{bundle from complx}. 

In fact fixing a point $\tilde{x}_0\in \pi^{-1}_{\XX}(x_0)$ we may identify $G$ as the group of deck transformations
on $\tilde{\XX}$. Then $G_1$ stabilizes the connected component of $\pi^{-1}(\XX_1)$ containing $\tilde{x}_0$. Since
$\XX_1\map \XX$ is $\pi_1$-injective this connected component, say $\tilde{\XX}_1$, is a universal cover of $\XX_1$.
We set $B_1=q(\tilde{\XX}_1)\cap B$ and $X_1=\pi^{-1}(B_1)$. The following proposition records these in a nutshell.

\begin{prop}\label{complex: pullback}
Suppose $\YY$ is a finite connected simplicial complex and $\boldG(\YY)$ is a developable complex of groups with qi condition and 
with fundamental group $G$ and suppose $\YY_1$ is a connected subcomplex of $\YY$. Suppose $G_1$ is the fundamental group of $\boldG(\YY_1)$.
Suppose the inclusion $\boldG(\YY_1)\map \boldG(\YY)$ induces injective homomorphism $G_1\map G$.
 
Then there is a metric graph bundle $\pi:X\map B$, a connected subgraph $B_1\subset B$ such that the following hold:

(1) $G$ acts on $X$ and on $B$ through simplicial maps. The map $\pi$ is $G$-equivariant.
The action is proper and cofinite on $X$ but it is only cofinite on $B$. Also, there is a simplicial $G$-equivariant map
$B\map \YY^{(1)}$ with trivial action on $\YY^{(1)}$ inducing an isomorphism of graphs $B/G\map \YY^{(1)}$. 
The group $G_b<G$ is a conjugate of $G_{\bar{b}}$ in $G$ where $\bar{b}$ is the image of $b$ under the map $B\map \YY^{(1)}$. 
Also the $G_b$-action on $F_b$ is proper and cofinite for all $b\in V(B)$.

(2) Let $X_1=\pi^{-1}(B_1)$. Then $G_1$ stabilizes $X_1$ and the $G_1$-action on $X_1$ is proper and cofinite. 
Also the restriction of the map $B/G\map \YY^{(1)}$  to $B_1/G_1$ is an isomorphism of graphs $B_1/G_1\map \YY^{(1)}_1$.

\end{prop}
Later on we shall work with rather special subcomplexes of groups as defined below.

\begin{defn}
Suppose $\YY$ is a finite connected simplicial complex and $(\GG,\YY)$ is a developable complexes of groups with qi 
condition over $\YY$. We shall call a connected subcomplex $\YY_1\subset \YY$ a {\em good} subcomplex if the following hold:

(1) The induced natural homomorphism $\pi_1(\GG,\YY_1)\map \pi_1(\GG,\YY)$ is injective. Suppose the image is $G_1$.

(2) If $\pi:X\map B$ is a metric graph bundle obtained as in Proposition \ref{bundle from complx} from $(\GG,\YY)$
and $B_1\subset B$ is as in Proposition \ref{complex: pullback}. Then the inclusion $B_1\subset B$ is a qi embedding.
\end{defn}

We note that $X$ is quasiisometric to $G$ and $X_1$ is quasiisometric to $G_1$. Thus it follows that $B$ is quasiisometric to the
`coned-off' space a la Farb(\cite{farb-relhyp}) obtained from $G$ by coning off the cosets of the various face groups of $(\GG,\YY)$.
Similarly $B_1$ is obtained by coning off various cosets of the face groups of $(\GG, \YY_1)$. Thus condition (2) of
the above definition is intrinsic and independent of the particular metric graph bundle obtained from $(\GG,\YY)$.


\section{Geometry of metric bundles} 
In this section, we recall some results from \cite{pranab-mahan} and also add a few of our own which are going to be useful for
the proof of our main theorem in the next section. Especially some of the results which were stated for geodesic metric
spaces in \cite{pranab-mahan} but whose proofs require little adjustments to hold true for length spaces are mentioned here.

\subsection{Metric graph bundles arising from metric bundles}

An analogue of the following result is proved in \cite{pranab-mahan}(see Lemma 1.17 through Lemma 1.21 in \cite{pranab-mahan}). 
We give an independent and relatively simpler proof here. We also construct an approximating metric graph bundle morphism 
starting with a given metric bundle morphism. However, one disadvantage of our construction is that the metric graphs so 
obtained are never proper.

\begin{prop}\label{bundle vs graph bundle}
Suppose  $\pi':X'\map B'$ is an $(\eta,c)$-metric bundle. Then there is a metric graph bundle $\pi:X\map B$ along 
with quasiisometries $\psi_B:B'\map B$ and $\psi_X: X'\map X$ such that (1) $\pi\circ \psi_X=\psi_B\circ \pi'$ and 
(2) for all $b\in B'$ the map $\psi_X$ restricted to $\pi'^{-1}(b)$ is a  $(1,1)$-quasiisometry onto $\pi^{-1}(\psi_B(b))$.

Moreover, the maps $\psi_X, \psi_B$ have coarse inverses $\phi_X$, $\phi_B$ respectively making the following diagram 
commutative: 

\begin{figure}[ht]
	\centering
	\begin{tikzpicture}[node distance=2cm,auto]
	\node (A) {$X'$}; 
	\node (B) [right of=A] {$X$};
	\node (C) [below=1cm of A] {$B'$};
	\node (D) [below=1cm of B] {$B$};
	\draw [transform canvas={yshift=0.5ex},->] (A) to node [above] {$\psi_X$} (B);
	\draw [transform canvas={yshift=-0.5ex},->,dashed] (B) to node [below] [swap] {$\phi_X$} (A);
	\draw [->] (A) to node [swap] {$\pi'$} (C);
	\draw [->] (B) to node {$\pi$} (D);
	\draw [transform canvas={yshift=0.5ex},->] (C) to node [above] [swap]{$\psi_B$}(D);
	\draw [transform canvas={yshift=-0.5ex},->,dashed] (D) to node [below] {$\phi_B$} (C);
	\end{tikzpicture}
	\caption{}\label{pic4}
\end{figure}
\end{prop}

\proof (1) For the proof we use the construction of Lemma \ref{length space qi to graph}.
We shall briefly recall the construction of the spaces. We define $\VV(B)=B'$ and $s,t\in \VV(B)$ are connected by an 
edge if and only if $s\neq t$ and $d_{B'}(s,t)\leq 1$. This defines the graph. We also have a natural map 
$\psi_B: B'\map B$ which is just the inclusion map when $B'$ is identified with the vertex set of $B$. 
To define $X$, we take $\VV(X)=X'$. Edges are of two types. 

{\bf Type 1 edges:} For all $s\in B'$, $x,y\in \pi'^{-1}(s)$ are connected by an edge if and only if $d_s(x,y)\leq 1$. 

{\bf Type 2 edges:} If $s\neq t\in B'$, $x\in \pi'^{-1}(s)$ and $y\in \pi'^{-1}(t)$ then $x,y$ are connected by
an edge if and only if $d_{B'}(s,t)\leq 1$ and $d_{X'}(x,y)\leq c$. 

The map $\psi_X:X'\map X$ is defined as before to be the inclusion map.
By Lemma \ref{length space qi to graph} $\psi_B$ is a qi. We also note that $\pi\circ \psi_X=\psi_B\circ \pi'$. We need to 
verify that $\psi_X$ is a qi. For that, it is enough to produce Lipschitz coarse inverses $\phi_X$, $\phi_B$ as claimed 
in the second part of the proposition and then apply Lemma \ref{elem-lemma1} since it is clear that $\psi_X$ is 
$1$-Lipschitz. We first choose a coarse inverse $\phi_B$ of $\psi_B$ as follows.
On $\VV(B)$ it is simply the identity map. The interior of each edge is then sent to one of its end points. 
The map $\phi_X$ on $\VV(X)$ is also defined as the identity map. The interior of a type 1 edge is sent to one of its end points. 
Then interior of each type 2 edge $e=[x,y]$ is sent to one of the end points $x$ or $y$ according as the edge
$\pi(e)$ is mapped by $\phi_B$ to $\pi(x)$ or $\pi(y)$ respectively. It follows that the diagram in Figure \ref{pic4} 
commutes. We just need to check that $\phi_X$ is coarsely Lipschitz, since $\phi_B, \phi_X$ are inverses of $\psi_B, \psi_X$ 
respectively on a $1$-dense subset, they will be coarse inverse automatically.
However, by Lemma \ref{proving Lipschitz} it is enough to show that edges are mapped to small diameter sets. This is 
again clear. In fact, the image of an edge has diameter at most $c$. This proves the first part of the proposition.

(2)  This is immediate from the definition of $\psi_X$ and the construction in
\lemref{length space qi to graph}.

(3) Finally, we need to check that $(X,B,\pi)$ is a metric graph bundle. 
Let $s \in B$ and $x, y \in \pi^{-1}(s)$ such that $d_{X}(x,y) \leq M$ for some $M>0$. Since $\phi_X$ is a quasiisometry, 
$d_{X'}(x,y) \leq M'$, where $M'>0$ depends on $M$ and $\phi_X$. Since $\pi'^{-1}(\phi_B(s))$ is properly embedded in 
$X'$ as measured by $\eta$, we have $d_{\phi_B(s)}(x,y) \leq \eta(M')$. Now, using the above fact that $\pi'^{-1}(\phi_B(s))$ 
is $(1,1)$-quasiisometric to $\pi^{-1}(s)$, we have $d_s(x,y) \leq \eta(M')+1$. Hence, $\pi^{-1}(s)$ is uniformly properly 
embedded in $X$. Next we check the condition $(2)$ of Definition \ref{defn-mgbdl}. 
Suppose $s,t \in \VV(B)$ are adjacent vertices. Then, $d_{B'}(s,t) \leq 1$. Let $\alpha$ be a path in $B'$ joining $s,t$ 
with $l_{B'}(\alpha) \leq 1$. Then, for any $x \in \pi'^{-1}(s)$, $\alpha$ can be lifted to a path of length at most $c$, 
joining $x$ to some $y \in \pi'^{-1}(t)$.  Then there exists an edge joining $x$ and $y$ in $X$, which is a lift of the edge 
joining $s$ and $t$ in $B$.\qed

\begin{rem}
We shall refer to the metric graph bundle $X$ obtained from $X'$ as the {\em canonical metric graph bundle} associated
to the bundle $X$. Since we are working with length metric spaces some of the machinery of \cite{pranab-mahan} may not
apply directly. The above proposition then comes to the rescue. We sometimes modify our definitions suitably to make 
things work. Consequently, all the results proved for metric graph bundles have their close analogs in metric bundles.
We shall make this precise for instance in Proposition \ref{existence-qi-section} and Definition \ref{defn-ladder}. 
\end{rem}

\smallskip
\noindent
{\bf Approximating a metric bundle morphism} 

Suppose $\pi':X'\map B'$ is a metric bundle and $g:A'\map B'$ is a Lipschitz map. 
Suppose $Y'$ is the pullback of the bundle under the map $g$ as constructed in the proof of Proposition \ref{pullback bundle prop},
i.e. $Y'$ is also the set theoretic pullback.
Let $g^{*}\pi':Y'\map A'$ be the corresponding bundle projection map and $f:Y'\map X'$ be the pullback map. 
Suppose we use the recipe of the above proposition to construct metric graph bundles $\pi_X:X\map B$, $\pi_Y:Y\map A$ with 
quasiisometries $\psi_A:A'\map A$, $\psi_B:B'\map B$, $\psi_Y:Y'\map Y$ and $\psi_X:X'\map X$ such that 
$\pi_Y\circ \psi_Y=\psi_A\circ g^*\pi'$ and $\pi_X\circ \psi_X=\psi_B\circ \pi'$.

Suppose $\phi_X, \phi_B, \phi_Y, \phi_A$ are the coarse inverses (as constructed in the proposition above) of $\psi_X$, 
$\psi_B$, $\psi_Y$, and $\psi_A$ respectively. We then have a commutative diagram: 

\begin{figure}[ht]
	\centering
	\begin{tikzpicture}[node distance=2cm,auto]
	\node (A) {$Y$}; 
	\node (B) [right of=A] {$Y'$};
	\node (C) [right of=B] {$X'$};
	\node (D) [right of=C] {$X$};
	\node (E) [below=1cm of A] {$A$};
	\node (F) [below=1cm of B] {$A'$};
	\node (G) [below=1cm of C] {$B'$};
	\node (H) [below=1cm of D] {$B$};
	\draw [transform canvas={yshift=0.5ex},->] (B) to node [above] {$\psi_Y$} (A);
	\draw [transform canvas={yshift=-0.5ex},->,dashed] (A) to node [below] [swap] {$\phi_Y$} (B);
	\draw [->] (B) to node [above] [swap] {$f$} (C);
	\draw [transform canvas={yshift=0.5ex},->] (C) to node [above] {$\psi_X$} (D);
	\draw [transform canvas={yshift=-0.5ex},->,dashed] (D) to node [below] [swap] {$\phi_X$} (C);
	\draw [->] (A) to node [swap] {$\pi_Y$} (E);
	\draw [->] (B) to node [swap] {$g^{*}\pi'$} (F);
	\draw [->] (C) to node [swap] {$\pi'$} (G);
	\draw [->] (D) to node [swap] {$\pi_X$} (H);
	\draw [transform canvas={yshift=0.5ex},->] (B) to node [above] {$\psi_Y$} (A);
	\draw [transform canvas={yshift=-0.5ex},->] (F) to node [below] {$\psi_A$} (E);
	\draw [transform canvas={yshift=0.5ex},->,dashed] (E) to node [above] [swap] {$\phi_A$} (F);
	\draw [->] (F) to node [below] [swap] {$g$} (G);
	\draw [transform canvas={yshift=-0.5ex},->] (G) to node [below] {$\psi_B$} (H);
	\draw [transform canvas={yshift=0.5ex},->,dashed] (H) to node [above] [swap] {$\phi_B$} (G);
	\end{tikzpicture}
	\caption{}\label{pic5}
\end{figure}

Let $\bar{f}, \bar{g}$ denote the restrictions of $\psi_X\circ f\circ \phi_Y$ and  $\psi_B\circ g\circ \phi_A$
on the vertex sets of $Y$ and $A$ respectively.

\begin{prop}\label{lem:approx morphism}
(1) The pair of maps $(\bar{f}, \bar{g})$ gives a morphism of metric graph bundles from $Y$ to $X$. 

Moreover, if $Y'$ is the pullback of $X'$ under $g$ and $f$ is the pullback map then $Y$ is the pullback of 
$X$ under $\bar{g}$ and $\bar{f}$ is the pullback map.

(2) In case, $X', Y'$ are hyperbolic then $f$ admits the CT map if and only if so does $\bar{f}$.
\end{prop}

\proof (1) Since all the maps in consideration, i.e. $\psi_X, f, \phi_Y, \psi_B, g, \phi_A$ are coarsely Lipschitz
the maps $\bar{f}, \bar{g}$ are also coarsely Lipschitz by Lemma \ref{qi composition}(1).
It also follows that $\pi_X\circ \bar{f}=\bar{g}\circ \pi_Y$. Thus $(\bar{f}, \bar{g})$ is a morphism.

\begin{figure}[ht]
	\centering
	\begin{tikzpicture}[node distance=1.8cm,auto]
	\node (A) {$Y$}; 
	\node (B) [right of=A] {$Y'$};
	\node (C) [right of=B] {$X'$};
	\node (D) [right of=C] {$X$};
	\node (E) [below=1cm of A] {$A$};
	\node (F) [below=1cm of B] {$A'$};
	\node (G) [below=1cm of C] {$B'$};
	\node (H) [below=1cm of D] {$B$};
    \node (I) [above=1cm of A] {$Y_1$};
	\draw [transform canvas={yshift=0.5ex},->] (B) to node [above] {$\psi_Y$} (A);
	\draw [transform canvas={yshift=-0.5ex},->,dashed] (A) to node [below] [swap] {$\phi_Y$} (B);
	\draw [->] (B) to node [above] [swap] {$f$} (C);
	\draw [transform canvas={yshift=0.5ex},->] (C) to node [above] {$\psi_X$} (D);
	\draw [transform canvas={yshift=-0.5ex},->,dashed] (D) to node [below] [swap] {$\phi_X$} (C);
	\draw [->] (A) to node {$\pi_Y$} (E);
	\draw [->] (B) to node {$g^{*}\pi'$} (F);
	\draw [->] (C) to node {$\pi'$} (G);
	\draw [->] (D) to node {$\pi_X$} (H);
	\draw [transform canvas={yshift=0.5ex},->] (B) to node [above] {$\psi_Y$} (A);
	\draw [transform canvas={yshift=-0.5ex},->] (F) to node [below] {$\psi_A$} (E);
	\draw [transform canvas={yshift=0.5ex},->,dashed] (E) to node [above] [swap] {$\phi_A$} (F);
	\draw [->] (F) to node [below] [swap] {$g$} (G);
	\draw [transform canvas={yshift=-0.5ex},->] (G) to node [below] {$\psi_B$} (H);
	\draw [transform canvas={yshift=0.5ex},->,dashed] (H) to node [above] [swap] {$\phi_B$} (G);
    \draw [->, bend left] (I) to node [swap] {$f_1$} (D);
    \draw [->, bend right] (I) to node [swap] {$\pi_1$} (E);
    \draw [->] (I) to node  {$f_2$} (B);
	\end{tikzpicture}
	\caption{}\label{pic6}
\end{figure}

Suppose $Y'$ is a the pullback of $X'$ under $g$.
To show that $Y$ is the pullback of $X$ we need to verify the universal property. Suppose $\pi_1:Y_1\map A$ is
any metric bundle and $f_1:\VV(Y_1)\map \VV(X)$ is a coarsely Lipschitz map such that the pair $(f_1, \bar{g})$
is a morphism of metric graph bundles from $Y_1$ to $X$. We note that 
$\pi'\circ (\phi_X \circ f_1)=g\circ (\phi_A\circ \pi_1)$. Since $Y'$ is a set theoretic pullback there is a
unique map $f_2:\VV(Y_1)\map Y'$ making the whole diagram below commutative. 

Now, by Lemma \ref{qi composition}(1) the maps $\phi_X \circ f_1$ and $\phi_A\circ \pi_1$ are coarsely Lipschitz.
Hence, it follows by Lemma \ref{pullback lemma} and Remark \ref{pullback lemma remark} that the map $f_2$ is coarsely Lipschitz. 
Let $h=\psi_Y\circ f_2$. Then $h$ is coarsely Lipschitz by Lemma \ref{qi composition}(1) and we have $\bar{f}\circ h=f_1$ and
$\pi_Y\circ h= \pi_1$. Hence, $(h, Id_A)$ is a morphism from $Y_1$ to $Y$. Finally coarse uniqueness of $h$ follows
from Lemma \ref{pullback lemma}.

(2) This is a simple application of \lemref{CT properties}. \qed

\begin{comment}
Suppose $X', Y'$ are hyperbolic. Then $X$ and $Y$ are also hyperbolic by Lemma \ref{qi vs gromov hyp}. Suppose $f$ admits a CT map 
$\partial{f}:\partial{Y'} \to \partial{X'}$. Since $\psi_X$ and $\phi_Y$ are quasiisometries, by
\lemref{CT properties}(4) they admit Cannon-Thurston maps $\partial{\psi_X} : \partial{X'} \to \partial{X}$ and 
$\partial{\phi_Y} : \partial{Y} \to \partial{Y'}$ respectively which are both homeomorphisms. 
Then by \lemref{CT properties}(1) the map $\bar{f}$ admits a Cannon-Thurston map.
Conversely, if $\psi_X\circ f\circ\phi_Y$ admits a CT map, since $\phi_X$ and $\psi_Y$ are coarse inverses 
of $\psi_X$ and $\phi_Y$ respectively and $f$ is coarsely Lipschitz, $\phi_X\circ \psi_X\circ f\circ\phi_Y\circ \psi_Y$ 
and $f$ are at finite distance. So, by  $(2)$ of \lemref{CT properties}, 
$\partial({\phi_X\circ \psi_X\circ f\circ\phi_Y\circ \psi_Y}) = \partial{f}$. Thus, $f$ admits a CT map 
$\partial{f} : \partial{Y'} \to \partial{X'}$.
\end{comment}


\medskip

\subsection{Metric bundles with hyperbolic fibers}\label{section 4.2}

For the rest of this section we shall assume that all our metric (graph) bundles $\pi:X\map B$ have the following property:

($\ast$) {\em Each of the fibers $F_b$ , $b\in B$ (resp. $b \in \mathcal{ V}(B)$) is a $\delta'$-hyperbolic metric space with respect to the
path metric $d_b$ induced from $X$.} 

We will refer to this by saying that the metric (graph) bundle has uniformly hyperbolic fibers. 
Moreover, the following property is crucial for the existence of (global) qi sections.

($\ast\ast$) {\em There is $N\geq 0$ such that for all $b\in B$ the barycenter map $\phi_b : \partial^3 F_b \rightarrow F_b$ is coarsely $N$-surjective. } (Recall that barycenter maps were defined right after
Lemma \ref{defn: bary map}.) 

\begin{prop}\textup{(\cite[Proposition 2.10, Proposition 2.12]{pranab-mahan})}
{\bf Global qi sections for metric (graph) bundles:}\label{existence-qi-section} For all $\delta',c\geq 0, N\geq 0$ and 
$\eta:[0,\infty) \rightarrow [0,\infty)$ there exists $K_0=K_0(c,\eta,\delta',N)$ such that the following holds.

Suppose $p : X' \rightarrow B'$ is an $(\eta,c)$-metric bundle or an $\eta$-metric graph bundle satisfying $(\ast)$ and $(\ast\ast)$.
Then there is a $K_0$-qi section over $B'$ through each point of $X'$ (where we assume $c=1$ for the metric graph bundle).
\end{prop}
\begin{comment}
\proof We shall briefly indicate a proof for the metric bundle case assuming the proposition for metric graph bundles.
Suppose $X'$ is a metric bundle over $B'$ with the properties mentioned in the proposition and suppose $X\map B$ is
the canonical metric graph bundle associated to $X'$.
Since any length space is uniformly quasiisometric to a metric graph by Lemma \ref{length space qi to graph} and quasiisometries 
induce bijection of the boundaries of hyperbolic spaces, by Lemma \ref{CT properties}(4), it follows 
the metric graph bundle satisfies the same properties (1) and (2), i.e. the fibers are uniformly hyperbolic and the
barycenter maps are uniformly coarsely surjective. Hence by the existence of qi sections in a metric graph bundle through
any point $x\in X$ there is a uniform qi section $\SSS$ over $B$. Now, clearly $\phi_X(\SSS)$ is a uniform qi section
through $x$ in $X'$ where $\phi_X:X\map X'$ is as in Proposition \ref{bundle vs graph bundle} \qed
\end{comment}
\begin{convention}
(1) With the notation of Proposition \ref{bundle vs graph bundle}, we note that for any qi section $\SSS$ in $X$ over $B$,
$\phi_X(\SSS)=\SSS$ since $\phi_X$ is the identity map when restricted to $\VV(X)$. We shall
refer to it as a qi section of the metric graph bundle {\em transported} to the metric bundle.

(2) Whenever we talk about a $K$-qi section in a metric bundle we shall mean that it is the transport of 
a $K$-qi section contained in the associated canonical metric graph bundle.
\end{convention}

\begin{defn}\label{defn-ladder} \textup{(\cite[Definition 2.13]{pranab-mahan})}
Suppose $\Sigma_1$ and $\Sigma_2$ are two $K$-qi sections of the metric graph bundle $X$.
For each $b\in \mathcal{ V}(B)$ we join the points $\SSS_1\cap F_b$, 
$\SSS_2\cap F_b$ by a geodesic in $F_b$. We denote the union of these 
geodesics by $\LL(\SSS_1,\SSS_2)$, and call it a $K$-{\bf ladder} (formed by the sections 
$\SSS_1$ and $\SSS_2$).

For a metric bundle by a ladder, we will mean one transported from the canonical metric graph bundle associated to it
(by the canonical map $\phi_X$ as in Proposition \ref{bundle vs graph bundle}.)
\end{defn}

\begin{comment}
\begin{lemma}(\cite[Lemma 3.3, Lemma 3.5, Remark 3.6]{pranab-mahan})\label{connected ladder}
The $2K$-neighborhood of a $K$-ladder in a metric (graph) ladder is connected and with respect to 
the induced path metric it is uniformly
properly embedded.
\end{lemma}
\end{comment}

The following are the most crucial properties of a ladder summarized from \cite{pranab-mahan}.
 
\begin{prop}\label{ladders are qi embedded}
Given $K\geq 0$, $\delta\geq 0$ there are $C=C_{\ref{ladders are qi embedded}}(K)\geq 0$,
$R= R_{\ref{ladders are qi embedded}}(K)\geq 0$ and $K_{\ref{ladders are qi embedded}}(\delta, K)\geq 0$ 
such that the following holds: 

Suppose $\pi:X\map B$ is an $\eta$-metric graph bundle satisfying $(\ast)$.
 Suppose $\SSS_1, \SSS_2$ are two $K$-qi sections in $X$ and $\LL=\LL(\SSS_1, \SSS_2)$
is the ladder formed by them. Then the following hold.

(1) {\bf (Ladders are coarse Lipschitz retracts)}
There is a coarsely $C$-Lipschitz retraction $\pi_{\LL}:X\map \LL$ defined as follows:

For all $x\in X$ we define $\pi_{\LL}(x)$ to be a nearest point projection of $x$ in $F_{\pi(x)}$ on
$\LL\cap F_{\pi(x)}$.

(2) Given a $k$-qi section $\gamma$ in $X$ over a geodesic in $B$, $\pi_{\LL}(\gamma)$ is a $(C+2kC)$-qi section in 
$X$ contained in $\LL$ over the same geodesic in $B$.

(3) {\bf (QI sections in ladders)} If $X$ also satisfies $(\ast\ast)$ 
then through any point of $\LL$ there is $(1+2K)C$-qi section contained in $\LL$.

(4) {\bf (Quasiconvexity of ladders)}
The $R$-neighborhood of $\LL$ is (i) connected and (ii) uniformly qi embedded in $X$.

In particular if $X$ is $\delta$-hyperbolic then $\LL$ is $K_{\ref{ladders are qi embedded}}(\delta, K)$-quasiconvex in $X$.
\end{prop}

\proof (1) is stated as Theorem 3.2 in \cite{pranab-mahan}. (2), (3) are immediate from
(1) or one can refer to Lemma 3.1 of \cite{pranab-mahan}. (4) is proved in Lemma 3.6
in \cite{pranab-mahan} assuming $(\ast\ast)$. However, we briefly indicate the argument here without assuming $(\ast\ast)$.

4(i) Suppose $b,b'\in B$, $d_B(b,b')= 1$. Let $x\in \LL\cap F_b$. Then there is a point
$x'\in F_{b'}$ such that $d(x,x')= 1$. Hence, $d(\pi_{\LL}(x), \pi_{\LL}(x'))=d(x,\pi_{\LL}(x'))\leq 2C$. 
If we define $R=2C$ then clearly the $R$-neighborhood of $\LL$ is connected. 

4(ii) We first claim that the $N_R(\LL)=Y$ say, is also properly embedded in $X$. Suppose $x',y'\in Y$ 
with $d_X(x',y')\leq N$. Let $x, y\in \LL$ be such that $d(x,x')\leq R, d(y,y')\leq R$. Then 
$d(x,y)\leq 2R+N$. Hence, $d_B(\pi(x), \pi(y))\leq 2R+N$. Let $\alpha$ be a geodesic in $B$ joining 
$\pi(x), \pi(y)$. Then by Lemma \ref{lifting geodesics} there is a geodesic lift $\tilde{\alpha}$
of $\alpha$ starting from $x$. It follows that for all adjacent vertices $b_1,b_2\in \alpha$ we have
$d(\pi_{\LL}(\tilde{\alpha})(b_1),\pi_{\LL}(\tilde{\alpha})(b_2))\leq 2C$. Hence, the length of
$\pi_{\LL}(\tilde{\alpha})$ is at most $2C(2R+N)$.  
Hence, $d(y, \pi_{\LL}(\tilde{\alpha}(\pi(y)))\leq d(x,y)+ d(x,\pi_{\LL}(\tilde{\alpha}(\pi(y))))\leq  2R+N+ l(\pi_{\LL}(\tilde{\alpha}))\leq
2R+N+2C(2R+N)$. Hence, $d_{\pi(y)}(y, \pi_{\LL}(\tilde{\alpha}(\pi(y))))\leq \eta(2R+N+4CR+2CN)$. Since $\pi_{\LL}(\tilde{\alpha})\subset Y$,
$d_Y(x,y)\leq d_{\pi(y)}(y, \pi_{\LL}(\tilde{\alpha}(\pi(y))))+l(\pi_{\LL}(\tilde{\alpha}))\leq \eta(2R+N+4CR+2CN)+ 4CR+2CN$.
Hence, $d_Y(x',y')\leq 4CR+2CN+\eta(2R+N+4CR+2CN)$.

Finally we prove the qi embedding. Let $f(N)=\eta(2R+N+4CR+2CN)+ 4CR+2CN$ for all $N\in \NN$.
Given $x,y\in \LL$, $d_X(x,y)=n$ and a geodesic $\gamma:[0,n]\map X$ joining them. By the proof of (4)(i)
we have $d_{Y}(\pi_{\LL}(\gamma(i)), \pi_{\LL}(\gamma(i+1))\leq f(2C)$ for all $0\leq i\leq n-1$ whence 
$d_{\LL}(x,y)\leq nf(2C)= f(2C)d_X(x,y)$. Clearly $d_X(x,y)\leq d_{\LL}(x,y)$. This proves the qi embedded part. 

It follows that for all $x, y\in \LL$ a geodesic joining $x,y$ in $Y$ is a $(f(2C),0)$-quasigeodesic
in $X$. Since $X$ is $\delta$-hyperbolic stability of quasigeodesics implies that $\LL$ is uniformly quasiconvex.
In fact, we can take 
$K_{\ref{ladders are qi embedded}}(\delta, K)=R+D_{\ref{cor: stab-qg}}(\delta, f(2C),0)$. \qed

\begin{rem}
Part (3) and (4) are clearly also true for metric bundles  which satisfy the 
properties
$(\ast)$ and $(\ast\ast)$.  
\end{rem}

The following corollary is immediate.
\begin{cor}{\bf (Ladders form subbundles)}\label{cor: ladders subbundles}
Suppose $\pi:X\map B$ is an $\eta$-metric graph bundle satisfying $(\ast)$ and $(\ast\ast)$. 
Let $C, R$ be as in the previous proposition.
Suppose $\LL=\LL(\SSS_1, \SSS_2)$ is a $K$-ladder. Consider the metric graph $Z$ obtained from $\LL$ 
by introducing some extra edges as follows: Suppose $b,b'\in B$ are adjacent vertices then for all $x\in \LL\cap F_b$, 
$x'\in \LL\cap F_{b'}$ we join $x,x'$ by an edge if and only if $d_X(x,x')\leq C+2KC$. Let $\pi_Z:Z\map B$ be the 
simplicial map such that $\pi=\pi_Z$ on $\VV(Z)$ and the extra edges are mapped isometrically to edges of $B$.

Then $Z$ is a metric graph bundle and the natural map $Z\map X$ gives a subbundle of $X$ which is also a (uniform)
qi onto $N_R(\LL)$ and hence a (uniform) qi embedding in $X$.
\end{cor}

In the next section of the paper, we will exclusively deal with bundles $\pi:X\map B$ which are hyperbolic
satisfying $(\ast)$ and $(\ast\ast)$ and we will need to understand geodesics in $X$. Since ladders are quasiconvex
we look for quasigeodesics contained in ladders. The lemma below is the last technical piece of
information needed for that purpose. However, we need the following definitions for stating the lemma.

\begin{defn}
Suppose $X$ is a metric graph bundle over $B$ and suppose $\SSS_1, \SSS_2$ are any two qi sections.

(1) {\bf Neck of ladders} {\em (\cite[Definition 2.16]{pranab-mahan})}.
Suppose $R\geq 0$. 
Then the set $U_R(\SSS_1,\SSS_2)=\{b\in B:\,d_b(\SSS_1\cap F_b,\SSS_2\cap F_b)\leq R\}$ is called the $R$-{\bf neck} of 
the ladder $\LL(\SSS_1,\SSS_2)$.

For a metric bundle the $R$-neck of a ladder will be defined to be the one transported from the canonical 
metric graph bundle associated to it, i.e. the image under $\phi_B$.

(2) {\bf Girth of ladders} {\em (\cite[Definition 2.15]{pranab-mahan})}.
The quantity $\min\{d_b(\SSS_1\cap F_b, \SSS_2\cap F_b): b\in B\}$ is called the {\em girth} of the
ladder $\LL(\SSS_1,\SSS_2)$ and it will be denoted by $d_h(\SSS_1,\SSS_2)$.
\end{defn}


\begin{defn}\label{defn-flare}{\em (\cite[Definition 1.12]{pranab-mahan})}{\em ({\bf Flaring for metric graph bundles})}
Suppose $\pi:X\rightarrow B$ is a metric graph bundle.
We say that it satisfies a {\bf flaring condition} if for all   $k \geq 1$, there exist
  $\nu_k>1$ and  $n_k,M_k\in \mathbb N$ such that
the following holds:\\
Let $\gamma:[-n_k,n_k]\rightarrow B$ be a geodesic and let
$\tilde{\gamma_1}$ and $\tilde{\gamma_2}$ be two
$k$-qi lifts of $\gamma$ in $X$.
If $d_{\gamma(0)}(\tilde{\gamma_1}(0),\tilde{\gamma_2}(0))\geq M_k$,
then we have
\[\mbox{
{\small $\nu_k.d_{\gamma(0)}(\tilde{\gamma_1}(0),\tilde{\gamma_2}(0))\leq \mbox{max}\{d_{\gamma(n_k)}(\tilde{\gamma_1}(n_k),\tilde{\gamma_2}(n_k)),d_{\gamma(-n_k)}(\tilde{\gamma_1}(-n_k),\tilde{\gamma_2}(-n_k))\}$}}.
\]
\end{defn}
We note that existence of flaring in a metric graph bundle implies the existence of three functions $\nu_k, n_k, M_k$ of $k$
with the said property in the above definition. This is independent of the hypotheses about metric graph bundles and 
the conditions $(\ast)$ and $(\ast\ast)$ mentioned in the beginning of this subsection. This notion is motivated from the 
hallway flaring condition of Bestvina-Feighn (\cite{BF}).

\begin{defn}\label{flaring defn}{\em ({\bf Flaring for metric bundles})}
We shall say that a metric bundle $\pi:X\rightarrow B$ satisfies a $(\nu_k, M_k,n_k)$-{\em flaring} condition if the canonical 
metric graph bundle associated to it satisfies a $(\nu_k, M_k,n_k)$-flaring condition.
\end{defn}

\begin{rem}
(1) Since the base for a metric bundle need not be a geodesic metric space, it is not reasonable to use
\cite[Definition 1.12]{pranab-mahan} of flaring for metric bundles. However, one can formulate
analogous flaring of qi sections over uniform quasigeodesics in the base and then show that this is indeed equivalent 
to Definition \ref{flaring defn}. Since this discussion is not directly related to the rest of the paper we move it 
to the end of the paper and we include it as an appendix. See Lemma \ref{append part 1} and Lemma \ref{append part 2}.

(2) This definition of flaring for metric bundles is equivalent to \cite[Definition 1.12]{pranab-mahan} in the case of
geodesic metric bundles. In fact 
it follows from Lemma \ref{append part 1} and Lemma \ref{append part 2} that a geodesic metric bundle satisfies flaring 
as per \cite[Definition 1.12]{pranab-mahan} iff the canonical metric graph bundle associated to it also satisfies flaring. 
\end{rem}
The following lemma will be crucial for the next section of the paper.

\begin{lemma}\label{qc-level-set-new}{\bf (Quasiconvexity of necks of ladders, \cite[Lemma 2.18]{pranab-mahan})}
Let $X$ be an $\eta$-metric graph bundle over $B$
satisfying {\em $(\nu_k,M_k, n_k)$-flaring condition} for all $k\geq 1$. 
Then for all $c_1\geq 1$ and $R>1 $ there are constants $D_{\ref{qc-level-set-new}}=D_{\ref{qc-level-set-new}}(c_1,R)$
and $K_{\ref{qc-level-set-new}}= K_{\ref{qc-level-set-new}}(c_1)$ such that the following holds:\\
Suppose $\SSS_1,\SSS_2$ are two $c_1$-qi sections of $B$ in $X$ and let $L\geq max \{ M_{c_1}, d_h(\SSS_1, \SSS_2) \}$.
\begin{enumerate}
\item Let $\gamma:[t_0,t_1]\rightarrow B$ be a geodesic, $t_0,t_1\in \mathbb Z$, such that\\
a)  $d_{\gamma(t_0)}(\SSS_1\cap F_{\gamma(t_0)},\SSS_2\cap F_{\gamma(t_0)})=LR$.\\
b) $\gamma(t_1)\in U_L:=U_L(\SSS_1,\SSS_2)$ but for all $t\in [t_0,t_1)\cap \mathbb Z$, $\gamma(t)\not \in U_L$.\\
Then the length of $\gamma$ is at most $D_{\ref{qc-level-set-new}}(c_1,R)$.

\item For any $b_1,b_2\in U_L$ and any geodesic $[b_1,b_2]$ joining them in $B$, we have
$[b_1,b_2]\subset N_{K_{\ref{qc-level-set-new}}}(U_L)$. In particular, if $B$ is hyperbolic then $U_L$ is $K_{\ref{qc-level-set-new}}$-quasiconvex in $B$.

\item If $d_h(\SSS_1,\SSS_2) \geq M_{c_1}$ then the diameter of the set $U_L$ is at most 
$D'_{\ref{qc-level-set-new}}=D'_{\ref{qc-level-set-new}}(c_1,L)$.
\end{enumerate}
\end{lemma}
 Part (2) of the above lemma is slightly different from that of \cite[Lemma 2.18]{pranab-mahan} but the proof there actually showed this.
However, ladders with short necks to which Lemma \ref{qc-level-set-new} applies are given a special name:
\begin{defn}
{\bf (Small girth ladders)} Given two $K$-qi sections $\SSS_1, \SSS_2$ in a metric graph bundle 
 satisfying a flaring condition the ladder $\LL(\SSS_1, \SSS_2)$ is called 
a small girth ladder if $U_L(\SSS_1, \SSS_2)\neq \emptyset$ where $L= M_K$.
\end{defn}

\begin{rem}
Suppose $X'\map B'$ is a metric bundle and $X\map B$ is the canonical metric graph bundle associated to it. 
Suppose a flaring condition holds for $X$. This is the case for instance when $X$ or equivalently $X'$ is hyperbolic.
In such a case, a small girth ladder in $X'$ for us will be, by definition, 
the transport of a small girth ladder from $X$ under $\phi_X$ (as in Proposition \ref{bundle vs graph bundle}).

\end{rem}

We end this section with two simple lemmas. We note that flaring condition is not needed for these to hold.

\begin{lemma}\label{distance from qi section}
Given $D\geq 0, K\geq 1$ there is $R=R_{\ref{distance from qi section}}(D, K)$ such that the following holds.

Suppose $\SSS$ is $K$-qi section in $X$ and $x\in X$. Let $b=\pi(x)$. 
Then $d(x, \SSS)\geq D$ if $d_b(x, \SSS\cap F_b)\geq R$. 
\end{lemma}
\proof Suppose $y\in \SSS$ a nearest point from $x$. Let $\alpha \subset \SSS$ be the lift of a
geodesic $[b,\pi(y)]$ joining $b$ to $\pi(y)$ joining $y$ to $\SSS\cap F_b$. We note that
$d_B(b, \pi(y))\leq d(x,y)$. Hence, $d(y, \alpha(b))\leq Kd(x,y)+K$. Therefore,
$d(x,\alpha(b))\leq d(x,y)+d(y,\alpha(b))\leq (K+1)d(x,y)+K$. This implies
$d(x,y)\geq \frac{1}{K+1}d(x,\alpha(b))$ since all distances are integers in this case.
Now fibers of $X$ are properly embedded as measured by $\eta$. Thus if 
$d_b(x,\alpha(b))\geq \eta((K+1)D)$ then $d(x,y)\geq D$. Hence, we can take
$R=\eta(KD+D)$. \qed

The corollary below gives a relation between the girth of a ladder $\LL(\SSS_1,\SSS_2)$ and $d(\SSS_1, \SSS_2)$.
\begin{cor}\label{girth vs section distance}
Given $D\geq 0, K\geq 1$ there is an $R=R_{\ref{girth vs section distance}}(D, K)$ such that the following holds.\\
Suppose $\SSS_1, \SSS_2$ are two $K$-qi sections in $X$. Then
$d(\SSS_1, \SSS_2)\geq D$ if $U_R(\SSS_1,\SSS_2)=\emptyset$.
\end{cor}

The next lemma is a generalization of Lemma \ref{distance from qi section}.
Nevertheless we keep both of them since they are used many times in the next section.
\begin{lemma}\label{distance from ladder}
Given $K, D$ there is $R=R_{\ref{distance from ladder}}(K, D)$ such that the following holds.

Suppose $\SSS_1, \SSS_2$ are two $K$-qi sections in $X$ and $\LL=\LL(\SSS_1, \SSS_2)$.
Suppose $x\in X$ and $\pi(x)=b$. Then $d(x, \LL)\geq D$ if $d_b(x, \LL\cap F_b)\geq R$.
\end{lemma}
\proof Suppose $y\in \LL$ is a nearest point from $x$. Let $\alpha$ be a geodesic lift of
any geodesic $[b,\pi(y)]$ joining $b$ to $\pi(y)$ such that $\alpha$ joins $y$ to
$F_b$. Now $\pi_{\LL}(\alpha)$ is a $2C$-qi lift of $[b,\pi(y)]$ where $C= C_{\ref{ladders are qi embedded}}(K)$.
Thus $d(y,\pi_{\LL}(\alpha)(b))\leq 2Cd_B(b,\pi(y)+2C\leq 2Cd(x,y)+2C$. Hence,
$d(x,\LL\cap F_b)\leq d(x,y)+d(y,\pi_{\LL}(\alpha)(b))\leq (2C+1)d(x,y)+2C$. Therefore,
$d(x,y)\geq \frac{1}{2C+1} d(x,\LL\cap F_b)$. Hence, we can take
$R=\eta((2C+1)D)$. \qed

\section{Cannon-Thurston maps for pull-back bundles}
In this section, we prove the main result of the paper. Here is the set-up. From now on 
we suppose that $\pi:X\rightarrow B$ is  an $(\eta,c)$-metric bundle or an $\eta$-metric graph bundle 
satisfying the following hypotheses.

\begin{itemize}
\item[{\bf(H1)}] $B$ is a $\delta_0$-hyperbolic metric space.
\item[{\bf(H2)}] Each of the fibers $F_b$, $b\in B$ is a $\delta_0$-hyperbolic 
metric space with respect to the path metric induced from $X$. 
\item[{\bf(H3)}] The barycenter maps $\partial^3F_b\rightarrow F_b$, $b\in B$ (resp. $b\in \mathcal{V}(B)$) are $N_0$-coarsely
surjective for some constant $N_0$.
\item[{\bf(H4)}] The $(\nu_k, M_k,n_k)$-flaring condition is satisfied for all $k\geq 1$.
\end{itemize} 

\smallskip
The following theorem is the main result of \cite{pranab-mahan}:

\begin{theorem}{\em (\cite[Theorem 4.3 and Proposition 5.8]{pranab-mahan})}\label{mbdl thm}
If $\pi:X\map B$ is a geodesic metric bundle or a metric graph bundle satisfying $H1, H2, H3$ then 
$X$ is a hyperbolic metric space if and only if $X$ satisfies a flaring condition. 
\end{theorem}

\subsection{Proof of the main theorem}
We are now ready to state and prove the main theorem of the paper.

\begin{theorem}{\em ({\bf Main Theorem})}\label{CT for mgbdl}
Suppose $\pi:X\map B$ is a metric (graph) bundle satisfying the hypotheses H1, H2, H3, and H4. Suppose 
$g: A\map B$ is a Lipschitz $k$-qi embedding and suppose $p:Y\map A$ is the pullback bundle. 
Let $f:Y\map X$ be the pullback map. 

Then $Y$ is a hyperbolic metric space and the CT map exists for $f:Y\map X$.
\end{theorem}

\proof We first note that $X$ is hyperbolic. This follows from Theorem \ref{mbdl thm} if $X$ is
a metric graph bundle (or a geodesic metric bundle). In case $X$ is a (length) metric bundle one may first pass to the canonical metric graph
bundle associated to it, and then verify the hypotheses of Theorem \ref{mbdl thm} for it. 
In fact, if any metric bundle satisfies (H1), (H2), and (H3) then the canonical metric graph bundle associated to it
also has these properties with possibly different values of the respective parameters. Flaring condition (H4) follows from
Definition \ref{flaring defn}. It then follows
that the metric graph bundle is hyperbolic. Consequently, $X$ is hyperbolic by Proposition \ref{bundle vs graph bundle}.
We shall assume that $X$ is $\delta$-hyperbolic.
We begin with the following reductions: (1) {\em It is enough to prove the theorem only for metric graph bundles.}
Indeed this follows from Proposition \ref{lem:approx morphism}(2). So for the rest of the proof we shall assume that $\pi:X\map B$
is a metric graph bundle satisfying (H1), (H2), (H3), (H4).

Since we work with graphs from now, for the rest of the section by hyperbolicity we shall mean 
Rips hyperbolicity.

(2) {\em We may moreover assume that
$A$ is a connected subgraph, $g:A\map B$ is the inclusion map and $Y$ is the restriction bundle for
that inclusion. In particular, $f:Y\map X$ is the inclusion map and $Y=\pi^{-1}(A)$.} 

\noindent
Since $g:A\map B$ is a $k$-qi embedding and $B$ is $\delta_0$-hyperbolic, $g(A)$ is 
$D_{\ref{cor: stab-qg}}(\delta_0, k,k)$-quasiconvex in $B$. Let $A'$ be the $D_{\ref{cor: stab-qg}}(\delta_0, k,k)$-neighborhood
of $g(A)$ in $B$. Then clearly $A'$ is connected subgraph of $B$ and $g:A\map A'$ is a quasiisometry with respect to the induced path
metric on $A'$ from $B$. Clearly $A'$ is $(1,4D_{\ref{cor: stab-qg}}(\delta_0, k,k))$-qi embedded.
Let $\pi':X'=\pi^{-1}(A')\map A'$ be the restriction of $\pi$ on $X'$. Then $\pi':X'\map A'$
is a metric graph bundle by Lemma \ref{restriction bundle}. Also, we note that $(f,g):Y\map X'$ is a morphism
of metric graph bundles. By Corollary \ref{pullback graph cor} the fiber maps of the morphism $f:Y\map X'$ are
uniform quasiisometries and hence by Theorem \ref{bundle isomorphism} we see that $f:Y\map X'$ is an isomorphism of
metric graph bundles. Since (Rips) hyperbolicity of graphs is a qi invariant, 
we are reduced to proving hyperbolicity of $X'$  and also by Lemma \ref{CT properties}(1) we are reduced to
proving the existence of the CT map for the inclusion $X'\map X$. 

{\bf Hyperbolicity of $Y$}\\
$Y$ is hyperbolic by Remark $4.4$ of \cite{pranab-mahan}. In fact, by Theorem \ref{mbdl thm} it is enough to check
that flaring holds for the bundle $Y\map A$. This is a consequence of flaring of the bundle $\pi:X\map B$ and bounded flaring.


\begin{rem}(1)
The sole purpose of $(H3)$ is to have global uniform qi sections through every point of $X$
which is guaranteed by Proposition \ref{existence-qi-section}. 
For the rest of this section, we shall also assume the following.

{\em 
{\bf (H3$\,'$)} Through any point of $X$ there is a global $K_0$-qi section.
}

(2) Clearly $Y$ is an $\eta$-metric graph bundle over $A$ satisfying H2, H3. We shall assume
that $A$ is $\delta'_0$-hyperbolic. We shall also assume the bundle $Y$ satisfies a 
$(\nu'_k,M'_k, n'_k)$-flaring condition for all $k\geq 1$.
\end{rem}

{\bf Existence of CT map}\\
{\em Outline of the proof:} To prove the existence of the CT map we use Lemma \ref{mitra-lemma}. 
The different steps used in the proof are as follows.
(1) Given $y,y'\in Y$ first we define a uniform quasigeodesic $c(y,y')$ in $X$ joining $y,y'$.
This is extracted from \cite{pranab-mahan}. (2) In the next step we modify $c(y,y')$ to obtain a path $\bar{c}(y,y')$ in $Y$. (3) We then
check that these paths are uniform quasigeodesics in $Y$. (4) Finally we verify the condition of Lemma \ref{mitra-lemma} for the paths
$c(y,y')$ and $\bar{c}(y,y')$. Since $X,Y$ are hyperbolic metric spaces, stability of quasigeodesics and Lemma \ref{mitra-lemma} finishes the proof.
To maintain modularity of the arguments we state intermediate observations as lemma, proposition etc.

\begin{rem}
Although we assumed that $y,y' \in Y$ as is necessary for our proof, $c(y,y')$ as defined below is a uniform quasigeodesic
for all $y,y'\in X$ as it will follow from the proof. 

However, we would like to note that description of uniform quasigeodesics in a metric graph bundle with the above properties
H1-H4 is already contained in \cite{pranab-mahan}, e.g. see Proposition 3.4, and Proposition 3.14 of \cite{pranab-mahan}.
We make it more explicit with the help of Proposition \ref{hamenstadt}. 
\end{rem}

\smallskip
{\bf Step 1: Descriptions of the uniform quasigeodesic $c(y,y')$.}\\
The description of the paths and the proof that they are uniform quasigeodesics in $X$ is broken up
into three further substeps.
\smallskip

{\bf Step 1(a):} {\bf Choosing a ladder containing $y,y'$.}
We begin by choosing any two $K_0$-qi sections
$\SSS, \SSS'$  in $X$ containing $y,y'$ respectively. Let $\LL(\SSS, \SSS')$ be the ladder formed by them.
Throughout Step 1 we shall work with these  qi sections and ladder. The path $c(y,y')$ that we shall construct in
Step 1(c) will be contained in this ladder.

\smallskip
{\bf Step 1(b)}: {\bf Decomposition of the ladder into small girth ladders.}\\
We next choose finitely many qi sections in $\LL(\SSS, \SSS')$ after \cite[Proposition 3.14]{pranab-mahan}
in a way suitable for using Proposition \ref{hamenstadt}. This requires a little preparation.
We start with the following.

\begin{lemma}\label{lem: ladders cobdd}
For all $K\geq 1$ there is $D_{\ref{lem: ladders cobdd}}(K)$ such that
the following holds in $X$.

Suppose $\SSS_1, \SSS_2$ are two $K$-qi sections and $d_h(\SSS_1,\SSS_2)\geq M_K$. 
Then $\SSS_1, \SSS_2$ are $D_{\ref{lem: ladders cobdd}}(K)$-cobounded.
\end{lemma}

\proof We note that $\SSS_1, \SSS_2$ are $K'=D_{\ref{cor: stab-qg}}(\delta, K,K)$-quasiconvex in $X$.
Suppose $P:X\map \SSS_1$ is an $1$-approximate nearest point projection map and
the diameter of $P(\SSS_2)$ is bigger than $D=D_{\ref{cor: lip proj}}(\delta, K',1)$.
Then $d(\SSS_1, \SSS_2)\leq R=R_{\ref{cor: lip proj}}(\delta, K',1)$. If $x\in \SSS_2$
such that $d(x, \SSS_1)\leq R$ and $b=\pi(x)$ then 
$d_b(x, \SSS_1\cap F_b)\leq R_{\ref{distance from qi section}}(R, K')=\bar{R}$, say. Hence,
$\pi(P(\SSS_2))\subset U_{\bar{R}}(\SSS_1, \SSS_2)$. However, by Lemma \ref{qc-level-set-new}
the diameter of $U_{\bar{R}}(\SSS_1, \SSS_2)$ is at most $D'_{\ref{qc-level-set-new}}(K',\bar{R})$.
It follows that the diameter of $P(\SSS_2)$ is at most
$K+KD'_{\ref{qc-level-set-new}}(K',\bar{R})$. 
Hence we may choose $D_{\ref{lem: ladders cobdd}}(K)=\max\{D_{\ref{cor: lip proj}}(\delta, K',1),
K+KD'_{\ref{qc-level-set-new}}(K',\bar{R})\}$. \qed

\begin{lemma}\label{ladder coarse bisection}
Suppose $\SSS_1, \SSS_2$ are two $K$-qi sections and $\SSS\subset \LL(\SSS_1, \SSS_2)$ is
$K$-qi section. Then $\SSS$ coarsely uniformly bisects $\LL(\SSS_1, \SSS_2)$ into the subladders
$\LL(\SSS_1, \SSS)$ and $\LL(\SSS, \SSS_2)$. 
\end{lemma}

\proof First of all any ladder formed by $K$-qi sections is $K_{\ref{ladders are qi embedded}}(\delta, K)$-quasiconvex.
Let $K'=K_{\ref{ladders are qi embedded}}(\delta, K)$.
Let $k\geq 1$, and $x_i\in \SSS_i$, $i=1,2$ be any points. Let $\gamma_{x_1x_2}: I\map X$ be a $k$-quasigeodesic
joining them where $I$ is an interval. Then there are points $t_1, t_2\in I$ with $|t_1-t_2|\leq 1$
such that $\gamma_{x_1x_2}(t_1)\in N_{K'}(\LL(\SSS_1, \SSS))$ and $\gamma_{x_1x_2}(t_2)\in N_{K'}(\LL(\SSS, \SSS_2))$.
Let $y_1\in \LL(\SSS_1, \SSS)$ and $y_2\in \LL(\SSS, \SSS_2)$ be such that $d(y_i, \gamma_{x_1x_2}(t_i))\leq K'$,
$i=1,2$. We note that $d(\gamma_{x_1x_2}(t_1), \gamma_{x_1x_2}(t_2))\leq 2k$. Hence,
$d(y_1,y_2)\leq 2K'+2k$. Let $b=\pi(y_1)$. Then 
$d_b(y_1, \LL(\SSS, \SSS_2)\cap F_b)\leq R_{\ref{distance from ladder}}(K, 2K'+2k)$. This implies
$d_b(y_1, \SSS\cap F_b)\leq R_{\ref{distance from ladder}}(K, 2K'+2k)$. Thus
$d(\gamma_{x_1x_2}(t_1), \SSS)\leq K'+ R_{\ref{distance from ladder}}(K, 2K'+2k)$.
This proves the lemma.
\qed

\begin{lemma}\label{retricting sections}
If $\QQQ$ is a $K$-qi section in $X$ then $\QQQ\cap Y$ is a $K_{\ref{retricting sections}}(K)$-qi 
section of $A$ in $Y$.
\end{lemma}
\proof Suppose $s:B\map X$ is the $K$-qi embedding such that $s(B)=\QQQ$. Let $s$ also denote the restriction on $A$.
Since the bundle map $Y\map A$ is $1$-Lipschitz we have $d_A(u,v)\leq d_Y(s(u), s(v))$ for all $u,v \in A$.
Thus it is enough to show that $s:A\map Y$ is uniformly coarsely Lipschitz. Suppose $u,v\in A$ are adjacent
vertices. Then $d_X(s(u),s(v))\leq 2K$. Now, there is a vertex $x\in F_v$ adjacent to $s(u)\in F_u$.
Hence, $d_X(s(v), x)\leq 1+2K$. Therefore, $d_v(s(v), x)\leq \eta(1+2K)$. Hence,
$d_Y(s(u), s(v))\leq 1+\eta(1+2K)$. It follows that for all $u,v\in A$ we have
$d_Y(s(u), s(v))\leq (1+\eta(1+2K))d_A(u,v)$. Hence, we can take $K_{\ref{retricting sections}}(K)=1+\eta(1+2K)$.
\qed

The following corollary is proved exactly as Lemma \ref{lem: ladders cobdd}. Hence we omit the proof.
\begin{cor}\label{Y ladders cobdd}
For all $K\geq 1$ there is 
$D_{\ref{Y ladders cobdd}}(K)\geq 0$ such that the following holds. 

Suppose $\SSS_1, \SSS_2$ are two $K$-qi sections in $X$ and $d_h(\SSS_1,\SSS_2)\geq M_K$. 
Then $\SSS_1\cap Y, \SSS_2\cap Y$ are $D_{\ref{Y ladders cobdd}}(K)$-cobounded in $Y$.
\end{cor}
\smallskip

Before describing the decomposition of ladders
the following conclusions and notation on qi sections and ladders will be useful to record.
\begin{convention}

{\bf (C0)} We recall that $A$ is $k$-qi embedded in $B$. We let $k_0=D_{\ref{stab-qg}}(\delta_0, k ,k)$
so that $A$ is $k_0$-quasiconvex in $B$. Finally we assume that $Y$ is $\delta'$ hyperbolic. 

{\bf (C1)} Let $K_{i+1}=(1+2K_0)C_{\ref{ladders are qi embedded}}(K_i)$ for all $i\in \NN$ where $K_0$ is as in (H$3'$).
Therefore, through any point of a $K_i$-ladder in $X$, there is a $K_{i+1}$-qi section contained in the ladder.
Let $K'_i=K_{\ref{retricting sections}}(K_i)$.

{\bf (C2)} We let
$\lambda_i=\max\{D_{\ref{cor: stab-qg}}(\delta, K_i, K_i), K_{\ref{ladders are qi embedded}}(\delta, K_i),
D_{\ref{cor: stab-qg}}(\delta', K'_i, K'_i), K_{\ref{ladders are qi embedded}}(\delta', K'_i)\}$ 
so that any $K_i$-qi section $\QQQ  \subset X$ and any ladder $\LL\subset X$ formed by two $K_i$-qi sections in $X$ are
$\lambda_i$-quasiconvex in $X$ and moreover $\QQQ\cap Y$ and $\LL\cap Y$ are $\lambda_i$-quasiconvex in $Y$.

{\bf (C3)} If $\SSS_1, \SSS_2$ are two $K_i$-qi sections in $X$ and $d_h(\SSS_1, \SSS_2)\geq M_{K_i}$ then
they are $D_i$-cobounded in $X$, as are $\SSS_1\cap Y, \SSS_2\cap Y$ in $Y$
where $D_i=\max\{D_{\ref{lem: ladders cobdd}}(K_i), D_{\ref{Y ladders cobdd}}(K_i)\}$. 

{\bf (C4)} For each pair of $K_i$-qi sections $\SSS_1, \SSS_2$ in $X$ with 
$d_h(\SSS_1, \SSS_2)>r_i= \max\{ R_{\ref{girth vs section distance}}(2\lambda_i+1, K_i), 
R_{\ref{girth vs section distance}}(2\lambda_i+1, K'_i)\}$ we have $d_X(\SSS_1, \SSS_2)>2\lambda_i+1$ and
$d_Y(\SSS_1\cap Y, \SSS_2\cap Y)>2\lambda_i+1$.
\end{convention}


The following proposition is extracted from Proposition 3.14 of \cite{pranab-mahan}.
The various parts of this proposition are contained in the different steps of the proof 
of \cite[Proposition 3.14]{pranab-mahan}.

Let us fix a point $b_0\in A$ once and for all. Suppose  $\alpha:[0,l]\map F_{b_0}\cap \LL(\SSS, \SSS')$ is an isometry 
such that $\alpha(0)=\SSS\cap F_{b_0}$ and $\SSS'\cap F_{b_0}=\alpha(l)$.

\begin{prop}{\em (See \cite[Corollary 3.13 and Proposition 3.14]{pranab-mahan})}\label{ps mj ladder}
There is a constant $L_0$ such that for all $L\geq L_0$
there is a partition $0=t_0<t_1<\cdots <t_n=l$ of $[0,l]$ and $K_1$-qi sections
$\SSS_i$ passing through $\alpha(t_i)$, $0\leq i\leq n$ inside $\LL(\SSS, \SSS')$ such that the following holds.
\begin{enumerate}
\item $\SSS_0=\SSS, \SSS_n=\SSS'$.

\item For $0\leq i\leq n-2$, $\SSS_{i+1}\subset \LL(\SSS_i, \SSS')$.

\item For $0\leq i\leq n-2$ either (I) $d_h(\SSS_i, \SSS_{i+1})=L$, or (II) $d_h(\SSS_i, \SSS_{i+1})>L$ and
there is a $K_2$-qi section $\SSS'_i$ through $\alpha(t_{i+1}-1)$ inside $\LL(\SSS_i, \SSS_{i+1})$ such that
$d_h(\SSS_i, \SSS'_i)<C+CL$ where $C=C_{\ref{ladders are qi embedded}}(K_1)$.

\item $d_h(\SSS_{n-1}, \SSS_n)\leq L$.
\end{enumerate}

\end{prop}

However, we will need a slightly different decomposition of $\LL(\SSS,\SSS')$ than what is described here.
It is derived as the following corollary to the Proposition \ref{ps mj ladder}.

\begin{convention}
We shall fix $L=L_0+M_{K_3}+r_3$ and denote it by $R_0$ for the rest of the paper. 
Also we shall define $R_1=C+CR_0$ where $C=C_{\ref{ladders are qi embedded}}(K_1)$.
Thus we have the following.
\end{convention}

\begin{cor}\label{ladder subdivision}{\em (Decomposition of $\LL(\SSS, \SSS')$)}
There is a partition $0=t_0<t_1<\cdots <t_n=l$ of $[0,l]$ and $K_1$-qi sections
$\SSS_i$ passing through $\alpha(t_i)$, $0\leq i\leq n$ inside $\LL(\SSS, \SSS')$ such that the following holds.
\begin{enumerate}
\item $\SSS_0=\SSS, \SSS_n=\SSS'$.

\item For $0\leq i\leq n-2$, $\SSS_{i+1}\subset \LL(\SSS_i, \SSS')$.

\item For $0\leq i\leq n-2$ either (I) $d_h(\SSS_i, \SSS_{i+1})=R_0$, or (II) $d_h(\SSS_i, \SSS_{i+1})>R_0$ and
there is a $K_2$-qi section $\SSS'_i$ through $\alpha(t_{i+1}-1)$ inside $\LL(\SSS_i, \SSS_{i+1})$ such that
$d_h(\SSS_i, \SSS'_i)<R_1$.

In either case $d_X(\SSS_i, \SSS_{i+1})>2\lambda_1+1$ and $\SSS_i, \SSS_{i+1}$ are $D_1$-cobounded in $X$.

\item $d_h(\SSS_{n-1}, \SSS_n)\leq R_0$.
\end{enumerate}
\end{cor}

We note that the second part of (3) follows from (C1), (C2), (C3) above. However, a subladder
$\LL(\SSS_i, \SSS_{i+1})$ of $\LL(\SSS, \SSS')$ will be referred to as a {\bf type (I) subladder} or a {\bf type (II) subladder}
according as $d_h(\SSS_i, \SSS_{i+1})=R_0$ or $d_h(\SSS_i, \SSS_{i+1})>R_0$ respectively. 

\begin{rem}\label{arbit rem}
(1) We note that by the choice of $R_0,R_1$ it follows that 
$d_Y(\SSS_i\cap Y, \SSS_{i+1}\cap Y)>2\lambda_1+1$ and $\SSS_i\cap Y, \SSS_{i+1}\cap Y$ are $D_1$-cobounded in $Y$
for $0\leq i\leq n-2$.

(2) {\bf We shall use
$\SSS_i$ to mean qi sections in $\LL(\SSS,\SSS')$ exactly as in the corollary above for the rest of this section.}

(3) Finally we note that $\SSS_n, \SSS_{n-1}$ need not be cobounded in general and the same remark applies to
$\SSS_n\cap Y, \SSS_{n-1}\cap Y$.
\end{rem}

\begin{lemma}\label{applying hamenstadt prop}
Let $\Pi: \LL(\SSS,\SSS')\map [0,n]$ be any map that sends $\SSS_i$ to $i\in [0,n]\cap \mathbb Z$
and sends any point of $\LL(\SSS_i, \SSS_{i+1})\setminus \{\SSS_i\cup \SSS_{i+1}\}$ to a point 
in $(i,i+1)$.  Then the hypotheses of Proposition \ref{hamenstadt} are verified
for both $\Pi$ and its restriction $\LL(\SSS,\SSS')\cap Y\map [0,n]$.
\end{lemma} 
\proof For both $\Pi$ and its restriction to $\LL(\SSS,\SSS')\cap Y$, $(P\,0), (P\, 1)$ follow from (C2), 
$(P\, 2)$ follows from Lemma \ref{ladder coarse bisection}, and $(P\,3)$ follows from (C4). 
$(P\, 4)$ for $\Pi$ follows from (C3) and for the restriction of $\Pi$ to $\LL(\SSS,\SSS')\cap Y$
from Remark \ref{arbit rem}(1).\qed



{\bf Step 1(c):} {\bf Joining $y,y'$ inside $\LL(\SSS, \SSS')$.}
We now inductively define a finite sequence of points $y_i\in \SSS_i$, $0\leq i\leq n+1$ with $y_0=y, y_{n+1}=y'$
such that each $y_i$, $1\leq i\leq n$, is a uniform approximate nearest point projection of $y_{i-1}$ on $\SSS_i$ in $X$. 
We also define uniform quasigeodesics $\gamma_i$ in $X$ joining $y_i, y_{i+1}$. The
concatenation of these $\gamma_i$'s then forms a uniform quasigeodesic in $X$ joining $y,y'$ by Proposition \ref{hamenstadt}
and Lemma \ref{applying hamenstadt prop}.

{\em We define $\gamma_n$ to be the lift of $[\pi(y_n),\pi(y_{n+1})]$ in $\SSS'$.}

\smallskip
Suppose $y_0,\ldots, y_i$ and $\gamma_0,\ldots, \gamma_{i-1}$ are already constructed, $0\leq i\leq n-2$.
We next explain how to define 
$y_{i+1}$ and $\gamma_i$.

{\bf Case I.} Suppose $\LL_i= \LL(\SSS_{i},\SSS_{i+1})$ is of type (I) i.e. $d_h(\SSS_i, \SSS_{i+1})=R_0$
or $i=n-1$. 
Then, $U_{R_0}(\SSS_{i},\SSS_{i+1})$ is non-empty.
Let $u_i$ be a nearest point projection of $\pi(y_i)$ on $U_{R_0}(\SSS_{i},\SSS_{i+1})$. 
We define $y_{i+1}=\SSS_{i+1}\cap F_{u_i}$. Let $\alpha_i$ be the lift of
$[\pi(y_i),u_i]$ in $\SSS_i$, and let $\sigma_i$ be the subsegment of $F_{u_i}\cap \LL_i$ joining 
$\alpha_i(u_i)$ and $y_{i+1}$.
We define $\gamma_i$ to be the concatenation of $\alpha_i$ and $\sigma_i$. 
Then clearly $\gamma_i$ is a $(K_1+R_0)$-quasigeodesic in $X$. That $y_{i+1}$
is a uniform approximate nearest point projection of $y_i$ on $\SSS_{i+1}$ follows from the following lemma.

\begin{lemma}\label{step 1(c)-1}
Given $K\geq 1$ and $R\geq M_K$ there are constants $\epsilon_{\ref{step 1(c)-1}}(K,R)$ and $\epsilon'_{\ref{step 1(c)-1}}(K,R)$
such that the following holds.

Suppose $\QQQ_1, \QQQ_2$ are two $K$-qi sections and $d_h(\QQQ_1, \QQQ_2)\leq R$.
Let $x\in \QQQ_1$ and let $U=U_R(\QQQ_1, \QQQ_2)$. Suppose $b$ is a nearest point projection of $\pi(x)$ on $U$.
Then $\QQQ_2\cap F_b$ is $\epsilon_{\ref{step 1(c)-1}}(K,R)$-approximate nearest point projection of $x$
on $\QQQ_2$.

If $d_h(\QQQ_1, \QQQ_2)\geq M_K$ then for any $b'\in U$ the point $\QQQ_2\cap F_{b'}$ is an
$\epsilon'_{\ref{step 1(c)-1}}(K,R)$-approximate nearest point projection of any point of $\QQQ_1$
on $\QQQ_2$.
\end{lemma}
This lemma follows from  Corollary 1.40 and Proposition 3.4 of \cite{pranab-mahan}
given that ladders are quasiconvex. However, we give an independent proof using the hyperbolicity of $X$.

\proof Suppose $\bar{x}$ is a nearest point projection of $x$ on $\QQQ_2$ and let $x'= \QQQ_2\cap F_b$.
Let $\gamma_{xx'}$ be the
concatenation of the lift in $\QQQ_1$ of any geodesic in $B$ joining $\pi(x)$ to $b$ and any geodesic
in $F_b$ joining $\QQQ_1\cap F_b$ to $\QQQ_2\cap F_b$. Clearly it is a $(K+R)$-quasigeodesic in $X$. Also
by Lemma \ref{subqc-elem} the concatenation of any $1$-quasigeodesics joining $x, \bar{x}$ and
$\bar{x}, x'$ is a $K_{\ref{subqc-elem}}(\delta, K,1,0)$-quasigeodesic. Hence,
by stability of quasigeodesics we have $\bar{x}\in N_D(\gamma_i)$ where 
$D=D_{\ref{cor: stab-qg}}(\delta, K', K')$ and $K'=\max\{K+R,K_{\ref{subqc-elem}}(\delta, K,1,0)\}$. 
This implies there is a point $z\in \gamma_{xx'}$ such that $d(z, \bar{x})\leq  D$. 
If $z\in F_b\cap \gamma_{xx'}$ then $d(\bar{x}, x')\leq D+R$ and hence $x'$ is a 
$(D+R)$-approximate nearest point projection of $x$ on $\QQQ_2$. 

Suppose $z \in \QQQ_1\cap \gamma_{xx'}$. Then 
$d_{\pi(z)}(z, \QQQ_2\cap F_{\pi(z)}) \leq R_{\ref{distance from qi section}}(D, K)$.
Hence, by Lemma \ref{qc-level-set-new} we have
$d_B(\pi(z), b)\leq D_{\ref{qc-level-set-new}}(K,R')$ where $R'=R_{\ref{distance from qi section}}(D, K)/R$.
Therefore, $d(\bar{x}, x')\leq d(\bar{x}, z)+d(z, \QQQ_1\cap F_b)+d(\QQQ_1\cap F_b, x')
\leq D+(K+KD_{\ref{qc-level-set-new}}(K,R'))+R$. Hence in this case $x'$
is a $(D+K+KD_{\ref{qc-level-set-new}}(K,R')+R)$-approximate nearest point projection of $x$
on $\QQQ_2$. We may set $\epsilon_{\ref{step 1(c)-1}}(K,R)=D+K+KD_{\ref{qc-level-set-new}}(K,R')+R$.

For the last part, we note that the diameter of $U$ is at most $D'_{\ref{qc-level-set-new}}(K,R)$.
Thus clearly $\epsilon'_{\ref{step 1(c)-1}}(K,R)=\epsilon_{\ref{step 1(c)-1}}(K,R)+K+KD'_{\ref{qc-level-set-new}}(K,R)$
works. \qed

{\bf Case II.} Suppose $\LL_i= \LL(\SSS_{i},\SSS_{i+1})$ is of type (II), i.e. $d_h(\SSS_{i},\SSS_{i+1})>R_0$.
In this case there exists a $K_2$-qi section $\SSS'_i$ inside $\LL_i = \LL(\SSS_i, \SSS_{i+1})$ 
passing through $\alpha(t_{i+1}-1)$ such that $d_h(\SSS_i, \SSS'_i)\leq R_1$.
We thus use Case (I) twice as follows. First we project
$y_i$ on $\SSS'_i$. Suppose the projection is $y'_i$. Then we project $y'_i$ on $\SSS_{i+1}$
which we call $y_{i+1}$ and so on. Here are the details involved. 

Let $v_i$ be a nearest point projection of $\pi(y_i)$ on $U_{R_1}(\SSS_i, \SSS'_i)$
and let $w_i$ be a nearest point projection $v_i$ on $U_{R_1}(\SSS'_i, \SSS_{i+1})$.
Then $y_{i+1}=\SSS_{i+1}\cap F_{w_i}$. In this case we let $\alpha_i$ denote the
lift of $[\pi(y_i),v_i]$ in $\SSS_i$ and let $\beta_i$ denote the lift of $[v_i,w_i]$ in $\SSS'_i$.
Then $\gamma_i$ is the concatenation of the paths $\alpha_i$, 
$[\SSS_i\cap F_{v_i},\SSS'_i\cap F_{v_i}]_{v_i}$, $\beta_i$ and 
$[\SSS'_i\cap F_{w_i},\SSS_{i+1}\cap F_{w_i}]_{w_i}$. That $y_{i+1}$ is a
uniform approximate nearest point projection of $y_i$ on $\SSS_{i+1}$ and that
$\gamma_i$ is a uniform quasigeodesic follow immediately from Lemma \ref{step 1(c)-1}
and the last part of Proposition \ref{hamenstadt}. 

\begin{rem}
We note that $\LL(\SSS, \SSS')\cap Y$ is a ladder in $Y$ formed by the qi sections $\SSS\cap Y$ and 
$\SSS'\cap Y$ defined over $A$. 
However, in this case 
the subladders $\LL(\SSS_i, \SSS_{i+1})\cap Y$ may not be of type (I) or (II). 
Therefore, we cannot directly use the above procedure to construct a uniform quasigeodesic in $Y$ joining $y,y'$.
\end{rem}

{\bf Step 2: Modification of the path $c(y,y')$.} 

In this step we shall construct a path $\bar{c}(y,y')$ in $Y$ joining $y,y'$ by modifying $c(y,y')$.
For $0\leq i \leq n$, let $b_i$ be a nearest point projection of $\pi(y_i)$ on $A$ and let $\bar{y}_i=F_{b_i}\cap \SSS_i$.
We define a path $\bar{\gamma}_i \subset Y$ joining the points $\bar{y}_i, \bar{y}_{i+1}$ for $0\leq i\leq n$. 
Finally the path $\bar{c}(y,y')$ is defined to be the concatenation of these paths. {\em The path $\bar{\gamma}_n$ is the
lift of $[\pi(y_{n+1}), \pi(\bar{y}_n)]_A$ in $\SSS'\cap Y$.} The definition of $\bar{\gamma}_i$, for $0\leq i\leq n-1$, 
depends on the type of the subladder $\LL_i=\LL(\SSS_i, \SSS_{i+1})$ given by Corollary \ref{ladder subdivision}(3).

\smallskip
{\bf Case 2(I):}  Suppose $\LL_i$ is of type (I) or $i=n-1$. 
Let $\bar{\alpha}_i$ denote the lift of $[b_{i},b_{i+1}]_A$ in $\SSS_{i}$ starting at 
$\bar{y}_{i}$. The path $\bar{\gamma}_i$ is
defined to be the concatenation of $\bar{\alpha}_i$ and the fiber geodesic $F_{b_{i+1}}\cap \LL(\SSS_{i},\SSS_{i+1})$. 

\smallskip
{\bf Case 2(II):} Suppose $\LL_i$ is of type (II). 
In this case, we apply Case 2(I) to each of the subladders $\LL(\SSS_i, \SSS'_i)$ and $\LL(\SSS'_i, \SSS_{i+1})$. 
Let  $y'_i$ be as defined in step 1(c). 
Let $b'_i\in A$ be a nearest point projection $\pi(y'_i)$ on $A$ and 
$\bar{y}'_i=\pi^{-1}({b}'_i)\cap \SSS'_i$.
Next we connect $\bar{y}_i, \bar{y}'_i$ and $\bar{y}'_i, \bar{y}_{i+1}$ as in Case 2(I) inside the ladders
$\LL(\SSS_i\cap Y,\SSS'_i\cap Y)$ and $\LL(\SSS'_i\cap Y, \SSS_{i+1}\cap Y)$ respectively. 
We shall denote by $\bar{\alpha}_i$ and $\bar{\beta}_i$ the lift of $[b_i,b'_i]_A$ in $\SSS_i\cap Y$
and $[b'_i,b_{i+1}]_A$ in $\SSS'_i\cap Y$ respectively. 
The concatenation of the paths $\bar{\alpha}_i$, 
$[\SSS_i\cap F_{b'_i}, \SSS'_i\cap F_{b'_i}]_{b'_i}\subset \LL(\SSS,\SSS')$, $\bar{\beta}_i$ and 
$[\SSS'_i\cap F_{b_{i+1}}, \SSS_{i+1}\cap F_{b_{i+1}}]_{b_{i+1}}\subset \LL(\SSS,\SSS')$ is 
defined to be $\bar{\gamma}_i$.

\smallskip

{\bf Step 3: Proving that $\bar{c}(y,y')$ is a uniform quasigeodesic in $Y$.}
To show that $\bar{c}(y,y')$ is a quasigeodesic it is enough, by Proposition \ref{hamenstadt}, to show
that the paths $\bar{\gamma}_i$ are all uniform quasigeodesics in $Y$ and that for $0\leq i\leq n-1$, 
$\bar{y}_{i+1}$ is an approximate nearest point projection of $\bar{y}_i$ in $\SSS_{i+1}\cap Y$. 
The proof of this is broken into three cases depending on the type of the ladder $\LL_i$.  
We start with the following lemma as a preparation for the proof.

The lemma below is true for any metric bundle that satisfies the hypotheses (H1)-(H4), (H$3'$)
although we are stating it for $X$ only. For instance, it is true for $Y$ too.

\begin{lemma}\label{Y-path 1.1}
Suppose $b\in B$, $x,y\in F_b$. Suppose for all $K\geq K_0$ and $R\geq M_K$
 there is a constant $D=D(K, R)\geq 0$ such that 
for all $x',y'\in [x,y]_b$ and any two $K$-qi sections $\mathcal Q_1$ and $\mathcal Q_2$ in $X$ passing through 
$x',y'$ respectively, either $U_{R}(\mathcal Q_1, \mathcal Q_2)=\emptyset$ or 
$d_B(b, U_{R}(\mathcal Q_1, \mathcal Q_2))\leq D$. Then the following hold:
\begin{enumerate}
\item $[x,y]_b$ is a $\lambda_{\ref{Y-path 1.1}}$-quasigeodesic in $X$ where 
$\lambda_{\ref{Y-path 1.1}}$ depends on the function $D$ (and the parameters of the metric bundle).
\item If $\mathcal Q$ and $\mathcal Q'$ are two $K$-qi sections passing through $x,y$ respectively 
then 
$x$ is a uniform approximate nearest point projection of $y$ 
on $\QQQ$ and $y$ is a uniform approximate nearest point projection of $x$ 
on $\QQQ'$. 
\end{enumerate}
\end{lemma}
\proof (1) Since the arc length parametrization of $[x,y]_b$ is a uniform proper embedding, 
by Lemma \ref{quasigeod criteria} it is enough to show that $[x,y]_b$ is uniformly close to a geodesic 
in $X$ joining $x,y$. 

{\bf Claim:} {\em Suppose $\SSS_x, \SSS_y$ are two $K_0$-qi sections passing through $x,y$ respectively.
Given any $z\in [x,y]_b$ and any $K_1$-qi section $\SSS_z$ passing through $z$ contained in the
ladder $\LL(\SSS_x, \SSS_y)$ the nearest point projection of $x$ on $\SSS_z$ is uniformly close to $z$.}

We note that once the claim is proved then applying Proposition \ref{hamenstadt} to the ladder 
$\LL(\SSS_x, \SSS_y)=\LL(\SSS_x, \SSS_z)\cup \LL(\SSS_z, \SSS_y)$ it follows that 
$z$ is uniformly close to a geodesic joining $x,y$. From this (1) follows immediately.

{\em Proof of the claim:} 
First suppose 
$U_{M_{K_1}}(\SSS_x, \SSS_z)\neq \emptyset$. Then we can find a uniform approximate 
nearest point projection of $x$ on $\SSS_z$ using Step 1(c), Case I and Lemma \ref{step 1(c)-1} above which 
is uniformly close to $z$ by hypothesis.

Now suppose $U_{M_{K_1}}(\SSS_x, \SSS_z)=\emptyset$. Let $\alpha_{zx}:[0,l]\map F_b$ be the unit speed 
parametrization of the geodesic $\LL(\SSS_x,\SSS_z)\cap F_b$ joining $z$ to $x$. 
By Corollary \ref{ladder subdivision} there is a $K_2$-qi section 
$\SSS_{z'}$ contained in the ladder $\LL(\SSS_x, \SSS_z)$ passing through $z'=\alpha_{zx}(t)$ 
for some $t\in [0,l]$ such that $\LL(\SSS_z, \SSS_{z'})$ is a $K_2$-ladder of type (I) or (II). 
Let $x'$ be a nearest point projection of $x$ on $\SSS_{z'}$. 
By the last part of Proposition \ref{hamenstadt} applied to $\LL(\SSS_x, \SSS_z)$,
it is enough to find a uniform approximate nearest point projection of $x'$ on $\SSS_z$ which is also uniformly
close to $z$. However, in this case $\SSS_z, \SSS_{z'}$ are $D_2$-cobounded. Hence
it is enough to find a uniform approximate nearest point projection of $z'$ on $\SSS_z$ which
is uniformly close to $z$. The proof of this is broken into two cases as follows.

(I) Suppose $d_h(\SSS_z, \SSS_{z'})=R_0$.  
By the last part of Lemma \ref{step 1(c)-1} if $v\in U_{R_0}(\SSS_z, \SSS_{z'})$ then
$F_v\cap \SSS_z$ is a uniform approximate nearest point projection of any point of $\SSS_{z'}$.
Since $d_A(b,v)$ is uniformly small by hypothesis, $d(z,F_v\cap \SSS_z)$ is also uniformly small. 

(II) Suppose $d_h(\SSS_z, \SSS_{z'})>R_0$. Then there is a $K_3$-qi section $\SSS_{z''}$ in $\LL(\SSS_z, \SSS_{z'})$
passing through $z''=\alpha_{zx}(t-1)$ such that $U_{R_0}(\SSS_z, \SSS_{z''}) \neq \emptyset$.
Let $v'$ be a nearest point projection of $b$ on $U_{R_0}(\SSS_z, \SSS_{z''})$. Then by hypothesis
$d(b, v')$ is uniformly small whence $d(z, F_{v'}\cap \SSS_z)$ is uniformly small.
 Also by Lemma \ref{step 1(c)-1} the point $\SSS_z\cap F_{v'}$ is
a uniform approximate nearest point projection of $z''$ on $\SSS_z$. It follows that $z$ is a uniform
approximate nearest point projection of $z''$ on $\SSS_z$. Finally,
since $d(z', z'')\leq 1$, $z$ is a uniform approximate
nearest point projection of $z'$. 

(2) We shall prove only the first statement since the proof of the second would be an
exact copy. Suppose $x_1\in \QQQ$ is a nearest point projection of $y$ on $\QQQ$.
Consider the $K$-qi section over $[b,\pi(x_1)]$ contained in $\QQQ$. This is a $K$-quasigeodesic
of $X$ joining $x, x_1$. Since $\QQQ$ is a $K$-qi section, by stability of quasigeodesics it is 
$D_{\ref{stab-qg}}(\delta, K,K)$-quasiconvex in $X$. Hence by Lemma \ref{subqc-elem}
the concatenation of this quasigeodesic with a geodesic in $X$ joining $y$ to $x_1$
is a $K_{\ref{subqc-elem}}(\delta, \tilde{K}, K,0)$-quasigeodesic where
$\tilde{K}=D_{\ref{stab-qg}}(\delta, K,K)$. Let 
$k'= \max \{\tilde{K},\lambda_{\ref{Y-path 1.1}}\}$.
Since $[x,y]_b$ is a $\lambda_{\ref{Y-path 1.1}}$-quasigeodesic, by stability of quasigeodesics we have
$x_1\in N_{2D'}([x,y]_b)$ where $D'=D_{\ref{stab-qg}}(\delta, k',k')$. Suppose $z\in [x,y]_b$ be such that
$d(x_1, z)\leq 2D'$. Then $d_B(\pi(x_1), \pi(z))=d_B(\pi(x_1), b)\leq 2D'$.
Hence, $d(x, x_1)\leq K+2D'K$. Thus $x$ is a $(K+2D'K)$-approximate nearest
point projection of $y$ on $\QQQ$. 
\qed

\begin{rem}
The proof of the first part of the above lemma uses the hypothesis for $K\leq K_3$ only whereas
the proof of the second part follows directly from the statement of the first part
and is independent of the hypotheses of the lemma.
\end{rem}

The following lemma is actually a trivial consequence of flaring ( Lemma \ref{distance from qi section})
and it is going to be used in the next two lemmas following it. 

\begin{lemma}\label{Y-path 1.2}
Given $R\geq 0,K,K'\geq 1$ and $R'\geq M_{K'}$ there is a constant $R_{\ref{Y-path 1.2}}(R,R',K,K')$ 
and $D_{\ref{Y-path 1.2}}(R,R',K,K')$  such that the following holds.

Suppose $u\in B$ and $P_A(u)=b$, where $P_A:B\map A$ is a nearest point projection map. 
Suppose $x,y\in F_b$ and let $\gamma_x, \gamma_y$ be two $K$-qi sections over $[u,b]$.
Let $\QQQ_1, \QQQ_2$ be two $K'$-qi sections over $A$ in $Y$ and $U=U_{R'}(\QQQ_1, \QQQ_2)$.
If $d_u(\gamma_x(u),\gamma_y(u))\leq R$ and $U\neq \emptyset$, then 
$d_b(x,y)\leq R_{\ref{Y-path 1.2}}(R,R',K,K')$ and $d_A(b, U)\leq D_{\ref{Y-path 1.2}}(R,R',K,K')$. 
\end{lemma}
\proof
Suppose $U \neq \emptyset$ and $d_u(\gamma_x(u),\gamma_y(y))\leq R$.
Let $b'\in U_{M_{K'}}(\QQQ_1,\QQQ_2)$ be any point and let $[b,b']$ denote a geodesic in $A$ 
joining $b,b'$.
Then the concatenation $[u,b]*[b,b']$ is a $K_{\ref{subqc-elem}}(\delta_0, k_0,k,0)$-quasigeodesic in $B$
by Lemma \ref{subqc-elem} since $A$ is $k$-qi embedded and $k_0$-quasiconvex.
Concatenation of $\gamma_x, \gamma_y$ with the qi sections over $[b,b']$ contained
in $\QQQ_1, \QQQ_2$ respectively defines $\max\{K, K'\}$-qi sections over $[u,b]*[b,b']$ passing 
through $x,y$ respectively. Let $k'=K_{\ref{subqc-elem}}(\delta_0, k_0,k,0)$
and $k''=\max\{K, K'\}$. Then by Lemma \ref{qi composition} these qi 
sections are $(k'k'',k''k'+k'')$-quasigeodesics in $X$. Since $X$ is $\delta$-hyperbolic
and $d(\gamma_x(u),\gamma_y(u))\leq R$ and $d(\QQQ_1\cap F_{b'}, \QQQ_2\cap F_{b'})\leq R'$,
by Corollary \ref{slim finite polygons} $x$ is contained in the 
$D':=(R+R'+2D_{\ref{slim iff gromov}}(\delta,k'k'',k'k''+k''))$-neighborhood of 
the qi section over $[u,b]*[b,b']$ passing through $y$. Applying Lemma \ref{distance from qi section}
to the restriction bundles over $[u,b]$ and $[b,b']$ we have $d_b(x,y)\leq R'_1$ where 
$R'_1= R_{\ref{distance from qi section}}(D', K)$.
Hence, we can take $R_{\ref{Y-path 1.2}}(R,R',K,K')=R'_1$.
Finally by Lemma \ref{qc-level-set-new} $d_A(b,U)\leq D_{\ref{qc-level-set-new}}(K',R'_1/M_{K'})$.
This completes the proof by taking 
$D_{\ref{Y-path 1.2}}(R,R',K,K')=D_{\ref{qc-level-set-new}}(K',R'_1/M_{K'})$. \qed

We recall that the paths $\bar{c}(y,y')$ were constructed from $c(y,y')$ by replacing parts of
$c(y,y')$ by some fiber geodesic segments. The main aim of the following three lemmas is to proving
that these fiber geodesic segments are uniform quasigeodesics in $Y$. Depending on how the corresponding
subladders of $X$ intersect $Y$ we have three scenarios and hence we divided the proof into three
lemmas.
\begin{lemma}\label{Y-path 1.3}
Given $K\geq K_0$ and $R\geq M_K$ there are constants $K_{\ref{Y-path 1.3}}=K_{\ref{Y-path 1.3}}(K,R)$,
$\epsilon_{\ref{Y-path 1.3}}=\epsilon_{\ref{Y-path 1.3}}(K,R)$ and
$D_{\ref{Y-path 1.3}}=D_{\ref{Y-path 1.3}}(K,R)$ such that the following holds.

Suppose $\QQQ, \QQQ'$ are two $K$-qi sections in $X$ and $d_h(\QQQ,\QQQ')\leq R$ in $X$.
Let $U=U_R(\QQQ, \QQQ')$. Suppose $d_h(\QQQ\cap Y, \QQQ'\cap Y)\geq R$ in $Y$. Then the following hold.

\begin{enumerate}
\item The projection of $U$ on $A$ is of diameter at most $D_{\ref{Y-path 1.3}}$.
\item For any $b\in P_A(U)$, $F_b\cap \LL(\QQQ,\QQQ')$ is a $K_{\ref{Y-path 1.3}}$-quasigeodesic in $Y$;
moreover, $F_b\cap \QQQ$ is an $\epsilon_{\ref{Y-path 1.3}}$-approximate nearest point projection
of any point of $\QQQ'$ on $\QQQ$ and vice versa.
\end{enumerate}
\end{lemma}

\proof (1) We know that $A$ is $k_0$-quasiconvex in $B$.
By Lemma \ref{qc-level-set-new}
$U$ is $K_{\ref{qc-level-set-new}}(K)$-quasiconvex in $B$. Let 
$\lambda'=\max\{k_0, K_{\ref{qc-level-set-new}}(K)\}$.
Suppose $P_A:B\map A$ is a nearest point projection map and $a,a'\in P_A(U)$
with $d_B(a,a')\geq D_{\ref{cor: lip proj}}(\delta, \lambda',0)$. Then there are $u,u'\in U$ such that
$d_B(a,u)\leq R_{\ref{cor: lip proj}}(\delta, \lambda',0)$ and $d_B(a',u')\leq R_{\ref{cor: lip proj}}(\delta, \lambda',0)$.
Let $D=R_{\ref{cor: lip proj}}(\delta, \lambda',0)$. We know $d_u(\QQQ\cap F_u,\QQQ'\cap F_u)\leq R$.
Hence by the bounded flaring condition we have $d_a(\QQQ\cap F_a,\QQQ'\cap F_a)\leq \mu_K(D)R$.
Similarly $d_{a'}(\QQQ\cap F_{a'},\QQQ'\cap F_{a'})\leq \mu_K(D)R$. Let $R_1= \mu_K(D)R$.
Thus, $a,a'\in U_{R_1}(\QQQ\cap Y, \QQQ'\cap Y)$. Since $R_1\geq M_K$, by Lemma \ref{qc-level-set-new}
we have $diam(U_{R_1}(\QQQ\cap Y, \QQQ'\cap Y))\leq D_{\ref{qc-level-set-new}}(K, R_1)$. This
proves (1). In fact, we can take
$D_{\ref{Y-path 1.3}}=\max \{D_{\ref{cor: lip proj}}(\delta, \lambda',0), D_{\ref{qc-level-set-new}}(K, R_1)\}$.

We derive (2) from Lemma \ref{Y-path 1.1} as follows. Let $u\in U$ be such that
$P_A(u)=b$ and let $x,y\in F_b\cap \LL(\QQQ, \QQQ')$. Suppose $\QQQ_1, \QQQ'_1$
are two $K'$-qi sections in $Y$ passing through $x,y$ respectively and 
$U'=U_{M_{K'}}(\QQQ_1,\QQQ'_1)$. Suppose $U'\neq \emptyset$.
Consider the restriction $Z$ of the bundle $X$ on $[u,b]\subset B$. In this bundle $\QQQ\cap Z, \QQQ'\cap Z$
are $K$-qi sections. By Proposition \ref{ladders are qi embedded}(3) there are
$(1+2K_0)C_{\ref{ladders are qi embedded}}(K)$-qi sections over $ub$ contained in the ladder
$\LL(\QQQ\cap Z, \QQQ'\cap Z)$ passing through $x,y$. Call them $\gamma_x, \gamma_y$ respectively.
We note that $d(\gamma_x(u), \gamma_y(u))\leq R$.
Now applying Lemma \ref{Y-path 1.2} we know that $d_B(b, U')$ is uniformly small. This verifies
the hypothesis of Lemma \ref{Y-path 1.1}. Thus $\QQQ\cap F_b$ is a uniform
approximate nearest point projection of $\QQQ'\cap F_b$ on $\QQQ$. Since 
$d_h(\QQQ\cap Y, \QQQ'\cap Y)\geq R\geq M_K$ the qi sections $\QQQ\cap Y, \QQQ'\cap Y$ are
uniformly cobounded by Lemma \ref{lem: ladders cobdd}. This shows that $\QQQ\cap F_b$ is
a uniform approximate nearest point projection of any point of $\QQQ'$ on $\QQQ$.
That $\QQQ'\cap F_b$ is a uniform approximate nearest point projection of any point of
$\QQQ$ on $\QQQ'$ is similar and hence we skip it.
\qed

\begin{lemma}\label{Y-path 1.4}
Given $D\geq 0$, $K\geq K_0$ and $R\geq M_K$ there are constants 
$K_{\ref{Y-path 1.4}}=K_{\ref{Y-path 1.4}}(D,K,R)$ 
$\epsilon_{\ref{Y-path 1.4}}=\epsilon_{\ref{Y-path 1.4}}(D,K,R)$ and
$D_{\ref{Y-path 1.4}}=D_{\ref{Y-path 1.4}}(D,K,R)$ such that the following holds.

Suppose $\QQQ, \QQQ'$ are two $K$-qi sections in $X$ and $d_h(\QQQ,\QQQ')\leq R$ in $X$.
Let $U=U_R(\QQQ, \QQQ')$. Suppose $U\neq \emptyset$ and $diam(U)\leq D$. Then the following holds.
\begin{enumerate}
\item $diam(P_A(U))\leq D_{\ref{Y-path 1.4}}$.
\item For any $b\in P_A(U)$, $F_b\cap \LL(\QQQ,\QQQ')$ is a $K_{\ref{Y-path 1.4}}$-quasigeodesic in $Y$.
\item $F_b\cap \QQQ$ is an $\epsilon_{\ref{Y-path 1.4}}$-approximate nearest point projection
of any point of $\QQQ'$ on $\QQQ$ and vice versa.
\end{enumerate}
\end{lemma}
\proof (1) Since $B$ is $\delta_0$-hyperbolic and $A$ is $k_0$-quasiconvex in $B$ any nearest
point projection map $P_A:B\map A$ is coarsely $L:=L_{\ref{cor: lip proj}}(\delta_0, k_0,0)$-Lipschitz.
Hence, $diam(P_A(U))\leq L+DL$.

We can derive (2), (3) from Lemma \ref{Y-path 1.1} and the hypotheses of Lemma \ref{Y-path 1.1}
can be verified using Lemma \ref{Y-path 1.2}. The proof is an exact copy of the proof of
Lemma \ref{Y-path 1.3}(2),(3). Hence we omit it. The only part that requires explanation
is why $\QQQ\cap Y$, $\QQQ'\cap Y$ are uniformly cobounded in $Y$.
If $d_h(\QQQ\cap Y,\QQQ'\cap Y)>R$ then we are done by Lemma \ref{lem: ladders cobdd}.
Suppose this is not the case. Then by the hypothesis 
$diam(U_R(\QQQ\cap Y,\QQQ'\cap Y))\leq k(k+D)$ since $A$ is $k$-qi embedded in $B$. Then
we are done by the first part of Lemma \ref{step 1(c)-1}. \qed

\begin{lemma}\label{Y-path 1.5}
Given $K\geq K_0$ and $R\geq M_K$ there is a constant 
$D_{\ref{Y-path 1.5}}=D_{\ref{Y-path 1.5}}(K,R)$  
such that the following holds.

Suppose $\QQQ, \QQQ'$ are two $K$-qi sections in $X$ and $d_h(\QQQ\cap Y,\QQQ'\cap Y)\leq R$.
Let $U=U_R(\QQQ, \QQQ')$. Then the following holds.
For any $b\in P_A(U)$, $d_b(\QQQ\cap F_b,\QQQ'\cap F_b)\leq D_{\ref{Y-path 1.5}}$.
\end{lemma}
\proof Suppose $u\in U$ and $P_A(u)=b$. If $u\in A$ then $b=u$ and
$d_b(\QQQ\cap F_b,\QQQ'\cap F_b)\leq R$. Suppose $u\not \in A$. We note that
$U\cap Y=U(\QQQ\cap Y,\QQQ'\cap Y)\neq \emptyset$. Let $v\in U(\QQQ\cap Y,\QQQ'\cap Y)$.
Then by Lemma \ref{subqc-elem} $[u,b]*[b,v]$ is a
$K_{\ref{subqc-elem}}(\delta_0, k_0, 1, 0)$-quasigeodesic in $B$.
Since $U$ is $K_{\ref{qc-level-set-new}}(K)$-quasiconvex in $B$. Let
$k'=K_{\ref{subqc-elem}}(\delta_0, k_0, 1, 0)$. 
Hence, by Lemma \ref{stab-qg}, $b\in N_D(U)$ where $D=D_{\ref{stab-qg}}(\delta_0, k',k')+K_{\ref{qc-level-set-new}}(K)$.
Finally by the bounded flaring
$d_b(\QQQ\cap F_b, \QQQ'\cap F_b)\leq R\max\{1, \mu_K(D)\}$. Hence we can take 
$D_{\ref{Y-path 1.5}}=R\max\{1, \mu_K(D)\}$.
\qed

\smallskip
Finally, we are ready to finish the proof of step 3.

\begin{lemma}\label{small-girth-ladder}
For $0\leq i \leq n-1$ we have the following.
\begin{enumerate}
\item $\bar{y}_{i+1}$ is a uniform approximate nearest point projection of $\bar{y}_i$
on $\SSS_{i+1}\cap Y$.
\item $\bar{\gamma}_i$ is a uniform quasigeodesic in $Y$.
\end{enumerate}
\end{lemma}
\proof The proof is broken into three cases depending on the type of $\LL_i$.

{\bf Case 1: $i\leq n-2$ and $\LL_i$ is of type (I):}
By Corollary \ref{qc-level-set-new} $U_{R_0}(\SSS_{i},\SSS_{i+1})$
has uniformly small diameter. Hence by Lemma \ref{Y-path 1.4}(2)
$[\SSS_i\cap F_{b_{i+1}}, \SSS_{i+1}\cap F_{b_{i+1}}]_{b_{i+1}}$ is a uniform quasigeodesic in $Y$.
By the part (3) of the same lemma $\SSS_{i+1}\cap F_{b_{i+1}}$ is a uniform approximate
nearest point projection of $\SSS_i\cap F_{b_{i+1}}$ on $\SSS_{i+1}\cap Y$ and $\SSS_{i}\cap F_{b_{i+1}}$ is a uniform approximate
nearest point projection of $\SSS_{i+1}\cap F_{b_{i+1}}$ on $\SSS_{i}\cap Y$ in $Y$.
Hence the second part of the lemma follows, in this case, by Lemma \ref{subqc-elem}.

{\bf Case 2: $i\leq n-2$ and $\LL_i$ is of type (II):}
Suppose $\LL_i$ is a ladder of type (II). In this case, it is enough, by Proposition \ref{hamenstadt}, 
to show the following two statements $(2')$ and $(2'')$:

$(2')$: {\em $\bar{y}'_i$ is a uniform approximate nearest point
projection of $\bar{y}_i$ on $\SSS'_i\cap Y$ in $Y$ and the concatenation of $\bar{\alpha}_i$ and the fiber
geodesic $[\SSS_i\cap F_{b'_i}, \SSS'_i\cap F_{b'_i}]_{b'_i}$ is a uniform quasigeodesic
in $Y$.} 

We know that $d_h(\SSS_i, \SSS'_i)\leq R_1$. Depending on the nature of
$d_h(\SSS_i\cap Y, \SSS'_{i}\cap Y)$ the proof of $(2')$ is broken into the following
two cases.


{\bf Case $(2')(i)$:} Suppose $d_h(\SSS_i\cap Y,\SSS'_i\cap Y)\leq R_1$. In this case 
$d_{b'_i}(\SSS_i\cap F_{b'_i}, \SSS'_i\cap F_{b'_i})$ is uniformly small by
Lemma \ref{Y-path 1.5}. By Lemma \ref{step 1(c)-1} if $b''_i$ is a nearest point projection
of $\pi(\bar{y}_i)$ on $U_{R_1}(\SSS_i\cap Y, \SSS'_i\cap Y)$ then $F_{b''_i}\cap \SSS'_i$ is
a uniform approximate nearest point projection of $\bar{y}_i$ on $\SSS'_i\cap Y$ in $Y$. Thus
it is enough to show that $d_B(b''_i, b'_i)$ uniformly bounded to prove that
$\bar{y}'_i$ is a uniform approximate nearest point projection of $\bar{y}_i$ on $\SSS'_i\cap Y$
in $Y$. Then since $\SSS_i\cap Y$ is $K'_1$-qi section in $Y$ and 
$d_{b'_i}(\SSS_i\cap F_{b'_i}, \SSS'_i\cap F_{b'_i})$ is uniformly small it will follow that 
the concatenation of $\bar{\alpha}_i$ and the fiber
geodesic $[\SSS_i\cap F_{b'_i}, \SSS'_i\cap F_{b'_i}]_{b'_i}$ is a uniform quasigeodesic
in $Y$.

That $d_B(b''_i, b'_i)$ uniformly bounded  is proved as follows.
Let $U=U_{R_1}(\SSS_i, \SSS'_i)$, $V=U\cap A=U_{R_1}(\SSS_i\cap Y, \SSS'_i\cap Y)$. 
Since $B$ is $\delta_0$-hyperbolic, $A$ is $k$-qi embedded in $B$ and $V$ is $\lambda_2$-quasiconvex
in $A$, $V$ is $K_{\ref{qc in subspace}}(\delta_0,k,\lambda_2)$-quasiconvex in $B$.
Let $k'=\max\{ \lambda_2, k_0, K_{\ref{qc-level-set-new}}(K_2), K_{\ref{qc in subspace}}(\delta_0,k,\lambda_2)\}$.
Then $A,U,V$ are all $k'$-quasiconvex in $B$.
By the definitions of $y_i$'s we know that $\pi(y'_i)$ is the nearest point projection of
$\pi(y_i)$ on $U$. Let $\bar{b}'_i$ be a nearest point projection of $\pi(y'_i)$ on $V$.
Also $b'_i=\pi(\bar{y}'_i)$ is the nearest point projection of $\pi(y'_i)$ on $A$. 
On the other hand, $b_i=\pi(\bar{y}_i)$ is a nearest point projection of $\pi(y_i)$ on $A$
and $b''_i$ is the nearest point projection of $b_i$ on $V$.
Therefore, $d_B(b''_i, \bar{b}'_i)\leq 2D_{\ref{nested qc sets}}(\delta_0, k', 0)$
by Corollary \ref{nested qc sets}.

Now, by Lemma \ref{Y-path 1.5}
$d_{b'_i}(\SSS_i\cap F_{b'_i}, \SSS'_i\cap F_{b'_i})\leq D_{\ref{Y-path 1.5}}(K_2,R_1)$.
Hence, by Lemma \ref{qc-level-set-new} 
$d_A(b'_i, V)\leq D_{\ref{qc-level-set-new}}(K_2,D_{\ref{Y-path 1.5}}(K_2,R_1)/R_1)=D_1$, say.
Let $v\in V$ be such that $d_A(b'_i,v)\leq D_1$. Then $d_B(b'_i,v)\leq kD_1+k$. Hence,
$Hd([\pi(y'_i), b'_i]_B, [\pi(y'_i),v]_B)\leq \delta_0+k+kD_1$. However, the concatenation
$[\pi(y'_i),\bar{b}'_i]_B*[\bar{b}'_i,v]_B$ is a $K_{\ref{subqc-elem}}(\delta_0,k',1,0)$-quasigeodesic.
Hence, there is a point $w\in [\pi(y'_i),v]_B$ such that 
$d_B(w,\bar{b}'_i)\leq D_{\ref{stab-qg}}(\delta_0, K_{\ref{subqc-elem}}(\delta_0,k',1,0), K_{\ref{subqc-elem}}(\delta_0,k',1,0))=D_2$,
say. Thus there is a point $w'\in [\pi(y'_i),b'_i]$ such that $d_B(w', \bar{b}'_i)\leq D_2+\delta_0+k+kD_1=D_3$, say.
But $b'_i$ is a nearest point projection of $\pi(y'_i)$ on $A$ and $\bar{b}'_i\in V\subset A$. Thus
$d_B(w', b'_i)\leq D_3$. Thus $d_B(\bar{b}'_i, b'_i)\leq 2D_3$. Hence,
$d_B(b'_i, b''_i)\leq d_B(b''_i, \bar{b}'_i)+d_B(\bar{b}'_i, b'_i)\leq 2D_{\ref{nested qc sets}}(\delta_0, k', 0)+2D_3$.

{\bf Case $(2')(ii)$:} Suppose $d_h(\SSS_i\cap Y,\SSS'_i\cap Y)\geq R_1$. In this case Lemma \ref{Y-path 1.3}
and Lemma \ref{subqc-elem} do the job. 

$(2'')$: {\em $\bar{y}_{i+1}$ is a uniform approximate nearest point projection
of $\bar{y}'_i$ on $\SSS_{i+1}\cap Y$ in $Y$ and the concatenation of $\bar{\beta}$ and the 
fiber geodesic $[\SSS'_i\cap F_{b_{i+1}}, \SSS_{i+1}\cap F_{b_{i+1}}]_{b_{i+1}}$ is a uniform 
quasigeodesic joining $\bar{y}'_i$ to $\bar{y}_{i+1}$ in $Y$.}

In this case $d_h(\SSS'_i\cap Y, \SSS_{i+1}\cap Y)\leq 1$ hence we are done as in Case $(2')(i)$.

{\bf Case 3: $i=n-1$: } The proof of this case is also analogous to that of the proof
of Case $(2')(i)$ since $d_h(\SSS_{n-1}, \SSS_n)\leq R_0$.
\qed

\begin{rem}
The conclusion of Lemma \ref{Y-path 1.4} is subsumed by Lemma \ref{Y-path 1.3} and Lemma \ref{Y-path 1.5}.
But we still keep Lemma \ref{Y-path 1.4} for the sake of ease of explanation.
\end{rem}

Thus by Lemma \ref{applying hamenstadt prop} and Lemma \ref{small-girth-ladder}, we have proved the following.
\begin{prop}\label{thm-1}
	Let $x,y \in Y$ and let $\SSS$ and $\SSS'$ be two $K_0$-qi sections in $X$ through $x$ and $y$ respectively. 
Let $c(x,y)$ be a uniform quasigeodesic in $X$ joining $x$ and $y$ which is contained in $\LL(\SSS,\SSS')$ as constructed
in step 1(c). Then the corresponding modified path $\tilde{c}(x,y)$, as constructed in step 2, is a uniform quasigeodesic in $Y$.
\end{prop}


{\bf Step 4. Verification of the hypothesis of Lemma \ref{mitra-lemma}.}

\begin{lemma}\label{prop-embed}{\bf (Proper embedding of the pullback $Y$)}
The pullback $Y$ is metrically properly embedded in $X$. In fact, the distortion function for $Y$ is the composition of
a linear function with $\eta$, the common distortion function for all the fibers of the bundle $X$.
\end{lemma}

\proof As was done in the proof of the main theorem, we shall assume that $g$ is the inclusion map
and $Y=\pi^{-1}(A)$ and $p$ is the restriction of $\pi$.
Let $x, y \in Y$ such that $d_{X}(x,y) \leq M$. Let $\pi(x) = b_1$ and $\pi(y) = b_2$. 
Then, $d_{B}(b_1,b_2) \leq M$ and hence $d_A(b_1,b_2) \leq k+kM$.
Let $[b_1,b_2]_A$ be a geodesic joining $b_1$ and $b_2$ in $A$. 
This is a quasigeodesic in $B$. By \lemref{lifting geodesics}, there exists an isometric section 
$\gamma$ over $[b_1,b_2]_A$, through $x$ in $Y$. Clearly, $\gamma$ is a qi lift in $X$, say 
$k'$-qi lift. We have, $l_{X}(\gamma) \leq k'(kM+k)+k'=:D(M)$. The concatenation of $\gamma$ 
and the fiber geodesic $[\gamma \cap F_{b_2},y]_{F_{b_2}}$ is a path, denoted by $\alpha$, 
joining $x$ and $y$ in $X$. So, 
$$d_{X}(\gamma\cap F_{b_2},y) \leq d_X(\gamma\cap F_{b_2},x)+d_X(x,y) \leq l_X(\gamma)+d_X(x,y) \leq D(M)+M.$$
Now, since $F_{b_2}$ is uniformly properly embedded as measured by $\eta$, we have, 
$d_{b_2}(\gamma\cap F_{b_2},y) \leq \eta(D(M)+M)$. Now, $\alpha$ lies in $Y$ and $l_{Y}({\gamma}) \leq kM+k$.
Then, 
$$ d_Y(x,y) \leq l_Y(\alpha) \leq l_Y({\gamma})+ d_{Y}({\gamma} \cap F_{b_2},y)
	\leq kM+k+ d_{b_2}(\tilde{\gamma} \cap F_{b_2},y).$$
Therefore, $d_{Y}(x,y) \leq kM+k+\eta(D(M)+M)$. 
Setting $\eta_0(M) := kM+k+\eta(D(M)+M)$, we have the following: for all $x,y \in Y$, $d(x,y)\leq M$ implies 
$d_{Y}(x,y) \leq \eta_0(M)$.
\qed

We recall that we fixed a vertex $b_0\in A$ to define the paths $c(y,y')$ in the last step. Let $y_0\in F_{b_0}$.
However, the following lemma completes the proof of Theorem \ref{CT for mgbdl}.
\begin{lemma}\label{prop-2}
Given $D>0$, there is $D_1>0$ such that the following holds.

If $d_{X}(y_0,c(y,y')) \leq D$ then $d_{Y}(y_0, \bar{c}(y,y')) \leq D_1$.
\end{lemma}
\proof
Let $x \in c(y,y')$ be such that $d_{X}(y_0,x) \leq D$. This implies that $d_{B}(\pi(x),b_0) \leq D$. 
We recall that the path $c(y,y')$ is a concatenation of $\gamma_j$, $j=0,1,\cdots, n$.
Suppose $x\in \gamma_i$, $0\leq i\leq n$. We claim that there is a point of $\bar{\gamma}_i$
uniformly close to $y_0$.  Now, $\gamma_i$ is either a lift of geodesic segments of $B$ in a 
$K_2$-qi section $\SSS_i$ or possibly $\SSS'_i$ or it is the concatenation of such a lift and a fiber geodesic 
of length at most $R_1$. Let $\QQQ$ denote the corresponding qi section and suppose $c(y,y')\cap \QQQ$ joins
the points $z\in \QQQ$ to $w\in \QQQ$. If $i=n$ then $\gamma_i$ is a qi lift of $[\pi(z), \pi(w)]_B$
in $\QQQ$ joining $z,w$. Otherwise there is a fiber geodesic $\sigma\subset c(y,y')\cap F_{\pi(w)}$
connecting $\QQQ$ to the next qi section $\QQQ'$, say. Then both the points $z$ and $\QQQ'\cap \sigma$ are 
one of the $y_i$'s or $y'_j$'s. Let $z'=\QQQ'\cap \sigma$ and $b'=\pi(z')$.
Let $b$ be the nearest point projection of $\pi(x)$ on $A$.  It follows that $d_B(\pi(x), b)\leq D$.

Suppose $x\in \sigma$. By the definition of $\bar{c}(y,y')$ we have $\QQQ'\cap F_b\in \bar{c}(y,y')$.
However, $d_B(b,b_0)\leq d_B(b,\pi(x))+d_B(b_0, \pi(x))\leq 2D$. Since $A$ is $k$-qi embedded in $B$ we have
$d_A(b,b_0)\leq k+2Dk$. Hence, $d_Y(\QQQ'\cap F_{b_0}, \QQQ'\cap F_{b})\leq K_2+ (k+2Dk).K_2$.
On the other hand in this case $\pi(x)=b'$ and $d_{b'}(z',x)\leq R_1$. Hence, $d_X(z',y_0)\leq R_1+D$. Thus
$d_X(y_0, \QQQ'\cap F_{b_0})\leq d_X(y_0, x)+d_X(x,z')+d_X(z',\QQQ'\cap F_{b_0})\leq D+R_1+K_2+DK_2$
since $d_B(b,b_0)\leq 2D$. 
Hence, $d_Y(y_0, \QQQ'\cap F_{b_0})\leq d_{b_0}(y_0, \QQQ'\cap F_{b_0})\leq \eta(D+R_1+K_2+DK_2)$.
Thus $d_Y(y_0, \QQQ'\cap F_b)\leq d_Y(y_0, \QQQ'\cap F_{b_0})+d_Y(\QQQ'\cap F_b, \QQQ'\cap F_{b_0})
\leq d_Y(y_0, \QQQ'\cap F_{b_0})+K_2+K_2d_A(b,b_0)\leq \eta(D+R_1+K_2+DK_2)+K_2+(k+2Dk)K_2$.
Hence, in this case $d_{Y}(y_0, \bar{c}(y,y'))\leq (1+k+2Dk)K_2+ \eta(D+R_1+K_2+DK_2)$.

\begin{comment}
Then $\pi(x)=b'$ and $d_{b'}(z',x)\leq R_1$ and therefore $d_X(z',y_0)\leq R_1+D$.
Also $d_B(b', b_0)\leq d_X(x,y_0)\leq D$. However, by the definition of the modified paths we have
$\QQQ'\cap F_b\in \bar{c}(y,y')$.
Hence, lifting $[b',b]$ in $\QQQ'$ we find a $K_2$-quasigeodesic in $X$ joining $z'$ to $\QQQ'\cap F_b$.
Thus $d_X(z', \QQQ'\cap F_b)\leq K_2d_B(b',b)+K_2\leq K_2d_B(b',b_0)+K_2\leq K_2+DK_2$. Hence,
$d_X(\QQQ'\cap F_b, y_0)\leq d_X(z',y_0)+ d_X(z', \QQQ'\cap F_b)\leq R_1+D+K_2+DK_2$. \end{comment}

Otherwise suppose $x$ is contained in the lift of $[\pi(z),\pi(w)]_B$ in $\QQQ$. 
We note that $\pi(x)\in [\pi(z),\pi(w)]_B$ and $d_B(\pi(x), A)\leq D$. Now $A$ is $k_0$-quasiconvex in $B$.
Hence, by Lemma \ref{trivial lemma} we have
$d_B(\pi(x), [\overline{\pi(z)}, \overline{\pi(w)}]_B)\leq D_{\ref{trivial lemma}}(D,k_0,\delta)$.
where $\overline{\pi(z)}$, $\overline{\pi(w)}$ are nearest point projections of $\pi(z), \pi(w)$
respectively on $A$. Since $A$ is $k$-qi embedded in $B$ by stability of quasigeodesics
$Hd([\overline{\pi(z)}, \overline{\pi(w)}]_B, [\overline{\pi(z)}, \overline{\pi(w)}]_A)\leq D_{\ref{stab-qg}}(\delta,k,k)$.
Hence, $d_B(\pi(x), [\overline{\pi(z)}, \overline{\pi(w)}]_A)\leq D_{\ref{trivial lemma}}(D,k_0,\delta)+D_{\ref{stab-qg}}(\delta,k,k)$.
Let $\alpha$ be the lift of $[\overline{\pi(z)}, \overline{\pi(w)}]_A$ in $\QQQ$. Then $\alpha\subset \bar{c}(y,y)$.
On the other hand, 
$d_X(x, \alpha)\leq K_2+K_2d_B(\pi(x), [\overline{\pi(z)}, \overline{\pi(w)}]_A)\leq K_2+K_2(D_{\ref{trivial lemma}}(D,k_0,\delta)+D_{\ref{stab-qg}}(\delta,k,k))=D_1$, say. Hence, $d_X(y_0, \alpha)\leq d_X(y_0,x)+d_X(x,\alpha)\leq D+D_1$. 
This implies that $d_Y(y_0, \alpha)$ is also bounded by a function of $D$ and the other parameters of the metric graph
bundles $X$ and $Y$, by Lemma \ref{prop-embed}.
\qed

\subsection{An example}
For the convenience of the reader,
we briefly illustrate a special case of our main theorem where $B=\RR, A=(-\infty, 0]$.
This discussion will also be used in the proof of the last proposition of the next section. 
We shall assume $b_0=0$ here.

As in the proof of Lemma \ref{prop-2} suppose $\QQQ, \QQQ'$ are two qi sections among the various 
$\SSS_i, \SSS'_j$'s and let $w'\in\QQQ', z,w\in \QQQ$ are points of $c(y,y')$ where $\pi(w')=\pi(w)$,
$d_{\pi(w)}(w, w')\leq R_1$ and the concatenation of the lift say $\alpha$, of $[\pi(z), \pi(w)]$ in $\QQQ$ and 
the vertical geodesic segment, say $\sigma$, in $F_{\pi(w)}$ is a part of $c(y,y')$. Following are the
possibilities.

{\bf Case 1.} If $w',z\in Y\cap c(y,y')$ then $\alpha*\sigma\subset Y$ and it is the corresponding
part of $\bar{c}(y,y')$.

{\bf Case 2.} $z\in Y, w'\not \in Y$. In this case, the modified segment is formed as the concatenation
of subsegment of $\alpha$ joining $z$ to $\QQQ\cap F_0$ and the fiber geodesic $[\QQQ\cap F_0, \QQQ'\cap F_0]_0$.

\begin{figure}[h]
\begin{tikzpicture}
    \draw (1,0) ellipse (.5cm and 1.5cm);
	\draw (5,0) ellipse (.5cm and 1.5cm);
	\draw (9,0) ellipse (.5cm and 1.5cm);
	\draw[black,thick,dotted](0,1.25) -- (1,1.2);
	\draw[black,thick,dotted](9,.8) -- (11,.7);
	\draw[black,thick,dotted](0,-1.25) -- (11,-.7);
	\draw[black,thick,dashed](5,1) node[above]{$\overline{w}$} -- (9,.8);
	\draw[black,thick](1,1.2) -- (5,1);
	\draw[black,thick,dotted](1,1.2) node[above]{$z$} -- (1,-1.2);
	\draw[black,thick](5,1) -- (5,-1) node[below]{$\overline{w'}$};
	\draw[black,thick,dashed](9,.8) node[above]{$w$} -- (9,-.8)node[below]{$w'$};
	\node at (1,-2) {$F_{\pi(z)}$};
	\node at (5,-2) {$F_0$};
	\node at (9,-2) {$F_{\pi(w')}$};
	\node at (11,.85) {$\mathcal{Q}$};
	\node at (11,-.85) {$\mathcal{Q}'$};
\end{tikzpicture}
\caption{Case 2}
\end{figure}

{\bf Case 3.} $w'\in Y, z\not \in Y$. In this case the modified segment is the concatenation of the segment of
$\alpha$ from $\QQQ\cap F_0$ to $w$ and the fiber geodesic segment $\sigma$.

\begin{figure}[h]
	\begin{tikzpicture}
	\draw (1,0) ellipse (.5cm and 1.5cm);
	\draw (5,0) ellipse (.5cm and 1.5cm);
	\draw (9,0) ellipse (.5cm and 1.5cm);
	\draw[black,thick,dotted](0,1.25) -- (1,1.2);
	\draw[black,thick,dotted](9,.8) -- (11,.7);
	\draw[black,thick,dotted](0,-1.25) -- (11,-.7);
	\draw[black,thick,dashed](5,1) node[above]{$\overline{z}$} -- (9,.8);
	\draw[black,thick](1,1.2) -- (5,1);
	\draw[black,thick](1,1.2) node[above]{$w=\overline{w}$} -- (1,-1.2) node[below]{$w'=\overline{w'}$};
	\draw[black,thick,dotted](5,1) -- (5,-1);
	\draw[black,thick,dotted](9,.8) node[above]{$z$} -- (9,-.8);
	\node at (1,-2) {$F_{\pi(w')}$};
	\node at (5,-2) {$F_0$};
	\node at (9,-2) {$F_{\pi(z)}$};
	\node at (11,.85) {$\mathcal{Q}$};
	\node at (11,-.85) {$\mathcal{Q}'$};
	\end{tikzpicture}
	\caption{Case 3}
\end{figure}

{\bf Case 4.} $z,w'\not \in Y$. In this case the modified segment is the fiber geodesic 
$[\QQQ\cap F_0, \QQQ'\cap F_0]_0$.

\begin{figure}[h]
	\begin{tikzpicture}
	\draw (1,0) ellipse (.5cm and 1.5cm);
	\draw (5,0) ellipse (.5cm and 1.5cm);
	\draw (9,0) ellipse (.5cm and 1.5cm);
	\draw[black,thick,dotted](0,1.25) -- (1,1.2);
	\draw[black,thick,dotted](9,.8) -- (11,.7);
	\draw[black,thick,dotted](0,-1.25) -- (11,-.7);
	\draw[black,thick,dashed](5,1) node[above]{$z$} -- (9,.8);
	\draw[black,thick, dotted](1,1.2) -- (5,1);
	\draw[black,thick](1,1.2) node[above]{$\overline{z}$} -- (1,-1.2) node[below]{$\overline{w'}$};
	\draw[black,thick,dotted](5,1) -- (5,-1);
	\draw[black,thick,dashed](9,.8) node[above]{$w$} -- (9,-.8)node[below]{$w'$};
	\node at (1,-2) {$F_{0}$};
	\node at (5,-2) {$F_{\pi(z)}$};
	\node at (9,-2) {$F_{\pi(w')}$};
	\node at (11,.85) {$\mathcal{Q}$};
	\node at (11,-.85) {$\mathcal{Q}'$};
	\end{tikzpicture}
	\caption{Case 4}
\end{figure}
Here, the dashed lines denote the portion of $c(y,y')$, the thick lines denote the portion of $\bar{c}(y,y')$
and dotted lines are portions of the qi sections $\QQQ, \QQQ'$.


\section{Applications, examples and related results}\label{sec:apps}

As the first application of our main theorem, we have the following.
Given a short exact sequence of finitely generated groups there is a natural way to associate a 
metric graph bundle to it as mentioned in Example 1.8 of \cite{pranab-mahan}. See also Example 
\ref{exact sequence example}. Having said that Theorem \ref{CT for mgbdl} gives the following as an immediate consequence.

\begin{theorem}\label{main application}
Suppose $1\map N\map G\stackrel{\pi}{\map} Q\map 1$ is a short exact sequence of hyperbolic groups
where $N$ is nonelementary hyperbolic. 
Suppose $Q_1$ is a finitely generated, qi embedded subgroup of $Q$ and $G_1=\pi^{-1}(Q_1)$.
Then the $G_1$ is hyperbolic and the inclusion $G_1\map G$ admits the CT map.
\end{theorem}

The next application is in the context of complexes of hyperbolic groups.
Suppose $\mathcal Y$ is a finite, connected simplicial complex and $ \GB(\mathcal Y)$
is a developable complex of nonelementary hyperbolic groups with qi condition defined over 
$\mathcal Y$ (see Section \ref{example: complex of groups})
such that the fundamental group $G$ of the complex of groups is hyperbolic. Suppose we have a good subcomplex
$\YY_1\subset \YY$ and $G_1$ is the image of $\pi_1(\GG,\YY_1)$ in $G$ under the natural homomorphism 
$\pi_1(\GG, \YY_1)\map \pi_1(\GG,\YY)$. Then we have the
following pullback diagram as obtained in Proposition \ref{complex: pullback} satisfying the properties of Theorem 
\ref{CT for mgbdl}.

\begin{figure}[h]
	\centering
	\begin{tikzpicture}[node distance=2cm,auto]
	\node (A) {$X_1$}; 
	\node (B) [right of= A] {$X$};
	\node (C) [below=1cm of A] {$B_1$};
	\node (D) [below=1.05cm of B] {$B$};
    \draw [->] (A) to node {$f$} (B);
	\draw [->] (A) to node {$\pi_1$} (C);
	\draw [->] (B) to node {$\pi$} (D);
	\draw [->] (C) to node {$i$} (D);
	\end{tikzpicture}
	\caption{}\label{complex of groups}
\end{figure}

\noindent
Thus we have:
\begin{theorem}\label{CT for complex of groups}
The group $G_1$ is hyperbolic and the inclusion $G_1\map G$ admits the CT map.
\end{theorem}

\begin{rem}
The rest of the paper is devoted to properties of the boundary of metric (graph) bundles and Cannon-Thurston maps.
We recall that qi sections, ladders etc for a metric bundle are defined as transport of the
same from the canonical metric graph bundle associated to it. {\em All the results in the rest of the section
are meant for metric bundles as well as metric graph bundles. However, using the dictionary provided 
by Proposition \ref{bundle vs graph bundle} it is enough to prove the results only for metric graph bundles.
Therefore, we shall state and prove results only for metric graph bundles in what follows starting with the
convention below.}
\end{rem}

\smallskip
\begin{convention} (1) For the rest of the paper we shall assume that 
$\pi: X\map B$ is a $\delta$-hyperbolic $\eta$-metric graph bundle over $B$
satisfying the hypothesis H1, H2, H3$'$ and H4 of section 5. 
(2) By Proposition \ref{visibility} any point of $\partial B$
can be joined to any point of $B\cup \partial B$ and any point of $\partial X$ can be joined
to $X\cup \partial X$ by a uniform quasigeodesic ray or line. We shall assume that these are $\kappa_0$-quasigeodesics.
(3) We shall assume that any geodesic in $B$ has a $c$-qi lift in $X$ using the path lifting lemma
for metric graph bundles. (4) We recall that through any point of $X$ there is a $K_0$-qi section over $B$.

\end{convention}

\subsection{Some properties of $\partial X$}
\begin{comment}For this subsection we let $\pi: X\map B$ be a $\delta$-hyperbolic $(\eta,c)$-metric bundle or 
$\eta$-metric graph bundle over $B$ satisfying the hypothesis H1 and H2 of section 5. For simplicity
paths in $B$ will be assumed to be continuous, arc length parametrized in the case of
a metric bundle and dotted edge paths for the case of a metric graph bundle. We shall
also assume that any $(1,1)$-quasigeodesic in $B$ has a $K_0$-qi lift through any point of $X$,
by dint of the path lifting lemma.
\end{comment}

\begin{lemma}\label{bdry 6.1}
Suppose $\alpha, \beta:[0,\infty)\map B$ are two $k$-quasigeodesic rays for some $k\geq 1$ with
$\alpha(\infty)=\beta(\infty)=\xi$. Suppose $\tilde{\beta}$ is a $K$-qi lift of $\beta$ for some
$K\geq 1$. Then there is a $K'$-qi lift $\tilde{\alpha}$ of $\alpha$ such that 
$\tilde{\alpha}(\infty)= \tilde{\beta}(\infty)$ where $K'$ depends on $k$, $K$, $d_B(\alpha(0), \beta(0))$
and the various parameters of the metric graph bundle.
\end{lemma}
\proof Suppose $\alpha, \beta:[0,\infty)\map B$ are two $k$-quasigeodesic rays for some $k\geq 1$ with
$\alpha(\infty)=\beta(\infty)=\xi$. This means $Hd(\alpha, \beta)<\infty$. Let $R=Hd(\alpha, \beta)$.
Then for all $s\in [0, \infty)$ there is $t=t(s)\in [0,\infty)$ such that $d_B(\alpha(s), \beta(t))\leq R$.
Let $\phi_{ts}:F_{\beta(t)}\map  F_{\alpha(s)}$ be fiber identification maps such that 
$d_X(x, \phi_{ts}(x))\leq 3c+3cR$ for all $x\in F_{\beta(t)}$, $t\in [0,\infty)$
where $c=1$ for metric graph bundles. (See Lemma \ref{fibers qi}.) 
Let $\tilde{\beta}$ be a $K$-qi lift of $\beta$. Now, for all $s\in [0,\infty)$ 
we define $\tilde{\alpha}(s)=\phi_{ts}(\tilde{\beta}(t))$. It is easy to verify that $\tilde{\alpha}$
thus defined is a uniform qi lift of $\alpha$. Also clearly $\tilde{\alpha}\subset N_{3c+3cR}(\tilde{\beta})$.
It follows that $\tilde{\alpha}(\infty)=\tilde{\beta}(\infty)$
\qed

\begin{cor}\label{cor: bundle bdry}
Let $\xi\in \partial B$ and let $\alpha$ be a quasigeodesic ray in $B$ joining $b$ to $\xi$.
Let $\partial^{\xi}_{\alpha} X:=\{\gamma(\infty):\gamma\, \mbox{is a qi lift of}\,\, \alpha \}$. 

Then $\partial^{\xi}_{\alpha} X$ is independent of $\alpha$; it is determined by $\xi$.
\end{cor}

Due to the above corollary, we shall use the notation $\partial^{\xi} X$ for all $\xi\in \partial B$ without
further explanation. The following proposition is motivated by a similar result proved
by Bowditch (\cite[Proposition 2.3.2]{bowditch-stacks}).

\begin{prop}\label{bundle boundary}
Let $b\in B$ be an arbitrary point and $F=F_b$. Then
we have $$\partial X=\Lambda(F)\cup(\coprod_{\xi\in \partial B} \partial^{\xi} X).$$
\end{prop}

\proof 
We first fix a point $x\in F$. Let $\gamma$ be a quasigeodesic ray in $X$ starting from $x$. 
Let $b_n=\pi(\gamma(n))$. Let $\alpha_n$ be a $(1,1)$-quasigeodesic in $B$ joining $b$ to $b_n$.
Let $\tilde{\alpha}_n$ be a $K_0$-qi lift of $\alpha_n$ joining $\gamma(n)$ to
$\tilde{\alpha_n}(b)=x_n\in F$. There are two possibilities.

Suppose $\{x_n\}$ has an unbounded subsequence say $\{x_{n_k}\}$. Then 
$d(x_{n_k},x)\map \infty$. We note that $\tilde{\alpha}_{n_k}$'s are uniform quasigeodesics in $X$ whose distance
from $x$ is going to infinity by Lemma \ref{distance from qi section}. Hence, by Lemma \ref{gromov product meaning} 
$x_{n_k}\map \gamma(\infty)$ and thus $\gamma(\infty)\in \Lambda(F)$. 

Otherwise, suppose $\{x_n\}$ is a bounded sequence. 

{\bf Claim:} In this case $\pi\circ\gamma$ is a quasigeodesic ray. 

{\em Proof of claim:} We note that by stability of quasigeodesics (Corollary \ref{cor: stab-qg})
and slimness of triangles (Lemma \ref{slim iff gromov}) $Hd(\tilde{\alpha}_n, \gamma|_{[0,n]})$
 is uniformly small for all $n$. 
This implies that $Hd(\alpha_n, (\pi\circ \gamma)|_{[0,n]})$ is uniformly small for all $n$; in particular
$d_B(b_m, \alpha_n)$ is uniformly small for all $n\geq m$. Next we note that $d_B(b,b_n)\map \infty$
for otherwise $d(\gamma(n), x)$ will be bounded. Then it follows that $\lim_{m,n\map \infty}(b_m.b_n)_b=\infty$.
Let $\xi=\lim_{n\map \infty} b_n$ and let $\alpha$ be a $\kappa_0$-quasigeodesic ray in $B$ joining $b$ to $\xi$.
Now, to show that $\pi\circ \gamma$ is a quasigeodesic it is enough to show by Lemma \ref{quasigeod criteria}
that $\pi\circ \gamma$ is (1) uniformly close to $\alpha$ and (2) properly embedded.

(1): Fix an arbitrary $m\in \NN$ and consider all $n\geq m$. Since $\lim_{n\map \infty}b_n=\alpha(\infty)=\xi$,
by Lemma \ref{convergence explained}(2) for any $\kappa_0$-quasigeodesic ray $\beta_n$ joining $b_n$ to $\xi$ we have
$d(b, \beta_n)\map \infty$. 
Since the triangles with vertices $b_n, b, \xi$ are uniformly slim by Lemma \ref{ideal triangles are slim}
and $d_B(b_m, \alpha_n)$ are uniformly small it follows that $b_m$ is uniformly close to $\alpha$. 
This shows (1). 

(2): Since $\pi$ is Lipschitz and $\gamma$ is a quasigeodesic it follows that $\pi\circ \gamma$ is coarsely Lipschitz.
Suppose $d_B(b_n, b_m)\leq D$ for some $D\geq 0$ and $m,n\in \NN$, $m\leq n$.
We claim that $d_X(\gamma(m), \gamma(n))$ is uniformly small. Note that this would then imply that
$n-m$ is uniformly small since $\gamma$ is quasigeodesic, and also that $\gamma$ is a qi lift of $\pi\circ \gamma$. 
We know that $Hd(\tilde{\alpha}_n, \gamma|_{[0,n]})\leq R$ for some constant $R$ independent of $n$. Hence,
$d_X(\gamma(m), \tilde{\alpha}_n)\leq R$. Let $y_{m, n}\in \tilde{\alpha}_n$ be such that $d_X(\gamma(m), y_{m,n})\leq R$.
Since $\pi$ is $1$-Lipschitz we have $d_X(b_m, \pi(y_{m,n}))\leq R$. Then 
$d_B(\pi(y_{m, n}), b_n)\leq d_B(\pi(y_{m, n}), b_m)+d_B(b_m,b_n)\leq R+D$. 
Since $\tilde{\alpha}_n$ is $K_0$-qi lift of $\alpha_n$ and $\pi\circ\tilde{\alpha}(n)=b_n$ it follows that
$d_X(y_{m,n}, \tilde{\alpha}(n))=d_X(y_{m,n},\gamma(n))\leq K_0(R+D)+K_0$. 
Hence, $d_X(\gamma(m), \gamma(n))\leq d_X(\gamma(m), y_{m,n})+ d_X(y_{m,n},\gamma(n))\leq R+K_0(R+D)+K_0$. 
Since $\gamma$ is quasigeodesic it follows that $(n-m)$ is uniformly small. This proves (2) and along with this the claim.
It follows that $\gamma(\infty)\in \partial^{\xi} X$.

 It remains to check that for all $\xi_1, \xi_2\in \partial B$,
$\partial^{\xi_1} X\cap \partial^{\xi_2} X\neq \emptyset$ implies $\xi_1=\xi_2$. 
Suppose $\gamma_i$ is a $\kappa_0$-quasigeodesic ray in $B$ joining $b$ to $\xi_i$, $i=1,2$. 
Suppose $\tilde{\gamma}_i$ is a qi lift of $\gamma_i$, $i=1,2$ such that $\tilde{\gamma}_1(\infty)= \tilde{\gamma}_2(\infty)$,
i.e. $Hd(\tilde{\gamma}_1,\tilde{\gamma}_2)<\infty $. Then $Hd(\gamma_1, \gamma_2)<\infty$ because 
$\pi:X\map B$ is $1$-Lipschitz.  Thus $\xi_1= \xi_2$. This finishes the proof.
\qed

\begin{cor}\label{cor: bundle boundary}
Suppose $F$ is a bounded metric space. Then $\partial X= \coprod_{\xi\in \partial B} \partial^{\xi} X$.
\end{cor}

For instance suppose $\SSS_1, \SSS_2$ are two qi sections and $\LL=\LL(\SSS_1, \SSS_2)$ then by Corollary
\ref{cor: ladders subbundles} there is a metric graph subbundle $\pi_Z:Z\map B$ of $X$ where 
the bundle map $Z\map X$ is a qi embedding onto a finite neighborhood of $\LL$. It follows that $Z$
is hyperbolic and fibers are uniformly quasiisometric to intervals. 
Therefore, the conclusion of Corollary  \ref{cor: bundle boundary} applies to the 
metric bundle $Z$ too. Hence, informally speaking we have the following.

\begin{cor} For any ladder $\LL=\LL(\SSS_1, \SSS_2)$ we have 
$$\partial \LL= \coprod_{\xi\in \partial B}\, \partial^{\xi} \LL.$$
\end{cor}

\begin{lemma}\label{lemma: ct to base}
Suppose $b\in B$ and $\alpha_n:[0,\infty)\map B$ is a sequence of uniform quasigeodesic rays starting 
from $b$. Suppose $\tilde{\alpha}_n$
is a uniform qi lift of $\alpha_n$ for all $n$ such that the set $\{\tilde{\alpha}_n(0)\}$ has finite
diameter. If $\tilde{\alpha}_n(\infty)\map z\in \partial X$ then $\lim_{n\map \infty}\alpha_n(\infty)$
exists. If $\xi=\lim_{n\map \infty}\alpha_n(\infty)$ and $\alpha:[0,\infty)\map B$ is a
$\kappa_0$-quasigeodesic ray joining $b$ to $\xi$ then there is a uniform qi lift 
$\tilde{\alpha}$ of $\alpha$ such that $\tilde{\alpha}(\infty)=z$.
\end{lemma}

\proof Since $\tilde{\alpha}_n(\infty)\map \xi$ there is a constant $D$ such that for all $M>0$ 
there is $N=N(M)>0$ with $Hd(\tilde{\alpha}_m|_{[0,M]}, \tilde{\alpha}_n|_{[0,M]})\leq D$ for all 
$m,n\geq N$ by Lemma \ref{convergence explained}(1). It follows that for all $M>0$,
$Hd(\alpha_m|_{[0,M]}, \alpha_n|_{[0,M]})\leq D$ for all $m,n\geq N$. Hence, again by
Lemma \ref{convergence explained}(1) $\alpha_n(\infty)$ converges to a point of $\xi\in\partial B$. 
Let $\alpha$ be a $\kappa_0$-quasigeodesic ray
in $B$ joining $b$ to $\xi$. We claim $z\in \partial^{\xi} X$. 
Given any $t\in [0,\infty)$ by Lemma \ref{convergence explained}(2) there is $N'=N'(t)\in \NN$ such that
$d(\alpha(t), \alpha_n)\leq D'$ for all $n\geq N'$ where $D'$ depends only on $\kappa_0$ and $\delta$. 
Let $N_0=\max\{N(t),N'(t)\}$. Let $t'$ be such that $d_X(\alpha(t), \alpha_{N_0}(t'))\leq D'$.
Define $\tilde{\alpha}(t)=\phi_{uv}(\tilde{\alpha}_{N_0}(t'))$ where $u=\alpha_{N_0}(t'), v=\alpha(t)$
and $\phi_{uv}$ is a fiber identification map. It is now easy to check that this 
defines a qi section over $\alpha$ and $z=\tilde{\alpha}(\infty)$.

\begin{cor}\label{cor: ct to base}
If  fibers of the metric (graph) bundle are of finite diameter then
the map $\partial X=\cup_{\xi\in \partial B} \partial^{\xi} X\map \partial B$ defined by
sending $\partial^{\xi} X$ to $\xi$ for all $\xi\in \partial B$ is continuous.
\end{cor}


\subsection{Cannon-Thurston lamination} \label{defn: ct lamination}

Suppose $b_0\in B$ is an arbitrary point and $F=F_{b_0}$. Then we know that the inclusion 
$i=i_{F,X}:F\hookrightarrow X$ admits the CT map
$\partial i: \partial F \map \partial X$. For any set $S$ we define
$$S^{(2)}=\{(a,b)\in S\times S: a\neq b\}.$$
Now, following Mitra(\cite{mitra-endlam}) we define the following.

\begin{defn} $(1)$  {\em (Cannon-Thurston lamination)}
Let $\Lam{X}(F)=\{(\alpha,\beta)\in \partial^{(2)}F:\partial i(\alpha)=\partial i (\beta)\}$.

$(2)$ Suppose $\xi\in \partial B$. Let $\Lam{\xi, X}(F)=\{(\alpha,\beta)\in \partial^{(2)}F:\partial i(\alpha)=\partial i (\beta)\in \partial^{\xi} X\}$. 
{\em We shall denote $\Lam{\xi, X}(F)$ simply by $\Lam{\xi}(F)$ when $X$ is understood.}
\end{defn}

In this subsection we are going to discuss the various properties of the CT lamination. First we need some 
definitions. We recall that for all $b,s\in B$ we have the fiber identification map $\phi_{bs}:F_b\map F_{s}$ 
which is a uniform quasiisometry depending on $d_B(b,s)$. This induces a bijection 
$\partial \phi_{bs}: \partial F_b\map \partial F_{s}$. Suppose $z\in \partial F_b$.
Let $z_s=\partial \phi_{bs}(z)$ for all $s\in B$. 

\begin{convention}
 For the rest of the subsection by `quasigeodesic rays' or `lines', we shall always mean 
$\kappa_0$-quasigeodesic rays and lines in the fibers of a metric (graph) bundle unless otherwise specified,
\end{convention}

\begin{defn}(1) {\bf (Semi-infinite ladders)} Suppose $\SSS_1$ is a qi section over $B$ in $X$.  
For all $s\in B$ let $\gamma_s\subset F_s$ be a (uniform) quasigeodesic ray joining $\SSS_1\cap F_s$ to 
$z_s=\partial \phi_{bs}(z)$. 
The union of all the rays will be denoted by $\LL(\SSS_1;z)$.

{\em This set is coarsely well-defined by Lemma \ref{ideal triangles are slim}.
 We shall refer to this as the {\em semi-infinite ladder} defined by $\SSS_1$ and $z$.}

(2) {\bf (Bi-infinite ladders)} Suppose $b\in B$ and $z,z'\in \partial F_b$, $z\neq z'$.
Now for all $s\in B$ join $z_s=\partial \phi_{bs}(z)$ to $z'_s=\partial \phi_{bs}(z')$ 
by a (uniform) quasigeodesic line in $F_s$. 
The union of all these lines will be denoted by $\LL(z;z')$. 

{\em As before, this set is coarsely well-defined by Lemma \ref{ideal triangles are slim}.
We shall refer to this as the {\em bi-infinite ladder} defined by $z$ and $z'$.}
\end{defn}

We shall refer to either of these ladders as an `infinite girth ladder'. 

\begin{lemma}{\bf (Properties of infinite girth ladders)}\label{infinite ladder property}
Suppose $\LL$ is an infinite girth ladder. 
\begin{enumerate}
\item {\em (Coarse retract)}  There is a uniformly coarsely Lipschitz retraction $\pi_{\LL}:X\map \LL$
such that for all $b\in B$ and $x\in F_b$, $\pi_{\LL}(x)$ is a (uniform approximate) nearest point
projection of $x$ in $F_b$ on $\LL\cap F_b$.

Consequently, infinite girth ladders are uniformly quasiconvex and their uniformly small neighborhoods are qi embedded in $X$.

\item {\em (QI sections in ladders)} Through any point of $\LL$, there exists a uniform qi section contained in $\LL$.

\item {\em (QI sections coarsely bisect ladders)}
 Any qi section in $\LL$ coarsely bisects it into two subladders.

\end{enumerate}
\end{lemma}
\proof We shall briefly indicate the proofs comparing with the proof of the analogous results for finite girth ladders.
(3) follows exactly as Lemma \ref{ladder coarse bisection}. (2) is immediate from (1). In fact given $x\in \LL$
one takes a $K_0$-qi section $\SSS$ in $X$ containing $x$ and then $\pi_{\LL}(\SSS)$ is the required qi section.
Therefore, we are left with proving (1). This is an exact analog of Proposition \ref{ladders are qi embedded}(1).
The reader is referred to \cite[Theorem 4.6]{mitra-endlam} for supporting arguments. \qed

\begin{convention}
All semi-infinite ladders $\LL(\SSS;z)$ are formed by $K_0$-qi section $\SSS$.
We shall assume that through any point of an infinite girth ladder there is a $\bar{K}_0$-qi section
contained in the ladder. Also, all infinite girth ladders are assumed to be $\bar{\lambda}_0$-quasiconvex.
\end{convention}

\subsubsection{Properties of the CT lamination $\Lam{X}(F)$}
In this subsection, we prove many properties of the CT lamination using coarse bisection of ladders by qi sections. 
These are motivated by analogous results proved in \cite{mitra-endlam} and \cite{bowditch-stacks}.
For the rest of the subsection, we will use the following set up.
Let $b_0\in B$ and $F=F_{b_0}$. Suppose $(z_1, z_2)\in \Lam{} = \Lam{X}(F)$ and $\LL=\LL(z_1;z_2)$. 
Let $\gamma:\RR\map F$ be a $\kappa_0$-quasigeodesic line in $F$ joining $z_1$ to $z_2$ such that $Im(\gamma)=\LL\cap F$.
Let $i_{F,X}:F\map X$ denote the inclusion map 
and $\partial i_{F,X}:\partial F\map \partial X$ denote the CT map.

\begin{lemma}\label{lamination lemma1}
 Suppose $\Sigma$ is any qi section contained in $\LL$. Then $\partial i_{F,X}(z_i)\in \Lambda(\SSS)$, $i=1,2$.
\end{lemma}
\proof Let $\SSS$ be a qi section contained in $\LL$. Then $\SSS$ coarsely separates $\LL$ in $X$
into $\LL_1=\LL(\SSS; z_1)$ and $\LL_2=\LL(\SSS; z_2)$. We note that 
$\partial i_{F,X}(z_1)=\partial i_{F,X}(z_2)\in \Lambda(\LL_1)\cap \Lambda(\LL_2)$.
Hence we are done by Lemma \ref{coarse separation vs limit set}. \qed

\begin{lemma}\label{CT lamination}
Suppose $(z_1,z_2)\in \Lam{X}(F)$ and $\LL=\LL(z_1;z_2)$. There is a unique $\xi\in \partial B$ such that $(z_1,z_2)\in \Lam{\xi, X}(F)$.
Moreover, for any $\kappa_0$-quasigeodesic $\beta:[0,\infty)\map B$ joining $b_0$
to $\xi$ and any qi section $\Sigma$ contained in $\LL$, if $\tilde{\beta}$ is the lift of $\beta$ in $\Sigma$ then 
$\tilde{\beta}(\infty)= \partial i_{F,X}(z_1)= \partial i_{F,X}(z_2)$.

In particular $\Lam{X}(F)=\coprod_{\xi\in \partial B} \Lam{\xi}(F)$.
\end{lemma}

\proof 
Let $\sigma: B\map X$ be a qi section with image $\Sigma$ contained in $\LL$. By Lemma \ref{lamination lemma1} 
$\partial i_{F,X}(z_1)\in \Lambda(\SSS)$. But $\Lambda(\SSS)=\partial \sigma(\partial B)$ by Lemma \ref{CT-limset}.
Hence, there is a $\kappa_0$-quasigeodesic ray $\beta:[0,\infty)\map B$ such that $\partial \sigma(\beta(\infty))= \partial i_{F,X}(z_1)$.
Let $\xi=\beta(\infty)$. If $\tilde{\beta}=\sigma\circ \beta$ then $\tilde{\beta}$ is a qi lift of $\beta$ 
and $\partial i_{F,X}(z_1)=\tilde{\beta}(\infty)\in \partial^{\xi} X$. Thus $(z_1,z_2)\in \Lam{\xi, X}(F)$.
This shows the existence of $\xi$. Thus we have $\Lam{X}(F)=\bigcup_{\xi\in \partial B} \Lam{\xi}(F)$.
Also for $\xi, \xi'\in \partial B$, $\xi\neq \xi'$ we have $\partial^{\xi_1} X\cap \partial^{\xi_2} X=\emptyset$ by Proposition
\ref{bundle boundary} which immediately implies $\Lam{\xi, X}(F)\cap \Lam{\xi', X}(F)=\emptyset$.
This shows that the point $\xi$ is independent of the chosen section $\Sigma$ in $\LL$. The last part of the lemma is immediate
from these observations. \qed

\begin{comment}
\begin{rem}
From the proof of Lemma \ref{CT lamination} it follows that the point $\xi$ is unique since $\partial \sigma$ is injective; 
it is also independent of the $\sigma$ chosen. In particular any two qi lifts of $\beta$ contained in two qi sections 
in $\LL$ are asymptotic.
\end{rem}\end{comment}
\smallskip

We next aim to show that the sets $\Lam{\xi,X}(F)$ are closed subsets of $\Lam{X}(F)$.
Let $\beta:[0,\infty)\map B$ be a continuous, arc length parameterized $\kappa_0$-quasigeodesic in $B$ 
with $\beta(0)=b_0$ and $\beta(\infty)=\xi$
as in the proof of Lemma \ref{CT lamination}. Let $A=\beta([0,\infty))$. Let $Y=\pi^{-1}(A)$ be the restriction of the bundle
$X$ over $A$. Let $i_{Y,X}:Y\map X$, $i_{F,Y}:F\map Y$ be inclusion maps.

\begin{lemma}\label{CT lamination revealed}
If $(z_1,z_2)\in \Lam{\xi, X}(F)$ then
$\partial i_{F,Y}(z_1)=\partial i_{F,Y}(z_2)$, i.e. $(z_1,z_2)\in \Lam{\xi, Y}(F)$.
\end{lemma}
\proof Let $\SSS_n$ be any qi section in $\LL$ over $B$ passing through $\gamma(n)$, $n\in \ZZ$.
Then by Lemma \ref{CT lamination}, $\SSS_m\cap Y$ and $\SSS_n\cap Y$ are asymptotic for all $m,n\in \ZZ$
in $X$. Since $Y$ is properly embedded in $X$ by Lemma \ref{prop-embed} they are still asymptotic in $Y$.
Clearly $d_Y(\gamma(0), \SSS_n\cap Y)\map \infty$ as $n\map \pm \infty$.
Thus by Lemma \ref{convergence explained}(1) 
$\lim_{n\map \pm \infty} \gamma(n)=\tilde{\beta}_0(\infty)$ in $Y$ 
where $\tilde{\beta}_0$ is the lift of $\beta$ in $\SSS_0$. This completes the proof.
\qed
\begin{comment}
\begin{cor}
Since $\LL\cap Y$ is qi embedded in $Y$
it follows that $z_1, z_2$ are identified under the CT map $\gamma\map \LL$.
\end{cor} \end{comment}

\begin{cor}\label{cor: asymptotic lifts}
Let $\tilde{\beta}$ be any qi lift of $\beta$ in $\LL$.
Then $\tilde{\beta}(\infty)= \partial i_{F,X}(z_1)$. In particular any two qi lifts of $\beta$ in $\LL$ are asymptotic.
\end{cor}
\proof 
We know that $\tilde{\beta}$ coarsely separates $\LL\cap Y$ into two semi-infinite ladders, $\LL^{+}$ and
$\LL^{-}$ in $Y$. It follows that $\Lambda(\LL^{+})\cap \Lambda(\LL^{-})=\Lambda(\tilde{\beta})=\tilde{\beta}(\infty)$.
It then follows that the limit of $\gamma(n)$ in $\partial \LL$ is $\tilde{\beta}(\infty)$.
\qed

\begin{cor}
(1) $\partial (\LL\cap Y)$ is a point. (2) $\Lambda_Y(\LL\cap Y)$ is a point. (3) $\Lambda_X(\LL\cap Y)$ is a point.
\end{cor}
\proof We know by Proposition \ref{infinite ladder property}(1) (see also Proposition \ref{ladders are qi embedded}(4))
that a small neighborhood, say $\LL'_Y =N_R(\LL\cap Y)$, of $\LL\cap Y$ in $Y$ is qi embedded in $Y$ and hence it
is a hyperbolic metric space by its own right. Also, this is a subbundle of $Y$ by Corollary \ref{cor: ladders subbundles}. 

(1) The first part is an informal way of saying that $\partial (\LL'_ Y)$ is a point. 
However, this is immediate from Proposition \ref{bundle boundary} and Corollary \ref{cor: asymptotic lifts}. 

(2)  By Lemma \ref{CT-limset} $\Lambda_Y(\LL'_Y)$ is the image of the CT map for the inclusion $\LL'_Y\map Y$
since $\LL'_Y $ is qi embedded in $Y$. But $\partial \LL'_Y$ is a point by the first part.
Thus $\Lambda_Y(\LL'_ Y)$ is a singleton. Finally, $\Lambda_Y(\LL'_ Y)=\Lambda_Y(\LL\cap Y)$ by Lemma
\ref{hausdorff limset}. Hence we are done.

(3) Lastly, it follows that $\LL\cap Y$ is quasiconvex in $X$ too since by Corollary \ref{cor: asymptotic lifts}
$\LL\cap Y$ is the union of qi lifts of $\beta$ contained in $\LL\cap Y$ all of which converge to the same point of $\partial X$. 
Hence $\LL'_Y$ is also quasiconvex in $X$. Since $Y$ is properly embedded in $X$ by Lemma \ref{prop-embed} and $\LL'_Y$ is qi embedded in $Y$
it follows that $\LL'_Y$ is properly embedded in $X$. Thus $\LL'_Y$ is qi embedded in $X$ by Lemma \ref{qc vs qi emb}(2).
As in (2) we are done by Lemma \ref{CT-limset}. 
\qed

\begin{cor}\label{cor: lamination closed}
We have $\Lam{Y}(F)=\Lam{\xi, Y}(F)=\Lam{\xi,X}(F)$. 

In particular, each $\Lam{\xi, X} (F)$ is a closed subset of $\partial^{(2)} F$.
\end{cor}

\proof  The first equality follows from Lemma \ref{CT lamination} applied to the metric bundle $Y$ over $A$.
We will now prove the second one. Since $\partial i_{F,X}=\partial i_{Y,X}\circ \partial i_{F,Y}$,
clearly $\Lam{\xi, Y}(F)\subset \Lam{\xi,X}(F)$. The opposite inclusion is an immediate
consequence of Lemma \ref{CT lamination revealed}. 

Since $\partial i_{F,Y}$ is continuous it follows that $\Lam{\xi, X} (F)$ is a closed subset of $\partial ^{(2)} F$.
One has to use the standard fact that the Gromov boundaries are Hausdorff spaces.
\qed

The following three results are motivated by similar results proved in \cite{mitra-endlam}.
The proof ideas are very similar. However, we get rid of the group actions that were there and in our setting
properness is never needed.

\begin{comment}
\begin{lemma}
Suppose $\xi_1\neq \xi_2\in \partial B$. If $(z_i,w_i)\in \Lam{\xi_i}$, $i=1,2$ 
then $\{z_1,w_1\}\cap \{z_2, w_2\}=\emptyset$. 
In particular, $\Lam{\xi_1}\cap\Lam{\xi_2}=\emptyset$.
\end{lemma}

\proof Suppose $\gamma_i$ is a $\kappa_0$-quasigeodesic ray in $B$ joining $b$ to $\xi_i$, $i=1,2$. 
Suppose $\tilde{\gamma}_i$ is a qi lift of $\gamma$ such that $\partial i(z_i)=\tilde{\gamma}_i(\infty)$. 
Since $\tilde{\gamma}_1(\infty)= \tilde{\gamma}_2(\infty)$ if and only if $\tilde{\gamma}_i$'s are asymptotic in
which case $\gamma_i$'s would also be asymptotic because $\pi:X\map B$ is $1$-Lipschitz. This would
be a contradiction since $\xi_1\neq \xi_2$.
\qed \end{comment}

\begin{defn}
Suppose $Z_1, Z_2$ are hyperbolic metric spaces. Suppose $f:Z_1\map Z_2$ is a metrically proper map that
admits the CT map. If $\gamma\subset Z_1$ is a quasigeodesic line such that $\partial f(\gamma(\infty))= \partial f(\gamma(-\infty))$
then we refer to $\gamma$ as a {\em leaf} of the CT lamination $\Lam{Z_2}(Z_1)$.
\end{defn}
We recall that in our context the quasigeodesic lines are assumed to be $\kappa_0$-quasigeodesic lines.

\begin{lemma}\label{transversality}
Suppose $\xi_1\neq \xi_2\in \partial B$. Given $D>0$  there exists $R=R_{\ref{transversality}}(D)>0$
such that the following holds:

Suppose $\gamma_1$ is a leaf of $\Lam{\xi_1, X}(F)$ and $\gamma_2$ is a leaf of $\Lam{\xi_2,X}(F)$. 
Then $\gamma_1\cap N_D(\gamma_2)$ has diameter less than $R$. 
\end{lemma}
\proof Let $\alpha$ be a $\kappa_0$-quasigeodesic line in $B$ joining $\xi_1, \xi_2$.
Let $b'_0\in \alpha$ be a nearest point projection of $b_0$ on $\alpha$.
Let $c$ be a geodesic in $B$ joining $b_0$ to $b'_0$. Let $\alpha_i$ be the concatenation of
$c$ with the portion of $\alpha$ joining $b'_0$ to $\xi_i$, $i=1,2$. 
We note that $\kappa_0$-quasigeodesics in $B$ are $D_{\ref{stab-qg}}(\delta_0, \kappa_0,\kappa_0)$-quasiconvex by stability of quasigeodesics.
Let $K=D_{\ref{stab-qg}}(\delta_0, \kappa_0,\kappa_0)$.
Hence, $\alpha_i$'s are $K_{\ref{subqc-elem}}(\delta_0,K,\kappa_0, 1)$-quasigeodesics by Lemma \ref{subqc-elem}(2). 
Let $k=K_{\ref{subqc-elem}}(\delta_0,K,\kappa_0, 1)$.

Next suppose $x_i, x'_i\in \gamma_i$, $i=1,2$ are such that $d_F(x_1, x_2)\leq D$ and $d_F(x'_1, x'_2)\leq D$.
Let $\SSS_i, \SSS'_i$ be two qi sections in each $\LL_i=\LL(\gamma_i(\infty), \gamma_i(-\infty))$ passing through
$x_i$ and $x'_i$ respectively, $i=1,2$.
Let  $\tilde{\alpha}_{i}$ and $\tilde{\alpha}'_{i}$ be lifts of $\alpha_i$ in $\LL_i$ through $x_i$
and $x'_i$ respectively for $i=1,2$. We now look at the quasigeodesic hexagon in $X$ with vertices
$x_i, x'_i, \xi_i$, $i=1,2$ where $\tilde{\alpha}_{i}$'s and $\tilde{\alpha}'_{i}$'s form four sides
and the other two sides are formed by geodesics joining $x_1$ to $x_2$ and $x'_1$ to $x'_2$ respectively.
We note that the infinite sides of this polygon are all $(k\bar{K}_0+k+\bar{K}_0)$-quasigeodesics. 
Let $\tilde{k}= k\bar{K}_0+k+\bar{K}_0$.
Hence, such a hexagon is $R_{\ref{ideal polygons are slim}}(\delta,\tilde{k}, 6)$-slim by 
Corollary \ref{ideal polygons are slim}. Let $R_1=R_{\ref{ideal polygons are slim}}(\delta,\tilde{k}, 6)$.
Let $b_2$ be a point on $\alpha_2$ such that $d_B(b_2, \alpha_1)= D+R_1+1=R$, say and let $y_2=\tilde{\alpha}_2(b_2)$.
Then $y_2\in N_{R_1}(\tilde{\alpha}'_2)$. In particular, $y_2\in N_R(\SSS'_2)$. Hence, by 
Lemma \ref{distance from qi section} 
$d_{b_2}(\SSS_2\cap F_{b_2}, \SSS'_2\cap F_{b_2})\leq R_{\ref{distance from qi section}}(\bar{K}_0, R)$.
It follows by bounded flaring that $d_{b_0}(x_2,x'_2)\leq \mu_{\tilde{k}}(R_{\ref{distance from qi section}}(\bar{K}_0, R))$.
\qed

\begin{lemma}\label{CT property final}
If $\xi_n\map \xi$ in $\partial B$, $(z_n, w_n)\in \Lam{\xi_n, X}(F)$ and $(z_n, w_n)\map (z,w)\in \partial^{(2)} F$.
Then $(z,w)\in \Lam{\xi, X}(F)$. 
\end{lemma}

\proof Since $\partial i_{F,X}(z_n)=\partial i_{F,X}(w_n)$ for all $n$ and $\partial i_{F,X}$ is
continuous it follows that $\partial i_{F,X}(z)=\partial i_{F,X}(w)$ whence $(z,w)\in \Lam{X} F$.
Let $[z_n,w_n], [z_n,z], [w_n,w]$ and $[z,w]$ denote $\kappa_0$-quasigeodesic lines in $F$
joining these pairs of points. Let $x\in [z,w]\cap F$ and let $\alpha$ be a $\kappa_0$-quasigeodesic 
ray in $B$ joining $b$ to $\xi$. 

{\bf Claim:} There is a uniform qi lift $\tilde{\alpha}$ of $\alpha$ through $x$ such that 
$\tilde{\alpha}(\infty)=\partial i_{F,X}(z)=\partial i_{F,X}(w)$.

\noindent
Since $z_n\map z$ and $w_n\map w$ by Lemma
\ref{convergence explained}(1), we have $d_{b_0}(x, [z_n, z])\map \infty$ and $d_{b_0}(x, [w_n, w])\map \infty$.
Hence, by Corollary \ref{ideal polygons are slim} there is $N\in \NN$ such that 
$d_{b_0}(x, [z_n,w_n])\leq R=R_{\ref{ideal polygons are slim}}(\delta_0, \kappa_0, 4)$ for all $n\geq N$.
Now, let $x_n\in [z_n, w_n]$ such that $d_{b_0}(x,x_n)\leq R$.
Let $\alpha_n$ be a $\kappa_0$-quasigeodesic ray in $B$ joining $b$ to $\xi_n$. 
Then by Corollary \ref{cor: asymptotic lifts} we know that there is a uniform qi lift $\tilde{\alpha}_n$ of each $\alpha_n$,
$n\geq N$ such that $\tilde{\alpha}_n(0)=x_n$ and  $\tilde{\alpha}_n(\infty)=\partial i_{F,X}(z_n)$.
Hence, by Lemma \ref{lemma: ct to base} and Lemma \ref{bdry 6.1} there is a qi lift $\tilde{\alpha}$
starting from $x$ such that $\tilde{\alpha}(\infty)=\partial i_{F,X}(z)=\partial i_{F,X}(w)$. This proves the claim.

However, this means that $\partial i_{F,X}(z)=\partial i_{F,X}(w)\in \partial^{\xi} X$. Therefore,
$(z,w)\in\Lam{\xi, X}(F)$. \qed

\subsubsection{Leaves of CT laminations for pullback bundles}
The following result is motivated by a similar result proved in \cite{ps-kap} for trees of hyperbolic 
spaces which in turn was suggested by Mahan Mj. We gratefully acknowledge the same.

{\em Suppose we have the hypotheses of Theorem \ref{CT for mgbdl}. 
We identify $Y$ as a subspace of $X$ and $A$ as a subspace of $B$. 
Similarly, $\partial A$ is identified as a subset of $\partial B$. 
With that in mind, we have the following:}

\begin{theorem}\label{CT leaf in Y}
Suppose we have a metric graph bundle satisfying the hypotheses of Theorem \ref{CT for mgbdl} such that
the fibers of the bundle are all proper metric spaces. Suppose $\gamma$ is a quasigeodesic line in $Y$ 
such that $(\gamma(\infty),\gamma(-\infty))\in \Lam{X}(Y)$. Let $F=F_b$ be any fiber of $Y$.

Then (1) $\gamma(\pm \infty)\in \partial i_{F,Y}(\partial F)$. 

(2) There is a point $\xi\in \partial B\setminus \partial A$ determined by $\gamma(\pm \infty)$ such that if $
z_{\pm}\in \partial F$ with $\partial_{F,Y}(z_{\pm})=\gamma(\pm)$ then $(z_+, z_-)\in \Lam{\xi, X}(F)$.

(3) $\pi(\gamma)$ is bounded. Moreover, $\gamma$ is within a finite Hausdorff distance 
from a $\kappa_0$-quasigeodesic line  $\sigma$ of $F$ so that $\partial i_{F,Y}(\sigma(\pm\infty))=\gamma(\pm\infty)$.
Also, $(\sigma(\infty),\sigma(-\infty))\in \Lam{\xi, X}(F)$ for some $\xi\in \partial B\setminus \partial A$.

(4) If $b$ is a nearest point projection of $\xi$ on $A$. Then $\sigma$ (as defined in (3))
is a uniform quasigeodesic line in $Y$.
\end{theorem}

\proof We have $\partial Y=\Lambda_Y(F)\cup(\cup_{\xi\in \partial A} \partial^{\xi} Y)$ by Proposition \ref{bundle boundary}.
Also since $F$ is a proper metric space, by Lemma \ref{CT-limset} $\Lambda_Y(F)=\partial i_{F,Y}(\partial F)$.
Thus $\partial Y=\partial i_{F,Y}(\partial F)\cup(\cup_{\xi\in \partial A} \partial^{\xi} Y)$. We shall use the following
observation a few times in the proof which are immediate from the fact that $A$ is qi embedded in $B$.

{\em Suppose $\alpha$ is a quasigeodesic ray in $A$ and $\tilde{\alpha}$ is a qi lift
of $\alpha$ in $Y$. Then $\tilde{\alpha}$ is a quasigeodesic ray in $Y$ as well as in $X$.
Also any pair of such rays are asymptotic in $Y$ if and only if they are asymptotic in $X$ since 
$Y$ is properly embedded in $X$.}

(1) The proof of this assertion is by elimination of the possibilities coming from the decomposition 
$\partial i_{F,Y}(\partial F)\cup(\cup_{\xi\in \partial A} \partial^{\xi} Y)$ of $\partial Y$.

Suppose $\gamma(\infty)\in \partial^{\xi_1} Y$ and $\gamma(-\infty)\in \partial^{\xi_2} Y$ for 
some $\xi_1, \xi_2\in \partial A$. However, this case is not possible due to the above observation.

Suppose $\gamma(\infty)\in \partial^{\xi} Y$ for some $\xi\in \partial A$ and 
$\gamma(-\infty)\in \partial i_{F,Y}(\partial F)\setminus \cup_{\xi\in \partial A} \partial^{\xi} Y$ or vice versa.
We show below that this case is also not possible.

Let $\alpha$ be a $\kappa_0$-quasigeodesic ray in $A$ joining $b$ to $\xi$ and let $\tilde{\alpha}$ be a $K_0$-qi lift
of $\alpha$ in $Y$ such that $\tilde{\alpha}(\infty)=\gamma(\infty)$. Also let $\beta$ be a $\kappa_0$-quasigeodesic ray
in $F$ such that $\partial i_{F,Y}(\beta(\infty))=\gamma(-\infty)$. 
Now, for all $n\in \NN$ let $\SSS_n$ be a $K_0$-qi section in $X$ passing through $\beta(n)$
and let $\LL_n=\LL(\SSS_n, \beta(\infty))$. Then $\LL_n$ is $\bar{\lambda}_0$-quasiconvex in $X$.
Clearly $\gamma(\infty)=\tilde{\alpha}(\infty)\in \Lambda_X(\LL_n)$. Hence, by Lemma \ref{limset lem}
$\tilde{\alpha}$ is asymptotic to $\LL_n$. 
It follows by Proposition \ref{ladders are qi embedded} and Lemma \ref{distance from ladder}
that $\pi_{\LL_n}(\tilde{\alpha})$ is a uniform qi lift of $\alpha$ and it is asymptotic to
$\tilde{\alpha}$. Since $Y$ properly embedded in $X$ by Lemma \ref{prop-embed}, it follows that these
qi lifts are asymptotic in $Y$ too. In particular, $\pi_{\LL_n}(\tilde{\alpha})(\infty)=\gamma(\infty)$.
Now, since $d_F(\beta(0), \beta(n))\map \infty$, by Lemma \ref{distance from ladder}
$d_Y(\beta(0),\pi_{\LL_n}(\tilde{\alpha}))\map \infty$. It follows from Lemma \ref{convergence explained}
that $\lim_{n\map \infty} \beta(n)=\gamma(\infty)$ in $\partial Y$. This gives a contradiction since
$\lim_{n\map \infty} \beta(n)=\gamma(-\infty)\neq \gamma(\infty)$.

Therefore, the only possibility is that 
$$\gamma(\pm \infty)\in \partial i_{F,Y}(\partial F) \setminus \cup_{\xi\in \partial A} \partial^{\xi} Y$$
proving part (1) of the theorem.

Let $z,z'\in \partial F$ be such that $\partial i_{F,Y}(z)=\gamma(\infty)$ and $\partial i_{F,Y}(z')=\gamma(-\infty)$.

(2) Since $\partial i_{F,X}=\partial i_{Y,X}\circ \partial i_{F,Y}$ by Lemma \ref{CT properties}(1), we have $(z,z')\in \Lam{X}(F)$ and hence $(z,z')\in \Lam{\xi, X}(F)$ for some $\xi\in \partial B$ by Lemma \ref{CT lamination}. From Corollary \ref{cor: lamination closed} 
it follows that $\xi\in \partial B\setminus \partial A$. This proves part (2) of the theorem.

(3) Let $\LL=\LL(z;z')$ be the bi-infinite ladder in $X$ formed by $z,z'$. Let $\sigma=\LL\cap F$ which is
an arc length parameterized $\kappa_0$-quasigeodesic line in $F$ joining $z,z'$.
Let $\alpha$ be a $\kappa_0$-quasigeodesic ray in $B$ joining $b$ to $\xi$.

Let $\Sigma_n$ be a $\bar{K}_0$-qi section in $\LL$ passing through $\sigma(n)$, $n\in \NN$. 
By Corollary \ref{cor: asymptotic lifts} qi lifts of $\alpha$ 
contained in these qi sections are asymptotic. Denote the qi section of $\alpha$ contained in $\SSS_n$
by $\tilde{\alpha}_n$. We note that these are $k=(\bar{K}_0\kappa_0+\bar{K}_0+\kappa_0)$-quasigeodesics by Lemma
\ref{qi composition}(2).
Hence, by Lemma \ref{ideal triangles are slim} given $m, n\in \NN$ we have
$\tilde{\alpha}_n(i)\in N_R(\tilde{\alpha}_{-m})$ (and $\tilde{\alpha}_{-m}(i)\in N_R(\tilde{\alpha}_n)$)
where $R=D_{\ref{ideal triangles are slim} }(\delta, k)$ as long as $\tilde{\alpha}_n(i)$ (resp.
$\tilde{\alpha}_{-m}(i)$) is not contained in the $R$-neighborhood of any $1$-quasigeodesic joining
$\sigma(-m), \sigma(n)$. In particular for such $i$ we have
$\tilde{\alpha}_n(i)\in N_R(\SSS_{-m})$, $\tilde{\alpha}_{-m}(i)\in N_R(\SSS_n)$.
Hence, by Lemma \ref{distance from qi section} we have 
$$d_{\alpha(i)}(\tilde{\alpha}_n(i), \tilde{\alpha}_{-m}(i))\leq R_1=R_{\ref{distance from qi section}}(R,\bar{K}_0)$$
for all such $i$. Let $R_2=\max\{R_1, M_{\bar{K}_0}\}$.
Thus for all $n\in \NN$, $U_{n}=U_{R_2}(\SSS_n, \SSS_{-n})\neq \emptyset$. Let $b_n\in U_n$ be a 
nearest point projection of $b$ on $U_n$ and let $b'_n$ be a nearest point projection of $b_n$ on $A$.
Then it follows from Lemma \ref{small-girth-ladder} 
that the concatenation of  
the segments of $\tilde{\alpha}_n, \tilde{\alpha}_{-n}$ over the portion of $\alpha$ joining $b, b'_n$
and the fiber geodesic segment $\LL\cap F_{b'_n}$ is a uniform quasigeodesic in $Y$ joining $\sigma(\pm n)$.
Call it $\gamma'_n$.
Since $\lim_{n\map \infty} \sigma(n)\neq \lim_{n\map \infty} \sigma(-n)$ in $Y$ there is a constant $D\geq 0$ 
such that $d_Y(\sigma(0),\gamma'_n)\leq D$ by Lemma \ref{gromov product meaning}. 
We claim that this means $d_B(b, b'_n)$ is bounded. In fact 
$d_Y(\sigma(0),\tilde{\alpha}_{\pm n})\map \infty$ by Lemma  \ref{distance from qi section}. Thus for all large
$n$ we have $d_Y(\sigma(0), \LL\cap F_{b'_n})\leq D$ whence $d_B(b, b'_n)\leq D$.
It follows from Proposition \ref{ladders are qi embedded}(3) that the Hausdorff distance of $\LL\cap F_{b'_n}$
and the segment of $\sigma$ between $\sigma(n)$ and $\sigma(-n)$ is at most
$(1+2K_0)C_{\ref{ladders are qi embedded}}(\bar{K}_0)$. Since $\sigma$ is a proper embedding in $Y$
it follows by Lemma \ref{quasigeod criteria} that $\sigma$ is a uniform quasigeodesic in $Y$ depending on $D$.
Let $K\geq 1$ be such that both $\sigma$ and $\gamma$ are $K$-quasigeodesics in $Y$. Then,
since $Y$ is $\delta'$-hyperbolic, $Hd(\sigma, \gamma)\leq R_{\ref{ideal polygons are slim}}(\delta', K, 2)$. Thus
$diam(\pi(\gamma))\leq R_{\ref{ideal polygons are slim}}(\delta', K, 2)$. 

We note here that $diam(\pi(\gamma))$ as well as the quasigeodesic constant of $\sigma$
depends only on $\max\{d_B(b,b'_n)\}$.

(4) Use shall use the notation of the proof of (3). Thus we know that there is $D_n\geq 0$ such that for all $i\geq D_n$
we have $d_{\alpha(i)}(\tilde{\alpha}_n(i), \tilde{\alpha}_{-n}(i))\leq R_1$ whence
$\alpha(i)\in U_n$ for all $i\geq D_n$. Also we know that the sets $U_n$ are 
$K_{\ref{qc-level-set-new}}(\bar{K}_0)$-quasiconvex in $B$ by Lemma \ref{qc-level-set-new}.
Let $t_n\geq \max\{D_n, d_B(b,b_n)\}$. Then $\alpha(t_n)\in U_n$. Thus $[b,b_n]_B*[b_n. \alpha(t_n)]_B$ is a
$K_{\ref{subqc-elem}}(\delta_0, K_{\ref{qc-level-set-new}}(\bar{K}_0),\kappa_0,\epsilon)$-quasigeodesic segment.
Let $K'=\max\{\kappa_0, K_{\ref{subqc-elem}}(\delta_0, K_{\ref{qc-level-set-new}}(\bar{K}_0),\kappa_0,\epsilon) \}$. 
Hence, by stability of quasigeodesics (Lemma \ref{stab-qg}) we get that 
$b_n\in N_{R'}(\alpha)$ where $R'=D_{\ref{stab-qg}}(\delta_0, K', K')$. 
We also note that $d_B(b,b_n)\map \infty$ by the bounded flaring condition (Lemma \ref{bdd-flaring})
since $d_b(\tilde{\alpha}_n(0), \tilde{\alpha}_{-n}(0))\map \infty$. This implies that $b_n\map \xi$.
Hence by Lemma \ref{qc last lemma} there exists $N>0$ such that $d(b'_n, b)\leq R_{\ref{qc last lemma}}(\delta_0, k_0)$
for all $n\geq N$ since $B$ is $\delta_0$-hyperbolic and $A$ is $k_0$-quasiconvex. Hence, we are done by
the note left at the end of the proof of (3). \qed

\smallskip
{\bf Surjectivity of the CT maps}

\smallskip
\begin{theorem}\label{thm:CT surjective}
Suppose we have the hypotheses of Theorem \ref{CT for mgbdl} such that the fibers of the bundle are proper metric spaces. 
Let $F$ be the fiber over a point $b\in A$. Suppose the CT map 
$\partial i_{F,X}:\partial F\map \partial X$
is surjective. Then the CT map $\partial i_{F,Y}:\partial F\map \partial Y$ is also surjective.

Conversely for any geodesic ray $\alpha:[0,\infty)\map B$ with $\alpha(0)=b$ let $Y_{\alpha}=\pi^{-1}(\alpha)$. 
If for all $z\in \partial B$ and for some (any) geodesic ray $\alpha$ joining $b$ to $z$ the CT map 
$\partial_{F,Y_{\alpha}}:\partial F\map \partial Y_{\alpha}$ is surjective then the CT map $\partial_{F,X}:\partial F\map \partial X$
is also surjective. 
\end{theorem}

\proof Let $\xi\in \partial Y$. We want to show that $\xi\in Im(\partial i_{F,Y})$.
Since $\partial i_{F,X}:\partial F\map \partial X$ is surjective 
there exists $z\in \partial F$ such that $\partial i_{F,X}(z)=\partial i_{Y,X}(\xi)$. If $\partial i_{F,Y}(z)= \xi$ we are done.
Suppose not. However, $\partial i_{F,X}=\partial i_{Y,X}  \circ \partial i_{F,Y}$. 
Hence, $\partial i_{Y,X}(\partial i_{F,Y}(z))=\partial i_{Y,X}(\xi)$. Then by Theorem \ref{CT leaf in Y}(3) we are done.

The converse part is a direct consequence of Corollary \ref{cor: bundle bdry} and Proposition \ref{bundle boundary}. \qed
\smallskip

\begin{cor}
Suppose $\pi:X\map B$ is a metric (graph) bundle such that $X,B$ are hyperbolic and the fibers are all proper, uniformly quasiisometric
to the hyperbolic plane $\mathbb H^2$. Then for all $b\in B$, the CT map $\partial_{F_b,X}:\partial F_b\map \partial X$ is surjective.
\end{cor}
\proof 
This is an immediate consequence of the second part of Theorem \ref{thm:CT surjective} and the following proposition of Bowditch.
\qed

\begin{prop}\textup{(\cite[Proposition 2.6.1]{bowditch-stacks})}
Suppose $\pi:X\map B$ is a metric (graph) bundle where $B=[0,\infty)$, $X$ is hyperbolic and the fibers are all uniformly quasiisometric
to the hyperbolic plane $\mathbb H^2$. Then for all $b\in B$, the CT map $\partial_{F_b,X}:\partial F_b\map \partial X$ is surjective.
\end{prop}

We would like to remark that Bowditch stated the above proposition in case the fibers are all isometric to the hyperbolic plane,
but the same proof goes through for fibers uniformly quasiisometric to the hyperbolic plane.

A special case of the following result was proved by E. Field (\cite[Theorem B]{field}).

\begin{theorem}
Suppose $1\map N\map G\stackrel{\pi}{\map} Q\map 1$ is a short exact sequence of infinite hyperbolic groups. Suppose $A\subset Q$ is qi embedded and $Y=\pi^{-1}(A)$.
Then the CT map $\partial N\map \partial Y$ is surjective.
\end{theorem}
\proof Since $N$ is a normal subgroup of the hyperbolic group $G$ it is a standard fact that $\Lambda(N)=\partial G$. 
Thus by Lemma \ref{CT-limset} the CT map $\partial N\map \partial G$ is surjective. Now we are done by
Corollary \ref{thm:CT surjective}. \qed

\smallskip
{\bf Fibers of the CT maps}

\smallskip

\begin{theorem}\label{main application2}
Suppose $X$ is a metric (graph) bundle over $B$ satisfying the hypotheses of Theorem \ref{CT for mgbdl} such
that $X$ is a proper metric space. Let $F=F_b$ where $b\in B$.
Suppose $\partial F$ is not homeomorphic to a dendrite and also the CT map $\partial F\map \partial X$
is surjective. 

Then for all $\xi\in \partial B$ we have $\Lam{\xi, X}(F)\neq \emptyset$.
\end{theorem}

\proof Suppose $\alpha$ is an arc length parameterized $\kappa_0$-quasigeodesic ray in $B$ joining $b$ to $\xi$. 
Let $Y=\pi^{-1}(\alpha)$. Since the CT map $\partial F\map \partial X$ is surjective, the map
$\partial i_{F,Y}: \partial F\map \partial Y$ is also surjective by Theorem \ref{thm:CT surjective}. 
Now, $\Lam{\xi, X}(F)=\Lam{\xi, Y}(F)$ by Corollary \ref{cor: lamination closed}. Hence, it is enough
to show that $\Lam{\xi, Y}(F)\neq \emptyset$. However, $\Lam{\xi, Y}(F)= \emptyset$ if and only if $\partial i_{F,Y}$
is injective. It follows that $\Lam{\xi, Y}(F)= \emptyset$ if and only if $\partial i_{F,Y}$ is bijective.
Since $X$ is proper, so are $F$ and $Y$. Hence, $\partial F$ and $\partial Y$ are compact metrizable spaces. 
(See \cite[Theorem Proposition 3.7, Proposition 3.21, Chapter III.H]{bridson-haefliger} for instance.)
Hence, $\partial i_{F,Y}$ is bijective implies $\partial i_{F,Y}$ is a homeomorphism between $\partial F$
and $\partial Y$. 
Since $\partial F$ is not a dendrite this is impossible due to the following result of Bowditch. 
Hence, $\Lam{\xi, Y}(F)\neq \emptyset$. \qed

\begin{theorem} \textup{(\cite[Proposition 2.5.2]{bowditch-stacks})}
Suppose $X$ is hyperbolic metric (graph) bundle over $B=[0,\infty)$ satisfying the hypotheses H1-H4 of section 5.
Suppose moreover that $X$ is a proper metric space. Then $\partial X$ is a dendrite.
\end{theorem}

We note that a special case of interest of Theorem \ref{main application2} is when the fibers are uniformly quasiisometric
to the hyperbolic plane. For instance, we have the following.

\begin{cor}
Suppose we have an exact sequence of infinite hyperbolic groups $1\map N\map G\map Q\map 1$ where $N$ is either the
fundamental group of an orientable closed surface of genus $g\geq 2$ or a free group $F_n$ on $n\geq 3$ generators. 
Then for all $\xi\in \partial Q$, $\Lam{\xi, G}(N)\neq \emptyset$.
\end{cor}

\begin{rem} We remark that much stronger results than the above corollary were already proved by 
Mj and Rafi in \cite{mj-rafi}. For instance, see Theorem 3.12, Theorem 5.7 and Proposition 5.8 there.
\end{rem}

Another context is that of complexes of groups where Theorem \ref{main application2} can be applied.
\begin{cor}\label{CT for cplx of gps}
Suppose $G$ is the fundamental group of a finite developable complexes of nonelementary hyperbolic groups $(\GG, \YY)$
with qi condition. Suppose $X$ is the metric bundle over $B$ obtained from this data as constructed in 
Example \ref{example: complex of groups}. Suppose $G$ is hyperbolic.

Then for all $\xi\in\partial B$ and any vertex group $G_v$, $v\in V(\YY)$ we have $\Lam{\xi, G}(G_v)\neq \emptyset$.
\end{cor}
\proof We need to check the hypotheses of Theorem \ref{main application2}. It is a standard fact that the
boundary of a hyperbolic group is not a dendrite. 
Since the fibers of the metric bundle under consideration are quasiisometric to nonelementary
hyperbolic groups $\partial F$ is not a dendrite for any fiber $F$. We also note that the metric
bundle satisfies H1-H4 of section 5. Finally, $G$ acts on $X$ and $B$ so that the map $\pi:X\map B$
is equivariant, the action of $G$ on $X$ is proper and cocompact and on $B$ is cocompact. Thus any
orbit map $G\map X$ is a qi by Milnor-Svarc lemma and therefore induces a homeomorphism $\partial X \map \partial G$.

Now, given any fiber $F$ and $g\in G$, $gF$ is another fiber of the metric bundle.
By Lemma \ref{fibers qi}(1) $Hd(F, gF)<\infty$. Hence, by Lemma \ref{hausdorff limset}
$\Lambda(F)=\Lambda(gF)=g\Lambda(F)$. It is a standard fact that the action of a nonelementary hyperbolic group
on its boundary is minimal, i.e. the only invariant closed subsets are the empty set and the whole set.
Hence, it follows that $\Lambda(F)= \partial X$. By Lemma \ref{CT-limset} we have $\Lambda(F)=\partial i_{F,X} (\partial F)$.
Thus the CT map $\partial i_{F,X}: \partial F\map \partial X$ is surjective. Finally, clearly $X$ is a proper metric
space. Hence, we have $\Lam{\xi, X}(F)\neq \emptyset$ by Theorem \ref{main application2}. Finally since
$G_v$ acts properly and cocompactly on $X_v$, any orbit map $G_v\map X_v$ is a quasiisometry. Hence, this
induces a homeomorphism $\partial G_v\map \partial X_v$. Therefore, taking $F=X_v$ we are done.\qed

\begin{defn}
Suppose $Z$ is any hyperbolic metric space and $S\subset Z$. Then a point $z\in \Lambda(S)\subset \partial Z$
will be called a {\em conical limit point} of $S$ if for some (any) quasigeodesic $\gamma$ converging to $z$ in $Z$ there is
a constant $D>0$ such that $N_D(\gamma)\cap S$ is a subset of infinite diameter in $Z$.
\end{defn}

\begin{prop}\label{conical-singleton}
Suppose we have the hypotheses of Theorem \ref{CT for mgbdl}.
Let $\partial i_{Y,X} : \partial Y \to \partial X$ be the CT map. If $\xi \in \partial X$ is a conical limit point of $Y$, 
then $|\partial i_{Y,X}^{-1}(\xi)| = 1$.
\end{prop} 
\proof 
Suppose $z\neq z'\in \partial Y$ such that $\partial i_{Y,X}(z)=\partial i_{Y,X}(z')=\xi$. 
Then by Theorem \ref{CT leaf in Y} there is $\xi_B\in \partial B\setminus \partial A$ 
and a qi lift of $\gamma$ of a quasigeodesic ray joining $b$ to $\xi_B$ such that 
$\xi=\gamma(\infty)$. Since $\xi_B\in \partial B\setminus \partial A$ and $A$ is quasiconvex
$\xi_B$ is not a limit point of $A$ in $\partial B$. Thus it is
clear that $\xi$ is not a conical limit point of $Y$. This gives a contradiction and
proves the proposition.\qed

\subsection{QI embedding fibers in a product of bundles}
The lemma below is the product of answering a question due to Misha Kapovich. 
\begin{lemma}
 Suppose $\pi: X \map \RR$ is a metric (graph) bundle satisfying the hypotheses of section 5 
and $X^{\pm}$ are the restrictions of it to $[0,\infty)$ and $(-\infty, 0]$ respectively. 
Then the diagonal embedding $f: F_0\map X^{+}\times X^{-}$ is a qi embedding
where the latter is given the $l_2$ metric.
\end{lemma}
\begin{proof}
Without loss of generality, we assume $(X,d)$ is a metric graph bundle.	
Let $d_{\pm}$ be the induced length metric on $X^{\pm}$ respectively.
Then the $l_2$ metric $d_Y$ on $Y:= X^+\times X^-$ is given by 
$d_{Y}((x_1,x_2), (y_1,y_2))^2 = d_{+}(x_1,y_1)^2+d_{-}(x_2,y_2)^2$
for all $x_1, y_1\in X^+$ and $x_2, y_2\in X^-$.
We note that the inclusion maps $F_0\map X^{\pm}$ are $1$-Lipschitz.

Let $x,y \in F_0$. Then, 
$d_Y(f(x),f(y))^2=d_Y((x,x),(y,y))^2 = d_{+}(x,y)^2+d_{-}(x,y)^2\leq d_{0}(x,y)^2+d_{0}(x,y)^2 = 2d_{0}(x,y)^2$, 
which implies that $d_Y(f(x),f(y))\leq \sqrt{2}d_{0}(x,y)$. A reverse inequality is obtained as follows.

Let $\SSS,\SSS'$ be a pair of $K_0$-qi sections in $X$ through $x,y$ respectively. 
Let $\LL=\LL(\SSS,\SSS')$ be the ladder formed by them. Let $\lambda = \LL\cap F_0$. 
This is a geodesic in $F_0$ joining $x,y$. Now, suppose $c(x,y)$ is a uniform quasigeodesic
in $X$ joining $x,y$ constructed as in section 5 by decomposing $\LL$ into subladders
using the the qi sections $\SSS_i$'s and $\SSS'_j$'s. 
Let $\bar{c}_+:=\bar{c}_{+}(x,y), \bar{c}_{-}:=\bar{c}_{-}(x,y)$ 
be the modified paths joining $x,y$ in $X^+,X^{-}$ respectively. By our main theorem in section 5,
$\bar{c}_+, \bar{c}_{-}$ are uniform quasigeodesics in $X^+,X^{-}$ respectively. Suppose these are
$K$-quasigeodesics. As in the discussion at the end of section 5, suppose $\QQQ, \QQQ'$ are consecutive qi
sections in the decomposition of $\LL=\LL(\SSS, \SSS')$ and $z,w\in \QQQ, w'\in \QQQ'$ with 
$b'=\pi(w)=\pi(w')$ are such that $\LL(\QQQ, \QQQ')\cap c(y, y')$ is made of the fiber geodesic
$[w,w']_{b'}$ and the lift of $[\pi(z),\pi(w)]_B$ in $\QQQ$. However, if $b'\in [0,\infty)$
then $\lambda \cap \LL(\QQQ, \QQQ') \subset \bar{c}_-$ and similarly if 
$b'\in (-\infty, 0]$ then $\lambda \cap \LL(\QQQ, \QQQ') \subset \bar{c}_+$.
Thus $\lambda \subset \bar{c}_+\cup \bar{c}_{-}$. Therefore we have, 
\begin{align*}
d_0(x,y) &\leq l_{+}(\tilde{c}_+) + l_{-}(\tilde{c}_-) \leq Kd_{+}(x,y)+K+Kd_{-}(x,y)+K\\
&= K(d_+(x,y)+d_-(x,y))+2K\\ 
&= 2Kd_{Y}((x,x),(y,y))+2K = 2Kd_{Y}(f(x),f(y))+2K.
\end{align*}
Thus, $-1+\frac{1}{2K} d_0(x,y)\leq d_Y(f(x),f(y))\leq \sqrt{2}d_{0}(x,y)$. Hence,
$f$ is $(2K, 1)$-qi embedding.
\end{proof}

In the same way, we obtain the following.

\begin{lemma}
If $v_0$ is a cut point of $B$ and removing it produces two quasiconvex subsets $A_1, A_2$ and $Y_1, Y_2$
are the restrictions of the bundle to $A_1, A_2$ respectively then the diagonal map $F_{v_0}\map Y_1\times Y_2$
is a qi embedding.  
\end{lemma}

\begin{cor}
If $v_0$ is a cut point of $B$ and removing it produces finitely many quasiconvex subsets $A_i$, $1\leq i\leq n$ and $Y_i$'s
are the restrictions of the bundle to $A_i$'s respectively then the diagonal map $F_{v_0}\map \Pi_i Y_i$
is a qi embedding.  
\end{cor}

\begin{rem}
In \cite{mitra-endlam} Mitra defined an ending lamination for an exact sequence of groups. 
Given any point $\xi\in \partial Q$ he defined a lamination $\Lambda_{\xi}$ and then
showed that $\Lambda_{\xi}=\Lam{\xi,X} (F)$. However, for formulating and proving these sorts of results one needs additional
structure on the bundle, e.g. action of a group on the bundle through morphisms which has uniformly bounded quotients
when restricted to the fibers. Results of this type are proved in \cite[Section 3]{mj-rafi}; see also 
\cite[Section 4.4]{bowditch-stacks}.
\end{rem}

\noindent {\bf Acknowledgements:} The authors gratefully acknowledge all the helpful comments, inputs, and suggestions 
received from Mahan Mj and Michael Kapovich. We are very thankful to the anonymous referee also for suggesting many
changes that helped to improve the exposition of the paper and for pointing out a number of gaps and inaccuracies in
an earlier version of the paper.
The second author was partially supported by DST INSPIRE grant DST/INSPIRE/04/2014/002236 and DST MATRICS grant 
MTR/2017/000485 of the Govt of India. Finally, we thank Sushil Bhunia for a careful reading of an earlier
draft of the paper and for making numerous helpful suggestions.


\appendix
\section{Flaring in metric bundle and its canonical metric graph bundle}

Suppose $\pi':X'\map B'$ is an $(\eta,c)$-metric bundle and $\pi:X\map B$ is the canonical metric graph bundle associated to it.
{\em We shall assume that $B'$ and $B$ are both $\delta$-hyperbolic. However, there will be no assumption about the 
fibers of the bundles.}
We shall freely use the notation from section 4 of the paper. The purpose of this appendix is to show that
a metric bundle satisfies a sort of 'generalized flaring property' (see property $(\dagger)$ below) iff
the associated canonical metric graph bundle satisfies a flaring condition.

{\bf Note:} If $b_0,b_1,\cdots, b_n$ are consecutive vertices on a geodesic in $B$ then $\alpha': i\mapsto b_i$
is a dotted $(1,3)$-quasigeodesic of $B'$ by Lemma \ref{length space qi to graph}. 
Thus there is a constant $D_0$ such that if $\beta'$ is any
$(1,1)$-quasigeodesic in $B'$ joining $b_0, b_n$ then $Hd(\alpha', \beta')\leq D_0$. We will preserve $D_0$
to denote this constant for the rest of this section.

 Suppose $b\in V(B)$ and  $p\in B'$ are such that $d_{B'}(p,b)\leq D_0$. Then for any $x\in \pi^{-1}(b)$
we can lift a $(1,1)$-quasigeodesic of $B'$ joining $b$ to $p$ to $X'$ which starts from $x$ and ends at $x'$, say.
This way we get a `fiber identification map' $V(\pi^{-1}(b))\map \pi'^{-1}(p)$. If we denote this map by
$f_{bp}$ then we have the following lemma. Since the proof is evident we skip it.

\begin{lemma}\label{append 1}
We have $-C_0+\frac{1}{C_0}d_b(x,y)\leq d'_p( f_{bp}(x),f_{bp}(y))\leq C_0+C_0d_b(x,y)$ for all $x,y\in \pi^{-1}(b)$ and
for some uniform constant $C_0$ where $d_b$ is the fiber distance in $\pi^{-1}(b)$ for the metric graph bundle $X$ and
$d'_p$ is the fiber distance in $\pi'^{-1}(p)$ for the metric bundle $X'$.
\end{lemma}

Suppose $\alpha$ is a geodesic in $B$ and $\tilde{\alpha}$ is a $C$-qi lift of $\alpha$ in $X$. 
Let $\alpha'$ be a $(1,1)$-quasigeodesic in $B'$ joining the end points of $\alpha$. Let $\sigma: \alpha\map \alpha'$
be any map such that $d_{B'}(b, \sigma(b))\leq D_0$ for all $b\in \alpha$. Let $\tilde{p}=f_{bp}(\tilde{\alpha}(b))\in \pi'^{-1}(p)$
for all $b\in \alpha$ where $p=\sigma(b)$. Now it is easy to find a uniform qi lift $\tilde{\alpha}'$ of $\alpha'$ such that 
$\tilde{\alpha}'(\sigma(b))=\tilde{p}$ where $p=\sigma(b)$ for all $b\in \alpha$. We record this as a lemma.

\begin{lemma}\label{append 2}
There is a constant $C'$ depending on $C$ and a $C'$-qi lift $\tilde{\alpha}'$ of $\alpha'$ such that
$\tilde{\alpha}'(\sigma(b))=\tilde{p}$ where $p=\sigma(b)$ for all $b\in \alpha$.
\end{lemma}

The following lemma roughly says that if two qi leaves start flaring in one direction then they keep on flaring
in the same direction. The proof follows immediately from the definition of flaring. 
One may also look up the proof of \cite[Lemma 2.17(1)]{pranab-mahan}.

\begin{lemma}{\em (Persistence of flaring in graph bundles)}\label{append graph 1}
Suppose the metric graph bundle satisfies $(\nu_k, M_k,n_k)$-flaring condition for all $k\geq 1$.
Suppose $\alpha:[-m,n]\rightarrow B$ is a geodesic where $m\geq n_k, n\geq n_k$ and $\tilde{\alpha}_1$ and $\tilde{\alpha}_2$ are two
$k$-qi lifts of $\alpha$ in $X$ with $d_{\alpha(0)}(\tilde{\alpha}_1(0),\tilde{\alpha}_2(0))\geq M_k$. Suppose
$$d_{\alpha(sn_k)}(\tilde{\alpha}_1(sn_k),\tilde{\alpha}_2(sn_k))\geq \nu_k d_{\alpha(0)}(\tilde{\alpha}_1(0),\tilde{\alpha}_2(0))$$
where $s$ is either $1$ or $-1$. Let $t$ be the largest integer smaller than $n/n_k$ or $m/n_k$ according as $s=1$ or $-1$.
Then for all integer $1\leq l\leq t$ we have  
$$d_{\alpha(lsn_k)}(\tilde{\alpha}_1(lsn_k),\tilde{\alpha}_2(lsn_k))\geq \nu^l_k d_{\alpha(0)}(\tilde{\alpha}_1(0),\tilde{\alpha}_2(0)).$$	
\end{lemma}

The same idea of proof gives the next lemma also. We will need a definition.

\smallskip
{\bf Property $(\dagger)$:} {\em We shall say that the metric bundle $X'$ has the property $(\dagger)$ if
 for any $k \geq 1$, there exist $\nu_k>1$ and  $n_k,M_k\in \mathbb N$ such that the following holds:
		
Suppose $\alpha':[-n_k,n_k]\rightarrow B'$ is a $1$-quasigeodesic and $\tilde{\alpha}'_1$ and $\tilde{\alpha}'_2$ 
are two $k$-qi lifts of $\alpha'$ in $X'$. If $d_{\gamma(0)}(\tilde{\gamma_1}(0),\tilde{\gamma_2}(0))\geq M_k$
then we have {\small
$$ \nu_k\cdot d_{\alpha'(0)}(\tilde{\alpha}'_1(0),\tilde{\alpha}'_2(0))\leq \max\{d_{\alpha'(n_k)}(\tilde{\alpha}'_1(n_k),\tilde{\alpha}'_2(n_k)),d_{\alpha'(-n_k)}(\tilde{\alpha}'_1(-n_k),\tilde{\alpha}'_2(-n_k))\}.
$$}
}

Note that one could define flaring condition for a length metric bundle using the property $(\dagger)$.

\begin{lemma}{\em (Persistence of flaring in metric bundles)}\label{append 3}
Suppose the metric bundle satisfies $(\dagger)$. Let $k\geq 1$. 
Suppose $\alpha':[-m,n]\rightarrow B$ is a geodesic where $m\geq n_k, n\geq n_k$ and $\tilde{\alpha}'_1$ and $\tilde{\alpha}'_2$ 
are two $k$-qi lifts of $\alpha'$ in $X$ with $d_{\alpha'(0)}(\tilde{\alpha}'_1(0),\tilde{\alpha}'_2(0))\geq M_k$. Suppose
$$d_{\alpha'(sn_k)}(\tilde{\alpha}'_1(sn_k),\tilde{\alpha}'_2(sn_k))\geq \nu_k 
d_{\alpha'(0)}(\tilde{\alpha}'_1(0),\tilde{\alpha}'_2(0))$$ where $s$ is either $1$ or $-1$. Let $t$ be the largest integer 
smaller than or equal to $n/n_k$ or $m/n_k$ according as $s=1$ or $-1$.
Then for all integer $l\leq t$ we have  
$$d_{\alpha'(lsn_k)}(\tilde{\alpha}'_1(lsn_k),\tilde{\alpha}'_2(lsn_k))\geq \nu^l_k 
d_{\alpha'(0)}(\tilde{\alpha}'_1(0),\tilde{\alpha}'_2(0)).$$	
\end{lemma}

Following is one of the main results of this appendix.
\begin{lemma}\label{append part 1}
Suppose the metric bundles $X'$ has the property $(\dagger)$.
Then the canonical metric graph bundle $\pi:X\map B$ associated to $X'$ satisfies a 
$(\hat{\nu}_k,\hat{M}_k,\hat{\lambda}_k)$-flaring condition.

In particular if a geodesic metric bundle satisfies flaring condition (see \cite[Definition 1.12]{pranab-mahan}) then
its canonical metric graph bundle satisfies flaring condition.
\end{lemma}
\proof
Suppose $\alpha:[-n,n]\map B$ is a geodesic and $\tilde{\alpha}, \tilde{\tilde{\alpha}}$ are two $k$-qi lifts
of $\alpha$ in $X$ where $n\in \NN$ and $k\geq 1$. Let $\alpha'$ be a $(1,1)$-quasigeodesic in $B'$ joining
$\alpha(n), \alpha(-n)$. Then there are $k'$-qi lifts $\tilde{\alpha}', \tilde{\tilde{\alpha}}'$ of $\alpha'$ respectively
as in Lemma \ref{append 2}. We shall choose a parametrization $\alpha':[-m',n']\map B'$ so that
$\alpha'(n')=\alpha(n), \alpha'(-m')=\alpha(-n)$ and $d_{B'}(\alpha(0), \alpha'(0))\leq D_0$. Note that
$d'_{\alpha'(0)}(\tilde{\alpha}'(0), \tilde{\tilde{\alpha}}'(0))\geq 
-C_0+\frac{1}{C_0}d_{\alpha(0)}(\tilde{\alpha}(0), \tilde{\tilde{\alpha}}(0))$ by Lemma \ref{append 1}.
Hence, if we assume $d_{\alpha(0)}(\tilde{\alpha}(0), \tilde{\tilde{\alpha}}(0))\geq C_0(C_0+M_{k'})$
then $d'_{\alpha'(0)}(\tilde{\alpha}'(0), \tilde{\tilde{\alpha}}'(0))\geq M_{k'}$.
Clearly, if we choose $n$ large enough then we have $n_{k'}<\min\{m',n'\}$. (In the course of the proof
we will be more precise.) Without loss of generality we shall assume  
$\nu_k\cdot d'_{\alpha'(0)}(\tilde{\alpha}'(0), \tilde{\tilde{\alpha}}'(0))
\leq d'_{\alpha'(n_{k'})}(\tilde{\alpha}'(n_{k'}), \tilde{\tilde{\alpha}}'(n_{k'}))$. Let
$l$ be the greatest integer less than or equal to $n'/n_{k'}$. Then by Lemma \ref{append 3}

\begin{equation}
\nu^l_k\cdot d'_{\alpha'(0)}(\tilde{\alpha}'(0), \tilde{\tilde{\alpha}}'(0))
\leq d'_{\alpha'(ln_{k'})}(\tilde{\alpha}'(ln_{k'}), \tilde{\tilde{\alpha}}'(ln_{k'})).
\end{equation} 
Note that $d_{B'}(\alpha'(ln_{k'}), \alpha'(n'))\leq 2n_{k'}$. Let $b'=\alpha'(ln_{k'})$ and
$b''=\alpha'(n')$. Then by Corollary \ref{fibers unif qi} the fiber identification map $\phi_{b'b''}$
referred to in that corollary is a $K_{\ref{fibers unif qi}}(2n_{k'})$-quasiisometry.
Let $K=K_{\ref{fibers unif qi}}(2n_{k'})$. Note that 
\begin{equation}
d_{X'}(\tilde{\alpha}'(ln_{k'}),\tilde{\alpha}'(n'))\leq k'+2n_{k'}k'
\end{equation}
since $\tilde{\alpha}'$ is a $k'$-qi section and $d_{B'}(\alpha'(ln_{k'}), \alpha'(n'))\leq 2n_{k'}$.
Also, by Corollary \ref{path lifting remark} we have 
\begin{equation}
d_{X'}(\tilde{\alpha}'(ln_{k'}), \phi_{b'b''}(\tilde{\alpha}'(ln_{k'})))\leq 3c+6cn_{k'}.
\end{equation}
Using the inequalities (2) and (3) we have 
$$d_{X'}(\tilde{\alpha}'(n'), \phi_{b'b''}(\tilde{\alpha}'(ln_{k'})))\leq 3c+6cn_{k'}+k'+2n_{k'}k'.$$
Since $X'$ is an $(\eta, c)$-metric bundle we have
$$d'_{\alpha'(n')}(\tilde{\alpha}'(n'), \phi_{b'b''}(\tilde{\alpha}'(ln_{k'})))\leq \eta(3c+6cn_{k'}+k'+2n_{k'}k').$$
In the same way we have 
$$d'_{\alpha'(n')}(\tilde{\tilde{\alpha}}'(n'), \phi_{b'b''}(\tilde{\tilde{\alpha}}'(ln_{k'})))\leq \eta(3c+6cn_{k'}+k'+2n_{k'}k').$$
Now using the fact that $\phi_{b'b''}$ is a $K$-quasiisometry and letting $R_1=2\eta(3c+6cn_{k'}+k'+2n_{k'}k')$
we have by triangle inequality
\begin{equation} d'_{\alpha'(ln_{k'})}(\tilde{\alpha}'(ln_{k'}), \tilde{\tilde{\alpha}}'(ln_{k'}))\leq
(K^2+2R_1K)+ Kd'_{\alpha'(n')}(\tilde{\alpha}'(n'), \tilde{\tilde{\alpha}}'(n')).
\end{equation}
However, by Lemma \ref{length space qi to graph} and Proposition \ref{bundle vs graph bundle}(2) we have
\begin{equation} 
d'_{\alpha'(n')}(\tilde{\alpha}'(n'), \tilde{\tilde{\alpha}}'(n'))\leq 3+d_{\alpha(n)}(\tilde{\alpha}(n), 
\tilde{\tilde{\alpha}}(n))
\end{equation}
 Then it follows from the inequalities (1), (4) and (5) that
\begin{equation}
\nu^l_k\cdot d'_{\alpha'(0)}(\tilde{\alpha}'(0), \tilde{\tilde{\alpha}}'(0)) \leq
RK+Kd_{\alpha(n)}(\tilde{\alpha}(n), \tilde{\tilde{\alpha}}(n))
\end{equation} where $R=3+K+2R_1$. Finally since
$d'_{\alpha'(0)}(\tilde{\alpha}'(0), \tilde{\tilde{\alpha}}'(0))\geq 
-C_0+\frac{1}{C_0}d_{\alpha(0)}(\tilde{\alpha}(0), \tilde{\tilde{\alpha}}(0))$
using (6) we have
\begin{equation}
\nu^l_{k'} (-C_0+\frac{1}{C_0}d_{\alpha(0)}(\tilde{\alpha}(0), \tilde{\tilde{\alpha}}(0))) \leq
RK+ Kd_{\alpha(n)}(\tilde{\alpha}(n), \tilde{\tilde{\alpha}}(n)).
\end{equation}
Recall that we assumed $d_{\alpha(0)}(\tilde{\alpha}(0), \tilde{\tilde{\alpha}}(0)))\geq C_0(C_0+M_{k'})$.
Hence, 
\begin{equation}-C_0+\frac{1}{C_0}d_{\alpha(0)}(\tilde{\alpha}(0), \tilde{\tilde{\alpha}}(0))\geq 
(\frac{1}{C_0}-\frac{1}{C_0+M_{k'}})d_{\alpha(0)}(\tilde{\alpha}(0), \tilde{\tilde{\alpha}}(0)).
\end{equation}
Let
$$\lambda=\frac{1}{K}(\frac{1}{C_0}-\frac{1}{C_0+M_{k'}}).$$
Then we have, using (7) and (8),
\begin{equation}
\nu^l_{k'}\cdot \lambda\cdot d_{\alpha(0)}(\tilde{\alpha}(0), \tilde{\tilde{\alpha}}(0)) \leq
R+d_{\alpha(n)}(\tilde{\alpha}(n),\tilde{\tilde{\alpha}}(n)).
\end{equation}
It is clear that 
$$-R+\nu^l_{k'}\cdot \lambda\cdot d_{\alpha(0)}(\tilde{\alpha}(0), \tilde{\tilde{\alpha}}(0))
\geq \frac{1}{2}\lambda \nu^l_{k'} d_{\alpha(0)}(\tilde{\alpha}(0), \tilde{\tilde{\alpha}}(0))$$
if $d_{\alpha(0)}(\tilde{\alpha}(0), \tilde{\tilde{\alpha}}(0))\geq \frac{2R}{\lambda\nu^l_{k'}}.$
In particular, since $\nu^l_{k'}>1$, we have
\begin{equation}
\frac{1}{2}\lambda\nu^l_{k'} d_{\alpha(0)}(\tilde{\alpha}(0), \tilde{\tilde{\alpha}}(0))\leq
d_{\alpha(n)}(\tilde{\alpha}(n),\tilde{\tilde{\alpha}}(n))
\end{equation} using (9)
if $d_{\alpha(0)}(\tilde{\alpha}(0), \tilde{\tilde{\alpha}}(0))\geq \frac{2R}{\lambda}$.
Thus it is enough to choose 
$$\hat{M}_k=\max\{\frac{2R}{\lambda}, C_0(C_0+M_{k'})\}, \,\hat{\lambda}_k=2$$ and to show that if
$n$ is sufficiently large then $l$ is so large that $\frac{1}{2}\lambda\nu^l_{k'}\geq 2$ which will give a choice
for $\hat{n}_k$. This is easy to verify and hence left to the reader.
\qed

Converse of Lemma \ref{append part 1} is also true and has an exactly similar proof. 
However, in this case one uses Lemma \ref{append graph 1} instead of Lemma \ref{append 2}.
We state it without proof to avoid repetition.

\begin{lemma}\label{append part 2}
Suppose the metric graph bundle $\pi:X\map B$ satisfies a $(\nu_k, M_k,n_k)$-flaring
condition for all $k\geq 1$. Then the metric bundle $\pi':X'\map B'$ satisfies the
condition $(\dagger)$ for three functions $\nu'_k, M'_k, n'_k$ of $k$.
\end{lemma}

\bibliography{pullback}
\bibliographystyle{amsalpha}

\end{document}